\documentclass{amsart}

\usepackage{amsmath,amssymb,amsthm,texdraw,color}
\usepackage{multicol}

\numberwithin{equation}{section}

\newsavebox{\tmpfiga}
\newsavebox{\tmpfigb}
\newsavebox{\tmpfigc}
\newsavebox{\tmpfigd}
\newsavebox{\tmpfige}
\newsavebox{\tmpfigf}
\newsavebox{\tmpfigg}
\newsavebox{\tmpfigh}
\newsavebox{\tmpfigi}
\newsavebox{\tmpfigj}
\newsavebox{\tmpfigk}
\newsavebox{\tmpfigl}
\newsavebox{\tmpfigm}
\newsavebox{\tmpfign}

\theoremstyle{plain}
\newtheorem{thm}{Theorem}[section]
\newtheorem{prop}[thm]{Proposition}
\newtheorem{lem}[thm]{Lemma}

\theoremstyle{definition}
\newtheorem{df}[thm]{Definition}

\newtheorem{example}[thm]{Example}
\theoremstyle{remark}
\newtheorem{remark}[thm]{Remark}

\newcommand{\nc}{\newcommand}

\nc{\eit}{\tilde{e}_i}
\nc{\fit}{\tilde{f}_i}
\nc{\ali}{\alpha_i}
\nc{\csa}{\mathfrak{h}}
\nc{\op}{\oplus}
\nc{\ot}{\otimes}
\nc{\C}{\mathbf{C}}
\nc{\Z}{\mathbf{Z}}
\nc{\g}{\mathfrak{g}}
\nc{\half}{\frac{1}{2}}
\nc{\slice}{\mathcal{S}^{(l)}}
\nc{\uq}{U_q}
\nc{\uqp}{U_q'}
\nc{\ftil}{\tilde{f}}
\nc{\etil}{\tilde{e}}
\nc{\newcrystal}{\mathcal{C}^{(l)}}
\nc{\anone}{A_n^{(1)}}
\nc{\cnone}{C_n^{(1)}}
\nc{\dnone}{D_n^{(1)}}
\nc{\atwontwo}{A_{2n}^{(2)}}
\nc{\dnponetwo}{D_{n+1}^{(2)}}
\nc{\bnone}{B_n^{(1)}}
\nc{\bthreeone}{B_3^{(1)}}
\nc{\atwonmonetwo}{A_{2n-1}^{(2)}}
\nc{\oldcrystal}{\mathcal{B}^{(l)}}
\nc{\defi}[1]{\emph{\textbf{#1}}}
\nc{\veps}{\varepsilon}
\nc{\vphi}{\varphi}
\nc{\La}{\Lambda}
\nc{\la}{\lambda}
\nc{\pathspace}{\mathcal{P}}
\nc{\pwspace}{\mathcal{Z}}
\nc{\rpwspace}{\mathcal{Y}}
\nc{\gspath}{\mathbf{b}}
\nc{\path}{\mathbf{p}}
\nc{\spaceiso}{\Psi}
\nc{\spaceosi}{\Phi}
\nc{\gswall}{\mathbf{Y}}
\nc{\wall}{\mathbf{Y}}
\nc{\hwc}{\mathcal{B}}
\nc{\bu}{\mathbf{u}}
\nc{\bv}{\mathbf{v}}
\nc{\bw}{\mathbf{w}}

\nc{\Na}{$(\overline{n})_{\underline{n\!-\!1}}$}
\nc{\Nb}{$(\underline{n\!-\!1})_{\overline{n}}$}

\nc{\Loo}{L11}
\nc{\Lot}{L12}
\nc{\Lto}{L21}
\nc{\Ltt}{L22}
\nc{\Roo}{R11}
\nc{\Rot}{R12}
\nc{\Rto}{R21}
\nc{\Rtt}{R22}

\nc{\Coo}{C11}
\nc{\Cot}{C12}
\nc{\Cto}{C21}
\nc{\Ctt}{C22}

\DeclareMathOperator{\wt}{wt}
\DeclareMathOperator{\cwt}{\overline{wt}}

\begin{document}

\title[Higher Level Affine Crystals and Young Walls]
{Higher level affine crystals and Young walls}

\author{Seok-Jin Kang$^*$}
\thanks{$^*$This research was supported in part by
        KOSEF Grant R01-2003-000-10012-0
        and KRF Grant 2003-070-C00001.}
\author{Hyeonmi Lee$^{**}$}
\thanks{$^{**}$This research was supported in part by
        KOSEF Grant R01-2003-000-10012-0}
\address{Korea Institute for Advanced Study\\
         207-43 Cheongryangri-dong, Dongdaemun-gu\\
         Seoul 130-722, Korea}
\email{sjkang@kias.re.kr, hmlee@kias.re.kr}

\begin{abstract}
Using combinatorics of Young walls, we give a new realization of
arbitrary level irreducible highest weight crystals
$\mathcal{B}(\la)$ for quantum affine algebras of type $\anone$,
$\bnone$, $\cnone$, $\atwonmonetwo$, $\atwontwo$, and
$\dnponetwo$. The irreducible highest weight crystals are realized
as the affine crystals consisting of reduced proper Young walls.
The notion of slices and splitting of blocks plays a crucial role
in the construction of crystals.
\end{abstract}

\maketitle
%\tableofcontents

%%%%%%%%%%%%%%%%%%%%%%%%%%%%%%%%%%%%%%%%%%%%%%%%%%%%%%%%%%%%%
%%%%%%%%%%%%%%%%%%%%%%%%%%%%%%%%%%%%%%%%%%%%%%%%%%%%%%%%%%%%%
\section{Introduction}\label{sec:01}

The {\it crystal bases}, introduced by Kashiwara in
\cite{MR93b:17045}, have many nice combinatorial features
reflecting the internal structure of integrable modules of quantum
groups. Moreover, it is known that the crystal bases are preserved
under the direct sum decomposition and have extremely simple
behavior with respect to taking the tensor product. Hence it is a
very natural and important problem to find explicit realizations
of crystal bases for irreducible highest weight modules over
quantum groups. (See \cite{MR1881971}, for example, and the
references there in.)

In \cite{Kang2000}, Kang introduced the notion of {\it Young
walls} as a new combinatorial scheme for realizing crystal bases
for quantum affine algebras. The Young walls consists of colored
blocks with various shapes, and can be viewed as generalizations
of colored Young diagrams.

For classical quantum affine algebras of type $\anone$, $\bnone$,
$\dnone$, $\atwonmonetwo$, $\atwontwo$, $\dnponetwo$, the rules
and patterns for building Young walls and the action of Kashiwara
operators are given explicitly in terms of combinatorics of Young
walls, which defines an affine crystal structure on the set of
{\it proper Young walls}. In particular, the level-$1$ irreducible
highest weight crystals are realized as the affine crystals
consisting of {\it reduced proper Young walls}.

However, in \cite{Kang2000}, the problem of Young wall realization
of crystal bases were left open for quantum affine algebras of
type $\cnone$ $(n \ge 2)$, because this case is more difficult to
deal with than the other classical quantum affine algebras. The
main difficulty lies in the fact that the level-$1$ perfect crystal
for this case are intrinsically of level-$2$.

This difficulty was resolved in \cite{CnOne} by introducing the
notion of {\it slices} and {\it splitting blocks}, which plays a
crucial role in constructing the desired realization of crystal
bases for quantum affine algebras of type $\cnone$: the level-$1$
irreducible highest weights crystals are realized as the affine
crystals consisting of reduced proper Young walls.

In this paper, we develop the combinatorics of higher level Young
walls for classical quantum affine algebras of type $\anone$,
$\bnone$, $\cnone$, $\atwonmonetwo$, $\atwontwo$, and
$\dnponetwo$. As in \cite{CnOne}, we first introduce the notion of
slices and splitting blocks, and give a new realization of higher
level perfect crystals as the equivalence classes of slices. We
then proceed to define the notion of higher level Young walls,
proper Young walls, reduced proper Young walls, and ground-state walls, etc.
Finally, we prove that the arbitrary level irreducible
highest weight crystals are realized as the affine crystals
consisting of reduced proper Young walls.

\vspace{2mm}
\noindent
{\bf Acknowledgments.} \ We would like to express our sincere
gratitude to Jin Hong for his interest in this work
and many valuable discussions.

\vskip 1cm

%%%%%%%%%%%%%%%%%%%%%%%%%%%%%%%%%%%%%%%%%%%%%%%%%%%%%%%%%%%%%
%%%%%%%%%%%%%%%%%%%%%%%%%%%%%%%%%%%%%%%%%%%%%%%%%%%%%%%%%%%%%
\section{Quantum affine algebras and perfect crystals}%
\label{sec:02}

In this section, we fix the notations for quantum affine algebras
and review some of the basic properties of perfect crystals. Let
$(A, P^{\vee}, \Pi^{\vee}, P, \Pi)$ be an {\it affine Cartan
datum}, where
\begin{itemize}
\item $I=\{0,1, \dots, n \}$ is the index set for the simple roots,
\item $A=(a_{ij})_{i, j\in I}$ is an affine generalized Cartan
      matrix of type $\anone$, $\bnone$, $\cnone$, $\atwonmonetwo$,
      $\atwontwo$, or $\dnponetwo$,
\item $P^{\vee} = \left(\bigoplus_{i\in I} \Z h_i\right) \op \Z d$ is
      the {\it dual weight lattice},
\item $\Pi^{\vee} = \{ h_i | i \in I \}$ is the set of {\it simple
      coroots},
\item $\csa = \C \ot_{\Z} P^{\vee}$ is the {\it Cartan subalgebra},
\item $P= \{ \la \in \csa^* | \la(P^{\vee}) \subset \Z \}$ is the
      {\it weight lattice},
\item $\Pi = \{ \ali | i \in I \}$ is the set of {\it simple roots}.
\end{itemize}

\noindent We denote the {\it null root} by $\delta$ and the {\it
fundamental weights} by $\La_i$ $(i=0,1, \dots, n)$, so that we
have
\begin{equation*}
P = \left( \bigoplus_{i\in I} \Z\La_i \right) \op \frac{1}{d_0} \Z\delta,
\end{equation*}
Here, $d_0$ is the coefficient of $\alpha_0$ in the null root
$\delta=d_0\alpha_0+\cdots+d_n\alpha_n$.
We also denote by $P^{+} = \{ \la \in P | \la(h_i) \in \Z_{\ge 0}
\ \ \text{for all} \ i \in I \}$ the set of {\it affine dominant
integral weights}.

Let $\uq(\g)$ be the {\it quantum affine algebra} associated with
the affine Cartan datum $(A, P^{\vee}, \Pi^{\vee}, P, \Pi)$, and
let $e_i$, $f_i$, $K_i^{\pm 1}$, $q^d$ $(i \in I)$ be the
generators of $\uq(\g)$. Let $\uq'(\g)$ be the subalgebra of
$\uq(\g)$ generated by $e_i$, $f_i$, $K_i^{\pm 1}$ $(i\in I)$.
Then $\uq'(\g)$ can be regarded as the quantum group associated
with the \emph{classical Cartan datum} $(A, \bar{P}^{\vee},
\Pi^{\vee}, \bar{P}, \Pi)$, where $\bar{P}^{\vee} =
\bigoplus_{i\in I} \Z h_i$ is the \emph{classical dual weight
lattice} and $\overline{P} = \bigoplus_{i\in I} \Z \La_i$ is the
\emph{classical weight lattice}.

\begin{df} \label{df:abstract crystal}
An \defi {affine crystal} (resp. \defi{classical crystal})
is a set $\mathcal{B}$ together with the maps
$\wt: \mathcal{B} \rightarrow P$\,(resp. $\wt: \mathcal{B} \rightarrow \bar{P}$),
$\veps_i, \vphi_i : \mathcal{B} \rightarrow \Z \cup \{- \infty \}$,
$\eit, \fit : \mathcal{B} \rightarrow \mathcal{B} \cup \{0\}$
satisfying the following conditions\, :\,
for all $i\in I$ and $b\in \mathcal{B}$,
\begin{enumerate}
\item $\vphi_i(b) = \veps_i(b) + \langle h_i, \wt(b) \rangle$;
\item $\wt(\eit b) = \wt(b) + \ali$
      if $\eit b \in \mathcal{B}$;
\item $\wt(\fit b) = \wt(b) - \ali$
      if $\fit b \in \mathcal{B}$;
\item if $\eit b \in \mathcal{B}$, then
  \begin{equation*}
  \veps_i(\eit b) = \veps_i(b) - 1, \quad \vphi_i(\eit b) =
  \vphi_i(b) + 1 ;
  \end{equation*}
\item if $\fit b \in \mathcal{B}$, then
  \begin{equation*}
  \veps_i(\fit b) = \veps_i(b) + 1, \quad \vphi_i(\fit b) =
  \vphi_i(b) - 1 ;
  \end{equation*}
\item $\fit b = b'$ if and only if $b = \eit b'$ for all $i\in I$
  and $b, b' \in \mathcal{B}$\,;
\item if $\veps_i(b)=-\infty$, then $\eit b = \fit b =0$.
\end{enumerate}
\end{df}

For example, for an affine dominant integral weight $\la \in
P^{+}$, the crystal graph $\hwc(\la)$ of the irreducible highest
weight module $V(\la)$ is an affine crystal, which will be called
the \defi{irreducible highest weight crystal}. An affine crystal
(resp. classical crystal) will also be called a
\defi {$\uq(\g)$-crystal} (resp. \defi{$\uq'(\g)$-crystal}). We
will often denote by $\overline{\wt}$ the weight function of a
classical crystal.

\begin{df} \label{df:crystal morphism}
Let $\mathcal{B}_1$ and $\mathcal{B}_2$ be (affine or classical) crystals. A
\defi{crystal morphism} $\psi: \mathcal{B}_1 \rightarrow \mathcal{B}_2$ is a map
$\psi: \mathcal{B} \cup \{0\} \rightarrow \mathcal{B} \cup \{0\}$ satisfying the
following conditions\,:
\begin{enumerate}
\item $\psi(0)=0$;
\item if $b \in \mathcal{B}_1$ and $\psi(b) \in \mathcal{B}_2$, then
  \begin{equation*}
  \wt(\psi(b))=\wt(b), \quad \veps_i(\psi(b))=\veps_i(b), \quad
  \vphi_i(\psi(b)) = \vphi_i(b) ;
  \end{equation*}
\item if $b, b' \in \mathcal{B}_1$, $\psi(b), \psi(b') \in \mathcal{B}_2$ and $\fit b
  = b'$, then
  \begin{equation*}
  \fit \psi(b) = \psi(b'), \quad \psi(b) = \eit\psi(b').
  \end{equation*}
\end{enumerate}
\end{df}

Let $\mathcal{B}$ be a classical crystal. For $b\in \mathcal{B}$, we define
\begin{equation*}
\veps(b)= \sum_{i \in I} \veps_i (b) \La_i, \qquad
\vphi(b)=\sum_{i\in I} \vphi_i(b) \La_i.
\end{equation*}

\begin{df}
For each positive integer $l>0$, a finite classical crystal $\mathcal{B}$ is
called a \defi{perfect crystal} of level-$l$ if

\begin{enumerate}
\item there is a finite dimensional $\uq'(\g)$-module with a
crystal basis whose crystal graph is isomorphic to $\mathcal{B}$,

\item $\mathcal{B} \ot \mathcal{B}$ is connected,

\item there exists some $\la_0 \in \bar{P}$ such that
\begin{equation*}
\wt(\mathcal{B}) \subset \la_0 + \frac{1}{d_0} \sum_{i \neq 0} \Z_{\le 0}
\ali, \qquad \#(\mathcal{B}_{\la_0})=1,
\end{equation*}

\item for any $b \in \mathcal{B}$, we have $\langle c, \veps(b) \rangle \ge
l$,

\item for each $\la \in \bar{P}^{+}$ with $\langle c, \la
\rangle =l$, there exist unique vectors $b^{\la} \in \mathcal{B}$ and $b_{\la} \in
\mathcal{B}$ such that
\begin{equation*}
\veps(b^{\la}) = \la, \qquad \vphi(b_{\la}) = \la.
\end{equation*}
\end{enumerate}
\end{df}
Here, $d_0$ is the coefficient of $\alpha_0$ in the null root
$\delta$. We recall the following fundamental crystal isomorphism
theorem proved in \cite{MR94a:17008}.

\begin{prop}\cite{MR94a:17008}\label{prop:KMN1}
Let $\mathcal{B}$ be a perfect crystal of level-$l>0$. Then for any dominant
integral weight $\la \in \bar{P}^{+}$ of level-$l$, there
exists a crystal isomorphism
\begin{equation*}
\Psi: \mathcal{B}(\la) \overset{\sim} \longrightarrow \mathcal{B}(\veps(b_{\la})) \ot
\mathcal{B} \qquad \text{given by} \qquad u_{\la} \longmapsto
u_{\veps(b_{\la})} \ot b_{\la},
\end{equation*}
where $b_{\la}$ is the unique element in $\mathcal{B}$ such that
$\vphi(b_{\la})=\la$ and $u_{\la}$ {\rm (}resp.
$u_{\veps(b_{\la})}${\rm )} is the highest weight vector of
$\mathcal{B}(\la)$ {\rm (}resp. $\mathcal{B}(\veps(b_{\la}))${\rm )}.
\end{prop}

For $k\ge 0$, set
\begin{equation*}
\la_0 = \la, \quad \la_{k+1} = \veps(b_{\la_k}), \qquad \text{and}
\qquad b_{0} = b_{\la}, \quad b_{k+1} = b_{\la_{k+1}}.
\end{equation*}
By taking the composition of crystal isomorphism given in
Proposition \ref{prop:KMN1}, we get a crystal isomorphism
\begin{equation*}
\Psi_{k} : \mathcal{B}(\la) \overset{\sim} \longrightarrow \mathcal{B}(\la_k) \ot
\mathcal{B}^{\ot k}
\end{equation*}
given by
\begin{equation*}
u_{\la} \longmapsto u_{\la_k} \ot b_{k-1} \ot \cdots \ot b_{1} \ot
b_{0}.
\end{equation*}

The sequence
\begin{equation*}
\path_{\la} = (b_k)_{k=0}^{\infty} = \cdots \ot
b_{k+1} \ot b_{k} \ot \cdots \ot b_{1} \ot b_{0}
\end{equation*}
is called the
\defi{ground-state path} of weight $\la$. A \defi{$\la$-path} in
$\mathcal{B}$ is a sequence
\begin{equation*}
\path = (\path(k))_{k=0}^{\infty} = \cdots \ot \path(k+1) \ot
\path(k) \ot \cdots \path(1) \ot \path(0)
\end{equation*}
in $\mathcal{B}$ such that
$\path(k)=b_{k}$ for all $k \gg 0$. Let $\pathspace(\la)$ denote
the set of all $\la$-paths. Then we can define a classical crystal
structure on $\pathspace(\la)$ by the tensor product rule, which
gives the \defi{path realization} of the irreducible highest
weight crystal $\hwc(\la)$.

\begin{prop} \cite{MR94a:17008} \label{prop:path realization}
There exists an isomorphism of classical crystals
\begin{equation*}
\Psi: \hwc(\la) \overset{\sim} \longrightarrow \pathspace(\la) \qquad
\text{given by} \quad u_{\la} \longmapsto \path_{\la}.
\end{equation*}
\end{prop}

\vskip 3mm For each of classical quantum affine algebras, it was
shown in \cite{MR95h:17016, MR94j:17013} that there exists a
\emph{coherent family of perfect crystals} $\{ \oldcrystal | l \in
\Z_{>0} \}$ . In the following, we will give an explicit
description of these perfect crystals.

\vskip 3mm

\begin{enumerate}
\item $\anone$ $(n\geq 1)$\\[1.5mm]
$\oldcrystal =
 \Big\{ (x_0,x_1,\dots,x_n) \,\Big\vert\, x_i \in \Z_{\geq0},\
\textstyle\sum_{i=0}^{n} x_i = \textnormal{$l$} \Big\}$.
\item $\bnone$ $(n\ge 3)$\\[1.5mm]
$\oldcrystal =
 \left\{
 (x_1,\dots,x_n |x_0| \bar{x}_n,\dots, \bar{x}_1) \, \Bigg\vert\,
 \begin{aligned}
 & x_0=0 \textup{ or } 1,\ x_i, \bar{x}_i \in \Z_{\geq0}, \\
 & x_0 + \textstyle\sum_{i=1}^{n} (x_i + \bar{x}_i) = \textnormal{$l$}
 \end{aligned}
 \right\}$.
\item $\cnone$ $(n \geq 2)$\\[1.5mm]
$\oldcrystal =
 \Big\{ (x_1,\dots,x_n | \bar{x}_n,\dots, \bar{x}_1) \,\Big\vert\,
x_i, \bar{x}_i \in \Z_{\geq0},\
\textnormal{$2l$} \geq \textstyle\sum_{i=1}^{n} (x_i +\bar{x}_i) \in 2\Z \Big\}$.
\item $\atwonmonetwo$ $(n\geq 3)$\\[1.5mm]
$\oldcrystal =
 \Big\{ (x_1,\dots,x_n | \bar{x}_n,\dots, \bar{x}_1) \,\Big\vert\,
x_i, \bar{x}_i \in \Z_{\geq0},\
\textstyle\sum_{i=1}^{n} (x_i + \bar{x}_i) = \textnormal{$l$} \Big\}$.
\item $\atwontwo$ $(n\geq 1)$\\[1.5mm]
$\oldcrystal =
 \Big\{ (x_1,\dots,x_n | \bar{x}_n,\dots, \bar{x}_1) \,\Big\vert\,
x_i, \bar{x}_i \in \Z_{\geq0},\
\textstyle\sum_{i=1}^{n} (x_i + \bar{x}_i) \leq \textnormal{$l$} \Big\}$.
\item $\dnponetwo$ $(n\geq 2)$\\[1.5mm]
$\oldcrystal =
 \left\{
 (x_1,\dots,x_n |x_0| \bar{x}_n,\dots, \bar{x}_1) \,\Bigg\vert\,
 \begin{aligned}
 & x_0 = 0 \textup{ or } 1,\ x_i, \bar{x}_i \in \Z_{\geq0},\\
 & x_0 + \textstyle\sum_{i=1}^{n} (x_i + \bar{x}_i) \leq \textnormal{$l$}
 \end{aligned}
 \right\}$.
\end{enumerate}

\vskip 3mm For the reader's convenience, we also give an explicit
description of the maps $\overline{\wt}$, $\eit$, $\fit$,
$\veps_i$ and $\vphi_i$ $(i\in I)$ for this coherent family of
perfect crystals over the quantum affine algebras of type
$\bnone$. For the other quantum affine algebras, see
\cite{MR95h:17016, MR99h:17008}.

Let  $b = (x_1,\dots,x_n|x_0 |\bar{x}_n,\dots, \bar{x}_1) \in
\oldcrystal$ so that $x_0=0$ or $x_0=1$, $x_i, \bar{x}_i \in
\Z_{\ge 0}$, and $x_0 + \sum x_i + \sum \bar{x}_i = l$.

For $i=0$ and $i=n$, the Kashiwara operators are given by
\begin{equation*}
\begin{aligned}
\tilde{e}_0 b & =
\begin{cases}
(x_1,x_2-1,x_3,\dots,x_n|x_0|\bar{x}_n,\dots,\bar{x}_2,\bar{x}_1+1)
&
 \textnormal{if $x_2 > \bar{x}_2$,}\\
(x_1-1,x_2,\dots,x_n|x_0|\bar{x}_n,\dots,\bar{x}_3,\bar{x}_2+1,\bar{x}_1)
&
 \textnormal{if $x_2 \le \bar{x}_2$,}
\end{cases}\\
\tilde{f}_0 b &=
\begin{cases}
(x_1,x_2+1,x_3,\dots,x_n|x_0|\bar{x}_n,\dots,\bar{x}_2,\bar{x}_1-1)
&
 \textnormal{if $x_2 \geq \bar{x}_2$,}\\
(x_1+1,x_2,\dots,x_n|x_0|\bar{x}_n,\dots,\bar{x}_3,\bar{x}_2-1,\bar{x}_1)
&
 \textnormal{if $x_2 < \bar{x}_2$,}
 \end{cases}\\
\tilde{e}_n b & =
\begin{cases}
(x_1,\dots,x_{n-1},x_n + 1|x_0-1|\bar{x}_n,\dots,\bar{x}_1) &
\textnormal{if $x_0=1$,}\\
(x_1,\dots,x_n|x_0+1|\bar{x}_n - 1,\bar{x}_{n-1},\dots,\bar{x}_1)
& \textnormal{if $x_0=0$,}
\end{cases}\\
\tilde{f}_n b & =
\begin{cases}
(x_1,\dots,x_{n-1},x_n - 1|x_0+1|\bar{x}_n,\dots,\bar{x}_1) &
\textnormal{if $x_0=0$,}\\
(x_1,\dots,x_n|x_0-1|\bar{x}_n + 1,\bar{x}_{n-1},\dots,\bar{x}_1)
& \textnormal{if $x_0=1$.}
\end{cases}
\end{aligned}
\end{equation*}

For $i=1,\dots,n-1$, we have
\begin{equation*}
\begin{aligned}
\eit b & =
\begin{cases}
(x_1,\dots,x_i+1,x_{i+1}-1,\dots,x_n|x_0|
\bar{x}_n,\dots,\bar{x}_1) &
\textnormal{if $x_{i+1} > \bar{x}_{i+1}$,}\\
(x_1,\dots,x_n|x_0| \bar{x}_n,\dots,\bar{x}_{i+1}+1,
\bar{x}_{i}-1,\dots,\bar{x}_1) & \textnormal{if $x_{i+1} \leq
\bar{x}_{i+1}$,}
\end{cases} \\
\fit b & =
\begin{cases}
(x_1,\dots,x_i-1,x_{i+1}+1,\dots,x_n|x_0|
\bar{x}_n,\dots,\bar{x}_1) &
\textnormal{if $x_{i+1} \geq \bar{x}_{i+1}$,}\\
(x_1,\dots,x_n|x_0| \bar{x}_n,\dots,,\bar{x}_{i+1}-1,
\bar{x}_{i}+1,\dots,\bar{x}_1) & \textnormal{if $x_{i+1} <
\bar{x}_{i+1}$.}
\end{cases}
\end{aligned}
\end{equation*}

The remaining maps are given below.
\begin{align*}
&\vphi_0 (b) =
  \bar{x}_1+(\bar{x}_2-x_2)_+,\\
&\vphi_i (b) =
  x_i + (\bar{x}_{i+1} - x_{i+1})_+ \quad (i=1,\dots,n-1),\\
&\vphi_n (b) =
  2x_n+x_0,\\
&\veps_0 (b) =
  x_1+(x_2-\bar{x}_2)_+,\\
&\veps_i (b) =
  \bar{x}_i+({x}_{i+1}-\bar{x}_{i+1})_+ \quad (i=1,\dots,n-1),\\
&\veps_n (b) =
  2\bar{x}_n+x_0,\\
&\cwt(b) =
  \sum_{i=0}^n (\vphi_i(b) - \veps_i(b))\La_i,
\end{align*}
where we use the notation $(x)_+ = \max(0,x)$.

\vskip 1cm

%%%%%%%%%%%%%%%%%%%%%%%%%%%%%%%%%%%%%%%%%%%%%%%%%%%%%%%%%%%%%%%%%%
%%%%%%%%%%%%%%%%%%%%%%%%%%%%%%%%%%%%%%%%%%%%%%%%%%%%%%%%%%%%%%%%%%
\section{Slices and splitting of blocks}
In this section, we introduce the notion of \emph{slices} and
\emph{splitting blocks}, and define a classical crystal structure
on the set $\newcrystal$ of equivalence classes of slices. In the
next section, we will show that the classical crystal
$\newcrystal$ is isomorphic to the level-$l$ perfect crystal
$\oldcrystal$.

To build a slice, we use the following \defi{colored blocks} of
three different types.\\
\savebox{\tmpfiga}
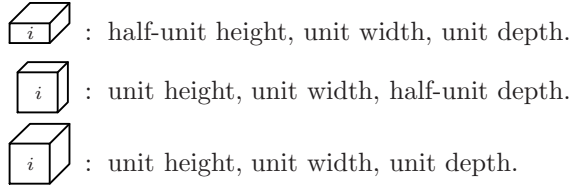
{
\begin{texdraw}
\fontsize{7}{7}\selectfont
\drawdim em
\textref h:C v:C
\setunitscale 1.9
\move(-1 0)\lvec(0 0)\lvec(0 0.5)\lvec(-1 0.5)\lvec(-1 0)
\move(0 0)\lvec(0.5 0.5)\lvec(0.5 1)\lvec(-0.5 1)\lvec(-1 0.5)
\move(0 0.5)\lvec(0.5 1)
\htext(-0.5 0.25){$i$}
\end{texdraw}}%
\savebox{\tmpfigb}{
\begin{texdraw}
\fontsize{7}{7}\selectfont
\drawdim em
\textref h:C v:C
\setunitscale 1.9
\move(-1 0)\lvec(0 0)\lvec(0 1)\lvec(-1 1)\lvec(-1 0)
\move(0 0)\lvec(0.25 0.25)\lvec(0.25 1.25)\lvec(-0.75 1.25)\lvec(-1 1)
\move(0 1)\lvec(0.25 1.25)
\htext(-0.5 0.5){$i$}
\end{texdraw}}
\savebox{\tmpfigc}{
\begin{texdraw}
\fontsize{7}{7}\selectfont
\drawdim em
\textref h:C v:C
\setunitscale 1.9
\move(-1 0)\lvec(0 0)\lvec(0 1)\lvec(-1 1)\lvec(-1 0)
\move(0 0)\lvec(0.5 0.5)\lvec(0.5 1.5)\lvec(-0.5 1.5)\lvec(-1 1)
\move(0 1.0)\lvec(0.5 1.5)
\htext(-0.5 0.5){$i$}
\end{texdraw}}%
\begin{align*}
\raisebox{-0.1em}{\usebox{\tmpfiga}}
 &\ : \text{\ half-unit height, unit width, unit depth.}\\
\raisebox{-0.5em}{\usebox{\tmpfigb}}
 &\ : \text{\ unit height, unit width, half-unit depth.}\\
\raisebox{-0.5em}{\usebox{\tmpfigc}}
 &\ : \text{\ unit height, unit width, unit depth.}\\
\end{align*}
The coloring of a block will be given differently according to the
types of blocks and the types of quantum affine algebras. For
simplicity, we will use the following notations\,:  \\[2mm]
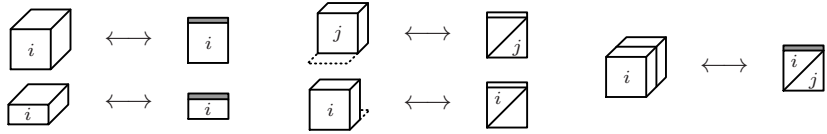
\begin{center}
\begin{tabular}{rcl}
\raisebox{-0.4\height}{
\begin{texdraw}
\fontsize{7}{7}\selectfont
\drawdim em
\textref h:C v:C
\setunitscale 1.9
\move(-1 0)\lvec(0 0)\lvec(0 1)\lvec(-1 1)\lvec(-1 0)
\move(0 0)\lvec(0.5 0.5)\lvec(0.5 1.5)\lvec(-0.5 1.5)\lvec(-1 1)
\move(0 1.0)\lvec(0.5 1.5)
\htext(-0.5 0.5){$i$}
\end{texdraw}
}
& $\longleftrightarrow$ &
\raisebox{-0.4\height}{
\begin{texdraw}
\fontsize{7}{7}\selectfont
\drawdim em
\textref h:C v:C
\setunitscale 1.9
\move(-1 1)\lvec(-1 1.133)\lvec(0 1.133)\lvec(0 1)\ifill f:0.6
\move(-1 0)\lvec(0 0)\lvec(0 1)\lvec(-1 1)\lvec(-1 0)
\move(-1 1)\lvec(-1 1.133)\lvec(0 1.133)\lvec(0 1)
\htext(-0.5 0.5){$i$}
\end{texdraw}
}\\[4mm]
\raisebox{-0.4\height}{
\begin{texdraw}
\fontsize{7}{7}\selectfont
\drawdim em
\textref h:C v:C
\setunitscale 1.9
\move(0 0)\lvec(1 0)\lvec(1 0.5)\lvec(0 0.5)\lvec(0 0)
\move(1 0)\lvec(1.5 0.5)\lvec(1.5 1)\lvec(0.5 1)\lvec(0 0.5)
\move(1 0.5)\lvec(1.5 1)
\htext(0.5 0.25){$i$}
\end{texdraw}
}
& $\longleftrightarrow$ &
\raisebox{-0.4\height}{
\begin{texdraw}
\fontsize{7}{7}\selectfont
\drawdim em
\textref h:C v:C
\setunitscale 1.9
\move(0 0.5)\lvec(0 0.633)\lvec(1 0.633)\lvec(1 0.5)\ifill f:0.6
\move(0 0.5)\lvec(0 0.633)\lvec(1 0.633)\lvec(1 0.5)
\move(0 0)\lvec(1 0)\lvec(1 0.5)\lvec(0 0.5)\lvec(0 0)
\htext(0.5 0.25){$i$}
\end{texdraw}
}
\end{tabular}
\quad
\begin{tabular}{rcl}
\raisebox{-0.4\height}{
\begin{texdraw}
\fontsize{7}{7}\selectfont
\drawdim em
\textref h:C v:C
\setunitscale 1.9
\move(-1 0)\lvec(0 0)\lvec(0 1)\lvec(-1 1)\lvec(-1 0)
\move(0 0)\lvec(0.25 0.25)\lvec(0.25 1.25)\lvec(-0.75 1.25)\lvec(-1 1)
\move(0 1)\lvec(0.25 1.25)
\lpatt(0.03 0.1)
\move(0 0)\lvec(-0.25 -0.25)\lvec(-1.25 -0.25)\lvec(-1 0)
\htext(-0.5 0.5){$j$}
\end{texdraw}
}
& $\longleftrightarrow$ &
\raisebox{-0.4\height}{
\begin{texdraw}
\fontsize{7}{7}\selectfont
\drawdim em
\textref h:C v:C
\setunitscale 1.9
\move(-1 1)\lvec(-1 1.133)\lvec(0 1.133)\lvec(0 1)
\move(-1 0)\lvec(0 0)\lvec(0 1)\lvec(-1 1)\lvec(-1 0)
\move(0 1)\lvec(-1 0)
\htext(-0.25 0.25){$j$}
\end{texdraw}
}\\[4mm]
\raisebox{-0.4\height}{
\begin{texdraw}
\fontsize{7}{7}\selectfont
\drawdim em
\textref h:C v:C
\setunitscale 1.9
\move(-1 0)\lvec(0 0)\lvec(0 1)\lvec(-1 1)\lvec(-1 0)
\move(0 0)\lvec(0.25 0.25)\lvec(0.25 1.25)\lvec(-0.75 1.25)\lvec(-1 1)
\move(0 1)\lvec(0.25 1.25)
\lpatt(0.03 0.1)
\move(0.25 0.25)\lvec(0.5 0.5)\lvec(0.25 0.5)
\htext(-0.5 0.5){$i$}
\end{texdraw}
}
& $\longleftrightarrow$ &
\raisebox{-0.4\height}{
\begin{texdraw}
\fontsize{7}{7}\selectfont
\drawdim em
\textref h:C v:C
\setunitscale 1.9
\move(-1 1)\lvec(-1 1.133)\lvec(0 1.133)\lvec(0 1)
\move(-1 0)\lvec(0 0)\lvec(0 1)\lvec(-1 1)\lvec(-1 0)
\move(0 1)\lvec(-1 0)
\htext(-0.75 0.75){$i$}
\end{texdraw}
}
\end{tabular}
\quad
\begin{tabular}{rcl}
\raisebox{-0.4\height}{
\begin{texdraw}
\fontsize{7}{7}\selectfont
\drawdim em
\textref h:C v:C
\setunitscale 1.9
\move(-1 0)\lvec(0 0)\lvec(0 1)\lvec(-1 1)\lvec(-1 0)
\move(0 0)\lvec(0.5 0.5)\lvec(0.5 1.5)\lvec(-0.5 1.5)\lvec(-1 1)
\move(0.25 0.25)\lvec(0.25 1.25)\lvec(-0.75 1.25)
\move(0 1.0)\lvec(0.5 1.5)
\htext(-0.5 0.5){$i$}
\end{texdraw}
}
& $\longleftrightarrow$ &
\raisebox{-0.4\height}{
\begin{texdraw}
\fontsize{7}{7}\selectfont
\drawdim em
\textref h:C v:C
\setunitscale 1.9
\move(-1 1)\lvec(-1 1.133)\lvec(0 1.133)\lvec(0 1)\ifill f:0.6
\move(-1 1)\lvec(-1 1.133)\lvec(0 1.133)\lvec(0 1)
\move(-1 0)\lvec(0 0)\lvec(0 1)\lvec(-1 1)\lvec(-1 0)
\move(0 1)\lvec(-1 0)
\htext(-0.75 0.75){$i$}
\htext(-0.25 0.25){$j$}
\end{texdraw}
}
\end{tabular}
\end{center}
\vspace{3mm}
The thin rectangle at the top of each notation is a reminder that
the block is stacked in a wall of unit depth.
The dark shading shows that it is of full unit depth.
Unshaded ones show that the unit depth has not been filled completely.
The following example is for a set of blocks stacked
in a wall of unit thickness.
\vspace{3mm}
\savebox{\tmpfiga}{\begin{texdraw}
\fontsize{7}{7}\selectfont
\textref h:C v:C
\drawdim em
\setunitscale 1.9
\move(0 0)
\bsegment
\move(0 0)\lvec(0 1)\lvec(1 1)\lvec(1 0)\lvec(0 0)
\move(0 1)\lvec(0.2 1.15)\lvec(1 1.15)
\esegment
\move(1 0)
\bsegment
\move(0 0)\lvec(0 1)\lvec(1 1)\lvec(1 0)\lvec(0 0)
\esegment
\move(2 0)
\bsegment
\move(0 0)\lvec(0 1)\lvec(1 1)\lvec(1 0)\lvec(0 0)
\esegment
\move(3 0)
\bsegment
\move(0 0)\lvec(0 1)\lvec(1 1)\lvec(1 0)\lvec(0 0)
\esegment
\move(4 0)
\bsegment
\move(0 0)\lvec(0 1)\lvec(1 1)\lvec(1 0)\lvec(0 0)
\move(0 1)\lvec(0.2 1.15)\lvec(1.2 1.15)
\move(0.2 1.15)\lvec(0.2 2)
\move(1 1)\lvec(1.4 1.3)
\move(1.2 1.15)\lvec(1.2 2.15)
\esegment
\move(5 0)
\bsegment
\move(0 0)\lvec(0 0.5)\lvec(1 0.5)\lvec(1 0)\lvec(0 0)
\move(1 0.5)\lvec(1.4 0.8)
\move(0 0.5)\lvec(0.4 0.8)
\move(1 0)\lvec(1.4 0.3)
\esegment
\move(1 1)
\bsegment
\move(0 0)\lvec(0 0.5)\lvec(1 0.5)\lvec(1 0)\lvec(0 0)
\esegment
\move(2 1)
\bsegment
\move(0 0)\lvec(0 1)\lvec(1 1)\lvec(1 0)\lvec(0 0)
\move(0 1)\lvec(0.2 1.15)\lvec(1 1.15)
\esegment
\move(3 1)
\bsegment
\move(0 0)\lvec(0 0.5)\lvec(1 0.5)\lvec(1 0)\lvec(0 0)
\esegment
\move(3 1.5)
\bsegment
\move(0 0)\lvec(0 0.5)\lvec(1 0.5)\lvec(1 0)\lvec(0 0)
\esegment
\move(1 1.5)
\bsegment
\move(0 0)\lvec(0 1)\lvec(1 1)\lvec(1 0)\lvec(0 0)
\move(0 1)\lvec(0.4 1.3)
\move(1 1)\lvec(1.4 1.3)
\move(0.4 1.3)\lvec(1.4 1.3)\lvec(1.4 0.65)
\esegment
\move(3 2)
\bsegment
\move(0 0)\lvec(0 0.5)\lvec(1 0.5)\lvec(1 0)\lvec(0 0)
\esegment
\move(3 2.5)
\bsegment
\move(0 0)\lvec(0 1)\lvec(1 1)\lvec(1 0)\lvec(0 0)
\esegment
\move(3 3.5)
\bsegment
\move(0 0)\lvec(0 0.5)\lvec(1 0.5)\lvec(1 0)\lvec(0 0)
\move(0 0.5)\lvec(0.4 0.8)
\esegment
\move(4 2)
\bsegment
\move(0 0)\lvec(0 1)\lvec(1 1)\lvec(1 0)\lvec(0 0)
\move(1 1)\lvec(1.4 1.3)
\move(1 0)\lvec(1.4 0.3)
\move(1.2 1.15)\lvec(1.2 2.15)\lvec(0.2 2.15)
\esegment
\move(4 3)
\bsegment
\move(0 0)\lvec(0 1)\lvec(1 1)\lvec(1 0)\lvec(0 0)
\move(1 1)\lvec(1.4 1.3)
\move(0 1)\lvec(0.4 1.3)
\esegment
\move(6.4 0.3)\lvec(6.4 0.8)\lvec(5.4 0.8)\lvec(5.4 4.3)
\lvec(3.4 4.3)
\htext(0.5 0.5){$1$}
\htext(1.5 0.5){$2$}
\htext(1.5 1.25){$0$}
\htext(1.5 2){$2$}
\htext(2.5 0.5){$1$}
\htext(2.5 1.5){$0$}
\htext(3.5 0.5){$3$}
\htext(3.5 1.25){$0$}
\htext(3.5 1.75){$0$}
\htext(3.5 2.25){$0$}
\htext(3.5 3){$2$}
\htext(3.5 3.75){$4$}
\htext(4.5 0.5){$2$}
\htext(4.7 1.65){$1$}
\htext(4.5 2.5){$3$}
\htext(4.5 3.5){$1$}
\htext(5.5 0.25){$0$}
\htext(5.3 3.7){$0$}
\end{texdraw}}%
\savebox{\tmpfigb}{\begin{texdraw}
\fontsize{7}{7}\selectfont
\textref h:C v:C
\drawdim em
\setunitscale 1.9
\move(0 0)
\bsegment
\move(0 0)\lvec(0 1)\lvec(1 1)\lvec(1 0)\lvec(0 0)
\move(1 1)\lvec(0 0)
\esegment
\move(1 0)
\bsegment
\move(0 0)\lvec(0 1)\lvec(1 1)\lvec(1 0)\lvec(0 0)
\esegment
\move(2 0)
\bsegment
\move(0 0)\lvec(0 1)\lvec(1 1)\lvec(1 0)\lvec(0 0)
\esegment
\move(3 0)
\bsegment
\move(0 0)\lvec(0 1)\lvec(1 1)\lvec(1 0)\lvec(0 0)
\esegment
\move(4 0)
\bsegment
\move(0 0)\lvec(0 1)\lvec(1 1)\lvec(1 0)\lvec(0 0)
\esegment
\move(5 0)
\bsegment
\move(0 0)\lvec(0 0.5)\lvec(1 0.5)\lvec(1 0)\lvec(0 0)
\esegment
\move(1 1)
\bsegment
\move(0 0)\lvec(0 0.5)\lvec(1 0.5)\lvec(1 0)\lvec(0 0)
\esegment
\move(2 1)
\bsegment
\move(0 0)\lvec(0 1)\lvec(1 1)\lvec(1 0)\lvec(0 0)
\move(1 1)\lvec(0 0)
\esegment
\move(3 1)
\bsegment
\move(0 0)\lvec(0 0.5)\lvec(1 0.5)\lvec(1 0)\lvec(0 0)
\esegment
\move(3 1.5)
\bsegment
\move(0 0)\lvec(0 0.5)\lvec(1 0.5)\lvec(1 0)\lvec(0 0)
\esegment
\move(4 1)
\bsegment
\move(0 0)\lvec(0 1)\lvec(1 1)\lvec(1 0)\lvec(0 0)
\move(1 1)\lvec(0 0)
\esegment
\move(1 1.5)
\bsegment
\move(0 0)\lvec(0 1)\lvec(1 1)\lvec(1 0)\lvec(0 0)
\esegment
\move(3 2)
\bsegment
\move(0 0)\lvec(0 0.5)\lvec(1 0.5)\lvec(1 0)\lvec(0 0)
\esegment
\move(3 2.5)
\bsegment
\move(0 0)\lvec(0 1)\lvec(1 1)\lvec(1 0)\lvec(0 0)
\esegment
\move(3 3.5)
\bsegment
\move(0 0)\lvec(0 0.5)\lvec(1 0.5)\lvec(1 0)\lvec(0 0)
\esegment
\move(4 2)
\bsegment
\move(0 0)\lvec(0 1)\lvec(1 1)\lvec(1 0)\lvec(0 0)
\esegment
\move(4 3)
\bsegment
\move(0 0)\lvec(0 1)\lvec(1 1)\lvec(1 0)\lvec(0 0)
\move(1 1)\lvec(0 0)
\esegment
\move(0 1)
\bsegment
\move(0 0)\lvec(0 0.133)\lvec(1 0.133)\lvec(1 0)\lvec(0 0)
\esegment
\move(1 2.5)
\bsegment
\move(0 0)\lvec(0 0.133)\lvec(1 0.133)\lvec(1 0)\ifill f:0.6
\move(0 0)\lvec(0 0.133)\lvec(1 0.133)\lvec(1 0)\lvec(0 0)
\esegment
\move(2 2)
\bsegment
\move(0 0)\lvec(0 0.133)\lvec(1 0.133)\lvec(1 0)\lvec(0 0)
\esegment
\move(3 4)
\bsegment
\move(0 0)\lvec(0 0.133)\lvec(1 0.133)\lvec(1 0)\ifill f:0.6
\move(0 0)\lvec(0 0.133)\lvec(1 0.133)\lvec(1 0)\lvec(0 0)
\esegment
\move(4 4)
\bsegment
\move(0 0)\lvec(0 0.133)\lvec(1 0.133)\lvec(1 0)\ifill f:0.6
\move(0 0)\lvec(0 0.133)\lvec(1 0.133)\lvec(1 0)\lvec(0 0)
\esegment
\move(5 0.5)
\bsegment
\move(0 0)\lvec(0 0.133)\lvec(1 0.133)\lvec(1 0)\ifill f:0.6
\move(0 0)\lvec(0 0.133)\lvec(1 0.133)\lvec(1 0)\lvec(0 0)
\esegment
\htext(0.25 0.75){$1$}
\htext(1.5 0.5){$2$}
\htext(1.5 1.25){$0$}
\htext(1.5 2){$2$}
\htext(2.5 0.5){$1$}
\htext(2.25 1.75){$0$}
\htext(3.5 0.5){$3$}
\htext(3.5 1.25){$0$}
\htext(3.5 1.75){$0$}
\htext(3.5 2.25){$0$}
\htext(3.5 3){$2$}
\htext(3.5 3.75){$4$}
\htext(4.5 0.5){$2$}
\htext(4.75 1.25){$1$}
\htext(4.5 2.5){$3$}
\htext(4.25 3.75){$1$}
\htext(4.75 3.25){$0$}
\htext(5.5 0.25){$0$}
\end{texdraw}}%
\begin{center}
\raisebox{-0.1em}{\usebox{\tmpfiga}}
\quad\raisebox{0.7em}{$\longleftrightarrow$}\quad\;
\usebox{\tmpfigb}
\end{center}
\vspace{4mm}

\begin{df}
If $\g \neq \cnone$, we define a \defi{level-$1$ slice} to be a set
of finitely many blocks of given type stacked in one column of unit depth
following the pattern given below.
In stacking the blocks, no block should be placed
on top of a column of half-unit depth.
If $\g = \cnone$, such a set of blocks will be called a
\defi{level-$\half$ slice}.
\end{df}

\vspace{3mm}

\savebox{\tmpfiga}{
\begin{texdraw}
\fontsize{7}{7}\selectfont
\drawdim em
\textref h:C v:C
\setunitscale 1.9
\move(-1 0)\lvec(0 0)\lvec(0 1)\lvec(-1 1)\lvec(-1 0)
\move(0 0)\lvec(0.5 0.5)\lvec(0.5 1.5)\lvec(-0.5 1.5)\lvec(-1 1)
\move(0 1.0)\lvec(0.5 1.5)
\htext(-0.5 0.5){$i$}
\end{texdraw}}
\savebox{\tmpfigb}{
\begin{texdraw}
\fontsize{7}{7}\selectfont
\drawdim em
\textref h:C v:C
\setunitscale 1.9
\move(-1 0)\lvec(0 0)\lvec(0 0.5)\lvec(-1 0.5)\lvec(-1 0)
\move(0 0)\lvec(0.5 0.5)\lvec(0.5 1)\lvec(-0.5 1)\lvec(-1 0.5)
\move(0 0.5)\lvec(0.5 1)
\htext(-0.5 0.25){$i$}
\end{texdraw}}%
\savebox{\tmpfigc}{
\begin{texdraw}
\fontsize{7}{7}\selectfont
\drawdim em
\textref h:C v:C
\setunitscale 1.9
\move(-1 0)\lvec(0 0)\lvec(0 1)\lvec(-1 1)\lvec(-1 0)
\move(0 0)\lvec(0.25 0.25)\lvec(0.25 1.25)\lvec(-0.75 1.25)\lvec(-1 1)
\move(0 1)\lvec(0.25 1.25)
\htext(-0.5 0.5){$i$}
\end{texdraw}}
\savebox{\tmpfigd}{
\begin{texdraw}
\fontsize{7}{7}\selectfont
\textref h:C v:C
\drawdim em
\setunitscale 1.9
\move(0 0)
\bsegment
\move(0 0)\rlvec(1 0)\move(0 1)\rlvec(1 0)
\move(0 2)\rlvec(1 0)\move(0 3.4)\rlvec(1 0)
\move(0 4.4)\rlvec(1 0)\move(0 5.4)\rlvec(1 0)
\move(0 6.4)\rlvec(1 0)\move(0 0)\rlvec(0 6.7)
\move(1 0)\rlvec(0 6.7)
\htext(0.5 0.5){$0$}\htext(0.5 1.5){$1$}
\vtext(0.5 2.7){$\cdots$}\htext(0.5 3.9){$n\!\!-\!\!1$}
\htext(0.5 4.9){$n$}\htext(0.5 5.9){$0$}
\esegment
\move(2.5 0)
\bsegment
\move(0 0)\rlvec(1 0)\move(0 1)\rlvec(1 0)
\move(0 2)\rlvec(1 0)\move(0 3.4)\rlvec(1 0)
\move(0 4.4)\rlvec(1 0)\move(0 5.4)\rlvec(1 0)
\move(0 6.4)\rlvec(1 0)\move(0 0)\rlvec(0 6.7)
\move(1 0)\rlvec(0 6.7)
\htext(0.5 0.5){$1$}\htext(0.5 1.5){$2$}
\vtext(0.5 2.7){$\cdots$}\htext(0.5 3.9){$n$}
\htext(0.5 4.9){$0$}\htext(0.5 5.9){$1$}
\esegment
\move(5 0)
\bsegment
\htext(0.5 0.5){$\cdots$}
\esegment
\move(7.5 0)
\bsegment
\move(0 0)\rlvec(1 0)\move(0 1)\rlvec(1 0)
\move(0 2)\rlvec(1 0)\move(0 3)\rlvec(1 0)
\move(0 4.4)\rlvec(1 0)\move(0 5.4)\rlvec(1 0)
\move(0 6.4)\rlvec(1 0)\move(0 0)\rlvec(0 6.7)
\move(1 0)\rlvec(0 6.7)
\htext(0.5 0.5){$n$}\htext(0.5 1.5){$0$}
\htext(0.5 2.5){$1$}\vtext(0.5 3.7){$\cdots$}
\htext(0.5 4.9){$n\!\!-\!\!1$}\htext(0.5 5.9){$n$}
\esegment
\end{texdraw}}
\savebox{\tmpfige}{
\begin{texdraw}
\fontsize{7}{7}\selectfont
\textref h:C v:C
\drawdim em
\setunitscale 1.9
\move(0 0)\rlvec(-1 0) \move(0 0.5)\rlvec(-1 0)
\move(0 1)\rlvec(-1 0) \move(0 2)\rlvec(-1 0)
\move(0 3.4)\rlvec(-1 0) \move(0 4.4)\rlvec(-1 0)
\move(0 5.8)\rlvec(-1 0) \move(0 6.8)\rlvec(-1 0)
\move(0 7.3)\rlvec(-1 0) \move(0 7.8)\rlvec(-1 0)
\move(0 8.8)\rlvec(-1 0) \move(0 0)\rlvec(0 9.1)
\move(-1 0)\rlvec(0 9.1)
\htext(-0.5 0.25){$0$} \htext(-0.5 0.75){$0$}
\htext(-0.5 1.5){$1$} \vtext(-0.5 2.8){$\cdots$}
\htext(-0.5 3.9){$n$} \vtext(-0.5 5.2){$\cdots$}
\htext(-0.5 6.3){$1$} \htext(-0.5 7.05){$0$}
\htext(-0.5 7.55){$0$} \htext(-0.5 8.3){$1$}
\textref h:L v:C
\htext(0.2 2.5){$\left.\rule{0pt}{4em}\right\}$ \defi{supporting} blocks}
\textref h:R v:C
\htext(-1.1 5.35){\defi{covering} blocks $\left\{\rule{0pt}{4em}\right.$}
\htext(-1.1 0.25){covering block $\rightarrow$}
\end{texdraw}}
\savebox{\tmpfigf}{
\begin{texdraw}
\fontsize{7}{7}\selectfont
\textref h:C v:C
\drawdim em
\setunitscale 1.9
\move(-1 0)
\bsegment
\move(0 0)\rlvec(1 0) \move(0 0.5)\rlvec(1 0)
\move(0 1)\rlvec(1 0) \move(0 2)\rlvec(1 0)
\move(0 3.4)\rlvec(1 0) \move(0 4.9)\rlvec(1 0)
\move(0 4.4)\rlvec(1 0) \move(0 5.4)\rlvec(1 0)
\move(0 6.4)\rlvec(1 0)
\move(0 8.8)\rlvec(1 0) \move(0 7.8)\rlvec(1 0)
\move(0 9.3)\rlvec(1 0) \move(0 9.8)\rlvec(1 0)
\move(0 10.8)\rlvec(1 0)
\move(0 0)\lvec(0 11.1) \move(1 0)\lvec(1 11.1)
\htext(0.5 0.25){$0$} \htext(0.5 0.75){$0$}
\htext(0.5 1.5){$1$} \vtext(0.5 2.8){$\cdots$}
\htext(0.5 3.9){$\!\!n\!\!-\!\!1\!\!$}
\htext(0.5 4.65){$n$} \htext(0.5 5.15){$n$}
\htext(0.5 5.9){$\!\!n\!\!-\!\!1\!\!$}
\vtext(0.5 7.2){$\cdots$}
\htext(0.5 8.3){$1$}
\htext(0.5 9.05){$0$} \htext(0.5 9.55){$0$}
\htext(0.5 10.3){$1$}
\textref h:L v:C
\htext(1.2 2.8){$\left.\rule{0pt}{4.3em}\right\}$ \defi{supporting} blocks}
\textref h:R v:C
\htext(-0.1 7.2){\defi{covering} blocks $\left\{\rule{0pt}{4.3em}\right.$}
\htext(-0.1 0.25){covering block $\rightarrow$}
\esegment
\move(1 0)
\end{texdraw}}
\savebox{\tmpfigg}{
\begin{texdraw}
\fontsize{7}{7}\selectfont \textref h:C v:C \drawdim em
\setunitscale 1.9 \move(0 0) \bsegment \move(0 0)\rlvec(1 0)
\move(0 0)\rlvec(1 1) \move(0 1)\rlvec(1 0) \move(0 2)\rlvec(1 0)
\move(0 3.4)\rlvec(1 0) \move(0 4.9)\rlvec(1 0) \move(0
4.4)\rlvec(1 0) \move(0 5.4)\rlvec(1 0) \move(0 6.4)\rlvec(1 0)
\move(0 8.8)\rlvec(1 0) \move(0 7.8)\rlvec(1 0) \move(0
8.8)\rlvec(1 1) \move(0 9.8)\rlvec(1 0) \move(0 10.8)\rlvec(1 0)
\move(0 0)\rlvec(0 11.1) \move(1 0)\rlvec(0 11.1) \htext(0.75
0.25){$0$} \htext(0.25 0.75){$1$} \htext(0.5 1.5){$2$} \vtext(0.5
2.8){$\cdots$} \htext(0.5 3.9){$\!\!n\!\!-\!\!1\!\!$} \htext(0.5
4.65){$n$} \htext(0.5 5.15){$n$} \htext(0.5
5.9){$\!\!n\!\!-\!\!1\!\!$} \vtext(0.5 7.2){$\cdots$} \htext(0.5
8.3){$2$} \htext(0.75 9.05){$0$} \htext(0.25 9.55){$1$} \htext(0.5
10.3){$2$}
\textref h:L v:C
\htext(1.1 2.55){$\left.\rule{0pt}{5em}\right\}$ \defi{supporting} blocks}
\textref h:R v:C
\htext(-0.1 7.45){\defi{covering} blocks $\left\{\rule{0pt}{5em}\right.$}
\htext(-0.1 0.5){covering blocks $\rightarrow$}
\esegment
\move(8 0)
\bsegment \move(0 0)\rlvec(1
0) \move(0 0)\rlvec(1 1) \move(0 1)\rlvec(1 0) \move(0 2)\rlvec(1
0) \move(0 3.4)\rlvec(1 0) \move(0 4.9)\rlvec(1 0) \move(0
4.4)\rlvec(1 0) \move(0 5.4)\rlvec(1 0) \move(0 6.4)\rlvec(1 0)
\move(0 8.8)\rlvec(1 0) \move(0 7.8)\rlvec(1 0) \move(0
8.8)\rlvec(1 1) \move(0 9.8)\rlvec(1 0) \move(0 10.8)\rlvec(1 0)
\move(0 0)\rlvec(0 11.1) \move(1 0)\rlvec(0 11.1) \htext(0.75
0.25){$1$} \htext(0.25 0.75){$0$} \htext(0.5 1.5){$2$} \vtext(0.5
2.8){$\cdots$} \htext(0.5 3.9){$\!\!n\!\!-\!\!1\!\!$} \htext(0.5
4.65){$n$} \htext(0.5 5.15){$n$} \htext(0.5
5.9){$\!\!n\!\!-\!\!1\!\!$} \vtext(0.5 7.2){$\cdots$} \htext(0.5
8.3){$2$} \htext(0.75 9.05){$1$} \htext(0.25 9.55){$0$} \htext(0.5
10.3){$2$}
\textref h:L v:C
\htext(1.1 2.55){$\left.\rule{0pt}{5.0em}\right\}$ \defi{supporting} blocks}
\textref h:R v:C
\htext(-0.1 7.45){\defi{covering} blocks $\left\{\rule{0pt}{5.0em}\right.$}
\htext(-0.1 0.5){covering blocks $\rightarrow$}
\esegment
\end{texdraw}}
\begin{enumerate}
\item $\anone$$(n\ge 1)$\\[3mm]
      \mbox{}\hspace{4mm} \raisebox{-0.5em}{\usebox{\tmpfiga}}
      \ $(i=0,\cdots,n)$\\[3mm]
      \mbox{}\hspace{4mm} \usebox{\tmpfigd}\\[2mm]
\item $\bnone$$(n\ge 3)$ and $\atwonmonetwo$$(n\ge 3)$\\[3mm]
      \mbox{}\hspace{4mm} \raisebox{-0.1em}{\usebox{\tmpfigb}}
      \ $(i=n)$
      \quad
      \raisebox{-0.5em}{\usebox{\tmpfigc}}
      \ $(i=0,1)$
      \quad
      \raisebox{-0.5em}{\usebox{\tmpfiga}}
      \ $(i=2,\cdots,n-1)$\\[3mm]
      \mbox{}\hspace{26mm} \usebox{\tmpfigg}\\[2mm]
\item $\cnone$$(n\ge 2)$ and $\atwontwo$$(n\ge 1)$\\[3mm]
      \mbox{}\hspace{4mm} \raisebox{-0.1em}{\usebox{\tmpfigb}}
      \ $(i=0)$
      \qquad
      \raisebox{-0.5em}{\usebox{\tmpfiga}}
      \ $(i=1,\cdots,n)$\\[3mm]
      \mbox{}\hspace{26mm} \usebox{\tmpfige}\\[2mm]
\item $\dnponetwo$$(n\ge 2)$\\[3mm]
      \mbox{}\hspace{4mm} \raisebox{-0.1em}{\usebox{\tmpfigb}}
      \ $(i=0,n)$
      \quad
      \raisebox{-0.5em}{\usebox{\tmpfiga}}
      \ $(i=1,\cdots,n-1)$\\[3mm]
      \mbox{}\hspace{26mm} \usebox{\tmpfigf}\\[2mm]
\end{enumerate}

As we can see in the figure, the blocks are stacked in a repeating
pattern, and, roughly speaking,
this pattern is symmetric with respect to the
$n$-block except for the $\anone$ type. We say that an $i$-block is a
\defi{covering $i$-block} (resp. \defi{supporting $i$-block}) if it is
closer to the $n$-block that sits below (resp. above) it than to
the $n$-block that sits above (resp. below) it. An $i$-block that
appears only once in each cycle is regarded as both a supporting
block and a covering block. An \defi{$i$-slot} is the top of a
level-1 (or level-$\half$) slice where one may add an $i$-block.
The notion of \defi{covering $i$-slot} or \defi{supporting
$i$-slot} is self-explanatory.

We define a \defi{$\delta$-column} to be a set of blocks appearing
in a cycle of the stacking pattern. For a level-1 (or
level-$\half$) slice $C$, we define $C + \delta$ (resp. $C -
\delta$) to be the level-1 (or level-$\half$) slice obtained from
$C$ by adding (resp. removing) a $\delta$-column.

\begin{example}
If $\g = B_3^{(1)}$, we have \\
\savebox{\tmpfiga}{
\begin{texdraw}
\fontsize{7}{7}\selectfont
\textref h:C v:C
\drawdim em
\setunitscale 1.9
\move(1 2)\lvec(1 2.133)\lvec(0 2.133)\lvec(0 2)\ifill f:0.6
\move(0 0)\lvec(1 0)\lvec(1 2)\lvec(0 2)\lvec(0 0)
\move(0 0)\lvec(1 1)
\move(0 1)\lvec(1 1)
\move(1 2)\lvec(1 2.133)\lvec(0 2.133)\lvec(0 2)
\htext(0.75 0.25){$0$}\htext(0.25 0.75){$1$}
\htext(0.5 1.5){$2$}
\end{texdraw}
\raisebox{0.5em}{$\;+\;\delta\ =\ $}
\begin{texdraw}
\fontsize{7}{7}\selectfont
\textref h:C v:C
\drawdim em
\setunitscale 1.9
\move(1 6)\lvec(1 6.133)\lvec(0 6.133)\lvec(0 6)\ifill f:0.6
\move(0 0)\lvec(1 0)\lvec(1 6)\lvec(0 6)\lvec(0 0)
\move(0 0)\lvec(1 1)
\move(0 1)\lvec(1 1)
\move(0 2)\lvec(1 2)
\move(0 2.5)\lvec(1 2.5)
\move(0 3)\lvec(1 3)
\move(0 4)\lvec(1 4)
\move(0 4)\lvec(1 5)
\move(0 5)\lvec(1 5)
\move(1 6)\lvec(1 6.133)\lvec(0 6.133)\lvec(0 6)
\htext(0.75 0.25){$0$}\htext(0.25 0.75){$1$}
\htext(0.5 1.5){$2$}\htext(0.5 2.25){$3$}
\htext(0.5 2.75){$3$}\htext(0.5 3.5){$2$}
\htext(0.75 4.25){$0$}\htext(0.25 4.75){$1$}
\htext(0.5 5.5){$2$}
\end{texdraw}}
\savebox{\tmpfigb}{
\begin{texdraw}
\fontsize{7}{7}\selectfont
\textref h:C v:C
\drawdim em
\setunitscale 1.9
\move(0 6.133)
\move(0 0)\lvec(1 0)\lvec(1 5)\lvec(0 5)\lvec(0 0)
\move(0 0)\lvec(1 1)
\move(0 1)\lvec(1 1)
\move(0 2)\lvec(1 2)
\move(0 2.5)\lvec(1 2.5)
\move(0 3)\lvec(1 3)
\move(0 4)\lvec(1 4)
\move(0 4)\lvec(1 5)
\move(1 5)\lvec(1 5.133)\lvec(0 5.133)\lvec(0 5)
\htext(0.75 0.25){$1$}\htext(0.25 0.75){$0$}
\htext(0.5 1.5){$2$}\htext(0.5 2.25){$3$}
\htext(0.5 2.75){$3$}\htext(0.5 3.5){$2$}
\htext(0.25 4.75){$0$}
\end{texdraw}
\raisebox{0.5em}{$\;-\;\delta\ =\ $}
\begin{texdraw}
\fontsize{7}{7}\selectfont
\textref h:C v:C
\drawdim em
\setunitscale 1.9
\move(0 0)\lvec(1 0)\lvec(1 1)\lvec(0 1)\lvec(0 0)
\move(0 0)\lvec(1 1)
\move(1 1)\lvec(1 1.133)\lvec(0 1.133)\lvec(0 1)
\htext(0.25 0.75){$0$}
\end{texdraw}}
\begin{center}
\begin{minipage}{0.3\textwidth}
\usebox{\tmpfiga}
\end{minipage}
\begin{minipage}{0.3\textwidth}
\usebox{\tmpfigb}
\end{minipage}
\end{center}
\end{example}
\vspace{1mm}

\begin{df}\hfill
\begin{enumerate}
\item
A \defi{level-$l$ slice} for $\g \neq \cnone$ is an ordered
$l$-tuple $C = (c_1, \dots, c_l)$ of level-$1$ slices satisfying
the following conditions\, :
\begin{enumerate}
\item
$\anone$, $\bnone$, and $\atwontwo$\\
$\bullet$ $c_1 \subset c_2 \subset \cdots \subset c_l \subset c_1
+ \delta$;
\item
$\atwonmonetwo$\\
$\bullet$ $c_1 \subset c_2 \subset \cdots \subset c_l \subset c_1 + \delta$,\\
$\bullet$ it contains an even number of $n$-blocks;
\item
$\dnponetwo$\\
$\bullet$ $c_1 \subset c_2 \subset \cdots \subset c_l \subset c_1 + \delta$,\\
$\bullet$ at most one of the top blocks is a supporting $n$-block.
\end{enumerate}
\item
A \defi{level-$l$ slice} for $\g=\cnone$ is an ordered $2l$-tuple
$C = (c_1, \dots, c_{2l})$ of level-$\half$ slices such that \\
$\bullet$ $c_1 \subset c_2 \subset \cdots \subset c_{2l} \subset c_1 + \delta$,\\
$\bullet$ it contains an even number of $0$-blocks.
\end{enumerate}
For affine types that allow more than one stacking pattern,
only one should be used in a level-$l$ slice.
Each level-$1$ (or level-$\half$) slice $c_i$ in $C$ is called the
\defi{$i$-th layer} of $C$.
For each affine type,
the set of all level-$l$ slices, that uses the same stacking pattern,
is denoted by $\slice$.
\end{df}

\begin{remark}
According to the above definition,
for affine types that allow more than one stacking pattern,
more than one set is being denoted by $\slice$.
But the contents of this section shall show that this will cause no confusion.
\end{remark}

We will often just say \emph{slice} for level-$l$ slice. A
level-$l$ slice can be viewed as the set of $l$ columns with the
$i$-th layer placed in front of the $(i+1)$-th layer. For
simplicity, we will often use the front-and-top view when
representing a slice. We will also use side views whenever it is
necessary.  We explain the two methods used in drawing a slice
with the following example for $B_3^{(1)}$-type.
\vspace{3mm}
\begin{center}
\raisebox{0.5em}{$(\ c_1\;=\ $}
\begin{texdraw}
\fontsize{7}{7}\selectfont
\textref h:C v:C
\drawdim em
\setunitscale 1.9
\move(1 2)\lvec(1 2.133)\lvec(0 2.133)\lvec(0 2)\ifill f:0.6
\move(0 0)\lvec(1 0)\lvec(1 2)\lvec(0 2)\lvec(0 0)
\move(0 0)\lvec(1 1)
\move(0 1)\lvec(1 1)
\move(1 2)\lvec(1 2.133)\lvec(0 2.133)\lvec(0 2)
\htext(0.75 0.25){$0$} \htext(0.25 0.75){$1$}
\htext(0.5 1.5){$2$}
\end{texdraw}\;,
\raisebox{0.5em}{$c_2\;=\ $}
\begin{texdraw}
\fontsize{7}{7}\selectfont
\textref h:C v:C
\drawdim em
\setunitscale 1.9
\move(1 2.5)\lvec(1 2.633)\lvec(0 2.633)\lvec(0 2.5)\ifill f:0.6
\move(0 0)\lvec(1 0)\lvec(1 2.5)\lvec(0 2.5)\lvec(0 0)
\move(0 1)\lvec(1 1) \move(0 2)\lvec(1 2)
\move(0 0)\lvec(1 1) \move(0 2.5)\lvec(1 2.5)
\move(1 2.5)\lvec(1 2.633)\lvec(0 2.633)\lvec(0 2.5)
\htext(0.75 0.25){$0$} \htext(0.25 0.75){$1$}
\htext(0.5 1.5){$2$} \htext(0.5 2.25){$3$}
\end{texdraw}\;,
\raisebox{0.5em}{$c_3\;=\ $}
\begin{texdraw}
\fontsize{7}{7}\selectfont
\textref h:C v:C
\drawdim em
\setunitscale 1.9
\move(1 2.5)\lvec(1 2.633)\lvec(0 2.633)\lvec(0 2.5)\ifill f:0.6
\move(0 0)\lvec(1 0)\lvec(1 2.5)\lvec(0 2.5)\lvec(0 0)
\move(0 1)\lvec(1 1) \move(0 2)\lvec(1 2)
\move(0 0)\lvec(1 1) \move(0 2.5)\lvec(1 2.5)
\move(1 2.5)\lvec(1 2.633)\lvec(0 2.633)\lvec(0 2.5)
\htext(0.75 0.25){$0$} \htext(0.25 0.75){$1$}
\htext(0.5 1.5){$2$} \htext(0.5 2.25){$3$}
\end{texdraw}\;,
\raisebox{0.5em}{$c_4\;=\ $}
\begin{texdraw}
\fontsize{7}{7}\selectfont
\textref h:C v:C
\drawdim em
\setunitscale 1.9
\move(0 0)\lvec(1 0)\lvec(1 5)\lvec(0 5)\lvec(0 0)
\move(0 1)\lvec(1 1) \move(0 2)\lvec(1 2)
\move(0 0)\lvec(1 1) \move(0 2.5)\lvec(1 2.5)
\move(0 3)\lvec(1 3) \move(0 4)\lvec(1 4)
\move(0 4)\lvec(1 5)
\move(1 5)\lvec(1 5.133)\lvec(0 5.133)\lvec(0 5)
\htext(0.75 0.25){$0$} \htext(0.25 0.75){$1$}
\htext(0.5 1.5){$2$} \htext(0.5 2.25){$3$}
\htext(0.5 2.75){$3$} \htext(0.5 3.5){$2$}
\htext(0.75 4.25){$0$}
\end{texdraw}\;,
\raisebox{0.5em}{$c_5\;=\ $}
\begin{texdraw}
\fontsize{7}{7}\selectfont
\textref h:C v:C
\drawdim em
\setunitscale 1.9
\move(0 0)\lvec(1 0)\lvec(1 5)\lvec(0 5)\lvec(0 0)
\move(0 1)\lvec(1 1) \move(0 2)\lvec(1 2)
\move(0 0)\lvec(1 1) \move(0 2.5)\lvec(1 2.5)
\move(0 3)\lvec(1 3) \move(0 4)\lvec(1 4)
\move(0 4)\lvec(1 5)
\move(1 5)\lvec(1 5.133)\lvec(0 5.133)\lvec(0 5)
\htext(0.75 0.25){$0$} \htext(0.25 0.75){$1$}
\htext(0.5 1.5){$2$} \htext(0.5 2.25){$3$}
\htext(0.5 2.75){$3$} \htext(0.5 3.5){$2$}
\htext(0.75 4.25){$0$}
\end{texdraw}\;,
\raisebox{0.5em}{$c_6\;=\ $}
\begin{texdraw}
\fontsize{7}{7}\selectfont
\textref h:C v:C
\drawdim em
\setunitscale 1.9
\move(1 5)\lvec(1 5.133)\lvec(0 5.133)\lvec(0 5)\ifill f:0.6
\move(0 0)\lvec(1 0)\lvec(1 5)\lvec(0 5)\lvec(0 0)
\move(0 1)\lvec(1 1) \move(0 2)\lvec(1 2)
\move(0 0)\lvec(1 1) \move(0 2.5)\lvec(1 2.5)
\move(0 3)\lvec(1 3) \move(0 4)\lvec(1 4)
\move(0 4)\lvec(1 5)
\move(1 5)\lvec(1 5.133)\lvec(0 5.133)\lvec(0 5)
\htext(0.75 0.25){$0$}\htext(0.25 0.75){$1$}
\htext(0.5 1.5){$2$}\htext(0.5 2.25){$3$}
\htext(0.5 2.75){$3$}\htext(0.5 3.5){$2$}
\htext(0.75 4.25){$0$}\htext(0.25 4.75){$1$}
\end{texdraw}
\raisebox{0.5em}{$\ )$}\\[3mm]
\raisebox{0.5em}{$\ =\ $}
\begin{texdraw}
\fontsize{7}{7}\selectfont
\textref h:C v:C
\drawdim em
\setunitscale 1.9
\move(-0.4 -0.3)
\bsegment
\move(0 0)\lvec(1 0)\lvec(1 2)\lvec(0 2)\lvec(0 0)
\move(0 1)\lvec(1 1)
\htext(0.5 0.5){$1$}
\htext(0.5 1.5){$2$}
\esegment
\move(0 0)
\bsegment
\move(1 0)\lvec(1 2.5)\lvec(0 2.5)\lvec(0 2)\lvec(1 2)
\htext(0.5 2.25){$3$}
\esegment
\move(0.4 0.3)
\bsegment
\move(1 0)\lvec(1 2.5)\lvec(0 2.5)
\esegment
\move(0.8 0.6)
\bsegment
\move(1 0)\lvec(1 2.5)\lvec(0 2.5)
\move(0 2.5)\lvec(0 4)\lvec(1 4)\lvec(1 2.5)
\move(0 3)\lvec(1 3)
\htext(0.5 2.75){$3$}
\htext(0.5 3.5){$2$}
\esegment
\move(1.0 0.75)
\bsegment
\move(1 0)\lvec(1 1)
\move(1 4)\lvec(1 5)\lvec(0 5)\lvec(0 4)\lvec(1 4)
\htext(0.5 4.5){$0$}
\esegment
\move(1.2 0.9)
\bsegment
\move(1 0)\lvec(1 5)\lvec(0 5)
\esegment
\move(1.4 1.05)
\bsegment
\move(1 0)\lvec(1 1)
\move(0.8 4)\lvec(1 4)\lvec(1 5)\lvec(0 5)\lvec(0 4.9)
\esegment
\move(1.6 1.2)
\bsegment
\move(1 0)\lvec(1 5)\lvec(0 5)
\esegment
\move(1.8 1.35)
\bsegment
\move(1 0)\lvec(1 1)
\move(1 4)\lvec(1 5)\lvec(0 5)
\esegment
\move(2.0 1.5)
\bsegment
\move(1 0)\lvec(1 5)\lvec(0 5)
\esegment
\move(0.6 -0.3)\lvec(3.0 1.5)
\move(0.6 0.7)\lvec(3.0 2.5)
\move(0.6 1.7)\lvec(3.0 3.5)
\move(-0.4 1.7)\lvec(0 2)
\move(1 2.5)\lvec(3 4)
\move(0 2.5)\lvec(0.8 3.1)
\move(1.8 3.6)\lvec(3 4.5)
\move(1.8 4.6)\lvec(3 5.5)
\move(0.8 4.6)\lvec(1 4.75)
\move(2 5.75)\lvec(2.2 5.9)
\move(1 5.75)\lvec(1.2 5.9)
\move(2.4 6.05)\lvec(3 6.5)
\move(1.4 6.05)\lvec(2 6.5)
\move(0.8 -0.15)\lvec(0.8 0.85)
\move(1.2 0.15)\lvec(1.2 1.15)
\move(1.6 0.45)\lvec(1.6 1.45)
\end{texdraw}
\quad
\raisebox{0.5em}{$\longleftrightarrow$}\quad
\raisebox{0.5em}{$C\;=\ $}
\begin{texdraw}
\fontsize{7}{7}\selectfont
\textref h:C v:C
\drawdim em
\setunitscale 1.9
\move(0 0)
\bsegment
\move(0 0)\lvec(1 0)\lvec(1 5)\lvec(0 5)\lvec(0 0)
\move(0 1)\lvec(1 1) \move(0 2)\lvec(1 2)
\move(0.5 0)\lvec(0.5 1) \move(0 2.5)\lvec(1 2.5)
\move(0 3)\lvec(1 3) \move(0 4)\lvec(1 4)
\move(0.5 4)\lvec(0.5 5)
\htext(0.25 0.5){$0$}\htext(0.75 0.5){$1$}
\htext(0.5 1.5){$2$}\htext(0.5 2.25){$3$}
\htext(0.5 2.75){$3$}\htext(0.5 3.5){$2$}
\htext(0.25 4.5){$0$}\htext(0.75 4.5){$1$}
\esegment
\move(1 0)
\bsegment
\move(0 0)\lvec(1 0)\lvec(1 4)\lvec(0.5 4)\lvec(0.5 5)\lvec(0 5)
\move(0 1)\lvec(1 1)
\move(0 2)\lvec(1 2) \move(0.5 0)\lvec(0.5 1)
\move(0 2.5)\lvec(1 2.5) \move(0 3)\lvec(1 3)
\move(0 4)\lvec(1 4) \move(0.5 4)\lvec(0.5 5)
\htext(0.25 0.5){$0$} \htext(0.75 0.5){$1$}
\htext(0.5 1.5){$2$} \htext(0.5 2.25){$3$}
\htext(0.5 2.75){$3$} \htext(0.5 3.5){$2$}
\htext(0.25 4.5){$0$}
\esegment
\move(2 0)
\bsegment
\move(0 0)\lvec(1 0)\lvec(1 4)\lvec(0.5 4)\lvec(0.5 5)\lvec(0 5)\lvec(0 4)
\move(0 1)\lvec(1 1) \move(0 2)\lvec(1 2)
\move(0.5 0)\lvec(0.5 1) \move(0 2.5)\lvec(1 2.5)
\move(0 3)\lvec(1 3) \move(0 4)\lvec(1 4)
\move(0.5 4)\lvec(0.5 5)
\htext(0.25 0.5){$0$}\htext(0.75 0.5){$1$}
\htext(0.5 1.5){$2$}\htext(0.5 2.25){$3$}
\htext(0.5 2.75){$3$}\htext(0.5 3.5){$2$}
\htext(0.25 4.5){$0$}
\esegment
\move(3 0)
\bsegment
\move(0 0)\lvec(1 0)\lvec(1 2.5)\lvec(0 2.5)\lvec(0 0)
\move(0 1)\lvec(1 1) \move(0 2)\lvec(1 2)
\move(0.5 0)\lvec(0.5 1) \move(0 2.5)\lvec(1 2.5)
\htext(0.25 0.5){$0$} \htext(0.75 0.5){$1$}
\htext(0.5 1.5){$2$}\htext(0.5 2.25){$3$}
\esegment
\move(4 0)
\bsegment
\move(0 0)\lvec(1 0)\lvec(1 2.5)\lvec(0 2.5)\lvec(0 0)
\move(0 1)\lvec(1 1) \move(0 2)\lvec(1 2)
\move(0.5 0)\lvec(0.5 1) \move(0 2.5)\lvec(1 2.5)
\htext(0.25 0.5){$0$} \htext(0.75 0.5){$1$}
\htext(0.5 1.5){$2$}\htext(0.5 2.25){$3$}
\esegment
\move(5 0)
\bsegment
\move(0 0)\lvec(1 0)\lvec(1 2)\lvec(0 2)\lvec(0 0)
\move(0.5 0)\lvec(0.5 1) \move(0 1)\lvec(1 1)
\htext(0.25 0.5){$0$} \htext(0.75 0.5){$1$}
\htext(0.5 1.5){$2$}
\esegment
\end{texdraw}
\raisebox{0.5em}{$\ =\ $}
\begin{texdraw}
\fontsize{7}{7}\selectfont
\textref h:C v:C
\drawdim em
\setunitscale 1.9
\move(0 2)\lvec(1 2)\lvec(1 2.133)\lvec(0 2.133)\lvec(0 2)\ifill f:0.6
\move(0 2.633)\lvec(1 2.633)\lvec(1 2.899)\lvec(0 2.899)\lvec(0 2.633)
\ifill f:0.6
\move(0 5.665)\lvec(1 5.665)\lvec(1 5.798)\lvec(0 5.798)\lvec(0 5.665)
\ifill f:0.6
\move(0 0)\lvec(1 0)\lvec(1 5.798)\lvec(0 5.798)\lvec(0 0)
\move(0 0)\lvec(1 1)
\move(0 1)\lvec(1 1)
\move(0 2)\lvec(1 2)
\move(0 2.133)\lvec(1 2.133)
\move(0 2.633)\lvec(1 2.633)
\move(0 2.766)\lvec(1 2.766)
\move(0 2.899)\lvec(1 2.899)
\move(0 3.399)\lvec(1 3.399)
\move(0 4.399)\lvec(1 4.399)
\move(0 5.399)\lvec(1 5.399)
\move(0 5.532)\lvec(1 5.532)
\move(0 5.665)\lvec(1 5.665)
\move(0 5.798)\lvec(1 5.798)
\move(0 4.399)\lvec(1 5.399)
\htext(0.75 0.25){$0$}
\htext(0.25 0.75){$1$}
\htext(0.5 1.5){$2$}
\htext(0.5 2.383){$3$}
\htext(0.5 3.149){$3$}
\htext(0.5 3.899){$2$}
\htext(0.75 4.649){$0$}
\end{texdraw}
\end{center}
\vspace{2mm}

Next, we will explain the notion of \defi{splitting an $i$-block}
in a level-$l$ slice. Let $C$ be a level-$l$ slice and fix an $i$
which may be chosen as follows\,:
\begin{itemize}
\item $i=01,2,\dots,n\!-\!1$ for $\g= \bnone$, $\atwonmonetwo$.
\item $i=1,\dots,n$ for $\g= \cnone$, $\atwontwo$.
\item $i=1,\dots,n\!-\!1$ for $\g= \dnponetwo$.
\end{itemize}
Here, the choice $i=01$ is \emph{not} a typographical error. The
01-block is the unit cube obtained by gluing a 0-block and a
1-block together. Note that in any fixed slice, there can be at
most two heights in which a covering or supporting $i$-block may
appear as the top block of a layer. Similarly, there can be at
most two heights in which a supporting or covering $i$-slot may
appear.

\vskip 3mm

\begin{df}\label{3.3} \hfill
\begin{enumerate}
\item
$\g=\bnone$, $\atwonmonetwo$.\\
Suppose that there is a layer whose top is a supporting $i$-block
and another layer whose top is a covering $i$-slot.
Recall that there can be at most two heights for such layers.
Among these layers, we choose the supporting $i$-block lying in the fore-front
layer (i.e., the one with the smallest layer index) among the ones
with the higher height, and the covering $i$-slot lying in the
very back layer (i.e., the one with the largest layer index) among
the ones with the lower height.
To \defi{split an $i$-block} means to break off the top half of the
chosen supporting $i$-block and to place it in the chosen covering
$i$-slot.

\vskip 2mm

\item
$\g= \cnone$, $\atwontwo$, $\dnponetwo$.\\
Suppose that there is a layer whose top is a covering $i$-block
and another layer whose top is a supporting $i$-slot. Among these
layers, we choose the covering $i$-block lying in the fore-front
layer among the ones with the higher height, and the supporting
$i$-slot lying in the very back layer among the ones with the
lower height. To \defi{split an $i$-block} means to break off the
top half of the chosen covering $i$-block and to place it in the
chosen supporting $i$-slot.

\end{enumerate}
\end{df}

%Simply put, splitting an $i$-block
%is breaking off the top half of a covering $i$-block
%or a supporting $i$-block,
%and placing it in a supporting $i$-slot
%or covering $i$-slot, respectively.

\vskip 2mm

\begin{remark}
There is no notion of splitting for $\g= \anone$.
\end{remark}

Dotted lines shall be used to denote broken blocks,
as seen in the next example.

\begin{example}
We illustrate the notion of splitting blocks for $\g=B_3^{(1)}$.

\savebox{\tmpfiga}{
\begin{texdraw}
\fontsize{7}{7}\selectfont
\textref h:C v:C
\drawdim em
\setunitscale 1.9
\move(0 0)
\bsegment
\move(-1 -0.2)\lvec(-1 5)\lvec(0 5)
\move(0.5 5)\lvec(1 5)\lvec(1 -0.2)
\move(-1 4)\lvec(2 4)\lvec(2 -0.2)
\move(1 1)\lvec(3 1)\lvec(3 -0.2)
\move(1 0)\lvec(3 0)
\move(1 3)\lvec(2 3)
\move(1 2)\lvec(2 2)
\move(1 2.5)\lvec(2 2.5)
\move(-0.5 5)\lvec(-0.5 4)
\move(0.5 5)\lvec(0.5 4)
\move(0 5)\lvec(0 4)
\move(1.5 0)\lvec(1.5 1)
\move(2.5 0)\lvec(2.5 1)
\htext(-0.25 4.5){$1$}
\htext(-0.75 4.5){$0$}
\htext(0.75 4.5){$1$}
\htext(1.25 0.5){$0$}
\htext(1.75 0.5){$1$}
\htext(2.25 0.5){$0$}
\htext(2.75 0.5){$1$}
\htext(1.5 3.5){$2$}
\htext(1.5 2.75){$3$}
\htext(1.5 2.25){$3$}
\htext(1.5 1.5){$2$}
\esegment
\end{texdraw}}
\savebox{\tmpfigb}{
\begin{texdraw}
\fontsize{7}{7}\selectfont
\textref h:C v:C
\drawdim em
\setunitscale 1.9
\move(0 0)
\bsegment
\move(-1 -0.2)\lvec(-1 4.5)
\move(1 4.5)\lvec(1 -0.2)
\move(-1 4)\lvec(2 4)\lvec(2 -0.2)
\move(1 1)\lvec(3 1)\lvec(3 -0.2)
\move(1 0)\lvec(3 0)
\move(1 3)\lvec(2 3)
\move(1 2)\lvec(2 2)
\move(1 2.5)\lvec(2 2.5)
\move(0.5 5)\lvec(0.5 4)
\move(1.5 4.5)\lvec(1.5 4)
\move(-0.5 4.5)\lvec(-0.5 4)
\move(1 5)\lvec(0.5 5)
\move(1 5)\lvec(1 4)
\move(0 4.5)\lvec(0 4)
\move(2 4.5)\lvec(2 4)
\move(1.5 0)\lvec(1.5 1)
\move(2.5 0)\lvec(2.5 1)
\htext(-0.25 4.25){$1$}
\htext(-0.75 4.25){$0$}
\htext(1.25 4.25){$0$}
\htext(1.75 4.25){$1$}
\htext(1.25 0.5){$0$}
\htext(1.75 0.5){$1$}
\htext(2.25 0.5){$0$}
\htext(2.75 0.5){$1$}
\htext(0.75 4.5){$1$}
\htext(1.5 3.5){$2$}
\htext(1.5 2.75){$3$}
\htext(1.5 2.25){$3$}
\htext(1.5 1.5){$2$}
\lpatt(0.05 0.15)
\move(-1 4.5)\lvec(0 4.5)
\move(1 4.5)\lvec(2 4.5)
\esegment
\end{texdraw}}
\savebox{\tmpfigc}{
\begin{texdraw}
\fontsize{7}{7}\selectfont
\textref h:C v:C
\drawdim em
\setunitscale 1.9
\move(0 -0.2)\lvec(0 5)\lvec(1 5)\lvec(1 -0.2)
\move(0 2)\lvec(2 2)\lvec(2 -0.2)
\move(1 1)\lvec(3 1)\lvec(3 -0.2)
\move(0 4)\lvec(1 4)
\move(0 3)\lvec(1 3)
\move(2 0)\lvec(3 0)
\move(0 3.5)\lvec(1 3.5)
\move(1.5 1)\lvec(1.5 2)
\htext(1.25 1.5){$0$}
\htext(1.75 1.5){$1$}
\htext(2.5 0.5){$2$}
\htext(0.5 4.5){$2$}
\htext(0.5 3.75){$3$}
\htext(0.5 3.25){$3$}
\htext(0.5 2.5){$2$}
\end{texdraw}}
\savebox{\tmpfigd}{
\begin{texdraw}
\fontsize{7}{7}\selectfont
\textref h:C v:C
\drawdim em
\setunitscale 1.9
\move(0 -0.2)\lvec(0 5)\lvec(1 5)\lvec(1 -0.2)
\move(2 1.5)\lvec(2 -0.2)
\move(0 1)\lvec(3 1)\lvec(3 -0.2)
\move(3 1.5)\lvec(3 -0.2)
\move(0.5 1)\lvec(0.5 2)
\move(0 2)\lvec(1 2)
\move(0 4)\lvec(1 4)
\move(0 3)\lvec(1 3)
\move(2 0)\lvec(3 0)
\move(0 3.5)\lvec(1 3.5)
\move(1.5 1)\lvec(1.5 1.5)
\move(2.5 1)\lvec(2.5 1.5)
\htext(1.25 1.25){$0$}
\htext(1.75 1.25){$1$}
\htext(2.25 1.25){$0$}
\htext(2.75 1.25){$1$}
\htext(0.25 1.5){$0$}
\htext(0.75 1.5){$1$}
\htext(2.5 0.5){$2$}
\htext(0.5 4.5){$2$}
\htext(0.5 3.75){$3$}
\htext(0.5 3.25){$3$}
\htext(0.5 2.5){$2$}
\lpatt(0.05 0.15)
\move(1 1.5)\lvec(3 1.5)
\end{texdraw}}
\savebox{\tmpfige}{
\begin{texdraw}
\fontsize{7}{7}\selectfont
\textref h:C v:C
\drawdim em
\setunitscale 1.9
\move(0 -0.2)\lvec(0 5)
\move(1 2.8)\lvec(1 5)
\move(2 1.8)\lvec(2 4)
\move(3 0.8)\lvec(3 3)
\move(4 -0.2)\lvec(4 2)
\move(5 -0.2)\lvec(5 2)
\move(6 -0.2)\lvec(6 1)
\move(1.5 3)\lvec(1.5 4)
\move(0 5)\lvec(1 5)
\move(0 4)\lvec(2 4)
\move(0.8 3)\lvec(3 3)
\move(1.8 2)\lvec(5 2)
\move(3 1.5)\lvec(5 1.5)
\move(2.8 1)\lvec(6 1)
\move(3.8 0)\lvec(6 0)
\htext(1.25 3.5){$0$}
\htext(1.75 3.5){$1$}
\htext(2.5 2.5){$2$}
\htext(0.5 4.5){$2$}
\htext(3.5 1.75){$3$}
\htext(3.5 1.25){$3$}
\htext(4.5 1.75){$3$}
\htext(4.5 1.25){$3$}
\htext(4.5 0.5){$2$}
\htext(5.5 0.5){$2$}
\end{texdraw}}
\savebox{\tmpfigf}{
\begin{texdraw}
\fontsize{7}{7}\selectfont
\textref h:C v:C
\drawdim em
\setunitscale 1.9
\move(0 -0.2)\lvec(0 4.5)
\move(1 2.8)\lvec(1 4.5)
\move(2 1.8)\lvec(2 4)
\move(3 0.8)\lvec(3 3)
\move(4 -0.2)\lvec(4 2.5)
\move(5 -0.2)\lvec(5 2)
\move(6 -0.2)\lvec(6 1)
\move(1.5 3)\lvec(1.5 4)
\move(0 4)\lvec(2 4)
\move(0.8 3)\lvec(3 3)
\move(1.8 2)\lvec(5 2)
\move(3 1.5)\lvec(5 1.5)
\move(2.8 1)\lvec(6 1)
\move(3.8 0)\lvec(6 0)
\htext(1.25 3.5){$0$}
\htext(1.75 3.5){$1$}
\htext(2.5 2.5){$2$}
\htext(0.5 4.25){$2$}
\htext(3.5 1.75){$3$}
\htext(3.5 1.25){$3$}
\htext(4.5 1.75){$3$}
\htext(4.5 1.25){$3$}
\htext(4.5 0.5){$2$}
\htext(5.5 0.5){$2$}
\htext(3.5 2.25){$2$}
\lpatt(0.05 0.15)
\move(0 4.5)\lvec(1 4.5)
\move(3 2.5)\lvec(4 2.5)
\end{texdraw}}
\savebox{\tmpfigg}{
\begin{texdraw}
\fontsize{7}{7}\selectfont
\textref h:C v:C
\drawdim em
\setunitscale 1.9
\move(0 -0.2)\lvec(0 4.5)
\move(1 2.8)\lvec(1 4.5)
\move(2 1.8)\lvec(2 4)
\move(3 0.8)\lvec(3 3)
\move(4 -0.2)\lvec(4 2.5)
\move(5 -0.2)\lvec(5 2.5)
\move(6 -0.2)\lvec(6 0.5)
\move(1.5 3)\lvec(1.5 4)
\move(0 4)\lvec(2 4)
\move(0.8 3)\lvec(3 3)
\move(1.8 2)\lvec(5 2)
\move(3 1.5)\lvec(5 1.5)
\move(2.8 1)\lvec(5 1)
\move(3.8 0)\lvec(6 0)
\htext(1.25 3.5){$0$}
\htext(1.75 3.5){$1$}
\htext(2.5 2.5){$2$}
\htext(0.5 4.25){$2$}
\htext(3.5 1.75){$3$}
\htext(3.5 1.25){$3$}
\htext(4.5 1.75){$3$}
\htext(4.5 1.25){$3$}
\htext(4.5 0.5){$2$}
\htext(5.5 0.25){$2$}
\htext(3.5 2.25){$2$}
\htext(4.5 2.25){$2$}
\lpatt(0.05 0.15)
\move(0 4.5)\lvec(1 4.5)
\move(3 2.5)\lvec(5 2.5)
\move(5 0.5)\lvec(6 0.5)
\end{texdraw}}
\begin{enumerate}
\item splitting an $01$-block\\
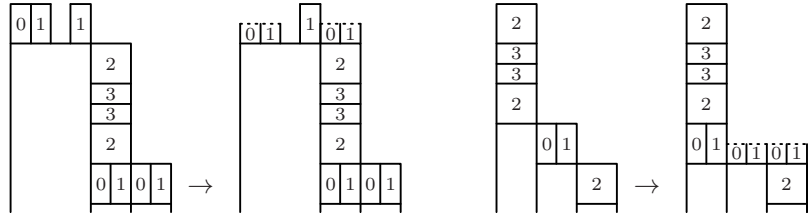
\begin{center}
\usebox{\tmpfiga}
\ \raisebox{0.7em}{$\rightarrow$} \
\usebox{\tmpfigb}
\quad\quad\quad
\usebox{\tmpfigc}
\ \raisebox{0.7em}{$\rightarrow$} \
\usebox{\tmpfigd}
\end{center}
\item splitting a $2$-block, twice\\
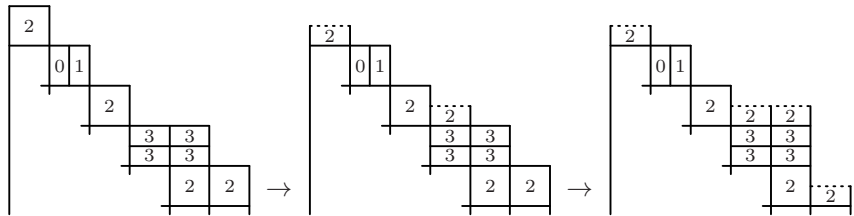
\begin{center}
\usebox{\tmpfige}
\ \raisebox{0.7em}{$\rightarrow$}
\usebox{\tmpfigf}
\ \raisebox{0.7em}{$\rightarrow$}
\usebox{\tmpfigg}
\end{center}
\end{enumerate}
\end{example}

\vskip 3mm

\begin{remark}
For a level-$l$ slice $C$, the result obtained after splitting all
possible $i$-blocks is not a level-$l$ slice. We will call it the
\defi{$i$-split form} of $C$. The \defi{split form} of $C$ is
defined to be the result obtained after splitting all possible
$i$-blocks for all $i$.
\end{remark}

We now define the action of Kashiwara operators on the set
$\slice$ of level-$l$ slices.
Let $C$ be a level-$l$ slice and fix an index $i\in I$.

\begin{enumerate}
\item $\g=\anone$\\
      The actions of $\eit$ and $\fit$ are defined by the rules
      $e\text{a})-e\text{b})$ and $f\text{a})-f\text{b})$, respectively.
\begin{enumerate}
\item[$e$a)] If there is no layer in $C$ whose top is an $i$-block,
             we define $\eit C =0$.
\item[$e$b)] If $C$ contains some layers whose top is an
             $i$-block, we remove an $i$-block from the (top of) fore-front
             layer among the ones with the higher height.
\item[$f$a)] If there is no layer whose top is an $i$-slot, we
             define $\fit C=0$.
\item[$f$b)] If $C$ contains some layer whose top is an
             $i$-slot, then we add an $i$-block on top of the very back layer
             among the ones with the lower height.
\end{enumerate}
\vskip 3mm
\item $\g=\bnone$\\
      For $i\neq 0, 1, 2$, let $C'$ be the $(i-1)$-split form of $C$.
      For $i=2$, let $C'$ be the $01$-split form of $C$, and for $i=0,
      1$, let $C'=C$. The actions of $\eit$ and $\fit$ are defined by
      the rules $e\text{a})-e\text{c})$ and $f\text{a})-f\text{c})$,
      respectively.
\begin{enumerate}
\item[$e$a)] If there is no layer in $C'$ whose top is an $i$-block, we define
             $\eit C=0$.
\item[$e$b)] If $C'$ contains some layer whose top is an $i$-block and
             all of these $i$-blocks are supporting blocks, then remove an
             $i$-block from the fore-front layer among the ones with the higher
             height.
\item[$e$c)] If $C'$ contains some layer whose top is an $i$-block and
             some of these $i$-blocks are covering blocks, then single out the
             layers with the higher height among the ones containing covering
             $i$-blocks. We remove an $i$-block from the fore-front layer among
             the chosen ones.
\item[$f$a)] If there is no layer in $C'$ whose top is an $i$-slot, we
             define $\fit C =0$.
\item[$f$b)] If $C'$ contains some layer whose top is an $i$-slot and all
             of these $i$-slots are covering slots, then add an $i$-block on
             top of the very back layer among the ones with the lower height.
\item[$f$c)] If $C'$ contains some layer whose top is an $i$-slot and
             some of these $i$-slots are supporting slots, then single out the
             layers with the lower height among the ones containing supporting
             $i$-slots. We add an $i$-block on top of the very back layer among
             the chosen ones.
\end{enumerate}
\vskip 3mm
\item $\g=\atwonmonetwo$\\
For $i=0,1, \dots, n-1$, we use the rules for $\g=\bnone$. For
$i=n$, let $C'$ be the $(n-1)$-split form of $C$. The actions of
$\eit$ and $\fit$ are defined by the rules $e\text{a})-e\text{d})$
and $f\text{a})-f\text{d})$, respectively.
\begin{enumerate}
\item[$e$a)] If there is no layer in $C'$ whose top is an $n$-block,
             we define $\eit C =0$.
\item[$e$b)] If $C'$ contains some layer whose top is an $n$-block
             and all of these $n$-blocks are supporting blocks, then the number
             of $n$-blocks must be even. Single out the layers with the higher
             height among the ones containing $n$-blocks.\\
             -- If there is only one such layer, remove an $n$-block from that
             layer and another $n$-block from the fore-front layer among
             the remaining ones with $n$-blocks.\\
             -- If there are more than one such layers, remove two $n$-blocks
             from the two front layers (i.e., one $n$-block from each layer)
             among the chosen ones.
\item[$e$c)] If $C'$ contains some layer whose top is a covering
             $n$-block and there is only one such layer,
             first of all, remove an $n$-blocks from that layer.
             Then the top of the layer from which the block was removed
             will be a supporting $n$-block.
             So this intermediate result contains
             at least one supporting $n$-block.
             Single out the layers with the higher
             height among the ones containing supporting $n$-blocks.
             We remove an $n$-block from the fore-front layer among
             the singled out layers.
\item[$e$d)] If $C'$ contains more than one layers whose top is a
             covering $n$-block, then single out the layers with the higher
             height among the ones containing covering $n$-blocks.\\
             -- If there is only one such layer, remove an $n$-block from that
             layer and remove another $n$-block from the fore-front layer among the
             remaining ones with covering $n$-blocks.\\
             -- If there are more than one such layers, remove two $n$-blocks
             from the two front layers among the chosen ones.
\item[$f$a)] If there is no layer in $C'$ whose top is an $n$-slot,
             we define $\fit C=0$.
\item[$f$b)] If $C'$ contains some layer whose top is an $n$-slot
             and all of these $n$-slots are covering slots, then the number of
             $n$-slots must be even. Single out the layers with the lower
             height among the ones containing $n$-slots.\\
             -- If there is only one such layer, add an $n$-block on top of
             that layer and another $n$-block on top of the very back layer
             among the remaining ones with $n$-slots.\\
             -- If there are more than one such layers, add two $n$-blocks on
             top of the two back layers (i.e., one $n$-block on top of each
             layer) among the chosen ones.
\item[$f$c)] If $C'$ contains some layer whose top is a supporting
             $n$-slot and there is only one such layer,
             first of all, add an $n$-block on top of that layer.
             Then the slot on top of the
             block just added will be a covering $n$-slot.
             So this intermediate result contains at least one covering
             $n$-slot.
             %(Now, there can be at most two height at which these covering
             %$n$-slots may appear.)
             Single out the layers with the lower
             height among the ones containing covering $n$-slots.
             We add an $n$-block to the very back of
             the singled out layers.
\item[$f$d)] If $C'$ contains more than one layers whose top is a
             supporting $n$-slot, then single out the layers with the lower
             height among the ones containing supporting $n$-slots.\\
             -- If there is only one such layer, add an $n$-block on top
             of that
             layer and add another $n$-block on top of the very back layer
             among the remaining ones with supporting $n$-slots.\\
             -- If there are more than one such layers, add two $n$-blocks on
             top of the two very back layers among the chosen ones.
\end{enumerate}
\vskip 3mm
\item $\g=\atwontwo$\\
      For $i\neq n$, let $C'$ be the $(i+1)$-split form of $C$. For
      $i=n$, let $C'=C$. Then we use the rules for $\g=\bnone$ with the
      following substitution of words\,:
\begin{itemize}
\item covering $\mapsto$ supporting,
\item supporting $\mapsto$ covering.
\end{itemize}
\vskip 3mm
\item $\g=\cnone$\\
      For $i=1, \dots, n$, we will use the rules for $\atwontwo$. For
      $i=0$, let $C'$ be the 1-split form of $C$, and use the rules for
      $i=n$ over $\atwonmonetwo$ with the following substitution of the
      words\,:
\begin{itemize}
\item $n$-block $\mapsto$ 0-block,
\item $n$-slot $\mapsto$ 0-slot,
\item covering $\mapsto$ supporting,
\item supporting $\mapsto$ covering.
\end{itemize}
\vskip 3mm
\item $\g=\dnponetwo$\\
For $i\neq n-1, n$, let $C'$ be the $(i+1)$-split form of $C$. For
$i=n-1, n$, let $C'=C$. Then we use the rules for $\g=\atwontwo$.
\end{enumerate}
\vspace{4mm}

Let $C= (c_1, \dots, c_l)$ (resp. $(c_1, \dots, c_{2l})$
for $\cnone$) be a level-$l$ slice.
We define the slices $C \pm \delta$ by
\begin{equation*}\label{eq:add a delta}
\begin{aligned}
C+\delta & = (c_2, \dots, c_l, c_1 + \delta)\
          (\text{resp}. \ (c_2, \dots, c_{2l}, c_1 + \delta) \text{ for } \cnone), \\
C-\delta & = (c_l - \delta, c_1, \cdots, c_{l-1})\
          (\text{resp}. \ (c_{2l} - \delta, c_1, \cdots, c_{2l-1}) \text{ for } \cnone).
\end{aligned}
\end{equation*}
We say that two slices $C$ and $C'$ are \defi{related}, denoted by
$C\sim C'$, if one of the two slices may be obtained from the
other by adding finitely many $\delta$'s. Let
\begin{equation*}
\newcrystal = \slice/\sim
\end{equation*}
be the set of equivalence classes of level-$l$ slices. We will use
the same symbol $C$ for the equivalence class containing the
level-$l$ slice $C$. By abuse of terminology, the equivalence
class containing a slice $C$ will be often referred to as the {\it
slice} $C$. Note that the map $C \mapsto C + \delta$ commutes with
the action of Kashiwara operators. Hence they are well-defined on
$\newcrystal$. We define
\begin{align*}
\veps_i ( C ) &= \max\{ n \mid \eit^n C \in \newcrystal \},\\
\vphi_i ( C ) &= \max\{ n \mid \fit^n C \in \newcrystal \},\\
\cwt( C ) &= \sum_i \big(\vphi_i(C) - \veps_i(C)\big) \La_i.
\end{align*}
Then it is lengthy but straightforward to prove the following
proposition.

\begin{prop}
The Kashiwara operators, together with the maps $\veps_i$,
$\vphi_i$ $(i\in I)$, $\cwt$, define a $\uq'(\g)$-crystal
structure on the set $\newcrystal$.
\end{prop}

\vskip 1cm

\section{New realization of perfect crystals}
In this section, we will show that the $\uq'(\g)$-crystal
$\newcrystal$ gives a new realization of the level-$l$ perfect
crystal $\oldcrystal$ described in Section~\ref{sec:02}. We first
define a \emph{canonical map} $\psi: \oldcrystal \rightarrow
\newcrystal$ as follows.

\vspace{3mm}
\begin{enumerate}
\item $\g=\anone$\\
Recall that every element $b \in \oldcrystal$ has the form
$b=(x_0, x_1, \dots, x_n)$ with $x_i \in \Z_{\ge 0}$, $\sum x_i = l$.
For ease of writing, we shall temporarily use the following notation :
\begin{center}
 \savebox{\tmpfiga}{\begin{texdraw}
 \fontsize{7}{7}\selectfont
 \textref h:C v:C
 \drawdim em
 \setunitscale 1.9
 \move(0 -0.2)\lvec(0 1)\lvec(1 1)\lvec(1 -0.2)
 \move(0 0)\lvec(1 0)
 \htext(0.5 0.5){$0$}
 \end{texdraw}}%
 \savebox{\tmpfigb}{\begin{texdraw}
 \fontsize{7}{7}\selectfont
 \textref h:C v:C
 \drawdim em
 \setunitscale 1.9
 \move(0 -0.2)\lvec(0 1)\lvec(1 1)\lvec(1 -0.2)
 \move(0 0)\rlvec(1 0)
 \htext(0.5 0.5){$i$}
 \end{texdraw}}%
 \savebox{\tmpfigc}{\begin{texdraw}
 \fontsize{7}{7}\selectfont
 \textref h:C v:C
 \drawdim em
 \setunitscale 1.9
 \move(0 -0.2)\lvec(0 1)\lvec(1 1)\lvec(1 -0.2)
 \move(0 0)\rlvec(1 0)
 \htext(0.5 0.5){$n$}
 \end{texdraw}}%
{\allowdisplaybreaks
\begin{align*}
%\bu_0 \quad&\longmapsto\quad
% \raisebox{-0.5em}{\usebox{\tmpfiga}}\\[0.3em]
\bu_{i} \ & = \
 \raisebox{-0.5em}{\usebox{\tmpfigb}}
 \qquad \text{ for each } i\in I.\\[0.3em]
%\bu_n \quad&\longmapsto\quad
% \raisebox{-0.5em}{\usebox{\tmpfigc}}
\end{align*}}% end of allowdisplaybreaks
\end{center}
The image of the above $b$ under the map $\psi$ is defined
to the equivalence class of a
level-$l$ slice obtained by pasting together $x_i$-many
$\bu_i$'s for each $i$.
Note that this equivalence class does not depend on the
way we have pasted the $\bu_i$'s together, as long as the pasted
result forms a slice.

We will use the notation $\psi(b)=[x_0,x_1,\cdots,x_n]$,
as needed.
%gathering together the number of $\bu_i$ used in the pasting.
%

\vskip 3mm

\item $\g=\atwontwo$ and $\cnone$\\
Let $b=(x_1, \dots, x_n | \bar x_n, \dots, \bar x_1) \in
\oldcrystal$ over $\atwontwo$ (resp. $\cnone$) with $x_i, \bar x_i
\in \Z_{\ge 0}$, $\sum (x_i + \bar x_i) = k \le l$ (resp. $\sum
(x_i + \bar x_i) = 2k \le 2l$).
We shall temporarily use the following notation :
\begin{center}
\savebox{\tmpfiga}{\begin{texdraw} \fontsize{7}{7}\selectfont
\textref h:C v:C \drawdim em \setunitscale 1.9 \move(0
-0.2)\lvec(0 1.5)\lvec(1 1.5)\lvec(1 -0.2) \move(0 1)\lvec(1 1)
\move(0 0)\lvec(1 0) \htext(0.5 1.25){$0$} \htext(0.5 0.5){$1$}
\end{texdraw}}%
\savebox{\tmpfigb}{\begin{texdraw} \fontsize{7}{7}\selectfont
\textref h:C v:C \drawdim em \setunitscale 1.9 \move(0
-0.2)\lvec(0 1)\lvec(1 1)\lvec(1 -0.2) \move(0 0)\rlvec(1 0)
\move(0 0.5)\rlvec(1 0) \htext(0.5 0.25){$0$} \htext(0.5
0.75){$0$}
\end{texdraw}}%
\savebox{\tmpfigc}{\begin{texdraw} \fontsize{7}{7}\selectfont
\textref h:C v:C \drawdim em \setunitscale 1.9 \move(0
-0.2)\lvec(0 2)\lvec(1 2)\lvec(1 -0.2) \move(0 0)\rlvec(1 0)
\move(0 1)\rlvec(1 0) \htext(0.5 0.5){$i\!\!-\!\!2$} \htext(0.5
1.5){$i\!\!-\!\!1$}
\end{texdraw}}%
\savebox{\tmpfigd}{\begin{texdraw} \fontsize{7}{7}\selectfont
\textref h:C v:C \drawdim em \setunitscale 1.9 \move(0
-0.2)\lvec(0 1)\lvec(1 1)\lvec(1 -0.2) \move(0 0)\rlvec(1 0)
\htext(0.5 0.5){$n$}
\end{texdraw}}%
\savebox{\tmpfige}{\begin{texdraw} \fontsize{7}{7}\selectfont
\textref h:C v:C \drawdim em \setunitscale 1.9 \move(0
-0.2)\lvec(0 2)\lvec(1 2)\lvec(1 -0.2) \move(0 0)\rlvec(1 0)
\move(0 1)\rlvec(1 0) \htext(0.5 0.5){$i\!\!+\!\!1$} \htext(0.5
1.5){$i$}
\end{texdraw}}%
{\allowdisplaybreaks
\begin{align*}
\bw_0 \ &= \
 \raisebox{-0.5em}{\usebox{\tmpfiga}}\\[0.3em]
\bu_1 \ &= \
 \raisebox{-0.5em}{\usebox{\tmpfigb}}\\[0.3em]
\bu_i \ &= \
 \raisebox{-0.5em}{\usebox{\tmpfigc}}
 \qquad \text{ for $i=2,\cdots,n$, }\\[0.3em]
\bv_n \ &= \
 \raisebox{-0.5em}{\usebox{\tmpfigd}}\\[0.3em]
\bv_i \ &= \
 \raisebox{-0.5em}{\usebox{\tmpfige}}
 \qquad \text{ for $i=1,\cdots,n-1$. }\\[0.3em]
\end{align*}}% end of allowdisplaybreaks
\end{center}
The image of $b$ under the map $\psi$ is defined
to the equivalence class of a
level-$l$ slice obtained by pasting together
$(l-k)$ (resp. $2(l-k)$) -many $\bw_0$'s, $x_i$-many
$\bu_i$'s, and $\bar{x}_i$-many $\bv_i$'s, for each $i$.

We will use the notation $\psi(b) = [t_0|x_1, \dots, x_n | \bar x_n, \dots,
\bar x_1]$, as needed.
Here, $t_0 = l-k$ (resp. $2(l-k)$).

\vskip 3mm

\item $\g= \dnponetwo$\\
Let $b=(x_1, \dots, x_n | x_0 | \bar x_n, \dots, \bar x_1) \in
\oldcrystal$ over $\dnponetwo$ with $x_0=0 \text{ or } 1$, $x_i, \bar x_i
\in \Z_{\ge 0}$, $x_0 + \sum (x_i + \bar x_i) = k \le l$.
We shall temporarily use the following notation :
\begin{center}
\savebox{\tmpfiga}{\begin{texdraw} \fontsize{7}{7}\selectfont
\textref h:C v:C \drawdim em \setunitscale 1.9 \move(0
-0.2)\lvec(0 1.5)\lvec(1 1.5)\lvec(1 -0.2) \move(0 1)\lvec(1 1)
\move(0 0)\lvec(1 0) \htext(0.5 1.25){$0$} \htext(0.5 0.5){$1$}
\end{texdraw}}%
\savebox{\tmpfigb}{\begin{texdraw} \fontsize{7}{7}\selectfont
\textref h:C v:C \drawdim em \setunitscale 1.9 \move(0
-0.2)\lvec(0 1)\lvec(1 1)\lvec(1 -0.2) \move(0 0)\rlvec(1 0)
\move(0 0.5)\rlvec(1 0) \htext(0.5 0.25){$0$} \htext(0.5
0.75){$0$}
\end{texdraw}}%
\savebox{\tmpfigc}{\begin{texdraw} \fontsize{7}{7}\selectfont
\textref h:C v:C \drawdim em \setunitscale 1.9 \move(0
-0.2)\lvec(0 2)\lvec(1 2)\lvec(1 -0.2) \move(0 0)\rlvec(1 0)
\move(0 1)\rlvec(1 0) \htext(0.5 0.5){$i\!\!-\!\!2$} \htext(0.5
1.5){$i\!\!-\!\!1$}
\end{texdraw}}%
\savebox{\tmpfigd}{\begin{texdraw} \fontsize{7}{7}\selectfont
\textref h:C v:C \drawdim em \setunitscale 1.9 \move(0
-0.2)\lvec(0 1.5)\lvec(1 1.5)\lvec(1 -0.2) \move(0 0)\rlvec(1 0)
\move(0 1)\rlvec(1 0) \htext(0.5 0.5){$n\!\!-\!\!1$} \htext(0.5
1.25){$n$}
\end{texdraw}}%
\savebox{\tmpfige}{\begin{texdraw} \fontsize{7}{7}\selectfont
\textref h:C v:C \drawdim em \setunitscale 1.9 \move(0
-0.2)\lvec(0 1)\lvec(1 1)\lvec(1 -0.2) \move(0 0)\rlvec(1 0)
\move(0 0.5)\rlvec(1 0) \htext(0.5 0.25){$n$} \htext(0.5
0.75){$n$}
\end{texdraw}}%
\savebox{\tmpfigf}{\begin{texdraw} \fontsize{7}{7}\selectfont
\textref h:C v:C \drawdim em \setunitscale 1.9 \move(0
-0.2)\lvec(0 2)\lvec(1 2)\lvec(1 -0.2) \move(0 0)\rlvec(1 0)
\move(0 1)\rlvec(1 0) \htext(0.5 0.5){$i\!\!+\!\!1$} \htext(0.5
1.5){$i$}
\end{texdraw}}%
{\allowdisplaybreaks
\begin{align*}
\bw_0 \ &= \
 \raisebox{-0.5em}{\usebox{\tmpfiga}}\\[0.3em]
\bu_1 \ &= \
 \raisebox{-0.5em}{\usebox{\tmpfigb}}\\[0.3em]
\bu_i \ &= \
 \raisebox{-0.5em}{\usebox{\tmpfigc}}
 \qquad \text{ for $i=2,\cdots,n$, }\\[0.3em]
\bu_0 \ &= \
 \raisebox{-0.5em}{\usebox{\tmpfigd}}\\[0.3em]
\bv_n \ &= \
 \raisebox{-0.5em}{\usebox{\tmpfige}}\\[0.3em]
\bv_i \ &= \
 \raisebox{-0.5em}{\usebox{\tmpfigf}}
 \qquad \text{ for $i=1,\cdots,n-1$. }\\[0.3em]
\end{align*}}% end of allowdisplaybreaks
\end{center}
The image of $b$ under the map $\psi$ is defined
to the equivalence class of a
level-$l$ slice obtained by pasting together $x_0$-many $\bu_0$,
$(l-k)$-many $\bw_0$'s, $x_i$-many
$\bu_i$'s, and $\bar{x}_i$-many $\bv_i$'s, for each $i$.

We will use the notation $\psi(b) = [t_0|x_1, \dots, x_n | x_0 | \bar x_n, \dots,
\bar x_1]$, as needed.
Here, $t_0 = l-k$.

\vskip 3mm

\item $\g= \bnone$\\
Let $b=(x_1, \dots, x_n | x_0 | \bar x_n, \dots, \bar x_1) \in
\oldcrystal$ over $\bnone$ with $x_0 = 0$ or $1$, $x_i, \bar x_i
\in \Z_{\ge 0}$, $x_0 + \sum (x_i + \bar x_i) = l$.
And set
\begin{equation*}\label{eq:bnone}
\begin{aligned}
& x'_1 = (x_1 - \bar x_1)_{+}, \\
& x'_2 = (x_2 - \bar x_2)_{+} + \text{min} \{x_1, \bar x_1 \}, \\
& \qquad \qquad \vdots \\
& x'_n = (x_n - \bar x_n)_{+} + \text{min} \{x_{n-1}, \bar x_{n-1}
\}, \\
& x'_0 = x_0 + 2 \text{min} \{x_n, \bar x_n \}, \\
& \bar x'_n = (\bar x_n - x_n)_{+} + \text{min} \{x_{n-1}, \bar
x_{n-1} \}, \\
& \qquad \qquad \vdots \\
& \bar x'_2 = (\bar x_2 - x_2)_{+} + \text{min} \{x_1, \bar x_1 \},
\\
& \bar x'_1 = (\bar x_1 - x_1)_{+}.
\end{aligned}
\end{equation*}
Here, we used the notation $(x)_+ = \max(0,x)$.
We shall temporarily use the following notation :
\begin{center}
\savebox{\tmpfiga}{\begin{texdraw} \fontsize{7}{7}\selectfont
\textref h:C v:C \drawdim em \setunitscale 1.9 \move(0
-0.2)\lvec(0 1)\lvec(1 1)\lvec(1 -0.2) \move(0 0)\lvec(1 0)
\move(0 0)\lvec(1 1) \move(0 0)\lvec(1 0) \htext(0.75 0.25){$0$}
\end{texdraw}}%
\savebox{\tmpfigf}{\begin{texdraw}
\fontsize{7}{7}\selectfont
\textref h:C v:C
\drawdim em
\setunitscale 1.9
\move(0 -0.2)\lvec(0 1)\lvec(1 1)\lvec(1 -0.2)
\move(0 0)\lvec(1 0)
\move(0 0)\lvec(1 1)
\move(0 0)\lvec(1 0)
\htext(0.25 0.75){$0$}
\end{texdraw}}%
\savebox{\tmpfigb}{\begin{texdraw}
\fontsize{7}{7}\selectfont
\textref h:C v:C
\drawdim em
\setunitscale 1.9
\move(0 -0.2)\lvec(0 2)\lvec(1 2)\lvec(1 -0.2)
\move(0 0)\rlvec(1 0)
\move(0 1)\rlvec(1 0)
\htext(0.5 0.5){$i\!\!-\!\!2$}
\htext(0.5 1.5){$i\!\!-\!\!1$}
\end{texdraw}}%
\savebox{\tmpfigc}{\begin{texdraw}
\fontsize{7}{7}\selectfont
\textref h:C v:C
\drawdim em
\setunitscale 1.9
\move(0 -0.2)\lvec(0 1.5)\lvec(1 1.5)\lvec(1 -0.2)
\move(0 0)\rlvec(1 0)
\move(0 1)\rlvec(1 0)
\htext(0.5 0.5){$n\!\!-\!\!1$}
\htext(0.5 1.25){$n$}
\end{texdraw}}%
\savebox{\tmpfigd}{\begin{texdraw}
\fontsize{7}{7}\selectfont
\textref h:C v:C
\drawdim em
\setunitscale 1.9
\move(0 -0.2)\lvec(0 1)\lvec(1 1)\lvec(1 -0.2)
\move(0 0)\rlvec(1 0)
\move(0 0.5)\rlvec(1 0)
\htext(0.5 0.25){$n$}
\htext(0.5 0.75){$n$}
\end{texdraw}}%
\savebox{\tmpfige}{\begin{texdraw}
\fontsize{7}{7}\selectfont
\textref h:C v:C
\drawdim em
\setunitscale 1.9
\move(0 -0.2)\lvec(0 2)\lvec(1 2)\lvec(1 -0.2)
\move(0 0)\rlvec(1 0)
\move(0 1)\rlvec(1 0)
\htext(0.5 0.5){$i\!\!+\!\!1$}
\htext(0.5 1.5){$i$}
\end{texdraw}}%
\savebox{\tmpfigg}{\begin{texdraw}
\fontsize{7}{7}\selectfont
\textref h:C v:C
\drawdim em
\setunitscale 1.9
\move(0 -0.2)\lvec(0 1)\lvec(1 1)\lvec(1 -0.2)
\move(0 0)\lvec(1 0)
\move(0 0)\lvec(1 1)
\move(0 0)\lvec(1 0)
\htext(0.25 0.75){$1$}
\end{texdraw}}%
\savebox{\tmpfigh}{\begin{texdraw}
\fontsize{7}{7}\selectfont
\textref h:C v:C
\drawdim em
\setunitscale 1.9
\move(0 -0.2)\lvec(0 1)\lvec(1 1)\lvec(1 -0.2)
\move(0 0)\lvec(1 0)
\move(0 0)\lvec(1 1)
\move(0 0)\lvec(1 0)
\htext(0.75 0.25){$1$}
\end{texdraw}}%
{\allowdisplaybreaks
\begin{align*}
\bu_1 \ &= \
 \raisebox{-0.5em}{\usebox{\tmpfiga}}
 \ \ \textup{or}\ \ \raisebox{-0.5em}{\usebox{\tmpfigf}}\\[0.3em]
\bu_i \ &= \
 \raisebox{-0.5em}{\usebox{\tmpfigb}}
 \qquad \text{ for $i=2,\cdots,n$, }\\[0.3em]
\bu_0 \ &= \
 \raisebox{-0.5em}{\usebox{\tmpfigc}}\\[0.3em]
\bv_n \ &= \
 \raisebox{-0.5em}{\usebox{\tmpfigd}}\\[0.3em]
\bv_i \ &= \
 \raisebox{-0.5em}{\usebox{\tmpfige}}
 \qquad \text{ for $i=2,\cdots,n-1$, }\\[0.3em]
\bv_1 \ &= \
 \raisebox{-0.5em}{\usebox{\tmpfigg}}
 \ \ \textup{or}\ \ \raisebox{-0.5em}{\usebox{\tmpfigh}}\\[0.3em]
\end{align*}}% end of allowdisplaybreaks
\end{center}
The image of $b$ under the map $\psi$ is defined
to the equivalence class of a
level-$l$ slice obtained by pasting together $x'_0$-many $\bu_0$'s,
$x'_i$-many $\bu_i$'s, and $\bar{x}'_i$-many $\bv_i$'s, for each $i$.

We will use the notation $\psi(b) = [x'_1, \dots, x'_n | x'_0 | \bar x'_n, \dots,
\bar x'_1]$, as needed.

\vskip 3mm

\item $\g=\atwonmonetwo$\\
Let $b=(x_1, \dots, x_n | \bar x_n, \dots, \bar x_1) \in
\oldcrystal$ over $\atwonmonetwo$ with $x_i, \bar x_i \in \Z_{\ge
0}$, $\sum (x_i + \bar x_i) = l$.
And set
\begin{equation*}\label{eq:atwonmonetwo}
\begin{aligned}
& x'_1 = (x_1 - \bar x_1)_{+}, \\
& x'_2 = (x_2 - \bar x_2)_{+} + \text{min} \{x_1, \bar x_1 \}, \\
& \qquad \qquad \vdots \\
& x'_n = (x_n - \bar x_n)_{+} + \text{min} \{x_{n-1}, \bar x_{n-1}
\}, \\
& x'_0 = 2 \text{min} \{x_n, \bar x_n \}, \\
& \bar x'_n = (\bar x_n - x_n)_{+} + \text{min} \{x_{n-1}, \bar
x_{n-1} \}, \\
& \qquad \qquad \vdots \\
& \bar x'_2 = (\bar x_2 - x_2)_{+} + \text{min} \{x_1, \bar x_1 \},
\\
& \bar x'_1 = (\bar x_1 - x_1)_{+}.
\end{aligned}
\end{equation*}
We shall temporarily use the following notation :
\begin{center}
\savebox{\tmpfiga}{\begin{texdraw}
\fontsize{7}{7}\selectfont
\textref h:C v:C
\drawdim em
\setunitscale 1.9
\move(0 -0.2)\lvec(0 1)\lvec(1 1)\lvec(1 -0.2)
\move(0 0)\lvec(1 0)
\move(0 0)\lvec(1 1)
\move(0 0)\lvec(1 0)
\htext(0.75 0.25){$0$}
\end{texdraw}}%
\savebox{\tmpfigf}{\begin{texdraw}
\fontsize{7}{7}\selectfont
\textref h:C v:C
\drawdim em
\setunitscale 1.9
\move(0 -0.2)\lvec(0 1)\lvec(1 1)\lvec(1 -0.2)
\move(0 0)\lvec(1 0)
\move(0 0)\lvec(1 1)
\move(0 0)\lvec(1 0)
\htext(0.25 0.75){$0$}
\end{texdraw}}%
\savebox{\tmpfigb}{\begin{texdraw}
\fontsize{7}{7}\selectfont
\textref h:C v:C
\drawdim em
\setunitscale 1.9
\move(0 -0.2)\lvec(0 2)\lvec(1 2)\lvec(1 -0.2)
\move(0 0)\rlvec(1 0)
\move(0 1)\rlvec(1 0)
\htext(0.5 0.5){$i\!\!-\!\!2$}
\htext(0.5 1.5){$i\!\!-\!\!1$}
\end{texdraw}}%
\savebox{\tmpfigc}{\begin{texdraw}
\fontsize{7}{7}\selectfont
\textref h:C v:C
\drawdim em
\setunitscale 1.9
\move(0 -0.2)\lvec(0 1.5)\lvec(1 1.5)\lvec(1 -0.2)
\move(0 0)\rlvec(1 0)
\move(0 1)\rlvec(1 0)
\htext(0.5 0.5){$n\!\!-\!\!1$}
\htext(0.5 1.25){$n$}
\end{texdraw}}%
\savebox{\tmpfigd}{\begin{texdraw}
\fontsize{7}{7}\selectfont
\textref h:C v:C
\drawdim em
\setunitscale 1.9
\move(0 -0.2)\lvec(0 1)\lvec(1 1)\lvec(1 -0.2)
\move(0 0)\rlvec(1 0)
\move(0 0.5)\rlvec(1 0)
\htext(0.5 0.25){$n$}
\htext(0.5 0.75){$n$}
\end{texdraw}}%
\savebox{\tmpfige}{\begin{texdraw}
\fontsize{7}{7}\selectfont
\textref h:C v:C
\drawdim em
\setunitscale 1.9
\move(0 -0.2)\lvec(0 2)\lvec(1 2)\lvec(1 -0.2)
\move(0 0)\rlvec(1 0)
\move(0 1)\rlvec(1 0)
\htext(0.5 0.5){$i\!\!+\!\!1$}
\htext(0.5 1.5){$i$}
\end{texdraw}}%
\savebox{\tmpfigg}{\begin{texdraw}
\fontsize{7}{7}\selectfont
\textref h:C v:C
\drawdim em
\setunitscale 1.9
\move(0 -0.2)\lvec(0 1)\lvec(1 1)\lvec(1 -0.2)
\move(0 0)\lvec(1 0)
\move(0 0)\lvec(1 1)
\move(0 0)\lvec(1 0)
\htext(0.25 0.75){$1$}
\end{texdraw}}%
\savebox{\tmpfigh}{\begin{texdraw}
\fontsize{7}{7}\selectfont
\textref h:C v:C
\drawdim em
\setunitscale 1.9
\move(0 -0.2)\lvec(0 1)\lvec(1 1)\lvec(1 -0.2)
\move(0 0)\lvec(1 0)
\move(0 0)\lvec(1 1)
\move(0 0)\lvec(1 0)
\htext(0.75 0.25){$1$}
\end{texdraw}}%
{\allowdisplaybreaks
\begin{align*}
\bu_1 \ &= \
 \raisebox{-0.5em}{\usebox{\tmpfiga}}
 \ \ \textup{or}\ \ \raisebox{-0.5em}{\usebox{\tmpfigf}}\\[0.3em]
\bu_i \ &= \
 \raisebox{-0.5em}{\usebox{\tmpfigb}}
 \qquad \text{ for $i=2,\cdots,n$, }\\[0.3em]
\bw_0 \ &= \
 \raisebox{-0.5em}{\usebox{\tmpfigc}}\\[0.3em]
\bv_n \ &= \
 \raisebox{-0.5em}{\usebox{\tmpfigd}}\\[0.3em]
\bv_i \ &= \
 \raisebox{-0.5em}{\usebox{\tmpfige}}
 \qquad \text{ for $i=2,\cdots,n$,}\\[0.3em]
\bv_1 \ &= \
 \raisebox{-0.5em}{\usebox{\tmpfigg}}
 \ \ \textup{or}\ \ \raisebox{-0.5em}{\usebox{\tmpfigh}}\\[0.3em]
\end{align*}}% end of allowdisplaybreaks
\end{center}
\end{enumerate}
The image of $b$ under the map $\psi$ is defined
to the equivalence class of a
level-$l$ slice obtained by pasting together $x'_0$-many $\bw_0$'s,
$x'_i$-many $\bu_i$'s, and $\bar{x}'_i$-many $\bv_i$'s, for each $i$.

We will use the notation $\psi(b) = [x'_1, \dots, x'_n | x'_0 | \bar x'_n, \dots,
\bar x'_1]$, as needed.

\vskip 3mm

When $\g = \anone$, $\atwontwo$, $\cnone$, $\dnponetwo$, it is
easy to see that the map $\psi$ is a bijection.

In the case $\g=\bnone$,
we define a new map $\phi: \newcrystal \rightarrow \oldcrystal$
to see that the map $\psi$ is a bijection.
Since any element of $\newcrystal$ may be obtained
by pasting together some number of $\bu_0$, $\bu_i$, and $\bv_i$,
we may denote an arbitrary element of $\newcrystal$ by
$C=[y_1, \dots, y_n | y_0 | \bar y_n, \dots, \bar y_1]$,
where $y_1 \bar{y}_1=0,\ y_0, y_i, \bar{y}_i \in \Z_{\geq0}, \text{ and }
y_0 + \textstyle\sum_{i=1}^{n} (y_i + \bar{y}_i) = \textnormal{$l$}$.
For such an element $C$,
we define $\phi(C)=(x_1,
\dots, x_n | x_0 | \bar x_n, \dots, \bar x_1)$, where
\begin{align*}
& x_1 = y_1 + \text{min} \{ y_2,\bar{y}_2 \},\\
& x_2 = (0,y_2-\bar{y}_2)_+ + \text{min} \{ y_3,\bar{y}_3 \},\\
& \qquad \qquad \vdots\\
& x_{n-1} = (0,y_{n-1}-\bar{y}_{n-1})_+
           + \text{min} \{ y_n,\bar{y}_n \},\\
& x_n = (0,y_{n}-\bar{y}_{n})_+
       + [\frac{y_0}{2}], \\
& x_0 = y_0-2[\frac{y_0}{2}],\\
& \bar{x}_n = (0,\bar{y}_{n}-y_{n})_+
       + [\frac{y_0}{2}], \\
& \bar{x}_{n-1} = (0,\bar{y}_{n-1}-y_{n-1})_+
           + \text{min} \{ y_n,\bar{y}_n \},\\
& \qquad \qquad \vdots\\
& \bar{x}_2 = (0,\bar{y}_2-y_2)_+ + \text{min} \{ y_3,\bar{y}_3 \},\\
& \bar{x}_1 = \bar{y}_1 + \text{min} \{ y_2,\bar{y}_2 \}.
\end{align*}
Similarly, in the case $\g=\atwonmonetwo$,
any element $C$ of $\newcrystal$ can be denoted by
$C=[y_1,
\dots, y_n | y_0 | \bar y_n, \dots, \bar y_1]$,
where
$y_1 \bar{y}_1=0,\ y_0\in 2\Z_{\ge 0},\ y_i, \bar{y}_i \in \Z_{\geq0},
\text{ and }
y_0+\textstyle\sum_{i=1}^{n} (y_i + \bar{y}_i) = \textnormal{$l$}$.
We define
$\phi(C)=(x_1, \dots, x_n | \bar x_n, \dots, \bar x_1)$,
where
\begin{align*}
& x_1 = y_1 + \text{min} \{ y_2,\bar{y}_2 \},\\
& x_2 = (0,y_2-\bar{y}_2)_+ + \text{min} \{ y_3,\bar{y}_3 \},\\
& \qquad \qquad \vdots\\
& x_{n-1} = (0,y_{n-1}-\bar{y}_{n-1})_+
           + \text{min} \{ y_n,\bar{y}_n \},\\
& x_n = (0,y_{n}-\bar{y}_{n})_+
       + \frac{y_0}{2}, \\
& \bar{x}_n = (0,\bar{y}_{n}-y_{n})_+
       + \frac{y_0}{2}, \\
& \bar{x}_{n-1} = (0,\bar{y}_{n-1}-y_{n-1})_+
           + \text{min} \{ y_n,\bar{y}_n \},\\
& \qquad \qquad \vdots\\
& \bar{x}_2 = (0,\bar{y}_2-y_2)_+ + \text{min} \{ y_3,\bar{y}_3 \},\\
& \bar{x}_1 = \bar{y}_1 + \text{min} \{ y_2,\bar{y}_2 \}.
\end{align*}
It is now easy to verify that the map $\phi$ is well-defined
and that it is actually the inverse of $\psi$.

It is
almost obvious that the maps $\cwt$, $\veps_i$, $\vphi_i$ are
preserved under $\psi$. It still remains to show that $\psi$
commutes with the Kashiwara operators. However, it is a lengthy
but a straightforward case-by-case verification. For instance, if
$\g=\bnone$ and $b=(x_1, \dots, x_n | x_0 | \bar x_n, \dots, \bar
x_1) \in \oldcrystal$ with $x_1 < \bar x_1$, $x_2 \ge \bar x_2$,
then we have
\begin{equation*}
\tilde f_0 b =(x_1, x_2 + 1, \dots, x_n | x_0 | \bar x_n, \dots,
\bar x_2, \bar x_1 -1),
\end{equation*}
and $\psi(b)=[y_1, \dots, y_n | y_0 | \bar y_n, \dots, \bar y_1]$,
where $y_1=0$, $y_2 = x_2 - \bar x_2 + x_1$, $\cdots$, $\bar y_2 =
x_1$, $\bar y_1 = \bar x_1 - x_1$. By applying $\tilde f_0$, we
get $\tilde f_0 \psi (b) = [y_1, y_2 +1, \dots, y_n | y_0 | \bar
y_n, \dots, \bar y_1 -1]$, which is the same as $\psi (\tilde f_0
b)$ as desired.

Therefore, we obtain a new realization of level-$l$ perfect
crystals as the sets of equivalence classes of level-$l$ slices.

\begin{thm}\label{thm:37}
The map $\psi: \oldcrystal \rightarrow \newcrystal$ defined above
is an isomorphism of $\uqp(\mathfrak{g})$-crystals.
\end{thm}

\vskip 3mm

\begin{example}
The following is a drawing of a portion of the level-$3$ perfect crystal
for type $B_3^{(1)}$.
The elements have been represented by elemnts of $\newcrystal$.
%We can compare this with Example~\ref{ex:02.3}.
%
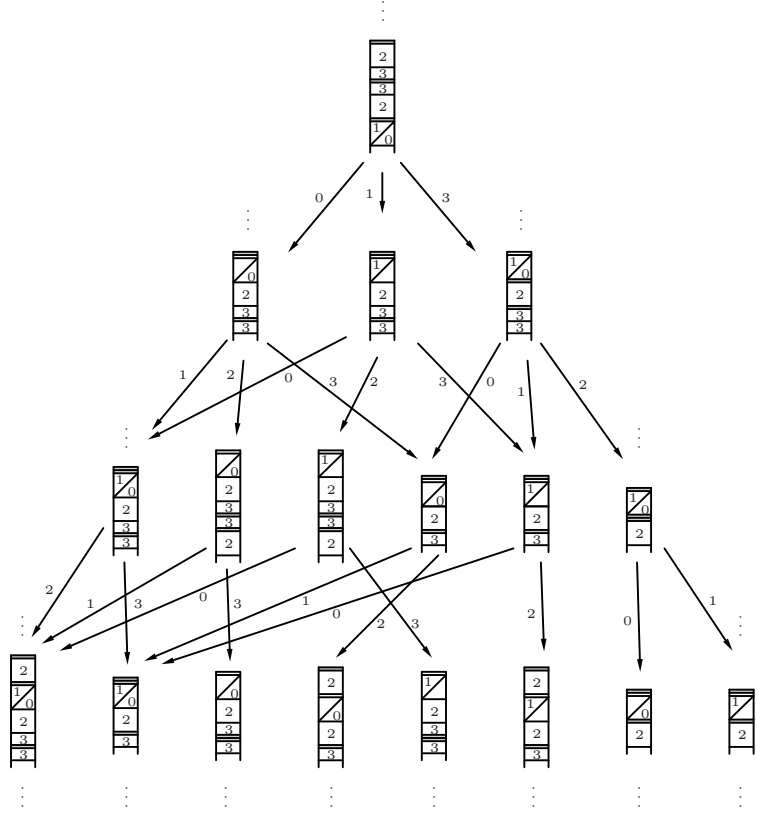
\vspace{3mm}
\begin{center}
\begin{texdraw}
\drawdim in
\arrowheadsize l:0.065 w:0.03
\arrowheadtype t:F
\fontsize{5}{5}\selectfont
\textref h:C v:C
\drawdim em
\setunitscale 1.9
\move(0 0)
\bsegment
\setsegscale 0.7
\move(1 1)\lvec(1 1.133)\lvec(0 1.133)\lvec(0 1)\ifill f:0.8
\move(1 2.766)\lvec(1 2.633)\lvec(0 2.633)\lvec(0 2.766)\ifill f:0.8
\move(1 4.266)\lvec(1 4.399)\lvec(0 4.399)\lvec(0 4.266)\ifill f:0.8
\move(0 -0.3)\lvec(0 4.399)
\move(1 -0.3)\lvec(1 4.399)
\move(0 0)\lvec(1 0)
\move(0 1)\lvec(1 1)
\move(1 1)\lvec(0 0)
\move(0 1.133)\lvec(1 1.133)
\move(0 2.133)\lvec(1 2.133)
\move(0 2.633)\lvec(1 2.633)
\move(0 2.766)\lvec(1 2.766)
\move(0 3.266)\lvec(1 3.266)
\move(0 4.266)\lvec(1 4.266)
\move(0 4.399)\lvec(1 4.399)
\htext(0.25 0.75){$1$}
\htext(0.75 0.25){$0$}
\htext(0.5 1.633){$2$}
\htext(0.5 2.383){$3$}
\htext(0.5 3.016){$3$}
\htext(0.5 3.766){$2$}
\esegment
\move(0 -6)
\bsegment
\move(-4 0.5)
\bsegment
\setsegscale 0.7
\move(1 0.5)\lvec(1 0.633)\lvec(0 0.633)\lvec(0 0.5)\ifill f:0.8
\move(1 3.266)\lvec(1 3.399)\lvec(0 3.399)\lvec(0 3.266)\ifill f:0.8
\move(0 -0.3)\lvec(0 3.399)
\move(1 -0.3)\lvec(1 3.399)
\move(0 0)\lvec(1 0)
\move(0 0.5)\lvec(1 0.5)
\move(0 0.633)\lvec(1 0.633)
\move(0 1.133)\lvec(1 1.133)
\move(0 2.133)\lvec(1 2.133)
\move(0 3.133)\lvec(1 3.133)
\move(0 3.266)\lvec(1 3.266)
\move(0 3.399)\lvec(1 3.399)
\move(0 2.133)\lvec(1 3.133)
\htext(0.5 0.25){$3$}
\htext(0.5 0.883){$3$}
\htext(0.5 1.633){$2$}
\htext(0.75 2.383){$0$}
\esegment
\move(4 0.5)
\bsegment
\setsegscale 0.7
\move(1 1)\lvec(1 1.133)\lvec(0 1.133)\lvec(0 1)\ifill f:0.8
\move(1 2.133)\lvec(1 2.266)\lvec(0 2.266)\lvec(0 2.133)\ifill f:0.8
\move(0 -0.3)\lvec(0 3.399)
\move(1 -0.3)\lvec(1 3.399)
\move(0 0)\lvec(1 0)
\move(0 0.5)\lvec(1 0.5)
\move(0 1)\lvec(1 1)
\move(0 1.133)\lvec(1 1.133)
\move(0 2.133)\lvec(1 2.133)
\move(0 2.266)\lvec(1 2.266)
\move(0 3.266)\lvec(1 3.266)
\move(0 3.399)\lvec(1 3.399)
\move(0 2.266)\lvec(1 3.266)
\htext(0.5 0.75){$3$}
\htext(0.5 0.25){$3$}
\htext(0.5 1.633){$2$}
\htext(0.75 2.516){$0$}
\htext(0.25 3.016){$1$}
\esegment
\move(0 0.5)
\bsegment
\setsegscale 0.7
\move(1 0.5)\lvec(1 0.633)\lvec(0 0.633)\lvec(0 0.5)\ifill f:0.8
\move(1 3.266)\lvec(1 3.399)\lvec(0 3.399)\lvec(0 3.266)\ifill f:0.8
\move(0 -0.3)\lvec(0 3.399)
\move(1 -0.3)\lvec(1 3.399)
\move(0 0)\lvec(1 0)
\move(0 0.5)\lvec(1 0.5)
\move(0 0.633)\lvec(1 0.633)
\move(0 1.133)\lvec(1 1.133)
\move(0 2.133)\lvec(1 2.133)
\move(0 3.133)\lvec(1 3.133)
\move(0 3.266)\lvec(1 3.266)
\move(0 3.399)\lvec(1 3.399)
\move(0 2.133)\lvec(1 3.133)
\htext(0.5 0.25){$3$}
\htext(0.5 0.883){$3$}
\htext(0.5 1.633){$2$}
\htext(0.25 2.883){$1$}
\esegment
\esegment
\move(0 -12)
\bsegment
\move(-7.5 0.2)
\bsegment
\setsegscale 0.7
\move(0 0.5)\lvec(1 0.5)\lvec(1 0.663)\lvec(0 0.663)\ifill f:0.8
\move(0 3.133)\lvec(1 3.133)\lvec(1 3.399)\lvec(0 3.399)\ifill f:0.8
\move(0 -0.3)\lvec(0 3.399)
\move(1 -0.3)\lvec(1 3.399)
\move(0 0)\lvec(1 0)
\move(0 0.5)\lvec(1 0.5)
\move(0 0.633)\lvec(1 0.633)
\move(0 1.133)\lvec(1 1.133)
\move(0 2.133)\lvec(1 2.133)
\move(0 3.133)\lvec(1 3.133)
\move(0 3.266)\lvec(1 3.266)
\move(0 3.399)\lvec(1 3.399)
\move(0 2.133)\lvec(1 3.133)
\htext(0.5 0.25){$3$}
\htext(0.5 0.883){$3$}
\htext(0.5 1.633){$2$}
\htext(0.75 2.383){$0$}
\htext(0.25 2.883){$1$}
\esegment
\move(-4.5 0)
\bsegment
\setsegscale 0.7
\move(1 1)\lvec(1 1.133)\lvec(0 1.133)\lvec(0 1)\ifill f:0.8
\move(1 1.766)\lvec(1 1.633)\lvec(0 1.633)\lvec(0 1.766)\ifill f:0.8
\move(0 -0.3)\lvec(0 4.399)
\move(1 -0.3)\lvec(1 4.399)
\move(0 0)\lvec(1 0)
\move(0 1)\lvec(1 1)
\move(0 1.133)\lvec(1 1.133)
\move(0 1.633)\lvec(1 1.633)
\move(0 1.766)\lvec(1 1.766)
\move(0 2.266)\lvec(1 2.266)
\move(0 3.266)\lvec(1 3.266)
\move(0 4.266)\lvec(1 4.266)
\move(0 4.399)\lvec(1 4.399)
\move(0 3.266)\lvec(1 4.266)
\htext(0.5 0.5){$2$}
\htext(0.5 1.383){$3$}
\htext(0.5 2.016){$3$}
\htext(0.5 2.766){$2$}
\htext(0.75 3.516){$0$}
\esegment
\move(-1.5 0)
\bsegment
\setsegscale 0.7
\move(1 1)\lvec(1 1.133)\lvec(0 1.133)\lvec(0 1)\ifill f:0.8
\move(1 1.766)\lvec(1 1.633)\lvec(0 1.633)\lvec(0 1.766)\ifill f:0.8
\move(0 -0.3)\lvec(0 4.399)
\move(1 -0.3)\lvec(1 4.399)
\move(0 0)\lvec(1 0)
\move(0 1)\lvec(1 1)
\move(0 1.133)\lvec(1 1.133)
\move(0 1.633)\lvec(1 1.633)
\move(0 1.766)\lvec(1 1.766)
\move(0 2.266)\lvec(1 2.266)
\move(0 3.266)\lvec(1 3.266)
\move(0 4.266)\lvec(1 4.266)
\move(0 4.399)\lvec(1 4.399)
\move(0 3.266)\lvec(1 4.266)
\htext(0.5 0.5){$2$}
\htext(0.5 1.383){$3$}
\htext(0.5 2.016){$3$}
\htext(0.5 2.766){$2$}
\htext(0.25 4.016){$1$}
\esegment
\move(1.5 0.3)
\bsegment
\setsegscale 0.7
\move(1 2.766)\lvec(1 2.899)\lvec(0 2.899)\lvec(0 2.766)\ifill f:0.8
\move(1 0.5)\lvec(1 0.633)\lvec(0 0.633)\lvec(0 0.5)\ifill f:0.8
\move(0 -0.3)\lvec(0 2.899)
\move(1 -0.3)\lvec(1 2.899)
\move(0 0)\lvec(1 0)
\move(0 0.5)\lvec(1 0.5)
\move(0 0.633)\lvec(1 0.633)
\move(0 1.633)\lvec(1 1.633)
\move(0 2.633)\lvec(1 2.633)
\move(0 2.766)\lvec(1 2.766)
\move(0 2.899)\lvec(1 2.899)
\move(0 1.633)\lvec(1 1.633)
\move(0 1.633)\lvec(1 2.633)
\htext(0.5 1.133){$2$}
\htext(0.5 0.25){$3$}
\htext(0.75 1.883){$0$}
\esegment
\move(4.5 0.3)
\bsegment
\setsegscale 0.7
\move(1 2.766)\lvec(1 2.899)\lvec(0 2.899)\lvec(0 2.766)\ifill f:0.8
\move(1 0.5)\lvec(1 0.633)\lvec(0 0.633)\lvec(0 0.5)\ifill f:0.8
\move(0 -0.3)\lvec(0 2.899)
\move(1 -0.3)\lvec(1 2.899)
\move(0 0)\lvec(1 0)
\move(0 0.5)\lvec(1 0.5)
\move(0 0.633)\lvec(1 0.633)
\move(0 1.633)\lvec(1 1.633)
\move(0 2.633)\lvec(1 2.633)
\move(0 2.766)\lvec(1 2.766)
\move(0 2.899)\lvec(1 2.899)
\move(0 1.633)\lvec(1 1.633)
\move(0 1.633)\lvec(1 2.633)
\htext(0.5 1.133){$2$}
\htext(0.5 0.25){$3$}
\htext(0.25 2.383){$1$}
\esegment
\move(7.5 0.3)
\bsegment
\setsegscale 0.7
\move(1 2.266)\lvec(1 2.399)\lvec(0 2.399)\lvec(0 2.266)\ifill f:0.8
\move(1 1)\lvec(1 1.266)\lvec(0 1.266)\lvec(0 1)\ifill f:0.8
\move(0 -0.3)\lvec(0 2.399)
\move(1 -0.3)\lvec(1 2.399)
\move(0 0)\lvec(1 0)
\move(0 1)\lvec(1 1)
\move(0 1.133)\lvec(1 1.133)
\move(0 1.266)\lvec(1 1.266)
\move(0 2.266)\lvec(1 2.266)
\move(0 2.399)\lvec(1 2.399)
\move(0 1.266)\lvec(1 2.266)
\htext(0.5 0.5){$2$}
\htext(0.25 2.016){$1$}
\htext(0.75 1.516){$0$}
\esegment
\esegment
\move(0 -18)
\bsegment
\move(-10.5 0)
\bsegment
\setsegscale 0.7
\move(0 0.5)\lvec(1 0.5)\lvec(1 0.663)\lvec(0 0.663)\ifill f:0.8
\move(0 3.133)\lvec(1 3.133)\lvec(1 3.266)\lvec(0 3.266)\ifill f:0.8
\move(0 4.266)\lvec(1 4.266)\lvec(1 4.399)\lvec(0 4.399)\ifill f:0.8
\move(0 -0.3)\lvec(0 4.399)
\move(1 -0.3)\lvec(1 4.399)
\move(0 0)\lvec(1 0)
\move(0 0.5)\lvec(1 0.5)
\move(0 0.633)\lvec(1 0.633)
\move(0 1.133)\lvec(1 1.133)
\move(0 2.133)\lvec(1 2.133)
\move(0 3.133)\lvec(1 3.133)
\move(0 3.266)\lvec(1 3.266)
\move(0 4.266)\lvec(1 4.266)
\move(0 4.399)\lvec(1 4.399)
\move(0 2.133)\lvec(1 3.133)
\htext(0.5 0.25){$3$}
\htext(0.5 0.883){$3$}
\htext(0.5 1.633){$2$}
\htext(0.75 2.383){$0$}
\htext(0.25 2.883){$1$}
\htext(0.5 3.766){$2$}
\esegment
\move(-7.5 0.4)
\bsegment
\setsegscale 0.7
\move(0 0.5)\lvec(1 0.5)\lvec(1 0.663)\lvec(0 0.663)\ifill f:0.8
\move(0 2.633)\lvec(1 2.633)\lvec(1 2.899)\lvec(0 2.899)\ifill f:0.8
\move(0 -0.3)\lvec(0 2.899)
\move(1 -0.3)\lvec(1 2.899)
\move(0 0)\lvec(1 0)
\move(0 0.5)\lvec(1 0.5)
\move(0 0.633)\lvec(1 0.633)
\move(0 1.633)\lvec(1 1.633)
\move(0 2.633)\lvec(1 2.633)
\move(0 2.766)\lvec(1 2.766)
\move(0 2.899)\lvec(1 2.899)
\move(0 1.663)\lvec(1 2.663)
\htext(0.5 0.25){$3$}
\htext(0.5 1.133){$2$}
\htext(0.75 1.883){$0$}
\htext(0.25 2.383){$1$}
\esegment
%\move(-1.5 0)
%\bsegment
%\setsegscale 0.7
%\move(1 1)\lvec(1 1.133)\lvec(0 1.133)\lvec(0 1)\ifill f:0.8
%\move(1 1.766)\lvec(1 1.633)\lvec(0 1.633)\lvec(0 1.766)\ifill f:0.8
%\move(1 4.266)\lvec(1 4.333)\lvec(0 4.333)\lvec(0 4.266)\ifill f:0.8
%\move(0 -0.3)\lvec(0 4.399)
%\move(1 -0.3)\lvec(1 4.399)
%\move(0 0)\lvec(1 0)
%\move(0 1)\lvec(1 1)
%\move(0 1.133)\lvec(1 1.133)
%\move(0 1.633)\lvec(1 1.633)
%\move(0 1.766)\lvec(1 1.766)
%\move(0 2.266)\lvec(1 2.266)
%\move(0 3.266)\lvec(1 3.266)
%\move(0 4.266)\lvec(1 4.266)
%\move(0 4.399)\lvec(1 4.399)
%\move(0 3.266)\lvec(1 4.266)
%\htext(0.5 0.5){$2$}
%\htext(0.5 1.383){$3$}
%\htext(0.5 2.016){$3$}
%\htext(0.5 2.766){$2$}
%\htext(0.75 3.516){$0$}
%\htext(0.25 4.016){$1$}
%\esegment
\move(-4.5 0.2)
\bsegment
\setsegscale 0.7
\move(1 0.5)\lvec(1 0.766)\lvec(0 0.766)\lvec(0 0.5)\ifill f:0.8
\move(0 -0.3)\lvec(0 3.399)
\move(1 -0.3)\lvec(1 3.399)
\move(0 0)\lvec(1 0)
\move(0 0.5)\lvec(1 0.5)
\move(0 0.633)\lvec(1 0.633)
\move(0 0.766)\lvec(1 0.766)
\move(0 1.266)\lvec(1 1.266)
\move(0 2.266)\lvec(1 2.266)
\move(0 3.266)\lvec(1 3.266)
\move(0 3.399)\lvec(1 3.399)
\move(0 2.266)\lvec(1 3.266)
\htext(0.5 0.25){$3$}
\htext(0.5 1.016){$3$}
\htext(0.5 1.766){$2$}
\htext(0.75 2.516){$0$}
\esegment
\move(1.5 0.2)
\bsegment
\setsegscale 0.7
\move(1 0.5)\lvec(1 0.766)\lvec(0 0.766)\lvec(0 0.5)\ifill f:0.8
\move(0 -0.3)\lvec(0 3.399)
\move(1 -0.3)\lvec(1 3.399)
\move(0 0)\lvec(1 0)
\move(0 0.5)\lvec(1 0.5)
\move(0 0.633)\lvec(1 0.633)
\move(0 0.766)\lvec(1 0.766)
\move(0 1.266)\lvec(1 1.266)
\move(0 2.266)\lvec(1 2.266)
\move(0 3.266)\lvec(1 3.266)
\move(0 3.399)\lvec(1 3.399)
\move(0 2.266)\lvec(1 3.266)
\htext(0.5 0.25){$3$}
\htext(0.5 1.016){$3$}
\htext(0.5 1.766){$2$}
\htext(0.25 3.016){$1$}
\esegment
\move(4.5 0)
\bsegment
\setsegscale 0.7
\move(1 3.766)\lvec(1 3.899)\lvec(0 3.899)\lvec(0 3.766)\ifill f:0.8
\move(1 0.5)\lvec(1 0.633)\lvec(0 0.633)\lvec(0 0.5)\ifill f:0.8
\move(0 -0.3)\lvec(0 3.899)
\move(1 -0.3)\lvec(1 3.899)
\move(0 0)\lvec(1 0)
\move(0 0.5)\lvec(1 0.5)
\move(0 0.633)\lvec(1 0.633)
\move(0 1.633)\lvec(1 1.633)
\move(0 2.633)\lvec(1 2.633)
\move(0 2.766)\lvec(1 2.766)
\move(0 3.766)\lvec(1 3.766)
\move(0 3.899)\lvec(1 3.899)
\move(0 1.633)\lvec(1 2.633)
\htext(0.5 1.133){$2$}
\htext(0.5 0.25){$3$}
\htext(0.25 2.383){$1$}
\htext(0.5 3.266){$2$}
\esegment
\move(-1.5 0)
\bsegment
\setsegscale 0.7
\move(1 3.766)\lvec(1 3.899)\lvec(0 3.899)\lvec(0 3.766)\ifill f:0.8
\move(1 0.5)\lvec(1 0.633)\lvec(0 0.633)\lvec(0 0.5)\ifill f:0.8
\move(0 -0.3)\lvec(0 3.899)
\move(1 -0.3)\lvec(1 3.899)
\move(0 0)\lvec(1 0)
\move(0 0.5)\lvec(1 0.5)
\move(0 0.633)\lvec(1 0.633)
\move(0 1.633)\lvec(1 1.633)
\move(0 2.633)\lvec(1 2.633)
\move(0 2.766)\lvec(1 2.766)
\move(0 3.766)\lvec(1 3.766)
\move(0 3.899)\lvec(1 3.899)
\move(0 1.633)\lvec(1 2.633)
\htext(0.5 1.133){$2$}
\htext(0.5 0.25){$3$}
\htext(0.75 1.883){$0$}
\htext(0.5 3.266){$2$}
\esegment
\move(7.5 0.4)
\bsegment
\setsegscale 0.7
\move(1 2.266)\lvec(1 2.399)\lvec(0 2.399)\lvec(0 2.266)\ifill f:0.8
\move(1 1)\lvec(1 1.133)\lvec(0 1.133)\lvec(0 1)\ifill f:0.8
\move(0 -0.3)\lvec(0 2.399)
\move(1 -0.3)\lvec(1 2.399)
\move(0 0)\lvec(1 0)
\move(0 1)\lvec(1 1)
\move(0 1.133)\lvec(1 1.133)
\move(0 2.133)\lvec(1 2.133)
\move(0 2.266)\lvec(1 2.266)
\move(0 2.399)\lvec(1 2.399)
\move(0 1.133)\lvec(1 2.133)
\htext(0.5 0.5){$2$}
\htext(0.75 1.383){$0$}
\esegment
\move(10.5 0.4)
\bsegment
\setsegscale 0.7
\move(1 2.266)\lvec(1 2.399)\lvec(0 2.399)\lvec(0 2.266)\ifill f:0.8
\move(1 1)\lvec(1 1.133)\lvec(0 1.133)\lvec(0 1)\ifill f:0.8
\move(0 -0.3)\lvec(0 2.399)
\move(1 -0.3)\lvec(1 2.399)
\move(0 0)\lvec(1 0)
\move(0 1)\lvec(1 1)
\move(0 1.133)\lvec(1 1.133)
\move(0 2.133)\lvec(1 2.133)
\move(0 2.266)\lvec(1 2.266)
\move(0 2.399)\lvec(1 2.399)
\move(0 1.133)\lvec(1 2.133)
\htext(0.5 0.5){$2$}
\htext(0.25 1.883){$1$}
\esegment
\esegment
\move(-0.2 -0.5)\ravec(-2.2 -2.62)
\move(0.35 -0.8)\ravec(0 -1.2)
\move(0.9 -0.5)\ravec(2.2 -2.6)
\htext(-1.5 -1.5){$0$}
\htext(-0.03 -1.4){$1$}
\htext(2.2 -1.5){$3$}
\move(-4.2 -5.7)\ravec(-2.2 -2.62)
\move(-3.7 -6.3)\ravec(-0.2 -2.2)
\move(-3 -5.8)\ravec(4.4 -3.4)
\move(-0.7 -5.6)\ravec(-5.8 -3)
\move(0.2 -6.2)\ravec(-1.1 -2.2)
\move(1.4 -5.8)\ravec(3.1 -3.2)
\move(3.8 -5.8)\ravec(-2 -3.4)
\move(4.6 -6.3)\ravec(0.2 -2.6)
\move(5 -5.8)\ravec(2.4 -3.4)
\htext(-5.5 -6.7){$1$}
\htext(-4.1 -6.7){$2$}
\htext(-1.1 -6.9){$3$}
\htext(-2.4 -6.8){$0$}
\htext(0.1 -6.9){$2$}
\htext(2.1 -6.9){$3$}
\htext(3.5 -6.9){$0$}
\htext(4.4 -7.2){$1$}
\htext(6.2 -7){$2$}
\move(-7.8 -11.2)\ravec(-2.1 -3.2)
\move(-7.2 -12.2)\ravec(0.1 -3)
\move(-4.8 -11.8)\ravec(-4.8 -2.8)
\move(-4.2 -12.4)\ravec(0.1 -2.7)
\move(-2.2 -11.8)\ravec(-6.9 -3)
\move(-0.6 -11.8)\ravec(2.4 -3.2)
\move(1.2 -11.8)\ravec(-7.8 -3.3)
\move(2 -12)\ravec(-3 -3)
\move(4.2 -11.8)\ravec(-10.3 -3.4)
\move(5 -12.2)\ravec(0.1 -2.6)
\move(7.8 -12.2)\ravec(0.1 -3.2)
\move(8.6 -11.8)\ravec(2 -3)
\htext(-9.4 -13){$2$}
\htext(-6.8 -13.4){$3$}
\htext(-8.2 -13.4){$1$}
\htext(-3.9 -13.4){$3$}
\htext(-4.9 -13.2){$0$}
\htext(1.3 -14){$3$}
\htext(-1.9 -13.3){$1$}
\htext(0.3 -14){$2$}
\htext(-1 -13.7){$0$}
\htext(4.7 -13.7){$2$}
\htext(7.5 -13.9){$0$}
\htext(10 -13.3){$1$}
\vtext(0.35 4){$\cdots$}
\vtext(-3.6 -2.1){$\cdots$}
\vtext(4.4 -2.1){$\cdots$}
\vtext(-7.15 -8.5){$\cdots$}
\vtext(7.85 -8.5){$\cdots$}
\vtext(-10.2 -14){$\cdots$}
\vtext(10.8 -14){$\cdots$}
\vtext(-10.2 -19){$\cdots$}
\vtext(-7.15 -19){$\cdots$}
\vtext(-4.25 -19){$\cdots$}
\vtext(-1.15 -19){$\cdots$}
\vtext(1.85 -19){$\cdots$}
\vtext(4.75 -19){$\cdots$}
\vtext(7.85 -19){$\cdots$}
\vtext(10.8 -19){$\cdots$}
\end{texdraw}
\end{center}
\end{example}

\vspace{2mm}

\vskip 1cm

%%%%%%%%%%%%%%%%%%%%%%%%%%%%%%%%%%%%%%%%%%%%%%%%%%%%%%%%%%%%%%%
%%%%%%%%%%%%%%%%%%%%%%%%%%%%%%%%%%%%%%%%%%%%%%%%%%%%%%%%%%%%%%%
\section{Combinatorics of higher level Young walls}

In this section, we introduce the notion of \emph{higher level
Young walls}. Roughly speaking, the level-$l$ Young walls are
constructed by lining up level-$l$ slices defined in the previous
section, and can be viewed as the $l$-tuples of level-1 Young
walls introduced in \cite{Kang2000} and \cite{CnOne}. The patterns
for building Young walls are given below.\\

\savebox{\tmpfiga}{\begin{texdraw} \fontsize{7}{7}\selectfont
\textref h:C v:C \drawdim em \setunitscale 1.9 \move(0
0)\rlvec(-5.8 0) \move(0 1)\rlvec(-5.8 0) \move(0 2)\rlvec(-5.8 0)
\move(0 3.5)\rlvec(-5.8 0) \move(0 4.5)\rlvec(-5.8 0) \move(0
5.5)\rlvec(-5.8 0) \move(0 6.5)\rlvec(-5.8 0) \move(0 0)\rlvec(0
6.8) \move(-1 0)\rlvec(0 6.8) \move(-2 0)\rlvec(0 6.8) \move(-3.5
0)\rlvec(0 6.8) \move(-4.5 0)\rlvec(0 6.8) \move(-5.5 0)\rlvec(0
6.8) \move(-0.5 0.5) \bsegment \htext(0 0){$0$} \htext(0 1){$1$}
\vtext(0 2.25){$\cdots$} \htext(0 3.5){$n$} \htext(0 4.5){$0$}
\htext(0 5.5){$1$} \htext(-1 0){$n$} \htext(-1 1){$0$} \vtext(-1
2.25){$\cdots$} \htext(-1 3.5){$n\!\!-\!\!1$} \htext(-1 4.5){$n$}
\htext(-1 5.5){$0$} \htext(-2.25 0){$\cdots$} \htext(-2.25
1){$\cdots$} \htext(-3.5 0){$1$} \htext(-3.5 1){$2$} \htext(-4.5
0){$0$} \htext(-4.5 1){$1$} \esegment \move(0 -2.5)
\end{texdraw}}%
\savebox{\tmpfigb}{\begin{texdraw}
\fontsize{7}{7}\selectfont
\textref h:C v:C
\drawdim em
\setunitscale 1.9
%\move(0 0)\lvec(-4 0)\lvec(-4 0.5)\lvec(0 0.5)\ifill f:0.7
\move(0 0)\rlvec(-4.3 0) \move(0 0.5)\rlvec(-4.3 0) \move(0 1)\rlvec(-4.3 0)
\move(0 2)\rlvec(-4.3 0) \move(0 3.5)\rlvec(-4.3 0) \move(0 4.5)\rlvec(-4.3 0)
\move(0 6)\rlvec(-4.3 0) \move(0 7)\rlvec(-4.3 0) \move(0 7.5)\rlvec(-4.3 0)
\move(0 8)\rlvec(-4.3 0) \move(0 9)\rlvec(-4.3 0) \move(0 0)\rlvec(0 9.3)
\move(-1 0)\rlvec(0 9.3) \move(-2 0)\rlvec(0 9.3) \move(-3 0)\rlvec(0 9.3)
\move(-4 0)\rlvec(0 9.3)
\htext(-0.5 0.25){$0$} \htext(-1.5 0.25){$0$} \htext(-2.5 0.25){$0$}
\htext(-3.5 0.25){$0$} \htext(-0.5 0.75){$0$} \htext(-1.5 0.75){$0$}
\htext(-2.5 0.75){$0$} \htext(-3.5 0.75){$0$}
\htext(-0.5 1.5){$1$} \htext(-1.5 1.5){$1$} \htext(-2.5 1.5){$1$}
\htext(-3.5 1.5){$1$}
\htext(-0.5 4){$n$} \htext(-1.5 4){$n$} \htext(-2.5 4){$n$}
\htext(-3.5 4){$n$}
\htext(-0.5 6.5){$1$} \htext(-1.5 6.5){$1$} \htext(-2.5 6.5){$1$}
\htext(-3.5 6.5){$1$}
\htext(-0.5 7.25){$0$} \htext(-1.5 7.25){$0$} \htext(-2.5 7.25){$0$}
\htext(-3.5 7.25){$0$} \htext(-0.5 7.75){$0$} \htext(-1.5 7.75){$0$}
\htext(-2.5 7.75){$0$} \htext(-3.5 7.75){$0$}
\htext(-0.5 8.5){$1$} \htext(-1.5 8.5){$1$} \htext(-2.5 8.5){$1$}
\htext(-3.5 8.5){$1$}
\vtext(-0.5 2.75){$\cdots$} \vtext(-1.5 2.75){$\cdots$}
\vtext(-2.5 2.75){$\cdots$} \vtext(-3.5 2.75){$\cdots$}
\vtext(-0.5 5.25){$\cdots$} \vtext(-1.5 5.25){$\cdots$}
\vtext(-2.5 5.25){$\cdots$} \vtext(-3.5 5.25){$\cdots$}
\end{texdraw}}%
\savebox{\tmpfigc}{\begin{texdraw}
\fontsize{7}{7}\selectfont
\textref h:C v:C
\drawdim em
\setunitscale 1.9
%\move(0 0)\lvec(-4 0)\lvec(-4 0.5)\lvec(0 0.5)\ifill f:0.7
\move(0 0)\rlvec(-4.3 0) \move(0 0.5)\rlvec(-4.3 0) \move(0 1)\rlvec(-4.3 0)
\move(0 2)\rlvec(-4.3 0) \move(0 3.5)\rlvec(-4.3 0) \move(0 4.5)\rlvec(-4.3 0)
\move(0 4.5)\rlvec(-4.3 0) \move(0 5)\rlvec(-4.3 0) \move(0 5.5)\rlvec(-4.3 0)
\move(0 6.5)\rlvec(-4.3 0) \move(0 8)\rlvec(-4.3 0) \move(0 9)\rlvec(-4.3 0)
\move(0 9.5)\rlvec(-4.3 0) \move(0 10)\rlvec(-4.3 0) \move(0 11)\rlvec(-4.3 0)
\move(0 0)\rlvec(0 11.3) \move(-1 0)\rlvec(0 11.3) \move(-2 0)\rlvec(0 11.3)
\move(-3 0)\rlvec(0 11.3) \move(-4 0)\rlvec(0 11.3)
\htext(-0.5 0.25){$0$} \htext(-1.5 0.25){$0$} \htext(-2.5 0.25){$0$}
\htext(-3.5 0.25){$0$} \htext(-0.5 0.75){$0$} \htext(-1.5 0.75){$0$}
\htext(-2.5 0.75){$0$} \htext(-3.5 0.75){$0$}
\htext(-0.5 1.5){$1$} \htext(-1.5 1.5){$1$} \htext(-2.5 1.5){$1$}
\htext(-3.5 1.5){$1$}
\vtext(-0.5 2.75){$\cdots$} \vtext(-1.5 2.75){$\cdots$}
\vtext(-2.5 2.75){$\cdots$} \vtext(-3.5 2.75){$\cdots$}
\htext(-0.5 4){$n\!\!-\!\!1$} \htext(-1.5 4){$n\!\!-\!\!1$}
\htext(-2.5 4){$n\!\!-\!\!1$} \htext(-3.5 4){$n\!\!-\!\!1$}
\htext(-0.5 4.75){$n$} \htext(-1.5 4.75){$n$} \htext(-2.5 4.75){$n$}
\htext(-3.5 4.75){$n$} \htext(-0.5 5.25){$n$} \htext(-1.5 5.25){$n$}
\htext(-2.5 5.25){$n$} \htext(-3.5 5.25){$n$}
\htext(-0.5 6){$n\!\!-\!\!1$} \htext(-1.5 6){$n\!\!-\!\!1$}
\htext(-2.5 6){$n\!\!-\!\!1$} \htext(-3.5 6){$n\!\!-\!\!1$}
\vtext(-0.5 7.25){$\cdots$} \vtext(-1.5 7.25){$\cdots$}
\vtext(-2.5 7.25){$\cdots$} \vtext(-3.5 7.25){$\cdots$}
\htext(-0.5 8.5){$1$} \htext(-1.5 8.5){$1$} \htext(-2.5 8.5){$1$}
\htext(-3.5 8.5){$1$}
\htext(-0.5 9.25){$0$} \htext(-1.5 9.25){$0$} \htext(-2.5 9.25){$0$}
\htext(-3.5 9.25){$0$} \htext(-0.5 9.75){$0$} \htext(-1.5 9.75){$0$}
\htext(-2.5 9.75){$0$} \htext(-3.5 9.75){$0$}
\htext(-0.5 10.5){$1$} \htext(-1.5 10.5){$1$} \htext(-2.5 10.5){$1$}
\htext(-3.5 10.5){$1$}
\end{texdraw}}%
\savebox{\tmpfigd}{\begin{texdraw}
\fontsize{7}{7}\selectfont
\textref h:C v:C
\drawdim em
\setunitscale 1.9
\nc{\dtri}{
\bsegment
%\move(-1 0)\lvec(0 1)\lvec(0 0)\lvec(-1 0)\ifill f:0.7
\esegment
}
\move(0 0)\dtri \move(-1 0)\dtri \move(-2 0)\dtri \move(-3 0)\dtri
\move(0 0)\rlvec(-4.3 0) \move(0 1)\rlvec(-4.3 0) \move(0 2)\rlvec(-4.3 0)
\move(0 3.5)\rlvec(-4.3 0) \move(0 4.5)\rlvec(-4.3 0)
\move(0 5.5)\rlvec(-4.3 0) \move(0 6.5)\rlvec(-4.3 0)
\move(0 8)\rlvec(-4.3 0) \move(0 9)\rlvec(-4.3 0)
\move(0 10)\rlvec(-4.3 0) \move(0 11)\rlvec(-4.3 0)
\move(0 0)\rlvec(0 11.3) \move(-1 0)\rlvec(0 11.3) \move(-2 0)\rlvec(0 11.3)
\move(-3 0)\rlvec(0 11.3) \move(-4 0)\rlvec(0 11.3)
\move(-1 0)\rlvec(1 1) \move(-2 0)\rlvec(1 1) \move(-3 0)\rlvec(1 1)
\move(-4 0)\rlvec(1 1) \move(-1 9)\rlvec(1 1) \move(-2 9)\rlvec(1 1)
\move(-3 9)\rlvec(1 1) \move(-4 9)\rlvec(1 1)
\move(0 5)\rlvec(-4.3 0)
\htext(-0.3 0.25){$0$} \htext(-0.75 0.75){$1$} \htext(-0.5 1.5){$2$}
\vtext(-0.5 2.75){$\cdots$} \htext(-0.5 4){$n\!\!-\!\!1$}
\htext(-0.5 6){$n\!\!-\!\!1$} \htext(-0.5 8.5){$2$}
\htext(-0.3 9.25){$0$} \htext(-0.75 9.75){$1$} \htext(-0.5 10.5){$2$}
\htext(-2.3 0.25){$0$} \htext(-2.75 0.75){$1$} \htext(-2.5 1.5){$2$}
\vtext(-2.5 2.75){$\cdots$} \htext(-2.5 4){$n\!\!-\!\!1$}
\htext(-2.5 6){$n\!\!-\!\!1$} \htext(-2.5 8.5){$2$}
\htext(-2.3 9.25){$0$} \htext(-2.75 9.75){$1$} \htext(-2.5 10.5){$2$}
\htext(-1.3 0.25){$1$} \htext(-1.75 0.75){$0$} \htext(-1.5 1.5){$2$}
\vtext(-1.5 2.75){$\cdots$} \htext(-1.5 4){$n\!\!-\!\!1$}
\htext(-1.5 6){$n\!\!-\!\!1$} \htext(-1.5 8.5){$2$}
\htext(-1.3 9.25){$1$} \htext(-1.75 9.75){$0$} \htext(-1.5 10.5){$2$}
\htext(-3.3 0.25){$1$} \htext(-3.75 0.75){$0$} \htext(-3.5 1.5){$2$}
\vtext(-3.5 2.75){$\cdots$} \htext(-3.5 4){$n\!\!-\!\!1$}
\htext(-3.5 6){$n\!\!-\!\!1$} \htext(-3.5 8.5){$2$}
\htext(-3.3 9.25){$1$} \htext(-3.75 9.75){$0$} \htext(-3.5 10.5){$2$}
\htext(-0.5 4.75){$n$} \htext(-2.5 4.75){$n$} \htext(-1.5 5.25){$n$}
\htext(-3.5 5.25){$n$} \htext(-1.5 4.75){$n$} \htext(-3.5 4.75){$n$}
\htext(-0.5 5.25){$n$} \htext(-2.5 5.25){$n$}
\vtext(-0.5 7.25){$\cdots$} \vtext(-1.5 7.25){$\cdots$}
\vtext(-2.5 7.25){$\cdots$} \vtext(-3.5 7.25){$\cdots$}
\end{texdraw}}%
%

%\item $A_n^{(1)}$ ($n\geq1$)
%      \raisebox{-\height}{\usebox{\tmpfiga}}
%\item $A_{2n}^{(2)}$ ($n\geq1$) and $\cnone$ ($n\geq2$)
%      \raisebox{-\height}{\usebox{\tmpfigb}}
%\item $D_{n+1}^{(2)}$ ($n\geq2$)
%      \raisebox{-\height}{\usebox{\tmpfigc}}
%\item $B_{n}^{(1)}$ ($n\geq3$) and $\atwonmonetwo$ ($n\geq3$)
%      \raisebox{-\height}{\usebox{\tmpfigd}}
%      \begin{center}
%      $D_{n+1}^{(2)}$ ($n\geq2$) : \raisebox{-\height}{\usebox{\tmpfigc}}
%      \quad
%      $B_{n}^{(1)}$ ($n\geq3$), $\atwonmonetwo$ ($n\geq3$)
%                                 : \raisebox{-\height}{\usebox{\tmpfigd}}
%      \end{center}

%
\begin{enumerate}
\item $A_n^{(1)}$($n\geq1$)\\[3mm]
  \mbox{} \hspace{22mm} \usebox{\tmpfiga}\\
  \vskip 2mm
\item $B_{n}^{(1)}$($n\geq3$) and $\atwonmonetwo$($n\geq3$)\\[3mm]
  \mbox{} \hspace{30mm} \usebox{\tmpfigd}\\
  \vskip 2mm
\item $\cnone$($n\geq2$) and $A_{2n}^{(2)}$($n\geq1$)\\[3mm]
  \mbox{} \hspace{30mm} \usebox{\tmpfigb}\\
  \vskip 2mm
\item $D_{n+1}^{(2)}$($n\geq2$)\\[3mm]
  \mbox{} \hspace{30mm} \usebox{\tmpfigc}
\end{enumerate}

\vskip 3mm

\begin{df}
A \defi{level-$1$ Young wall of type $\anone$, $\bnone$,
$\atwonmonetwo$, $\atwontwo$, $\dnponetwo$} (or a
\defi{level-$\half$ Young wall of type $\cnone$}) is a set of
blocks stacked in a wall of unit depth satisfying the following
conditions\,:
\begin{enumerate}
\item The colored blocks are stacked in the pattern given above.
\item Except for the rightmost column, there is no free
      space to the right of any block.
\item No block can be placed on top of a column of half-unit depth.
\end{enumerate}
\end{df}

Note that every column of a level-$1$ Young wall (resp. a
level-$\half$ Young wall for $\cnone$)is a level-$1$ slice (resp.
level-$\half$ slice for $\cnone$) defined in the previous section.
For a level-$1$ Young wall (or a level-$\half$ Young wall for
$\cnone$) $\mathbf{Y}$, we define $\mathbf{Y}+\delta$ (resp.
$\mathbf{Y}-\delta$) to be the Young wall obtained by \emph{adding
a $\delta$} to (resp. \emph{removing a $\delta$} from) each and
every column of $\mathbf{Y}$.

\begin{df}\hfill
\begin{enumerate}
\item A \defi{level-$l$ Young wall of type $\anone$, $\bnone$,
$\atwontwo$} is an ordered $l$-tuple of level-1 Young walls
$\mathbf{Y} = (Y_1, \dots, Y_l)$ such that
\begin{itemize}
\item $Y_1 \subset \cdots \subset Y_l \subset Y_1 + \delta$.
\end{itemize}

\item A \defi{level-$l$ Young wall of type $\atwonmonetwo$} is an
ordered $l$-tuple of level-1 Young walls $\mathbf{Y} = (Y_1,
\dots, Y_l)$ such that
\begin{itemize}
\item $Y_1 \subset \cdots \subset Y_l \subset Y_1 + \delta$,
\item each column contains an even number of $n$-blocks.
\end{itemize}

\item A \defi{level-$l$ Young wall of type $\dnponetwo$} is an
ordered $l$-tuple of level-1 Young walls $\mathbf{Y} = (Y_1,
\dots, Y_l)$ such that
\begin{itemize}
\item $Y_1 \subset \cdots \subset Y_l \subset Y_1 + \delta$,
\item in each column, at most one of the top blocks is a supporting $n$-block.
\end{itemize}

\item A \defi{level-$l$ Young wall of type $\cnone$} is an
ordered $2l$-tuple of level-$\half$ Young walls $\mathbf{Y} =
(Y_1, \dots, Y_{2l})$ such that
\begin{itemize}
\item $Y_1 \subset \cdots \subset Y_{2l-1} \subset Y_{2l} \subset Y_1 + \delta$,
\item each column contains an even number of $0$-blocks.
\end{itemize}
\end{enumerate}
\end{df}

The level-$1$ Young wall (or level-$\half$ Young wall for
$\cnone$) $Y_i$ in a level-$l$ Young wall $\mathbf{Y}=(Y_1, \dots,
Y_l)$ (or $\mathbf{Y}=(Y_1, \dots, Y_{2l})$ for $\cnone$) is
called the
\defi{$i$-th layer} of the Young wall $\mathbf{Y}$.

We will write $\mathbf{Y}=(\mathbf{Y}(k))_{k=0}^{\infty}$, where
$\mathbf{Y}(k)$ denotes the $k$-th column of $\mathbf{Y}$ (reading
from right to left). We will denote by $Y_{j}(k)$ the $k$-th
column of the $j$-th layer of $\mathbf{Y}$, which is, of course,
the same as the $j$-th layer of the $k$-th column of $\mathbf{Y}$.

Let $\mathbf{Y}$ be a level-$l$ Young wall and let $\mathbf{Y}'$
be the \defi{split form} of $\mathbf{Y}$, which is obtained from
$\mathbf{Y}$ by splitting every possible block in every column of
$\mathbf{Y}$.

\begin{df} \hfill
\begin{enumerate}
\item A level-$l$ Young wall $\mathbf{Y}$ is said to be
      \defi{proper} if it satisfies the following conditions\,:
      \begin{itemize}
      \item
      in the split form $\mathbf{Y}'$, except for the rightmost column,
      there is no free space to the right of any block (or broken halves
      of blocks);
      \item
      for each layer in the split form $\mathbf{Y}'$,
      none of the columns which is of integer height
      and whose top is of unit depth have the same height.
      \end{itemize}
\item A column in a level-$l$ proper Young wall is said to \defi{contain
a removable $\delta$} if one may remove a $\delta$ from that
column to get a proper Young wall.

\item A level-$l$ proper Young wall is said to be \defi{reduced}
if none of its columns contain a removable $\delta$.

\end{enumerate}
\end{df}

We will denote by $\pwspace$ and $\rpwspace$ the set of all proper
Young walls and reduced proper Young walls, respectively.

Let $\mathbf{Y}$ be a level-$l$ proper Young wall and let $C$ be a
column of $\mathbf{Y}$. For each $i\in I$, we define
$\bar{\veps}_i(C)$ (resp. $\bar{\vphi}_i(C)$) to be the
largest integer $k \ge 0$ such that $\eit^k (C) \neq 0$ (resp.
$\fit^k (C) \neq 0$) and that $k$ is the maximal number of times
we may act $\eit$ (resp. $\fit$) to the column $C$ while the resulting
wall remains a proper Young wall.
We will often say that $C$ is
\emph{$\bar{\veps}_i$-times $i$-removable} and
\emph{$\bar{\vphi}_i$-times $i$-admissible}.

\vskip 3mm

\begin{remark}
Here, we would like to give a warning to the readers. Viewed a
level-$l$ slice in the $\uq'(\g)$-crystal $\newcrystal$, we have
already assigned the nonnegative integers $\veps_i(C)$ and
$\vphi_i(C)$. However, these numbers are not necessarily the same
as $\bar{\veps}_i(C)$ and $\bar{\vphi}_i(C)$. Note that
the number $\bar{\veps}_i(C)$ (resp. $\bar{\vphi}_i(C)$)
depends on the column lying to the left (resp. right) of $C$.
Hence we have $\bar{\veps}_i(C) \le \veps_i(C)$ and
$\bar{\vphi}_i(C) \le \vphi_i(C)$ for all $i \in I$.
\end{remark}

\vskip 3mm

We now define the action of Kashiwara operators $\eit$, $\fit$
$(i\in I)$ on $\mathbf{Y}$ as follows.

\begin{enumerate}
\item For each column $C$ of $\mathbf{Y}$, we write
$\bar{\veps}_i(C)$-many $1$'s followed by
$\bar{\vphi}_i(C)$-many $0$'s under $C$. This sequence is
called the \defi{$i$-signature of $C$}.

\item From this sequence of 1's and 0's,
cancel out each (0,1)-pair to obtain a finite sequence of 1's
followed by 0's (reading from left to right). This sequence is
called the \defi{$i$-signature of $\mathbf{Y}$}.

\item We define $\eit \mathbf{Y}$ to be the proper Young wall
obtained from $\mathbf{Y}$ by replacing the column $C$
corresponding the rightmost 1 in the $i$-signature of $\mathbf{Y}$
with the column $\eit C$.

\item We define $\fit \mathbf{Y}$ to be the proper Young wall
obtained from $\mathbf{Y}$ by replacing the column $C$
corresponding the leftmost 0 in the $i$-signature of $\mathbf{Y}$
with the column $\fit C$.

\item If there is no 1 (resp. 0) in the $i$-signature of
$\mathbf{Y}$, we define $\eit \mathbf{Y} = 0$ (resp. $\fit
\mathbf{Y} = 0$).
\end{enumerate}

We define the maps $\cwt: \pwspace \rightarrow \bar P$, $\veps_i,
\vphi_i : \pwspace \rightarrow \Z$ by
\begin{equation*}
\begin{aligned}
\cwt(\mathbf{Y}) & = \sum_{i\in I}
(\vphi_i(\mathbf{Y})-\veps_i(\mathbf{Y})) \La_i, \\
\veps_i(\mathbf{Y}) & = \text{the number of 1's in the
$i$-signature of $\mathbf{Y}$}, \\
\vphi_i(\mathbf{Y}) & = \text{the number of 0's in the
$i$-signature of $\mathbf{Y}$}.
\end{aligned}
\end{equation*}

Then it is straightforward to verify that the following theorem
holds.

\begin{thm}
The set $\pwspace$ of all level-$l$ proper Young walls, together
with the maps $\cwt$, $\eit$, $\fit$, $\veps_i$, $\vphi_i$ $(i\in
I)$ , forms a $\uq'(\g)$-crystal.
\end{thm}

\vskip 1cm

%%%%%%%%%%%%%%%%%%%%%%%%%%%%%%%%%%%%%%%%%%%%%%%%%%%%%%%%%%%%%
\section{Young wall realization of $\hwc(\la)$}

In this section, we give a new realization of higher level
irreducible highest weight crystals in terms of reduced proper
Young walls. Let $\la = a_0 \La_0 + a_1 \La_1 + \cdots + a_n
\La_n$ be a dominant integral weight of level-$l$. We define the
\defi{ground-state wall $\mathbf{Y}_{\la}$ of weight $\la$}
as the following reduced proper Young walls.

\savebox{\tmpfiga}{\begin{texdraw}
\fontsize{7}{7}\selectfont
\textref h:C v:C
\drawdim em
\setunitscale 1.9
\move(0 1)\lvec(-5.6 1)\lvec(-5.6 1.665)\lvec(0 1.665)\lvec(0 1)
\ifill f:0.8
\move(0 2.665)\lvec(-5.6 2.665)\lvec(-5.6 3.330)\lvec(0 3.330)
\lvec(0 2.665)\ifill f:0.8
\move(0 5.830)\lvec(-5.6 5.830)\lvec(-5.6 6.495)\lvec(0 6.495)
\lvec(0 5.830)\ifill f:0.8
\move(0 7.495)\lvec(-5.6 7.495)\lvec(-5.6 8.160)\lvec(0 8.160)
\lvec(0 7.495)\ifill f:0.8
\move(0.5 1)\lvec(0.8 1)
\move(0.5 1.665)\lvec(0.8 1.665)
\move(0.65 1)\lvec(0.65 1.665)
\move(0.5 2.665)\lvec(0.8 2.665)
\move(0.5 3.330)\lvec(0.8 3.330)
\move(0.65 2.665)\lvec(0.65 3.330)
\move(0.5 5.830)\lvec(0.8 5.830)
\move(0.5 6.495)\lvec(0.8 6.495)
\move(0.65 5.830)\lvec(0.65 6.495)
\move(0.5 7.495)\lvec(0.8 7.495)
\move(0.5 8.160)\lvec(0.8 8.160)
\move(0.65 7.495)\lvec(0.65 8.160)
\move(0 0)\lvec(-5.7 0)
\move(0 1)\lvec(-5.7 1)
\move(0 1.133)\lvec(-5.7 1.133)
\move(0 1.532)\lvec(-5.7 1.532)
\move(0 1.665)\lvec(-5.7 1.665)
\move(0 2.665)\lvec(-5.7 2.665)
\move(0 2.798)\lvec(-5.7 2.798)
\move(0 3.197)\lvec(-5.7 3.197)
\move(0 3.330)\lvec(-5.7 3.330)
\move(0 4.830)\lvec(-5.7 4.830)
\move(0 5.830)\lvec(-5.7 5.830)
\move(0 5.963)\lvec(-5.7 5.963)
\move(0 6.362)\lvec(-5.7 6.362)
\move(0 6.495)\lvec(-5.7 6.495)
\move(0 7.495)\lvec(-5.7 7.495)
\move(0 7.628)\lvec(-5.7 7.628)
\move(0 8.027)\lvec(-5.7 8.027)
\move(0 8.160)\lvec(-5.7 8.160)
\move(0 0)\lvec(0 1.266)
\move(0 1.399)\lvec(0 2.931)
%\move(0 3.064)\lvec(0 3.830)
%\move(0 4.330)\lvec(0 6.096)
\move(0 3.064)\lvec(0 6.096)%
\vtext(-0.5 4.080){$\cdots$}%
\move(0 6.229)\lvec(0 7.751)
\move(0 7.884)\lvec(0 8.160)
\move(-1 0)\lvec(-1 1.266)
\move(-1 1.399)\lvec(-1 2.931)
%\move(-1 3.064)\lvec(-1 3.830)
%\move(-1 4.330)\lvec(-1 6.096)
\move(-1 3.064)\lvec(-1 6.096)%
\vtext(-1.5 4.080){$\cdots$}%
\move(-1 6.229)\lvec(-1 7.751)
\move(-1 7.884)\lvec(-1 8.160)
\move(-2 0)\lvec(-2 1.266)
\move(-2 1.399)\lvec(-2 2.931)
%\move(-2 3.064)\lvec(-2 3.830)
%\move(-2 4.330)\lvec(-2 6.096)
\move(-2 3.064)\lvec(-2 6.096)%
%\vtext(-2.5 4.080){$\cdots$}%
\move(-2 6.229)\lvec(-2 7.751)
\move(-2 7.884)\lvec(-2 8.160)
\move(-3.4 0)\lvec(-3.4 1.266)
\move(-3.4 1.399)\lvec(-3.4 2.931)
%\move(-3.4 3.064)\lvec(-3.4 3.830)
%\move(-3.4 4.330)\lvec(-3.4 6.096)
\move(-3.4 3.064)\lvec(-3.4 6.096)%
\vtext(-3.9 4.080){$\cdots$}%
\move(-3.4 6.229)\lvec(-3.4 7.751)
\move(-3.4 7.884)\lvec(-3.4 8.160)
\move(-4.4 0)\lvec(-4.4 1.266)
\move(-4.4 1.399)\lvec(-4.4 2.931)
%\move(-4.4 3.064)\lvec(-4.4 3.830)
%\move(-4.4 4.330)\lvec(-4.4 6.096)
\move(-4.4 3.064)\lvec(-4.4 6.096)%
\vtext(-4.9 4.080){$\cdots$}%
\move(-4.4 6.229)\lvec(-4.4 7.751)
\move(-4.4 7.884)\lvec(-4.4 8.160)
\move(-5.4 0)\lvec(-5.4 1.266)
\move(-5.4 1.399)\lvec(-5.4 2.931)
%\move(-5.4 3.064)\lvec(-5.4 3.830)
%\move(-5.4 4.330)\lvec(-5.4 6.096)
\move(-5.4 3.064)\lvec(-5.4 6.096)%
\move(-5.4 6.229)\lvec(-5.4 7.751)
\move(-5.4 7.884)\lvec(-5.4 8.160)
\htext(-3.9 0.5){$1$}
\htext(-1.5 0.5){$n$}
\htext(-0.5 0.5){$0$}
\htext(-4.9 0.5){$0$}
\htext(-2.7 0.5){$\cdots$}
\htext(-3.9 2.165){$2$}
\htext(-1.5 2.165){$0$}
\htext(-0.5 2.165){$1$}
\htext(-4.9 2.165){$1$}
\htext(-2.7 2.165){$\cdots$}
\htext(-3.9 5.330){$n$}
\htext(-1.5 5.330){$n\!\!-\!\!2$}
\htext(-0.5 5.330){$n\!\!-\!\!1$}
\htext(-4.9 5.330){$n\!\!-\!\!1$}
\htext(-2.7 5.330){$\cdots$}
\htext(-3.9 6.995){$0$}
\htext(-1.5 6.995){$n\!\!-\!\!1$}
\htext(-0.5 6.995){$n$}
\htext(-4.9 6.995){$n$}
\htext(-2.7 6.995){$\cdots$}
\textref h:L v:C
\htext(1.2 1.3325){$a_1$}
\htext(1.2 2.9975){$a_2$}
\htext(1.2 6.1625){$a_n$}
\htext(1.2 7.8275){$a_0$}
\end{texdraw}}%
\savebox{\tmpfigh}{\begin{texdraw}
\fontsize{7}{7}\selectfont
\drawdim em
\setunitscale 1.9
\textref h:L v:C
\htext(1.0 1){$\left.\rule{0pt}{0.3em}\right\}$ $a_1$}
\end{texdraw}}%
\savebox{\tmpfigb}{\begin{texdraw}
\fontsize{7}{7}\selectfont
\textref h:C v:C
\drawdim em
\setunitscale 1.9
\move(0 0)
\bsegment
\move(0 0.5)\lvec(0 1.165)\lvec(-4.2 1.165)\lvec(-4.2 0.5)
\lvec(0 0.5)\ifill f:0.8
\move(0 1.665)\lvec(-4.2 1.665)\lvec(-4.2 2.330)\lvec(0 2.330)
\lvec(0 1.665)\ifill f:0.8
\move(0 3.330)\lvec(-4.2 3.330)\lvec(-4.2 3.995)\lvec(0 3.995)
\lvec(0 3.330)\ifill f:0.8
\move(0.5 0.5)\lvec(0.8 0.5)
\move(0.5 1.165)\lvec(0.8 1.165)
\move(0.65 0.5)\lvec(0.65 1.165)%
\move(0.5 1.665)\lvec(0.8 1.665)
\move(0.5 2.330)\lvec(0.8 2.330)
\move(0.65 1.665)\lvec(0.65 2.330)%
\move(0.5 3.330)\lvec(0.8 3.330)
\move(0.5 3.995)\lvec(0.8 3.995)
\move(0.65 3.330)\lvec(0.65 3.995)%
\move(0 0)\lvec(-4.3 0)
\move(0 0.5)\lvec(-4.3 0.5)
\move(0 0.633)\lvec(-4.3 0.633)
\move(0 1.032)\lvec(-4.3 1.032)
\move(0 1.165)\lvec(-4.3 1.165)
\move(0 1.665)\lvec(-4.3 1.665)
\move(0 1.798)\lvec(-4.3 1.798)
\move(0 2.197)\lvec(-4.3 2.197)
\move(0 2.330)\lvec(-4.3 2.330)
\move(0 3.330)\lvec(-4.3 3.330)
\move(0 3.463)\lvec(-4.3 3.463)
\move(0 3.862)\lvec(-4.3 3.862)
\move(0 3.995)\lvec(-4.3 3.995)
\move(0 0)\lvec(0 0.766)
\move(0 0.899)\lvec(0 1.931)
\move(0 2.064)\lvec(0 3.596)
%\move(0 3.729)\lvec(0 4.495)
\move(-1 0)\lvec(-1 0.766)
\move(-1 0.899)\lvec(-1 1.931)
\move(-1 2.064)\lvec(-1 3.596)
%\move(-1 3.729)\lvec(-1 4.495)
\move(-2 0)\lvec(-2 0.766)
\move(-2 0.899)\lvec(-2 1.931)
\move(-2 2.064)\lvec(-2 3.596)
%\move(-2 3.729)\lvec(-2 4.495)
\move(-3 0)\lvec(-3 0.766)
\move(-3 0.899)\lvec(-3 1.931)
\move(-3 2.064)\lvec(-3 3.596)
%\move(-3 3.729)\lvec(-3 4.495)
\move(-4 0)\lvec(-4 0.766)
\move(-4 0.899)\lvec(-4 1.931)
\move(-4 2.064)\lvec(-4 3.596)
%\move(-4 3.729)\lvec(-4 4.495)
\htext(-0.5 0.25){$0$}
\htext(-1.5 0.25){$0$}
\htext(-2.5 0.25){$0$}
\htext(-3.5 0.25){$0$}
\htext(-0.5 1.415){$0$}
\htext(-1.5 1.415){$0$}
\htext(-2.5 1.415){$0$}
\htext(-3.5 1.415){$0$}
\htext(-0.5 2.830){$1$}
\htext(-1.5 2.830){$1$}
\htext(-2.5 2.830){$1$}
\htext(-3.5 2.830){$1$}
\textref h:L v:C
\htext(1.4 0.8325){$2a_0$}
\htext(1.2 1.9975){$a_1$}
\htext(1.2 3.6625){$a_2$}
\esegment
\move(0 6.495)\lvec(0 7.160)\lvec(-4.2 7.160)\lvec(-4.2 6.495)
\lvec(0 6.495)\ifill f:0.8
\move(0 8.160)\lvec(-4.2 8.160)\lvec(-4.2 8.825)\lvec(0 8.825)
\lvec(0 8.293)\ifill f:0.8
\move(0 9.852)\lvec(-4.2 9.852)\lvec(-4.2 10.490)\lvec(0 10.490)
\lvec(0 9.958)\ifill f:0.8
\move(0.5 6.495)\lvec(0.8 6.495)
\move(0.5 7.160)\lvec(0.8 7.160)
\move(0.65 6.495)\lvec(0.65 7.160)%
\move(0.5 8.160)\lvec(0.8 8.160)
\move(0.5 8.825)\lvec(0.8 8.825)
\move(0.65 8.160)\lvec(0.65 8.825)%
\move(0.5 9.852)\lvec(0.8 9.852)
\move(0.5 10.490)\lvec(0.8 10.490)
\move(0.65 9.852)\lvec(0.65 10.490)%
%\move(0 4.995)\lvec(0 6.761)
\move(0 3.729)\lvec(0 6.761)
\vtext(-0.5 4.745){$\cdots$}
%\move(-1 4.995)\lvec(-1 6.761)
\move(-1 3.729)\lvec(-1 6.761)
\vtext(-1.5 4.745){$\cdots$}
%\move(-2 4.995)\lvec(-2 6.761)
\move(-2 3.729)\lvec(-2 6.761)
\vtext(-2.5 4.745){$\cdots$}
%\move(-3 4.995)\lvec(-3 6.761)
\move(-3 3.729)\lvec(-3 6.761)
\vtext(-3.5 4.745){$\cdots$}
%\move(-4 4.995)\lvec(-4 6.761)
\move(-4 3.729)\lvec(-4 6.761)
\move(0 5.495)\lvec(-4.3 5.495)
\move(0 6.495)\lvec(-4.3 6.495)
\move(0 6.628)\lvec(-4.3 6.628)
\move(0 7.027)\lvec(-4.3 7.027)
\move(0 7.160)\lvec(-4.3 7.160)
\move(0 8.160)\lvec(-4.3 8.160)
\move(0 8.293)\lvec(-4.3 8.293)
\move(0 8.692)\lvec(-4.3 8.692)
\move(0 8.825)\lvec(-4.3 8.825)
\move(0 9.825)\lvec(-4.3 9.825)
\move(0 9.958)\lvec(-4.3 9.958)
\move(0 10.357)\lvec(-4.3 10.357)
\move(0 10.490)\lvec(-4.3 10.490)
\move(0 6.894)\lvec(0 8.426)
\move(-1 6.894)\lvec(-1 8.426)
\move(-2 6.894)\lvec(-2 8.426)
\move(-3 6.894)\lvec(-3 8.426)
\move(-4 6.894)\lvec(-4 8.426)
\move(0 8.559)\lvec(0 10.091)
\move(-1 8.559)\lvec(-1 10.091)
\move(-2 8.559)\lvec(-2 10.091)
\move(-3 8.559)\lvec(-3 10.091)
\move(-4 8.559)\lvec(-4 10.091)
%\move(0 10.224)\lvec(0 10.990)
%\move(-1 10.224)\lvec(-1 10.990)
%\move(-2 10.224)\lvec(-2 10.990)
%\move(-3 10.224)\lvec(-3 10.990)
%\move(-4 10.224)\lvec(-4 10.990)
\htext(-0.5 5.995){$n\!\!-\!\!1$}
\htext(-1.5 5.995){$n\!\!-\!\!1$}
\htext(-2.5 5.995){$n\!\!-\!\!1$}
\htext(-3.5 5.995){$n\!\!-\!\!1$}
\htext(-0.5 7.660){$n$}
\htext(-1.5 7.660){$n$}
\htext(-2.5 7.660){$n$}
\htext(-3.5 7.660){$n$}
\htext(-0.5 9.325){$n\!\!-\!\!1$}
\htext(-1.5 9.325){$n\!\!-\!\!1$}
\htext(-2.5 9.325){$n\!\!-\!\!1$}
\htext(-3.5 9.325){$n\!\!-\!\!1$}
\move(0 -0.133)
\bsegment
\move(0 13.123)\lvec(0 13.788)\lvec(-4.2 13.788)\lvec(-4.2 13.123)
\lvec(0 13.123)\ifill f:0.8
\move(0.5 13.123)\lvec(0.8 13.123)
\move(0.5 13.788)\lvec(0.8 13.788)
\move(0.65 13.123)\lvec(0.65 13.788)%
\move(0 12.123)\lvec(-4.3 12.123)
\move(0 13.123)\lvec(-4.3 13.123)
\move(0 13.256)\lvec(-4.3 13.256)
\move(0 13.655)\lvec(-4.3 13.655)
\move(0 13.788)\lvec(-4.3 13.788)
%\move(0 11.623)\lvec(0 13.389)
%\move(-1 11.623)\lvec(-1 13.389)
%\move(-2 11.623)\lvec(-2 13.389)
%\move(-3 11.623)\lvec(-3 13.389)
%\move(-4 11.623)\lvec(-4 13.389)
\move(0 10.357)\lvec(0 13.389)
\move(-1 10.357)\lvec(-1 13.389)
\move(-2 10.357)\lvec(-2 13.389)
\move(-3 10.357)\lvec(-3 13.389)
\move(-4 10.357)\lvec(-4 13.389)
\vtext(-0.5 11.375){$\cdots$}
\vtext(-1.5 11.375){$\cdots$}
\vtext(-2.5 11.375){$\cdots$}
\vtext(-3.5 11.375){$\cdots$}
\move(0 13.522)\lvec(0 13.788)
\move(-1 13.522)\lvec(-1 13.788)
\move(-2 13.522)\lvec(-2 13.788)
\move(-3 13.522)\lvec(-3 13.788)
\move(-4 13.522)\lvec(-4 13.788)
\htext(-0.5 12.623){$1$}
\htext(-1.5 12.623){$1$}
\htext(-2.5 12.623){$1$}
\htext(-3.5 12.623){$1$}
\esegment
\textref h:L v:C
\htext(1.2 6.8275){$a_n$}
\htext(1.2 8.4925){$a_n$}
\htext(1.4 10.1575){$a_{n\!-\!1}$}
\htext(1.2 13.3225){$a_1$}
\end{texdraw}}%
\savebox{\tmpfigc}{\begin{texdraw}
\fontsize{7}{7}\selectfont
\textref h:C v:C
\drawdim em
\setunitscale 1.9
\move(0 0)
\bsegment
\move(0 0.5)\lvec(0 1.165)\lvec(-4.2 1.165)\lvec(-4.2 0.5)
\lvec(0 0.5)\ifill f:0.8
\move(0 1.665)\lvec(-4.2 1.665)\lvec(-4.2 2.330)\lvec(0 2.330)
\lvec(0 1.665)\ifill f:0.8
\move(0 3.330)\lvec(-4.2 3.330)\lvec(-4.2 3.995)\lvec(0 3.995)
\lvec(0 3.330)\ifill f:0.8
\move(0.5 0.5)\lvec(0.8 0.5)
\move(0.5 1.165)\lvec(0.8 1.165)
\move(0.65 0.5)\lvec(0.65 1.165)%
\move(0.5 1.665)\lvec(0.8 1.665)
\move(0.5 2.330)\lvec(0.8 2.330)
\move(0.65 1.665)\lvec(0.65 2.330)%
\move(0.5 3.330)\lvec(0.8 3.330)
\move(0.5 3.995)\lvec(0.8 3.995)
\move(0.65 3.330)\lvec(0.65 3.995)%
\move(0 0)\lvec(-4.3 0)
\move(0 0.5)\lvec(-4.3 0.5)
\move(0 0.633)\lvec(-4.3 0.633)
\move(0 1.032)\lvec(-4.3 1.032)
\move(0 1.165)\lvec(-4.3 1.165)
\move(0 1.665)\lvec(-4.3 1.665)
\move(0 1.798)\lvec(-4.3 1.798)
\move(0 2.197)\lvec(-4.3 2.197)
\move(0 2.330)\lvec(-4.3 2.330)
\move(0 3.330)\lvec(-4.3 3.330)
\move(0 3.463)\lvec(-4.3 3.463)
\move(0 3.862)\lvec(-4.3 3.862)
\move(0 3.995)\lvec(-4.3 3.995)
\move(0 0)\lvec(0 0.766)
\move(0 0.899)\lvec(0 1.931)
\move(0 2.064)\lvec(0 3.596)
%\move(0 3.729)\lvec(0 4.495)
\move(-1 0)\lvec(-1 0.766)
\move(-1 0.899)\lvec(-1 1.931)
\move(-1 2.064)\lvec(-1 3.596)
\move(-1 3.729)\lvec(-1 4.495)
\move(-2 0)\lvec(-2 0.766)
\move(-2 0.899)\lvec(-2 1.931)
\move(-2 2.064)\lvec(-2 3.596)
\move(-2 3.729)\lvec(-2 4.495)
\move(-3 0)\lvec(-3 0.766)
\move(-3 0.899)\lvec(-3 1.931)
\move(-3 2.064)\lvec(-3 3.596)
\move(-3 3.729)\lvec(-3 4.495)
\move(-4 0)\lvec(-4 0.766)
\move(-4 0.899)\lvec(-4 1.931)
\move(-4 2.064)\lvec(-4 3.596)
\move(-4 3.729)\lvec(-4 4.495)
\htext(-0.5 0.25){$0$}
\htext(-1.5 0.25){$0$}
\htext(-2.5 0.25){$0$}
\htext(-3.5 0.25){$0$}
\htext(-0.5 1.415){$0$}
\htext(-1.5 1.415){$0$}
\htext(-2.5 1.415){$0$}
\htext(-3.5 1.415){$0$}
\htext(-0.5 2.830){$1$}
\htext(-1.5 2.830){$1$}
\htext(-2.5 2.830){$1$}
\htext(-3.5 2.830){$1$}
\textref h:L v:C
\htext(1.2 0.8325){$a_0$}
\htext(1.2 1.9975){$a_1$}
\htext(1.2 3.6625){$a_2$}
\esegment
\move(0 6.495)\lvec(0 7.160)\lvec(-4.2 7.160)\lvec(-4.2 6.495)
\lvec(0 6.495)\ifill f:0.8
\move(0 8.160)\lvec(-4.2 8.160)\lvec(-4.2 8.825)\lvec(0 8.825)
\lvec(0 8.160)\ifill f:0.8
\move(0 9.852)\lvec(-4.2 9.852)\lvec(-4.2 10.490)\lvec(0 10.490)
\lvec(0 9.958)\ifill f:0.8
\move(0.5 6.495)\lvec(0.8 6.495)
\move(0.5 7.160)\lvec(0.8 7.160)
\move(0.65 6.495)\lvec(0.65 7.160)%
\move(0.5 8.160)\lvec(0.8 8.160)
\move(0.5 8.825)\lvec(0.8 8.825)
\move(0.65 8.160)\lvec(0.65 8.825)%
\move(0.5 9.852)\lvec(0.8 9.852)
\move(0.5 10.490)\lvec(0.8 10.490)
\move(0.65 9.852)\lvec(0.65 10.490)%
%\move(0 4.995)\lvec(0 6.761)
%\move(-1 4.995)\lvec(-1 6.761)
%\move(-2 4.995)\lvec(-2 6.761)
%\move(-3 4.995)\lvec(-3 6.761)
%\move(-4 4.995)\lvec(-4 6.761)
\move(0 3.729)\lvec(0 6.761)
\move(-1 3.729)\lvec(-1 6.761)
\move(-2 3.729)\lvec(-2 6.761)
\move(-3 3.729)\lvec(-3 6.761)
\move(-4 3.729)\lvec(-4 6.761)
\vtext(-0.5 4.745){$\cdots$}
\vtext(-1.5 4.745){$\cdots$}
\vtext(-2.5 4.745){$\cdots$}
\vtext(-3.5 4.745){$\cdots$}
\move(0 5.495)\lvec(-4.3 5.495)
\move(0 6.495)\lvec(-4.3 6.495)
\move(0 6.628)\lvec(-4.3 6.628)
\move(0 7.027)\lvec(-4.3 7.027)
\move(0 7.160)\lvec(-4.3 7.160)
\move(0 8.160)\lvec(-4.3 8.160)
\move(0 8.293)\lvec(-4.3 8.293)
\move(0 8.692)\lvec(-4.3 8.692)
\move(0 8.825)\lvec(-4.3 8.825)
\move(0 9.825)\lvec(-4.3 9.825)
\move(0 9.958)\lvec(-4.3 9.958)
\move(0 10.357)\lvec(-4.3 10.357)
\move(0 10.490)\lvec(-4.3 10.490)
\move(0 6.894)\lvec(0 8.426)
\move(-1 6.894)\lvec(-1 8.426)
\move(-2 6.894)\lvec(-2 8.426)
\move(-3 6.894)\lvec(-3 8.426)
\move(-4 6.894)\lvec(-4 8.426)
\move(0 8.559)\lvec(0 10.091)
\move(-1 8.559)\lvec(-1 10.091)
\move(-2 8.559)\lvec(-2 10.091)
\move(-3 8.559)\lvec(-3 10.091)
\move(-4 8.559)\lvec(-4 10.091)
%\move(0 10.224)\lvec(0 10.990)
%\move(-1 10.224)\lvec(-1 10.990)
%\move(-2 10.224)\lvec(-2 10.990)
%\move(-3 10.224)\lvec(-3 10.990)
%\move(-4 10.224)\lvec(-4 10.990)
\htext(-0.5 5.995){$n\!\!-\!\!1$}
\htext(-1.5 5.995){$n\!\!-\!\!1$}
\htext(-2.5 5.995){$n\!\!-\!\!1$}
\htext(-3.5 5.995){$n\!\!-\!\!1$}
\htext(-0.5 7.660){$n$}
\htext(-1.5 7.660){$n$}
\htext(-2.5 7.660){$n$}
\htext(-3.5 7.660){$n$}
\htext(-0.5 9.325){$n\!\!-\!\!1$}
\htext(-1.5 9.325){$n\!\!-\!\!1$}
\htext(-2.5 9.325){$n\!\!-\!\!1$}
\htext(-3.5 9.325){$n\!\!-\!\!1$}
\move(0 -0.133)
\bsegment
\move(0 13.123)\lvec(0 13.788)\lvec(-4.2 13.788)\lvec(-4.2 13.123)
\lvec(0 13.123)\ifill f:0.8
\move(0.5 13.123)\lvec(0.8 13.123)
\move(0.5 13.788)\lvec(0.8 13.788)
\move(0.65 13.123)\lvec(0.65 13.788)%
\move(0 12.123)\lvec(-4.3 12.123)
\move(0 13.123)\lvec(-4.3 13.123)
\move(0 13.256)\lvec(-4.3 13.256)
\move(0 13.655)\lvec(-4.3 13.655)
\move(0 13.788)\lvec(-4.3 13.788)
%\move(0 11.623)\lvec(0 13.389)
%\move(-1 11.623)\lvec(-1 13.389)
%\move(-2 11.623)\lvec(-2 13.389)
%\move(-3 11.623)\lvec(-3 13.389)
%\move(-4 11.623)\lvec(-4 13.389)
\move(0 10.357)\lvec(0 13.389)
\move(-1 10.357)\lvec(-1 13.389)
\move(-2 10.357)\lvec(-2 13.389)
\move(-3 10.357)\lvec(-3 13.389)
\move(-4 10.357)\lvec(-4 13.389)
\vtext(-0.5 11.373){$\cdots$}
\vtext(-1.5 11.373){$\cdots$}
\vtext(-2.5 11.373){$\cdots$}
\vtext(-3.5 11.373){$\cdots$}
\move(0 13.522)\lvec(0 13.788)
\move(-1 13.522)\lvec(-1 13.788)
\move(-2 13.522)\lvec(-2 13.788)
\move(-3 13.522)\lvec(-3 13.788)
\move(-4 13.522)\lvec(-4 13.788)
\htext(-0.5 12.623){$1$}
\htext(-1.5 12.623){$1$}
\htext(-2.5 12.623){$1$}
\htext(-3.5 12.623){$1$}
\esegment
\textref h:L v:C
\htext(1.2 6.8275){$a_n$}
\htext(1.2 8.4925){$a_n$}
\htext(1.2 10.1575){$a_{n\!-\!1}$}
\htext(1.2 13.3225){$a_1$}
\end{texdraw}}%
\savebox{\tmpfigd}{\begin{texdraw}
\fontsize{7}{7}\selectfont
\textref h:C v:C
\drawdim em
\setunitscale 1.9
\move(0 -1.665)
\bsegment
\move(0 0.5)\lvec(0 1.165)\lvec(-4.2 1.165)\lvec(-4.2 0.5)
\lvec(0 0.5)\ifill f:0.8
\move(0 1.665)\lvec(-4.2 1.665)\lvec(-4.2 2.330)\lvec(0 2.330)
\lvec(0 1.665)\ifill f:0.8
\move(0 3.330)\lvec(-4.2 3.330)\lvec(-4.2 3.995)\lvec(0 3.995)
\lvec(0 3.330)\ifill f:0.8
\move(0.5 0.5)\lvec(0.8 0.5)
\move(0.5 1.165)\lvec(0.8 1.165)
\move(0.65 0.5)\lvec(0.65 1.165)%
\move(0.5 1.665)\lvec(0.8 1.665)
\move(0.5 2.330)\lvec(0.8 2.330)
\move(0.65 1.665)\lvec(0.65 2.330)%
\move(0.5 3.330)\lvec(0.8 3.330)
\move(0.5 3.995)\lvec(0.8 3.995)
\move(0.65 3.330)\lvec(0.65 3.995)%
\move(0 0)\lvec(-4.3 0)
\move(0 0.5)\lvec(-4.3 0.5)
\move(0 0.633)\lvec(-4.3 0.633)
\move(0 1.032)\lvec(-4.3 1.032)
\move(0 1.165)\lvec(-4.3 1.165)
\move(0 1.665)\lvec(-4.3 1.665)
\move(0 1.798)\lvec(-4.3 1.798)
\move(0 2.197)\lvec(-4.3 2.197)
\move(0 2.330)\lvec(-4.3 2.330)
\move(0 3.330)\lvec(-4.3 3.330)
\move(0 3.463)\lvec(-4.3 3.463)
\move(0 3.862)\lvec(-4.3 3.862)
\move(0 3.995)\lvec(-4.3 3.995)
\move(0 0)\lvec(0 0.766)
\move(0 0.899)\lvec(0 1.931)
\move(0 2.064)\lvec(0 3.596)
\move(0 3.729)\lvec(0 4.495)
\move(-1 0)\lvec(-1 0.766)
\move(-1 0.899)\lvec(-1 1.931)
\move(-1 2.064)\lvec(-1 3.596)
\move(-1 3.729)\lvec(-1 4.495)
\move(-1 4.995)\lvec(-1 6.761)
\move(-1 6.894)\lvec(-1 8.559)
\move(-2 0)\lvec(-2 0.766)
\move(-2 0.899)\lvec(-2 1.931)
\move(-2 2.064)\lvec(-2 3.596)
\move(-2 3.729)\lvec(-2 4.495)
\move(-2 4.995)\lvec(-2 6.761)
\move(-2 6.894)\lvec(-2 8.559)
\move(-3 0)\lvec(-3 0.766)
\move(-3 0.899)\lvec(-3 1.931)
\move(-3 2.064)\lvec(-3 3.596)
\move(-3 3.729)\lvec(-3 4.495)
\move(-3 4.995)\lvec(-3 6.761)
\move(-3 6.894)\lvec(-3 8.559)
\move(-4 0)\lvec(-4 0.766)
\move(-4 0.899)\lvec(-4 1.931)
\move(-4 2.064)\lvec(-4 3.596)
\move(-4 3.729)\lvec(-4 4.495)
\move(-4 4.995)\lvec(-4 6.761)
\move(-4 6.894)\lvec(-4 8.559)
\htext(-0.5 0.25){$0$}
\htext(-1.5 0.25){$0$}
\htext(-2.5 0.25){$0$}
\htext(-3.5 0.25){$0$}
\htext(-0.5 1.415){$0$}
\htext(-1.5 1.415){$0$}
\htext(-2.5 1.415){$0$}
\htext(-3.5 1.415){$0$}
\htext(-0.5 2.830){$1$}
\htext(-1.5 2.830){$1$}
\htext(-2.5 2.830){$1$}
\htext(-3.5 2.830){$1$}
\esegment
\move(0 4.830)\lvec(0 5.495)\lvec(-4.2 5.495)\lvec(-4.2 4.830)
\lvec(0 4.830)\ifill f:0.8
\move(0 6.495)\lvec(0 7.160)\lvec(-4.2 7.160)\lvec(-4.2 6.495)
\lvec(0 6.495)\ifill f:0.8
\move(0 8.160)\lvec(-4.2 8.160)\lvec(-4.2 8.825)\lvec(0 8.825)\ifill f:0.8
\move(0 9.852)\lvec(-4.2 9.852)\lvec(-4.2 10.490)\lvec(0 10.490)\ifill f:0.8
\move(0.5 4.830)\lvec(0.8 4.830)
\move(0.5 5.495)\lvec(0.8 5.495)
\move(0.65 4.830)\lvec(0.65 5.495)%
\move(0.5 6.495)\lvec(0.8 6.495)
\move(0.5 7.160)\lvec(0.8 7.160)
\move(0.65 6.495)\lvec(0.65 7.160)%
\move(0.5 8.160)\lvec(0.8 8.160)
\move(0.5 8.825)\lvec(0.8 8.825)
\move(0.65 8.160)\lvec(0.65 8.825)%
\move(0.5 9.852)\lvec(0.8 9.852)
\move(0.5 10.490)\lvec(0.8 10.490)
\move(0.65 9.825)\lvec(0.65 10.490)%
\move(0 3.830)\lvec(-4.3 3.830)
\move(0 4.830)\lvec(-4.3 4.830)
\move(0 4.963)\lvec(-4.3 4.963)
\move(0 5.362)\lvec(-4.3 5.362)
\move(0 5.495)\lvec(-4.3 5.495)
\move(0 3.330)\lvec(0 5.063)
\move(-1 3.330)\lvec(-1 5.063)
\move(-2 3.330)\lvec(-2 5.063)
\move(-3 3.330)\lvec(-3 5.063)
\move(-4 3.330)\lvec(-4 5.063)
\move(0 5.196)\lvec(0 6.761)
\move(-1 5.196)\lvec(-1 6.761)
\move(-2 5.196)\lvec(-2 6.761)
\move(-3 5.196)\lvec(-3 6.761)
\move(-4 5.196)\lvec(-4 6.761)
\move(0 5.495)\lvec(-4.3 5.495)
\move(0 6.495)\lvec(-4.3 6.495)
\move(0 6.628)\lvec(-4.3 6.628)
\move(0 7.027)\lvec(-4.3 7.027)
\move(0 7.160)\lvec(-4.3 7.160)
\move(0 7.660)\lvec(-4.3 7.660)
\move(0 8.160)\lvec(-4.3 8.160)
\move(0 8.293)\lvec(-4.3 8.293)
\move(0 8.692)\lvec(-4.3 8.692)
\move(0 8.825)\lvec(-4.3 8.825)
\move(0 9.825)\lvec(-4.3 9.825)
\move(0 9.958)\lvec(-4.3 9.958)
\move(0 10.357)\lvec(-4.3 10.357)
\move(0 10.490)\lvec(-4.3 10.490)
\move(0 4.995)\lvec(0 6.761)
\move(0 6.894)\lvec(0 8.426)
\move(-1 6.894)\lvec(-1 8.426)
\move(-2 6.894)\lvec(-2 8.426)
\move(-3 6.894)\lvec(-3 8.426)
\move(-4 6.894)\lvec(-4 8.426)
\move(0 8.559)\lvec(0 10.091)
\move(-1 8.559)\lvec(-1 10.091)
\move(-2 8.559)\lvec(-2 10.091)
\move(-3 8.559)\lvec(-3 10.091)
\move(-4 8.559)\lvec(-4 10.091)
\move(0 10.224)\lvec(0 10.990)
\move(-1 10.224)\lvec(-1 10.990)
\move(-2 10.224)\lvec(-2 10.990)
\move(-3 10.224)\lvec(-3 10.990)
\move(-4 10.224)\lvec(-4 10.990)
\htext(-0.5 5.995){$n\!\!-\!\!1$}
\htext(-1.5 5.995){$n\!\!-\!\!1$}
\htext(-2.5 5.995){$n\!\!-\!\!1$}
\htext(-3.5 5.995){$n\!\!-\!\!1$}
\htext(-0.5 7.410){$n$}
\htext(-1.5 7.410){$n$}
\htext(-2.5 7.410){$n$}
\htext(-3.5 7.410){$n$}
\htext(-0.5 7.910){$n$}
\htext(-1.5 7.910){$n$}
\htext(-2.5 7.910){$n$}
\htext(-3.5 7.910){$n$}
\htext(-0.5 4.330){$n\!\!-\!\!2$}
\htext(-1.5 4.330){$n\!\!-\!\!2$}
\htext(-2.5 4.330){$n\!\!-\!\!2$}
\htext(-3.5 4.330){$n\!\!-\!\!2$}
\htext(-0.5 9.325){$n\!\!-\!\!1$}
\htext(-1.5 9.325){$n\!\!-\!\!1$}
\htext(-2.5 9.325){$n\!\!-\!\!1$}
\htext(-3.5 9.325){$n\!\!-\!\!1$}
\move(0 -0.133)
\bsegment
\move(0 13.123)\lvec(0 13.788)\lvec(-4.2 13.788)\lvec(-4.2 13.123)
\lvec(0 13.123)\ifill f:0.8
\move(0.5 13.123)\lvec(0.8 13.123)
\move(0.5 13.788)\lvec(0.8 13.788)
\move(0.65 13.123)\lvec(0.65 13.788)%
\move(0 12.123)\lvec(-4.3 12.123)
\move(0 13.123)\lvec(-4.3 13.123)
\move(0 13.256)\lvec(-4.3 13.256)
\move(0 13.655)\lvec(-4.3 13.655)
\move(0 13.788)\lvec(-4.3 13.788)
\move(0 11.623)\lvec(0 13.389)
\move(-1 11.623)\lvec(-1 13.389)
\move(-2 11.623)\lvec(-2 13.389)
\move(-3 11.623)\lvec(-3 13.389)
\move(-4 11.623)\lvec(-4 13.389)
\move(0 13.522)\lvec(0 13.788)
\move(-1 13.522)\lvec(-1 13.788)
\move(-2 13.522)\lvec(-2 13.788)
\move(-3 13.522)\lvec(-3 13.788)
\move(-4 13.522)\lvec(-4 13.788)
\htext(-0.5 12.623){$1$}
\htext(-1.5 12.623){$1$}
\htext(-2.5 12.623){$1$}
\htext(-3.5 12.623){$1$}
\esegment
\move(0 -1.665)
\bsegment
\textref h:L v:C
\htext(1.2 0.8325){$a_0$}
\htext(1.2 1.9975){$a_1$}
\htext(1.2 3.6625){$a_2$}
\esegment
\textref h:L v:C
\htext(1.2 5.1625){$a_{n\!-\!1}$}
\htext(1.2 6.8275){$\frac{a_n}{2}$}
\htext(1.2 8.4925){$\frac{a_n}{2}$}
\htext(1.2 10.1575){$a_{n\!-\!1}$}
\htext(1.2 13.3225){$a_1$}
\end{texdraw}}%
\savebox{\tmpfige}{\begin{texdraw}
\fontsize{7}{7}\selectfont
\textref h:C v:C
\drawdim em
\setunitscale 1.9
\move(0 -1.665)
\bsegment
\move(0 0.5)\lvec(0 1.165)\lvec(-4.2 1.165)\lvec(-4.2 0.5)
\lvec(0 0.5)\ifill f:0.8
\move(0 1.665)\lvec(-4.2 1.665)\lvec(-4.2 2.330)\lvec(0 2.330)
\lvec(0 1.665)\ifill f:0.8
\move(0 3.330)\lvec(-4.2 3.330)\lvec(-4.2 3.995)\lvec(0 3.995)
\lvec(0 3.330)\ifill f:0.8
\move(0.5 0.5)\lvec(0.8 0.5)
\move(0.5 1.165)\lvec(0.8 1.165)
\move(0.65 0.5)\lvec(0.65 1.165)%
\move(0.5 1.665)\lvec(0.8 1.665)
\move(0.5 2.330)\lvec(0.8 2.330)
\move(0.65 1.665)\lvec(0.65 2.330)%
\move(0.5 3.330)\lvec(0.8 3.330)
\move(0.5 3.995)\lvec(0.8 3.995)
\move(0.65 3.330)\lvec(0.65 3.995)%
\move(0 0)\lvec(-4.3 0)
\move(0 0.5)\lvec(-4.3 0.5)
\move(0 0.633)\lvec(-4.3 0.633)
\move(0 1.032)\lvec(-4.3 1.032)
\move(0 1.165)\lvec(-4.3 1.165)
\move(0 1.665)\lvec(-4.3 1.665)
\move(0 1.798)\lvec(-4.3 1.798)
\move(0 2.197)\lvec(-4.3 2.197)
\move(0 2.330)\lvec(-4.3 2.330)
\move(0 3.330)\lvec(-4.3 3.330)
\move(0 3.463)\lvec(-4.3 3.463)
\move(0 3.862)\lvec(-4.3 3.862)
\move(0 3.995)\lvec(-4.3 3.995)
\move(0 0)\lvec(0 0.766)
\move(0 0.899)\lvec(0 1.931)
\move(0 2.064)\lvec(0 3.596)
%\move(0 3.729)\lvec(0 4.495)
\move(-1 0)\lvec(-1 0.766)
\move(-1 0.899)\lvec(-1 1.931)
\move(-1 2.064)\lvec(-1 3.596)
%\move(-1 3.729)\lvec(-1 4.495)
\move(-1 4.995)\lvec(-1 6.761)
\move(-1 6.894)\lvec(-1 8.559)
\move(-2 0)\lvec(-2 0.766)
\move(-2 0.899)\lvec(-2 1.931)
\move(-2 2.064)\lvec(-2 3.596)
%\move(-2 3.729)\lvec(-2 4.495)
\move(-2 4.995)\lvec(-2 6.761)
\move(-2 6.894)\lvec(-2 8.559)
\move(-3 0)\lvec(-3 0.766)
\move(-3 0.899)\lvec(-3 1.931)
\move(-3 2.064)\lvec(-3 3.596)
%\move(-3 3.729)\lvec(-3 4.495)
\move(-3 4.995)\lvec(-3 6.761)
\move(-3 6.894)\lvec(-3 8.559)
\move(-4 0)\lvec(-4 0.766)
\move(-4 0.899)\lvec(-4 1.931)
\move(-4 2.064)\lvec(-4 3.596)
%\move(-4 3.729)\lvec(-4 4.495)
\move(-4 4.995)\lvec(-4 6.761)
\move(-4 6.894)\lvec(-4 8.559)
\htext(-0.5 0.25){$0$}
\htext(-1.5 0.25){$0$}
\htext(-2.5 0.25){$0$}
\htext(-3.5 0.25){$0$}
\htext(-0.5 1.415){$0$}
\htext(-1.5 1.415){$0$}
\htext(-2.5 1.415){$0$}
\htext(-3.5 1.415){$0$}
\htext(-0.5 2.830){$1$}
\htext(-1.5 2.830){$1$}
\htext(-2.5 2.830){$1$}
\htext(-3.5 2.830){$1$}
\esegment
\move(0 4.830)\lvec(0 5.495)\lvec(-4.2 5.495)\lvec(-4.2 4.830)
\lvec(0 4.830)\ifill f:0.8
\move(0 6.495)\lvec(0 7.160)\lvec(-4.2 7.160)\lvec(-4.2 6.495)
\lvec(0 6.495)\ifill f:0.8
\move(0 7.660)\lvec(0 7.793)\lvec(-4.2 7.793)\lvec(-4.2 7.660)
\lvec(0 7.660)\ifill f:0.8
\move(0 8.293)\lvec(-4.2 8.293)\lvec(-4.2 8.958)\lvec(0 8.958)
\lvec(0 8.293)\ifill f:0.8
\move(0 9.958)\lvec(-4.2 9.958)\lvec(-4.2 10.623)\lvec(0 10.623)
\lvec(0 9.958)\ifill f:0.8
\move(0.5 4.830)\lvec(0.8 4.830)
\move(0.5 5.495)\lvec(0.8 5.495)
\move(0.65 4.830)\lvec(0.65 5.495)%
\move(0.5 6.495)\lvec(0.8 6.495)
\move(0.5 7.160)\lvec(0.8 7.160)
\move(0.65 6.495)\lvec(0.65 7.160)%
\move(0.5 7.660)\lvec(0.8 7.660)
\move(0.5 7.793)\lvec(0.8 7.793)
\move(0.65 7.660)\lvec(0.65 7.793)%
\move(0.5 8.293)\lvec(0.8 8.293)
\move(0.5 8.958)\lvec(0.8 8.958)
\move(0.65 8.293)\lvec(0.65 8.958)%
\move(0.5 9.958)\lvec(0.8 9.958)
\move(0.5 10.623)\lvec(0.8 10.623)
\move(0.65 9.958)\lvec(0.65 10.623)%
\move(0 3.830)\lvec(-4.3 3.830)
\move(0 4.830)\lvec(-4.3 4.830)
\move(0 4.963)\lvec(-4.3 4.963)
\move(0 5.362)\lvec(-4.3 5.362)
\move(0 5.495)\lvec(-4.3 5.495)
\move(0 6.495)\lvec(-4.3 6.495)
\move(0 6.628)\lvec(-4.3 6.628)
\move(0 7.027)\lvec(-4.3 7.027)
\move(0 7.160)\lvec(-4.3 7.160)
\move(0 7.660)\lvec(-4.3 7.660)
\move(0 7.793)\lvec(-4.3 7.793)
\move(0 8.293)\lvec(-4.3 8.293)
\move(0 8.426)\lvec(-4.3 8.426)
\move(0 8.825)\lvec(-4.3 8.825)
\move(0 8.958)\lvec(-4.3 8.958)
\move(0 9.958)\lvec(-4.3 9.958)
\move(0 10.091)\lvec(-4.3 10.091)
\move(0 10.490)\lvec(-4.3 10.490)
\move(0 10.623)\lvec(-4.3 10.623)
%\move(0 3.330)\lvec(0 5.063)
%\move(-1 3.330)\lvec(-1 5.063)
%\move(-2 3.330)\lvec(-2 5.063)
%\move(-3 3.330)\lvec(-3 5.063)
%\move(-4 3.330)\lvec(-4 5.063)
\move(0 2.030)\lvec(0 5.063)
\move(-1 2.030)\lvec(-1 5.063)
\move(-2 2.030)\lvec(-2 5.063)
\move(-3 2.030)\lvec(-3 5.063)
\move(-4 2.030)\lvec(-4 5.063)
\vtext(-0.5 3.080){$\cdots$}
\vtext(-1.5 3.080){$\cdots$}
\vtext(-2.5 3.080){$\cdots$}
\vtext(-3.5 3.080){$\cdots$}
\move(0 5.196)\lvec(0 6.761)
\move(-1 5.196)\lvec(-1 6.761)
\move(-2 5.196)\lvec(-2 6.761)
\move(-3 5.196)\lvec(-3 6.761)
\move(-4 5.196)\lvec(-4 6.761)
\move(0 6.894)\lvec(0 8.559)
\move(-1 6.894)\lvec(-1 8.559)
\move(-2 6.894)\lvec(-2 8.559)
\move(-3 6.894)\lvec(-3 8.559)
\move(-4 6.894)\lvec(-4 8.559)
\move(0 8.692)\lvec(0 10.224)
\move(-1 8.692)\lvec(-1 10.224)
\move(-2 8.692)\lvec(-2 10.224)
\move(-3 8.692)\lvec(-3 10.224)
\move(-4 8.692)\lvec(-4 10.224)
%\move(0 10.357)\lvec(0 11.123)
%\move(-1 10.357)\lvec(-1 11.123)
%\move(-2 10.357)\lvec(-2 11.123)
%\move(-3 10.357)\lvec(-3 11.123)
%\move(-4 10.357)\lvec(-4 11.123)
\htext(-0.5 4.330){$n\!\!-\!\!2$}
\htext(-1.5 4.330){$n\!\!-\!\!2$}
\htext(-2.5 4.330){$n\!\!-\!\!2$}
\htext(-3.5 4.330){$n\!\!-\!\!2$}
\htext(-0.5 5.995){$n\!\!-\!\!1$}
\htext(-1.5 5.995){$n\!\!-\!\!1$}
\htext(-2.5 5.995){$n\!\!-\!\!1$}
\htext(-3.5 5.995){$n\!\!-\!\!1$}
\htext(-0.5 7.410){$n$}
\htext(-1.5 7.410){$n$}
\htext(-2.5 7.410){$n$}
\htext(-3.5 7.410){$n$}
\htext(-0.5 8.043){$n$}
\htext(-1.5 8.043){$n$}
\htext(-2.5 8.043){$n$}
\htext(-3.5 8.043){$n$}
\htext(-0.5 9.458){$n\!\!-\!\!1$}
\htext(-1.5 9.458){$n\!\!-\!\!1$}
\htext(-2.5 9.458){$n\!\!-\!\!1$}
\htext(-3.5 9.458){$n\!\!-\!\!1$}
\move(0 13.123)\lvec(0 13.788)\lvec(-4.2 13.788)\lvec(-4.2 13.123)
\lvec(0 13.123)\ifill f:0.8
\move(0.5 13.123)\lvec(0.8 13.123)
\move(0.5 13.788)\lvec(0.8 13.788)
\move(0.65 13.123)\lvec(0.65 13.788)%
\move(0 12.123)\lvec(-4.3 12.123)
\move(0 13.123)\lvec(-4.3 13.123)
\move(0 13.256)\lvec(-4.3 13.256)
\move(0 13.655)\lvec(-4.3 13.655)
\move(0 13.788)\lvec(-4.3 13.788)
%\move(0 11.623)\lvec(0 13.389)
%\move(-1 11.623)\lvec(-1 13.389)
%\move(-2 11.623)\lvec(-2 13.389)
%\move(-3 11.623)\lvec(-3 13.389)
%\move(-4 11.623)\lvec(-4 13.389)
\move(0 10.357)\lvec(0 13.389)
\move(-1 10.357)\lvec(-1 13.389)
\move(-2 10.357)\lvec(-2 13.389)
\move(-3 10.357)\lvec(-3 13.389)
\move(-4 10.357)\lvec(-4 13.389)
\vtext(-0.5 11.373){$\cdots$}
\vtext(-1.5 11.373){$\cdots$}
\vtext(-2.5 11.373){$\cdots$}
\vtext(-3.5 11.373){$\cdots$}
\move(0 13.522)\lvec(0 13.788)
\move(-1 13.522)\lvec(-1 13.788)
\move(-2 13.522)\lvec(-2 13.788)
\move(-3 13.522)\lvec(-3 13.788)
\move(-4 13.522)\lvec(-4 13.788)
\htext(-0.5 12.623){$1$}
\htext(-1.5 12.623){$1$}
\htext(-2.5 12.623){$1$}
\htext(-3.5 12.623){$1$}
\textref h:L v:C
\move(0 -1.665)
\bsegment
\htext(1.2 0.8325){$a_0$}
\htext(1.2 1.9975){$a_1$}
\htext(1.2 3.6625){$a_2$}
\esegment
\htext(1.2 5.1625){$a_{n\!-\!1}$}
%\htext(1.2 6.8275){$\frac{a_{n}\!-\!1}{2}$}
%\htext(1.2 8.6255){$\frac{a_{n}\!-\!1}{2}$}
\htext(1.2 10.2905){$a_{n\!-\!1}$}
\htext(1.2 13.4555){$a_1$}
\htext(1.2 6.8275){$[\frac{a_{n}}{2}]$}
\htext(1.2 8.6255){$[\frac{a_{n}}{2}]$}
\htext(1.2 7.7261){$a_n-2[\frac{a_{n}}{2}]$}
\end{texdraw}}%
\savebox{\tmpfigf}{\begin{texdraw}
\fontsize{7}{7}\selectfont
\textref h:C v:C
\drawdim em
\setunitscale 1.9
\move(0 1.665)\lvec(-4.2 1.665)\lvec(-4.2 2.330)\lvec(0 2.330)
\lvec(0 1.665)\ifill f:0.8
\move(0 3.330)\lvec(-4.2 3.330)\lvec(-4.2 3.995)\lvec(0 3.995)
\lvec(0 3.330)\ifill f:0.8
\move(0.5 1.665)\lvec(0.8 1.665)
\move(0.5 2.330)\lvec(0.8 2.330)
\move(0.65 1.665)\lvec(0.65 2.330)%
\move(0.5 3.330)\lvec(0.8 3.330)
\move(0.5 3.995)\lvec(0.8 3.995)
\move(0.65 3.330)\lvec(0.65 3.995)%
\move(0.5 1)\lvec(0.8 1)
\move(0.5 1.665)\lvec(0.8 1.665)
\move(0.65 1)\lvec(0.65 1.665)%
\move(0 0)\lvec(-4.3 0)
\move(0 1)\lvec(-4.3 1)
\move(0 1.133)\lvec(-4.3 1.133)
\move(0 1.532)\lvec(-4.3 1.532)
\move(0 1.665)\lvec(-4.3 1.665)
\move(0 1.798)\lvec(-4.3 1.798)
\move(0 2.197)\lvec(-4.3 2.197)
\move(0 2.330)\lvec(-4.3 2.330)
\move(0 3.330)\lvec(-4.3 3.330)
\move(0 3.463)\lvec(-4.3 3.463)
\move(0 3.862)\lvec(-4.3 3.862)
\move(0 3.995)\lvec(-4.3 3.995)
\move(0 0)\lvec(0 1.266)
\move(-1 0)\lvec(-1 1.266)
\move(-2 0)\lvec(-2 1.266)
\move(-3 0)\lvec(-3 1.266)
\move(-4 0)\lvec(-4 1.266)
\move(0 1.399)\lvec(0 1.931)
\move(-1 1.399)\lvec(-1 1.931)
\move(-2 1.399)\lvec(-2 1.931)
\move(-3 1.399)\lvec(-3 1.931)
\move(-4 1.399)\lvec(-4 1.931)
\move(0 2.064)\lvec(0 3.596)
\move(-1 2.064)\lvec(-1 3.596)
\move(-2 2.064)\lvec(-2 3.596)
\move(-3 2.064)\lvec(-3 3.596)
\move(-4 2.064)\lvec(-4 3.596)
%\move(0 3.729)\lvec(0 4.495)
%\move(-1 3.729)\lvec(-1 4.495)
%\move(-2 3.729)\lvec(-2 4.495)
%\move(-3 3.729)\lvec(-3 4.495)
%\move(-4 3.729)\lvec(-4 4.495)
\move(0 1)\lvec(-1 0)
\move(-1 1)\lvec(-2 0)
\move(-2 1)\lvec(-3 0)
\move(-3 1)\lvec(-4 0)
\htext(-1.25 0.25){$1$}
\htext(-2.25 0.25){$0$}
\htext(-3.25 0.25){$1$}
\htext(-0.25 0.25){$0$}
\htext(-0.5 2.830){$2$}
\htext(-3.5 2.830){$2$}
\htext(-2.5 2.830){$2$}
\htext(-1.5 2.830){$2$}
\move(0 6.495)\lvec(-4.2 6.495)\lvec(-4.2 7.160)\lvec(0 7.160)
\lvec(0 6.495)\ifill f:0.8
\move(0 7.660)\lvec(-4.2 7.660)\lvec(-4.2 8.325)\lvec(0 8.325)
\lvec(0 7.660)\ifill f:0.8
\move(0 8.825)\lvec(-4.2 8.825)\lvec(-4.2 9.490)\lvec(0 9.490)
\lvec(0 8.825)\ifill f:0.8
\move(0 10.490)\lvec(-4.2 10.490)\lvec(-4.2 11.156)\lvec(0 11.156)
\lvec(0 10.490)\ifill f:0.8
\move(0.5 6.495)\lvec(0.8 6.495)
\move(0.5 7.160)\lvec(0.8 7.160)
\move(0.65 6.495)\lvec(0.65 7.160)%
\move(0.5 7.660)\lvec(0.8 7.660)
\move(0.5 8.325)\lvec(0.8 8.325)
\move(0.65 7.660)\lvec(0.65 8.325)%
\move(0.5 8.825)\lvec(0.8 8.825)
\move(0.5 9.490)\lvec(0.8 9.490)
\move(0.65 8.825)\lvec(0.65 9.490)%
\move(0.5 10.490)\lvec(0.8 10.490)
\move(0.5 11.156)\lvec(0.8 11.156)
\move(0.65 10.490)\lvec(0.65 11.156)%
\move(0 5.495)\lvec(-4.3 5.495)
\move(0 6.495)\lvec(-4.3 6.495)
\move(0 6.628)\lvec(-4.3 6.628)
\move(0 7.027)\lvec(-4.3 7.027)
\move(0 7.160)\lvec(-4.3 7.160)
\move(0 7.660)\lvec(-4.3 7.660)
\move(0 7.793)\lvec(-4.3 7.793)
\move(0 8.192)\lvec(-4.3 8.192)
\move(0 8.325)\lvec(-4.3 8.325)
\move(0 8.825)\lvec(-4.3 8.825)
\move(0 8.958)\lvec(-4.3 8.958)
\move(0 9.357)\lvec(-4.3 9.357)
\move(0 9.490)\lvec(-4.3 9.490)
\move(0 10.490)\lvec(-4.3 10.490)
\move(0 10.623)\lvec(-4.3 10.623)
\move(0 11.023)\lvec(-4.3 11.023)
\move(0 11.156)\lvec(-4.3 11.156)
%\move(0 4.995)\lvec(0 6.761)
%\move(-1 4.995)\lvec(-1 6.761)
%\move(-2 4.995)\lvec(-2 6.761)
%\move(-3 4.995)\lvec(-3 6.761)
%\move(-4 4.995)\lvec(-4 6.761)
\move(0 3.729)\lvec(0 6.761)%
\move(-1 3.729)\lvec(-1 6.761)%
\move(-2 3.729)\lvec(-2 6.761)%
\move(-3 3.729)\lvec(-3 6.761)%
\move(-4 3.729)\lvec(-4 6.761)%
\vtext(-0.5 4.745){$\cdots$}
\vtext(-1.5 4.745){$\cdots$}
\vtext(-2.5 4.745){$\cdots$}
\vtext(-3.5 4.745){$\cdots$}
\move(0 6.894)\lvec(0 7.926)
\move(-1 6.894)\lvec(-1 7.926)
\move(-2 6.894)\lvec(-2 7.926)
\move(-3 6.894)\lvec(-3 7.926)
\move(-4 6.894)\lvec(-4 7.926)
\move(-4 8.059)\lvec(-4 9.091)
\move(-1 8.059)\lvec(-1 9.091)
\move(-2 8.059)\lvec(-2 9.091)
\move(-3 8.059)\lvec(-3 9.091)
\move(0 8.059)\lvec(0 9.091)
\move(0 9.224)\lvec(0 10.756)
\move(-1 9.224)\lvec(-1 10.756)
\move(-2 9.224)\lvec(-2 10.756)
\move(-3 9.224)\lvec(-3 10.756)
\move(-4 9.224)\lvec(-4 10.756)
%\move(-4 10.889)\lvec(-4 11.656)
%\move(-3 10.889)\lvec(-3 11.656)
%\move(-2 10.889)\lvec(-2 11.656)
%\move(-1 10.889)\lvec(-1 11.656)
%\move(0 10.889)\lvec(0 11.656)
\htext(-0.5 5.995){$n\!\!-\!\!1$}
\htext(-1.5 5.995){$n\!\!-\!\!1$}
\htext(-2.5 5.995){$n\!\!-\!\!1$}
\htext(-3.5 5.995){$n\!\!-\!\!1$}
\htext(-0.5 7.410){$n$}
\htext(-1.5 7.410){$n$}
\htext(-2.5 7.410){$n$}
\htext(-3.5 7.410){$n$}
\htext(-3.5 8.575){$n$}
\htext(-2.5 8.575){$n$}
\htext(-1.5 8.575){$n$}
\htext(-0.5 8.575){$n$}
\htext(-0.5 9.990){$n\!\!-\!\!1$}
\htext(-1.5 9.990){$n\!\!-\!\!1$}
\htext(-2.5 9.990){$n\!\!-\!\!1$}
\htext(-3.5 9.990){$n\!\!-\!\!1$}
\move(0 13.656)\lvec(-4.2 13.656)\lvec(-4.2 14.321)\lvec(0 14.321)
\ifill f:0.8
\move(0.5 13.656)\lvec(0.8 13.656)
\move(0.5 14.321)\lvec(0.8 14.321)
\move(0.65 13.656)\lvec(0.65 14.321)%
\move(0 12.656)\lvec(-4.3 12.656)
\move(0 13.656)\lvec(-4.3 13.656)
\move(0 13.789)\lvec(-4.3 13.789)
\move(0 14.188)\lvec(-4.3 14.188)
\move(0 14.321)\lvec(-4.3 14.321)
%\move(0 12.156)\lvec(0 13.922)
%\move(-1 12.156)\lvec(-1 13.922)
%\move(-2 12.156)\lvec(-2 13.922)
%\move(-3 12.156)\lvec(-3 13.922)
%\move(-4 12.156)\lvec(-4 13.922)
\move(-4 10.889)\lvec(-4 13.922)
\move(-3 10.889)\lvec(-3 13.922)
\move(-2 10.889)\lvec(-2 13.922)
\move(-1 10.889)\lvec(-1 13.922)
\move(0 10.889)\lvec(0 13.922)
\vtext(-0.5 11.906){$\cdots$}
\vtext(-1.5 11.906){$\cdots$}
\vtext(-2.5 11.906){$\cdots$}
\vtext(-3.5 11.906){$\cdots$}
\move(0 14.055)\lvec(0 14.321)
\move(-1 14.055)\lvec(-1 14.321)
\move(-2 14.055)\lvec(-2 14.321)
\move(-3 14.055)\lvec(-3 14.321)
\move(-4 14.055)\lvec(-4 14.321)
\htext(-0.5 13.156){$2$}
\htext(-1.5 13.156){$2$}
\htext(-2.5 13.156){$2$}
\htext(-3.5 13.156){$2$}
\textref h:L v:C
\htext(1.2 1.3325){$a_1\!-\!a_0$}
\htext(1.2 1.9975){$a_0$}
\htext(1.2 3.6625){$a_2$}
\htext(1.2 6.8275){$a_{n\!-\!1}$}
\htext(1.2 7.9925){$a_n$}
\htext(1.2 9.1575){$a_{n\!-\!1}$}
\htext(1.2 10.823){$a_{n\!-\!2}$}
\htext(1.2 13.9885){$a_0$}
\end{texdraw}}%
\savebox{\tmpfigh}{\begin{texdraw}
\fontsize{7}{7}\selectfont
\textref h:C v:C
\drawdim em
\setunitscale 1.9
\move(0 1.665)\lvec(-4.2 1.665)\lvec(-4.2 2.330)\lvec(0 2.330)
\lvec(0 1.665)\ifill f:0.8
\move(0 3.330)\lvec(-4.2 3.330)\lvec(-4.2 3.995)\lvec(0 3.995)
\lvec(0 3.330)\ifill f:0.8
\move(0.5 1.665)\lvec(0.8 1.665)
\move(0.5 2.330)\lvec(0.8 2.330)
\move(0.65 1.665)\lvec(0.65 2.330)%
\move(0.5 3.330)\lvec(0.8 3.330)
\move(0.5 3.995)\lvec(0.8 3.995)
\move(0.65 3.330)\lvec(0.65 3.995)%
\move(0.5 1)\lvec(0.8 1)
\move(0.5 1.665)\lvec(0.8 1.665)
\move(0.65 1)\lvec(0.65 1.665)%
\move(0 0)\lvec(-4.3 0)
\move(0 1)\lvec(-4.3 1)
\move(0 1.133)\lvec(-4.3 1.133)
\move(0 1.532)\lvec(-4.3 1.532)
\move(0 1.665)\lvec(-4.3 1.665)
\move(0 1.798)\lvec(-4.3 1.798)
\move(0 2.197)\lvec(-4.3 2.197)
\move(0 2.330)\lvec(-4.3 2.330)
\move(0 3.330)\lvec(-4.3 3.330)
\move(0 3.463)\lvec(-4.3 3.463)
\move(0 3.862)\lvec(-4.3 3.862)
\move(0 3.995)\lvec(-4.3 3.995)
\move(0 0)\lvec(0 1.266)
\move(-1 0)\lvec(-1 1.266)
\move(-2 0)\lvec(-2 1.266)
\move(-3 0)\lvec(-3 1.266)
\move(-4 0)\lvec(-4 1.266)
\move(0 1.399)\lvec(0 1.931)
\move(-1 1.399)\lvec(-1 1.931)
\move(-2 1.399)\lvec(-2 1.931)
\move(-3 1.399)\lvec(-3 1.931)
\move(-4 1.399)\lvec(-4 1.931)
\move(0 2.064)\lvec(0 3.596)
\move(-1 2.064)\lvec(-1 3.596)
\move(-2 2.064)\lvec(-2 3.596)
\move(-3 2.064)\lvec(-3 3.596)
\move(-4 2.064)\lvec(-4 3.596)
%\move(0 3.729)\lvec(0 4.495)
%\move(-1 3.729)\lvec(-1 4.495)
%\move(-2 3.729)\lvec(-2 4.495)
%\move(-3 3.729)\lvec(-3 4.495)
%\move(-4 3.729)\lvec(-4 4.495)
\move(0 1)\lvec(-1 0)
\move(-1 1)\lvec(-2 0)
\move(-2 1)\lvec(-3 0)
\move(-3 1)\lvec(-4 0)
\htext(-1.75 0.75){$0$}
\htext(-2.75 0.75){$1$}
\htext(-3.75 0.75){$0$}
\htext(-0.75 0.75){$1$}
\htext(-0.5 2.830){$2$}
\htext(-3.5 2.830){$2$}
\htext(-2.5 2.830){$2$}
\htext(-1.5 2.830){$2$}
\move(0 6.495)\lvec(-4.2 6.495)\lvec(-4.2 7.160)\lvec(0 7.160)
\lvec(0 6.495)\ifill f:0.8
\move(0 7.660)\lvec(-4.2 7.660)\lvec(-4.2 8.325)\lvec(0 8.325)
\lvec(0 7.660)\ifill f:0.8
\move(0 8.825)\lvec(-4.2 8.825)\lvec(-4.2 9.490)\lvec(0 9.490)
\lvec(0 8.825)\ifill f:0.8
\move(0 10.490)\lvec(-4.2 10.490)\lvec(-4.2 11.156)\lvec(0 11.156)
\lvec(0 10.490)\ifill f:0.8
\move(0.5 6.495)\lvec(0.8 6.495)
\move(0.5 7.160)\lvec(0.8 7.160)
\move(0.65 6.495)\lvec(0.65 7.160)%
\move(0.5 7.660)\lvec(0.8 7.660)
\move(0.5 8.325)\lvec(0.8 8.325)
\move(0.65 7.660)\lvec(0.65 8.325)%
\move(0.5 8.825)\lvec(0.8 8.825)
\move(0.5 9.490)\lvec(0.8 9.490)
\move(0.65 8.825)\lvec(0.65 9.490)%
\move(0.5 10.490)\lvec(0.8 10.490)
\move(0.5 11.156)\lvec(0.8 11.156)
\move(0.65 10.490)\lvec(0.65 11.156)%
\move(0 5.495)\lvec(-4.3 5.495)
\move(0 6.495)\lvec(-4.3 6.495)
\move(0 6.628)\lvec(-4.3 6.628)
\move(0 7.027)\lvec(-4.3 7.027)
\move(0 7.160)\lvec(-4.3 7.160)
\move(0 7.660)\lvec(-4.3 7.660)
\move(0 7.793)\lvec(-4.3 7.793)
\move(0 8.192)\lvec(-4.3 8.192)
\move(0 8.325)\lvec(-4.3 8.325)
\move(0 8.825)\lvec(-4.3 8.825)
\move(0 8.958)\lvec(-4.3 8.958)
\move(0 9.357)\lvec(-4.3 9.357)
\move(0 9.490)\lvec(-4.3 9.490)
\move(0 10.490)\lvec(-4.3 10.490)
\move(0 10.623)\lvec(-4.3 10.623)
\move(0 11.023)\lvec(-4.3 11.023)
\move(0 11.156)\lvec(-4.3 11.156)
%\move(0 4.995)\lvec(0 6.761)
%\move(-1 4.995)\lvec(-1 6.761)
%\move(-2 4.995)\lvec(-2 6.761)
%\move(-3 4.995)\lvec(-3 6.761)
%\move(-4 4.995)\lvec(-4 6.761)
\move(0 3.729)\lvec(0 6.761)%
\move(-1 3.729)\lvec(-1 6.761)%
\move(-2 3.729)\lvec(-2 6.761)%
\move(-3 3.729)\lvec(-3 6.761)%
\move(-4 3.729)\lvec(-4 6.761)%
\vtext(-0.5 4.745){$\cdots$}
\vtext(-1.5 4.745){$\cdots$}
\vtext(-2.5 4.745){$\cdots$}
\vtext(-3.5 4.745){$\cdots$}
\move(0 6.894)\lvec(0 7.926)
\move(-1 6.894)\lvec(-1 7.926)
\move(-2 6.894)\lvec(-2 7.926)
\move(-3 6.894)\lvec(-3 7.926)
\move(-4 6.894)\lvec(-4 7.926)
\move(-4 8.059)\lvec(-4 9.091)
\move(-1 8.059)\lvec(-1 9.091)
\move(-2 8.059)\lvec(-2 9.091)
\move(-3 8.059)\lvec(-3 9.091)
\move(0 8.059)\lvec(0 9.091)
\move(0 9.224)\lvec(0 10.756)
\move(-1 9.224)\lvec(-1 10.756)
\move(-2 9.224)\lvec(-2 10.756)
\move(-3 9.224)\lvec(-3 10.756)
\move(-4 9.224)\lvec(-4 10.756)
%\move(-4 10.889)\lvec(-4 11.656)
%\move(-3 10.889)\lvec(-3 11.656)
%\move(-2 10.889)\lvec(-2 11.656)
%\move(-1 10.889)\lvec(-1 11.656)
%\move(0 10.889)\lvec(0 11.656)
\htext(-0.5 5.995){$n\!\!-\!\!1$}
\htext(-1.5 5.995){$n\!\!-\!\!1$}
\htext(-2.5 5.995){$n\!\!-\!\!1$}
\htext(-3.5 5.995){$n\!\!-\!\!1$}
\htext(-0.5 7.410){$n$}
\htext(-1.5 7.410){$n$}
\htext(-2.5 7.410){$n$}
\htext(-3.5 7.410){$n$}
\htext(-3.5 8.575){$n$}
\htext(-2.5 8.575){$n$}
\htext(-1.5 8.575){$n$}
\htext(-0.5 8.575){$n$}
\htext(-0.5 9.990){$n\!\!-\!\!1$}
\htext(-1.5 9.990){$n\!\!-\!\!1$}
\htext(-2.5 9.990){$n\!\!-\!\!1$}
\htext(-3.5 9.990){$n\!\!-\!\!1$}
\move(0 13.656)\lvec(-4.2 13.656)\lvec(-4.2 14.321)\lvec(0 14.321)
\ifill f:0.8
\move(0.5 13.656)\lvec(0.8 13.656)
\move(0.5 14.321)\lvec(0.8 14.321)
\move(0.65 13.656)\lvec(0.65 14.321)%
\move(0 12.656)\lvec(-4.3 12.656)
\move(0 13.656)\lvec(-4.3 13.656)
\move(0 13.789)\lvec(-4.3 13.789)
\move(0 14.188)\lvec(-4.3 14.188)
\move(0 14.321)\lvec(-4.3 14.321)
%\move(0 12.156)\lvec(0 13.922)
%\move(-1 12.156)\lvec(-1 13.922)
%\move(-2 12.156)\lvec(-2 13.922)
%\move(-3 12.156)\lvec(-3 13.922)
%\move(-4 12.156)\lvec(-4 13.922)
\move(-4 10.889)\lvec(-4 13.922)
\move(-3 10.889)\lvec(-3 13.922)
\move(-2 10.889)\lvec(-2 13.922)
\move(-1 10.889)\lvec(-1 13.922)
\move(0 10.889)\lvec(0 13.922)
\vtext(-0.5 11.906){$\cdots$}
\vtext(-1.5 11.906){$\cdots$}
\vtext(-2.5 11.906){$\cdots$}
\vtext(-3.5 11.906){$\cdots$}
\move(0 14.055)\lvec(0 14.321)
\move(-1 14.055)\lvec(-1 14.321)
\move(-2 14.055)\lvec(-2 14.321)
\move(-3 14.055)\lvec(-3 14.321)
\move(-4 14.055)\lvec(-4 14.321)
\htext(-0.5 13.156){$2$}
\htext(-1.5 13.156){$2$}
\htext(-2.5 13.156){$2$}
\htext(-3.5 13.156){$2$}
\textref h:L v:C
\htext(1.2 1.3325){$a_0\!-\!a_1$}
\htext(1.2 1.9975){$a_1$}
\htext(1.2 3.6625){$a_2$}
\htext(1.2 6.8275){$a_{n\!-\!1}$}
\htext(1.2 7.9925){$a_n$}
\htext(1.2 9.1575){$a_{n\!-\!1}$}
\htext(1.2 10.823){$a_{n\!-\!2}$}
\htext(1.2 13.9885){$a_1$}
\end{texdraw}}%
\savebox{\tmpfigg}{\begin{texdraw}
\fontsize{7}{7}\selectfont
\textref h:C v:C
\drawdim em
\setunitscale 1.9
\move(0 1.665)\lvec(-4.2 1.665)\lvec(-4.2 2.330)\lvec(0 2.330)
\lvec(0 1.665)\ifill f:0.8
\move(0 3.330)\lvec(-4.2 3.330)\lvec(-4.2 3.995)\lvec(0 3.995)
\lvec(0 3.330)\ifill f:0.8
\move(0.5 1.665)\lvec(0.8 1.665)
\move(0.5 2.330)\lvec(0.8 2.330)
\move(0.65 1.665)\lvec(0.65 2.330)%
\move(0.5 3.330)\lvec(0.8 3.330)
\move(0.5 3.995)\lvec(0.8 3.995)
\move(0.65 3.330)\lvec(0.65 3.995)%
\move(0.5 1)\lvec(0.8 1)
\move(0.5 1.665)\lvec(0.8 1.665)
\move(0.65 1)\lvec(0.65 1.665)%
\move(0 0)\lvec(-4.3 0)
\move(0 1)\lvec(-4.3 1)
\move(0 1.133)\lvec(-4.3 1.133)
\move(0 1.532)\lvec(-4.3 1.532)
\move(0 1.665)\lvec(-4.3 1.665)
\move(0 1.798)\lvec(-4.3 1.798)
\move(0 2.197)\lvec(-4.3 2.197)
\move(0 2.330)\lvec(-4.3 2.330)
\move(0 3.330)\lvec(-4.3 3.330)
\move(0 3.463)\lvec(-4.3 3.463)
\move(0 3.862)\lvec(-4.3 3.862)
\move(0 3.995)\lvec(-4.3 3.995)
\move(0 0)\lvec(0 1.266)
\move(-1 0)\lvec(-1 1.266)
\move(-2 0)\lvec(-2 1.266)
\move(-3 0)\lvec(-3 1.266)
\move(-4 0)\lvec(-4 1.266)
\move(0 1.399)\lvec(0 1.931)
\move(-1 1.399)\lvec(-1 1.931)
\move(-2 1.399)\lvec(-2 1.931)
\move(-3 1.399)\lvec(-3 1.931)
\move(-4 1.399)\lvec(-4 1.931)
\move(0 2.064)\lvec(0 3.596)
\move(-1 2.064)\lvec(-1 3.596)
\move(-2 2.064)\lvec(-2 3.596)
\move(-3 2.064)\lvec(-3 3.596)
\move(-4 2.064)\lvec(-4 3.596)
%\move(0 3.729)\lvec(0 4.495)
%\move(-1 3.729)\lvec(-1 4.495)
%\move(-2 3.729)\lvec(-2 4.495)
%\move(-3 3.729)\lvec(-3 4.495)
%\move(-4 3.729)\lvec(-4 4.495)
\move(0 1)\lvec(-1 0)
\move(-1 1)\lvec(-2 0)
\move(-2 1)\lvec(-3 0)
\move(-3 1)\lvec(-4 0)
\htext(-1.75 0.75){$0$}
\htext(-2.75 0.75){$1$}
\htext(-3.75 0.75){$0$}
\htext(-0.75 0.75){$1$}
\htext(-0.5 2.830){$2$}
\htext(-3.5 2.830){$2$}
\htext(-2.5 2.830){$2$}
\htext(-1.5 2.830){$2$}
\move(0 6.495)\lvec(-4.2 6.495)\lvec(-4.2 7.160)\lvec(0 7.160)
\lvec(0 6.495)\ifill f:0.8
\move(0 7.660)\lvec(-4.2 7.660)\lvec(-4.2 8.325)\lvec(0 8.325)
\lvec(0 7.660)\ifill f:0.8
\move(0 8.825)\lvec(-4.2 8.825)\lvec(-4.2 9.490)\lvec(0 9.490)
\lvec(0 8.825)\ifill f:0.8
\move(0 10.490)\lvec(-4.2 10.490)\lvec(-4.2 11.156)\lvec(0 11.156)
\lvec(0 10.490)\ifill f:0.8
\move(0.5 6.495)\lvec(0.8 6.495)
\move(0.5 7.160)\lvec(0.8 7.160)
\move(0.65 6.495)\lvec(0.65 7.160)%
\move(0.5 7.660)\lvec(0.8 7.660)
\move(0.5 8.325)\lvec(0.8 8.325)
\move(0.65 7.660)\lvec(0.65 8.325)%
\move(0.5 8.825)\lvec(0.8 8.825)
\move(0.5 9.490)\lvec(0.8 9.490)
\move(0.65 8.825)\lvec(0.65 9.490)%
\move(0.5 10.490)\lvec(0.8 10.490)
\move(0.5 11.156)\lvec(0.8 11.156)
\move(0.65 10.490)\lvec(0.65 11.156)%
\move(0 5.495)\lvec(-4.3 5.495)
\move(0 6.495)\lvec(-4.3 6.495)
\move(0 6.628)\lvec(-4.3 6.628)
\move(0 7.027)\lvec(-4.3 7.027)
\move(0 7.160)\lvec(-4.3 7.160)
\move(0 7.660)\lvec(-4.3 7.660)
\move(0 7.793)\lvec(-4.3 7.793)
\move(0 8.192)\lvec(-4.3 8.192)
\move(0 8.325)\lvec(-4.3 8.325)
\move(0 8.825)\lvec(-4.3 8.825)
\move(0 8.958)\lvec(-4.3 8.958)
\move(0 9.357)\lvec(-4.3 9.357)
\move(0 9.490)\lvec(-4.3 9.490)
\move(0 10.490)\lvec(-4.3 10.490)
\move(0 10.623)\lvec(-4.3 10.623)
\move(0 11.023)\lvec(-4.3 11.023)
\move(0 11.156)\lvec(-4.3 11.156)
%\move(0 4.995)\lvec(0 6.761)
%\move(-1 4.995)\lvec(-1 6.761)
%\move(-2 4.995)\lvec(-2 6.761)
%\move(-3 4.995)\lvec(-3 6.761)
%\move(-4 4.995)\lvec(-4 6.761)
\move(0 3.729)\lvec(0 6.761)
\move(-1 3.729)\lvec(-1 6.761)
\move(-2 3.729)\lvec(-2 6.761)
\move(-3 3.729)\lvec(-3 6.761)
\move(-4 3.729)\lvec(-4 6.761)
\vtext(-0.5 4.745){$\cdots$}
\vtext(-1.5 4.745){$\cdots$}
\vtext(-2.5 4.745){$\cdots$}
\vtext(-3.5 4.745){$\cdots$}
\move(0 6.894)\lvec(0 7.926)
\move(-1 6.894)\lvec(-1 7.926)
\move(-2 6.894)\lvec(-2 7.926)
\move(-3 6.894)\lvec(-3 7.926)
\move(-4 6.894)\lvec(-4 7.926)
\move(-4 8.059)\lvec(-4 9.091)
\move(-1 8.059)\lvec(-1 9.091)
\move(-2 8.059)\lvec(-2 9.091)
\move(-3 8.059)\lvec(-3 9.091)
\move(0 8.059)\lvec(0 9.091)
\move(0 9.224)\lvec(0 10.756)
\move(-1 9.224)\lvec(-1 10.756)
\move(-2 9.224)\lvec(-2 10.756)
\move(-3 9.224)\lvec(-3 10.756)
\move(-4 9.224)\lvec(-4 10.756)
%\move(-4 10.889)\lvec(-4 11.656)
%\move(-3 10.889)\lvec(-3 11.656)
%\move(-2 10.889)\lvec(-2 11.656)
%\move(-1 10.889)\lvec(-1 11.656)
%\move(0 10.889)\lvec(0 11.656)
\htext(-0.5 5.995){$n\!\!-\!\!1$}
\htext(-1.5 5.995){$n\!\!-\!\!1$}
\htext(-2.5 5.995){$n\!\!-\!\!1$}
\htext(-3.5 5.995){$n\!\!-\!\!1$}
\htext(-0.5 7.410){$n$}
\htext(-1.5 7.410){$n$}
\htext(-2.5 7.410){$n$}
\htext(-3.5 7.410){$n$}
\htext(-3.5 8.575){$n$}
\htext(-2.5 8.575){$n$}
\htext(-1.5 8.575){$n$}
\htext(-0.5 8.575){$n$}
\htext(-0.5 9.990){$n\!\!-\!\!1$}
\htext(-1.5 9.990){$n\!\!-\!\!1$}
\htext(-2.5 9.990){$n\!\!-\!\!1$}
\htext(-3.5 9.990){$n\!\!-\!\!1$}
\move(0 13.656)\lvec(-4.2 13.656)\lvec(-4.2 14.321)\lvec(0 14.321)
\ifill f:0.8
\move(0.5 13.656)\lvec(0.8 13.656)
\move(0.5 14.321)\lvec(0.8 14.321)
\move(0.65 13.656)\lvec(0.65 14.321)%
\move(0 12.656)\lvec(-4.3 12.656)
\move(0 13.656)\lvec(-4.3 13.656)
\move(0 13.789)\lvec(-4.3 13.789)
\move(0 14.188)\lvec(-4.3 14.188)
\move(0 14.321)\lvec(-4.3 14.321)
%\move(0 12.156)\lvec(0 13.922)
%\move(-1 12.156)\lvec(-1 13.922)
%\move(-2 12.156)\lvec(-2 13.922)
%\move(-3 12.156)\lvec(-3 13.922)
%\move(-4 12.156)\lvec(-4 13.922)
\move(-4 10.889)\lvec(-4 13.922)
\move(-3 10.889)\lvec(-3 13.922)
\move(-2 10.889)\lvec(-2 13.922)
\move(-1 10.889)\lvec(-1 13.922)
\move(0 10.889)\lvec(0 13.922)
\vtext(-0.5 11.906){$\cdots$}
\vtext(-1.5 11.906){$\cdots$}
\vtext(-2.5 11.906){$\cdots$}
\vtext(-3.5 11.906){$\cdots$}
\move(0 14.055)\lvec(0 14.321)
\move(-1 14.055)\lvec(-1 14.321)
\move(-2 14.055)\lvec(-2 14.321)
\move(-3 14.055)\lvec(-3 14.321)
\move(-4 14.055)\lvec(-4 14.321)
\htext(-0.5 13.156){$2$}
\htext(-1.5 13.156){$2$}
\htext(-2.5 13.156){$2$}
\htext(-3.5 13.156){$2$}
\textref h:L v:C
\htext(1.2 1.3325){$a_0\!-\!a_1$}
\htext(1.2 1.9975){$a_1$}
\htext(1.2 3.6625){$a_2$}
\htext(1.2 6.8275){$a_{n\!-\!1}$}
\htext(1.2 7.9925){$2a_n$}
\htext(1.2 9.1575){$a_{n\!-\!1}$}
\htext(1.2 10.823){$a_{n\!-\!2}$}
\htext(1.2 13.9885){$a_1$}
\end{texdraw}}%
\savebox{\tmpfigi}{\begin{texdraw}
\fontsize{7}{7}\selectfont
\textref h:C v:C
\drawdim em
\setunitscale 1.9
\move(0 1.665)\lvec(-4.2 1.665)\lvec(-4.2 2.330)\lvec(0 2.330)
\lvec(0 1.665)\ifill f:0.8
\move(0 3.330)\lvec(-4.2 3.330)\lvec(-4.2 3.995)\lvec(0 3.995)
\lvec(0 3.330)\ifill f:0.8
\move(0.5 1.665)\lvec(0.8 1.665)
\move(0.5 2.330)\lvec(0.8 2.330)
\move(0.65 1.665)\lvec(0.65 2.330)%
\move(0.5 3.330)\lvec(0.8 3.330)
\move(0.5 3.995)\lvec(0.8 3.995)
\move(0.65 3.330)\lvec(0.65 3.995)%
\move(0.5 1)\lvec(0.8 1)
\move(0.5 1.665)\lvec(0.8 1.665)
\move(0.65 1)\lvec(0.65 1.665)%
\move(0 0)\lvec(-4.3 0)
\move(0 1)\lvec(-4.3 1)
\move(0 1.133)\lvec(-4.3 1.133)
\move(0 1.532)\lvec(-4.3 1.532)
\move(0 1.665)\lvec(-4.3 1.665)
\move(0 1.798)\lvec(-4.3 1.798)
\move(0 2.197)\lvec(-4.3 2.197)
\move(0 2.330)\lvec(-4.3 2.330)
\move(0 3.330)\lvec(-4.3 3.330)
\move(0 3.463)\lvec(-4.3 3.463)
\move(0 3.862)\lvec(-4.3 3.862)
\move(0 3.995)\lvec(-4.3 3.995)
\move(0 0)\lvec(0 1.266)
\move(-1 0)\lvec(-1 1.266)
\move(-2 0)\lvec(-2 1.266)
\move(-3 0)\lvec(-3 1.266)
\move(-4 0)\lvec(-4 1.266)
\move(0 1.399)\lvec(0 1.931)
\move(-1 1.399)\lvec(-1 1.931)
\move(-2 1.399)\lvec(-2 1.931)
\move(-3 1.399)\lvec(-3 1.931)
\move(-4 1.399)\lvec(-4 1.931)
\move(0 2.064)\lvec(0 3.596)
\move(-1 2.064)\lvec(-1 3.596)
\move(-2 2.064)\lvec(-2 3.596)
\move(-3 2.064)\lvec(-3 3.596)
\move(-4 2.064)\lvec(-4 3.596)
%\move(0 3.729)\lvec(0 4.495)
%\move(-1 3.729)\lvec(-1 4.495)
%\move(-2 3.729)\lvec(-2 4.495)
%\move(-3 3.729)\lvec(-3 4.495)
%\move(-4 3.729)\lvec(-4 4.495)
\move(0 1)\lvec(-1 0)
\move(-1 1)\lvec(-2 0)
\move(-2 1)\lvec(-3 0)
\move(-3 1)\lvec(-4 0)
\htext(-1.25 0.25){$1$}
\htext(-2.25 0.25){$0$}
\htext(-3.25 0.25){$1$}
\htext(-0.25 0.25){$0$}
\htext(-0.5 2.830){$2$}
\htext(-3.5 2.830){$2$}
\htext(-2.5 2.830){$2$}
\htext(-1.5 2.830){$2$}
\move(0 6.495)\lvec(-4.2 6.495)\lvec(-4.2 7.160)\lvec(0 7.160)
\lvec(0 6.495)\ifill f:0.8
\move(0 7.660)\lvec(-4.2 7.660)\lvec(-4.2 8.325)\lvec(0 8.325)
\lvec(0 7.660)\ifill f:0.8
\move(0 8.825)\lvec(-4.2 8.825)\lvec(-4.2 9.490)\lvec(0 9.490)
\lvec(0 8.825)\ifill f:0.8
\move(0 10.490)\lvec(-4.2 10.490)\lvec(-4.2 11.156)\lvec(0 11.156)
\lvec(0 10.490)\ifill f:0.8
\move(0.5 6.495)\lvec(0.8 6.495)
\move(0.5 7.160)\lvec(0.8 7.160)
\move(0.65 6.495)\lvec(0.65 7.160)%
\move(0.5 7.660)\lvec(0.8 7.660)
\move(0.5 8.325)\lvec(0.8 8.325)
\move(0.65 7.660)\lvec(0.65 8.325)%
\move(0.5 8.825)\lvec(0.8 8.825)
\move(0.5 9.490)\lvec(0.8 9.490)
\move(0.65 8.825)\lvec(0.65 9.490)%
\move(0.5 10.490)\lvec(0.8 10.490)
\move(0.5 11.156)\lvec(0.8 11.156)
\move(0.65 10.490)\lvec(0.65 11.156)%
\move(0 5.495)\lvec(-4.3 5.495)
\move(0 6.495)\lvec(-4.3 6.495)
\move(0 6.628)\lvec(-4.3 6.628)
\move(0 7.027)\lvec(-4.3 7.027)
\move(0 7.160)\lvec(-4.3 7.160)
\move(0 7.660)\lvec(-4.3 7.660)
\move(0 7.793)\lvec(-4.3 7.793)
\move(0 8.192)\lvec(-4.3 8.192)
\move(0 8.325)\lvec(-4.3 8.325)
\move(0 8.825)\lvec(-4.3 8.825)
\move(0 8.958)\lvec(-4.3 8.958)
\move(0 9.357)\lvec(-4.3 9.357)
\move(0 9.490)\lvec(-4.3 9.490)
\move(0 10.490)\lvec(-4.3 10.490)
\move(0 10.623)\lvec(-4.3 10.623)
\move(0 11.023)\lvec(-4.3 11.023)
\move(0 11.156)\lvec(-4.3 11.156)
%\move(0 4.995)\lvec(0 6.761)
%\move(-1 4.995)\lvec(-1 6.761)
%\move(-2 4.995)\lvec(-2 6.761)
%\move(-3 4.995)\lvec(-3 6.761)
%\move(-4 4.995)\lvec(-4 6.761)
\move(0 3.729)\lvec(0 6.761)
\move(-1 3.729)\lvec(-1 6.761)
\move(-2 3.729)\lvec(-2 6.761)
\move(-3 3.729)\lvec(-3 6.761)
\move(-4 3.729)\lvec(-4 6.761)
\vtext(-0.5 4.745){$\cdots$}
\vtext(-1.5 4.745){$\cdots$}
\vtext(-2.5 4.745){$\cdots$}
\vtext(-3.5 4.745){$\cdots$}
\move(0 6.894)\lvec(0 7.926)
\move(-1 6.894)\lvec(-1 7.926)
\move(-2 6.894)\lvec(-2 7.926)
\move(-3 6.894)\lvec(-3 7.926)
\move(-4 6.894)\lvec(-4 7.926)
\move(-4 8.059)\lvec(-4 9.091)
\move(-1 8.059)\lvec(-1 9.091)
\move(-2 8.059)\lvec(-2 9.091)
\move(-3 8.059)\lvec(-3 9.091)
\move(0 8.059)\lvec(0 9.091)
\move(0 9.224)\lvec(0 10.756)
\move(-1 9.224)\lvec(-1 10.756)
\move(-2 9.224)\lvec(-2 10.756)
\move(-3 9.224)\lvec(-3 10.756)
\move(-4 9.224)\lvec(-4 10.756)
%\move(-4 10.889)\lvec(-4 11.656)
%\move(-3 10.889)\lvec(-3 11.656)
%\move(-2 10.889)\lvec(-2 11.656)
%\move(-1 10.889)\lvec(-1 11.656)
%\move(0 10.889)\lvec(0 11.656)
\htext(-0.5 5.995){$n\!\!-\!\!1$}
\htext(-1.5 5.995){$n\!\!-\!\!1$}
\htext(-2.5 5.995){$n\!\!-\!\!1$}
\htext(-3.5 5.995){$n\!\!-\!\!1$}
\htext(-0.5 7.410){$n$}
\htext(-1.5 7.410){$n$}
\htext(-2.5 7.410){$n$}
\htext(-3.5 7.410){$n$}
\htext(-3.5 8.575){$n$}
\htext(-2.5 8.575){$n$}
\htext(-1.5 8.575){$n$}
\htext(-0.5 8.575){$n$}
\htext(-0.5 9.990){$n\!\!-\!\!1$}
\htext(-1.5 9.990){$n\!\!-\!\!1$}
\htext(-2.5 9.990){$n\!\!-\!\!1$}
\htext(-3.5 9.990){$n\!\!-\!\!1$}
\move(0 13.656)\lvec(-4.2 13.656)\lvec(-4.2 14.321)\lvec(0 14.321)
\ifill f:0.8
\move(0.5 13.656)\lvec(0.8 13.656)
\move(0.5 14.321)\lvec(0.8 14.321)
\move(0.65 13.656)\lvec(0.65 14.321)%
\move(0 12.656)\lvec(-4.3 12.656)
\move(0 13.656)\lvec(-4.3 13.656)
\move(0 13.789)\lvec(-4.3 13.789)
\move(0 14.188)\lvec(-4.3 14.188)
\move(0 14.321)\lvec(-4.3 14.321)
%\move(0 12.156)\lvec(0 13.922)
%\move(-1 12.156)\lvec(-1 13.922)
%\move(-2 12.156)\lvec(-2 13.922)
%\move(-3 12.156)\lvec(-3 13.922)
%\move(-4 12.156)\lvec(-4 13.922)
\move(-4 10.889)\lvec(-4 13.922)
\move(-3 10.889)\lvec(-3 13.922)
\move(-2 10.889)\lvec(-2 13.922)
\move(-1 10.889)\lvec(-1 13.922)
\move(0 10.889)\lvec(0 13.922)
\vtext(-0.5 11.906){$\cdots$}
\vtext(-1.5 11.906){$\cdots$}
\vtext(-2.5 11.906){$\cdots$}
\vtext(-3.5 11.906){$\cdots$}
\move(0 14.055)\lvec(0 14.321)
\move(-1 14.055)\lvec(-1 14.321)
\move(-2 14.055)\lvec(-2 14.321)
\move(-3 14.055)\lvec(-3 14.321)
\move(-4 14.055)\lvec(-4 14.321)
\htext(-0.5 13.156){$2$}
\htext(-1.5 13.156){$2$}
\htext(-2.5 13.156){$2$}
\htext(-3.5 13.156){$2$}
\textref h:L v:C
\htext(1.2 1.3325){$a_1\!-\!a_0$}
\htext(1.2 1.9975){$a_0$}
\htext(1.2 3.6625){$a_2$}
\htext(1.2 6.8275){$a_{n\!-\!1}$}
\htext(1.2 7.9925){$2a_n$}
\htext(1.2 9.1575){$a_{n\!-\!1}$}
\htext(1.2 10.823){$a_{n\!-\!2}$}
\htext(1.2 13.9885){$a_0$}
\end{texdraw}}%
\begin{enumerate}
\item $A_n^{(1)}$ $(n\ge 1)$\\[2mm]
  \mbox{} \hspace{15mm} \usebox{\tmpfiga}\\
\item $B_{n}^{(1)}$ $(n\ge 3)$\\[2mm]
  \mbox{} \hspace{15mm} \usebox{\tmpfigf} \qquad or \qquad \usebox{\tmpfigh}\\[2mm]
The choice between the two depends on
whether $a_1\geq a_0$ or $a_1\le a_0$.\\
\item $C_n^{(1)}$ $(n\ge 2)$\\[2mm]
  \mbox{} \hspace{15mm} \usebox{\tmpfigb}\\
\item $A_{2n-1}^{(2)}$$(n\ge 3)$\\[2mm]
  \mbox{} \hspace{15mm} \usebox{\tmpfigi} \qquad or \qquad \usebox{\tmpfigg}\\[2mm]
The choice between the two depends on
whether $a_1\geq a_0$ or $a_1\le a_0$.\\
\item $A_{2n}^{(2)}$ $(n\ge 1)$\\[2mm]
  \mbox{} \hspace{15mm} \usebox{\tmpfigc}\\
\item $D_{n+1}^{(2)}$ $(n\ge 2)$\\[2mm]
  \mbox{} \hspace{15mm} \usebox{\tmpfige}\\
\end{enumerate}

\vspace{2mm}

A level-$l$ proper Young wall obtained by adding finitely many
blocks to the ground-state wall $\gswall_{\la}$ is said to have
been \defi{built on $\gswall_{\la}$}. We denote by $\pwspace(\la)$
(resp. $\rpwspace(\la)$) the set of all proper Young walls (resp.
reduced proper Young walls) built on $\gswall_{\la}$.

For each $\mathbf{Y} \in \pwspace(\la)$, we define its
\defi{affine weight} to be
\begin{equation*}
\wt(\mathbf{Y}) = \la - \sum_{i=0}^{n} k_i \ali,
\end{equation*}
where $k_i$ is the number of $i$-blocks
(or half the number when dealing with $0$-blocks of
$\cnone$ type and $n$-blocks of $\atwonmonetwo$ type)
that have been added to
$\gswall_{\la}$. Then we obtain\,:

\begin{prop}
The set $\pwspace(\la)$ of all level-$l$ proper Young walls built
on $\gswall_{\la}$, together with the maps $\eit$, $\fit$,
$\veps_i$, $\vphi_i$ $(i\in I)$, $\wt$, forms a $\uq(\g)$-crystal.
\end{prop}

Recall that the set $\pathspace(\la)$ of all $\la$-paths gives a
realization of the irreducible highest weight crystal $\hwc(\la)$
(see Section 2). We will show that there exists a
$\uq(\g)$-crystal isomorphism $\Phi: \rpwspace(\la)
\overset{\sim}\longrightarrow \pathspace(\la)$.

Let $\mathbf{Y} = (\mathbf{Y}(k))_{k=0}^{\infty}$ be a reduced
proper Young wall built on $\mathbf{Y}_{\la}$. Using the
$\uq'(\g)$-crystal isomorphism $\psi: \oldcrystal \overset{\sim}
\longrightarrow \newcrystal$ given in Theorem~\ref{thm:37} , we get a map
$\Phi: \rpwspace(\la) \longrightarrow \pathspace(\la)$ which is
defined by
\begin{equation} \label{eq:main}
\Phi(\mathbf{Y}) = (\psi^{-1}(\mathbf{Y}(k)))_{k=0}^{\infty}.
\end{equation}
Note that the ground-state wall $\gswall_{\la}$ is mapped onto the
ground-state path $\path_{\la}$.

Conversely, to each $\la$-path $\path = (\path(k))_{k=0}^{\infty}
\in \pathspace(\la)$, we can associate a wall
$\mathbf{Y}=(\mathbf{Y}(k))_{k=0}^{\infty}$ built on
$\mathbf{Y}_{\la}$ such that $\psi(\path(k)) = \mathbf{Y}(k)$ for
all $k \ge 0$. By adding or removing an appropriate number of
$\delta$'s (from left to right), one can easily see that there
exists a unique reduced proper Young wall
$\mathbf{Y}=(\mathbf{Y}(k))_{k=0}^{\infty}$ with this property,
which shows that $\Phi: \rpwspace(\la) \longrightarrow
\pathspace(\la)$ is a bijection.

Our main result is the following realization theorem.

\begin{thm}\label{thm:main thm}
The bijection $\Phi: \rpwspace(\la) \longrightarrow
\pathspace(\la)$ defined by \eqref{eq:main} is a $\uq(\g)$-crystal
isomorphism.
Therefore, we have a $\uq(\g)$-crystal isomorphism
\begin{equation*}
\rpwspace(\la) \overset{\sim} \longrightarrow \pathspace(\la)
\overset {\sim} \longrightarrow \hwc(\la).
\end{equation*}
\end{thm}

We shall focus
our efforts on showing that
the set $\rpwspace(\la)$ is a $\uq(\mathfrak{g})$-subcrystal
of $\pwspace(\la)$ and that
the map $\Phi$ commutes with the Kashiwara
operators $\fit$ and $\eit$. Other parts of the proof are similar or easy.

Now, recalling the fact that the map $\Phi$ is
determined by the $\uq'(\g)$-crystal isomorphism $\psi:
\oldcrystal \overset{\sim} \longrightarrow \newcrystal$,
we find that, to prove these two points,
it suffices to prove the following two lemmas.

\begin{lem}\label{lem:55}
Let $\mathbf{Y}=(\mathbf{Y}(k))_{k=0}^{\infty}$ be a reduced
proper Young wall in $\rpwspace(\la)$ and let $\path =
\Phi(\mathbf{Y})$ be the corresponding $\la$-path in
$\pathspace(\la)$. Then we have
\begin{enumerate}
\item The Kashiwara operator $\eit$ \textup{(}resp. $\fit$\textup{)}
acts on the $j$-th column of
$\mathbf{Y}$ if and only if $\eit$ \textup{(}resp. $\fit$\textup{)} acts on the
$j$-th component of $\path$.
\item $\eit \mathbf{Y}=0$ \textup{(}resp. $\fit \mathbf{Y} = 0$\textup{)}
if and only if $\eit
\path=0$ \textup{(}resp. $\fit \path =0$\textup{)}.
\end{enumerate}
\end{lem}

\begin{lem}\label{lem:56}
The action of Kashiwara operators on $\rpwspace(\la)$
satisfies the following properties \textup{:}
\begin{equation*}
\eit \rpwspace(\la) \subset \rpwspace(\la) \cup \{0\}, \qquad \fit
\rpwspace(\la) \subset \rpwspace(\la) \cup \{0\} \qquad \text{for
all} \ \ i\in I.
\end{equation*}
\end{lem}
\begin{proof}
Suppose that there exists some $\mathbf{Y}\in\rpwspace(\la)$ for which
$\fit \mathbf{Y} \not\in \rpwspace(\la) \cup \{0\}$. We assume that $\fit$
has acted on the $j$-th column of $\mathbf{Y}$ and set $\path =
\spaceosi(\mathbf{Y})$. Then, by Lemma~\ref{lem:55} the action of
$\fit$ on $\path$ would also have been on the $j$-th tensor
component of $\path$.

Since $\fit \mathbf{Y} \in \pwspace(\la)$,
we may remove finitely many $\delta$'s
from $\fit \mathbf{Y}$ to obtain a reduce proper Young wall
$\mathbf{Y}'$. The number of $\delta$'s removed is nonzero since $\fit \mathbf{Y}
\not\in \rpwspace(\la)$.
Since every columns of $\fit \mathbf{Y}$ is related to
the corresponding column of $\mathbf{Y}'$ under the previously defined
equivalence relation, we have
\begin{equation*}
\fit(\path) = \spaceosi(\mathbf{Y}').
\end{equation*}

Let us apply $\eit$ to both $\mathbf{Y}'$ and $\fit(\path)$. We have $\path
= \eit(\fit(\path))$ and the action of $\eit$ on $\fit(\path)$
would have been on the $j$-th tensor component. By
Lemma~\ref{lem:55}, the action of $\eit$ of $\mathbf{Y}'$ will also be on
the $j$-th column.
Hence
the proper Young wall $\eit \mathbf{Y}'$ may be obtained from the
reduced proper Young wall $\mathbf{Y}$ by removing finitely many
$\delta$'s.
We may now remove finitely many $\delta$'s from $\eit
\mathbf{Y}'$ to obtain a reduced proper Young wall $\mathbf{Y}''$
which also
corresponds to $\path$ under the map
$\spaceosi$.

Recall that we started out with a reduced proper Young wall $\mathbf{Y}$,
added an $i$-block to the $j$-th column of $\mathbf{Y}$ to obtain $\fit \mathbf{Y}$,
removed finitely many $\delta$'s from $\fit \mathbf{Y}$
to obtain a reduced proper Young wall $\mathbf{Y}'$, and removed an $i$-block
from the $j$-th column of $\mathbf{Y}'$ to obtain $\eit \mathbf{Y}'$,
and finally, removed finitely
many $\delta$'s from $\eit \mathbf{Y}'$ to obtain a reduced
proper Young wall $\mathbf{Y}''$ which corresponds to $\path$ under the map
$\spaceosi$.
Therefore, we have
$\Phi(\mathbf{Y})=\path=\Phi(\mathbf{Y''})$, but $\mathbf{Y} \neq
\mathbf{Y''}$, which is a contradiction. Hence $\fit \mathbf{Y}$
must be reduced.

Similarly, one can show that $\eit \rpwspace(\la) \subset
\rpwspace(\la) \cup \{0\}$ for all $i\in I$, which completes the
proof of our claim.
\end{proof}

The rest of this section is devoted to proving Lemma \ref{lem:55}.
We first fix some notations. Let us consider two consecutive
columns that form a part of a Young wall, or two consecutive
perfect crystal elements from a path. We shall denote by
left-$\vphi$, the number of $0$'s to be written under the left
column or crystal element when preparing for the Kashiwara
operators. Similarly, left-$\veps$ denotes the number of $1$'s to
be written under the left item. Right-$\vphi$ and right-$\veps$
are defined in a similar manner.
In the case of Young walls, left-$\veps$ depends on the column
that sits to the left of the left column and right-$\vphi$ depends
on the column that comes to the right of the right column. So it
is not possible to find the values left-$\veps$ and right-$\vphi$
from the two Young wall columns. Hence, a straightforward
comparison of the signatures between a Young wall and the
corresponding path after cancellations of $(0,1)$ pairs is not
possible. But still, we can verify that what is left of the
left-$\vphi$ and right-$\veps$ signatures, after the $(0,1)$-pair
cancellation, is the same for the path and the Young wall.

\begin{proof}[Proof of Lemma~\ref{lem:55}]
We give a proof of our claim only for the $\bnone$ case. Other
cases may be proved in a similar manner and are less complicated.

It suffices to check that, for all possible left-right pairs of
perfect crystal elements and their corresponding Young wall
columns, what remains after $(0,1)$-cancellation of left-$\vphi$
and right-$\veps$ signatures agrees. (The right-most column may be
dealt with in a similar way.)

We shall deal with the $i=n$ case first, and the remaining cases
will be covered in the Appendices.

Let us write the general level-$l$ Young wall column(or, slice) in the
following form.\\
%We have centered our attention on the $n$-blocks.
%Notice that any level-$l$ Young wall column is
%\emph{related} to this slice.
%
\savebox{\tmpfiga}{
\begin{texdraw}
\fontsize{6}{6}\selectfont
\textref h:C v:C
\drawdim em
\setunitscale 1.7
\move(0 0)
\bsegment
\move(0 -0.75)\lvec(0 3)\lvec(1 3)\lvec(1 2)
\move(0 2.5)\lvec(1 2.5)
\move(0 2)\lvec(1 2)
\htext(0.5 2.25){$n$}
\htext(0.5 2.75){$n$}
\esegment
\move(1 0)
\bsegment
\move(0 2)\lvec(0 3)\lvec(1.4 3)\lvec(1.4 2)
\move(0 2.5)\lvec(1.4 2.5)
\move(0 2)\lvec(1.4 2)
\htext(0.7 2.25){$\cdots$}
\htext(0.7 2.75){$\cdots$}
\esegment
\move(2.4 0)
\bsegment
\move(0 2)\lvec(0 3)\lvec(1 3)\lvec(1 -0.75)
\move(0 2.5)\lvec(1 2.5)
\move(0 2)\lvec(1 2)
\htext(0.5 2.25){$n$}
\htext(0.5 2.75){$n$}
\esegment
\move(3.4 0)
\bsegment
\move(0 -0.75)\lvec(0 2.5)\lvec(1 2.5)\lvec(1 1)
\move(0 1)\lvec(1 1)
\move(0 2)\lvec(1 2)
\htext(0.5 2.25){$n$}
\htext(0.5 1.5){$n\!\!-\!\!1$}
\esegment
\move(4.4 0)
\bsegment
\move(0 1)\lvec(0 2.5)\lvec(1.4 2.5)\lvec(1.4 1)
\move(0 1)\lvec(1.4 1) \move(0 2)\lvec(1.4 2)
\htext(0.7 2.25){$\cdots$}
\htext(0.7 1.5){$\cdots$}
\esegment
\move(5.8 0)
\bsegment
\move(0 1)\lvec(0 2.5)\lvec(1 2.5)\lvec(1 -0.75)
\move(0 1)\lvec(1 1)
\move(0 2)\lvec(1 2)
\htext(0.5 2.25){$n$}
\htext(0.5 1.5){$n\!\!-\!\!1$}
\esegment
\move(6.8 0)
\bsegment
\move(0 -0.75)\lvec(0 2)\lvec(1 2)\lvec(1 0)
\move(0 1)\lvec(1 1)
\move(0 0)\lvec(1 0)
\htext(0.5 0.5){$n\!\!-\!\!2$}
\htext(0.5 1.5){$n\!\!-\!\!1$}
\esegment
\move(7.8 0)
\bsegment
\move(0 0)\lvec(0 2)\lvec(1.4 2)\lvec(1.4 0)
\move(0 1)\lvec(1.4 1)
\move(0 0)\lvec(1.4 0)
\htext(0.7 0.5){$\cdots$}
\htext(0.7 1.5){$\cdots$}
\esegment
\move(9.2 0)
\bsegment
\move(0 0)\lvec(0 2)\lvec(1 2)\lvec(1 -0.75)
\move(0 1)\lvec(1 1) \move(0 0)\lvec(1 0)
\htext(0.5 0.5){$n\!\!-\!\!2$}
\htext(0.5 1.5){$n\!\!-\!\!1$}
\esegment
\move(10.2 0)
\bsegment
\move(0 -0.75)\lvec(0 1)\lvec(3.4 1)\lvec(3.4 -0.75)
\esegment
\htext(1.7 -1.25){$\underbrace{\rule{5.7em}{0em}}$}
\htext(5.1 -1.25){$\underbrace{\rule{5.7em}{0em}}$}
\htext(8.5 -1.25){$\underbrace{\rule{5.7em}{0em}}$}
\htext(11.9 -1.25){$\underbrace{\rule{5.7em}{0em}}$}
\htext(1.7 -2){$a_1$}
\htext(5.1 -2){$a_2$}
\htext(8.5 -2){$a_3$}
\htext(11.9 -2){$a_4$}
\move(0 -2.3)
%\drawbb
\end{texdraw}}
\begin{center}
\begin{minipage}{0.5\textwidth}
\usebox{\tmpfiga}
\end{minipage}
\quad
\begin{minipage}{0.15\textwidth}
\setlength{\tabcolsep}{0.8mm}
%\begin{tabular}{rcl}
%$\sum_{i=1}^{4} a_i$ &:& $l$\\
%$\sum_{i=1}^{4} b_i$ &:& $l$\\
%notation &:& L
%\end{tabular}
$l = \sum_{i=1}^{4} a_i$
\end{minipage}
\end{center}
Here, $a_1, a_2, a_3$ denote the number of layers having top
blocks of the form given in the figure, and $a_4$ denotes the
number of all other layers; that is, any layer having a top block
that comes between the supporting $(n-2)$-block and the covering
$(n-1)$-block (inclusive) is counted in $a_4$.
%It is clear that $\sum_i a_i = l$.
%Columns that are \emph{related} (Definition~\ref{df:35})
%are denoted by the same notation.

We break this into two cases and fix the notations for the two
cases.
\begin{itemize}
\item $a_1 \geq a_3$ : \Na
\item $a_1 \leq a_3$ : \Nb
\end{itemize}
The same notation will be used to denote any other layer-rotations
of these slices. Also the perfect crystal elements corresponding
to these slices will be denoted by the same notations.

If we split every block possible from these slices, the
result will be of the following two shapes.\\[2mm]
\savebox{\tmpfiga}
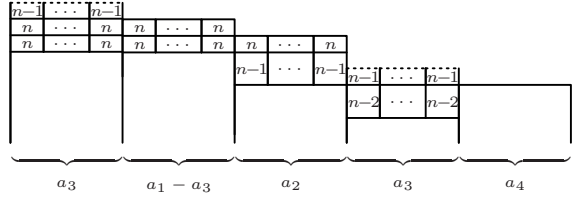
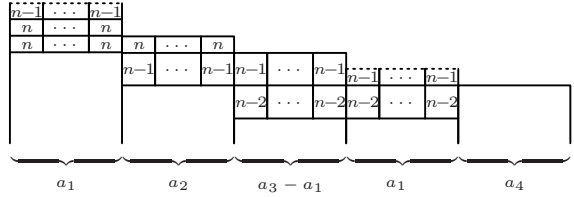
{
\begin{texdraw}
\fontsize{6}{6}\selectfont
\textref h:C v:C
\drawdim em
\setunitscale 1.7
\move(-3.4 5)
\bsegment
\move(0 -5.75)\lvec(0 -1.5)
\move(0 -5.75)\lvec(0 -1.5)
\move(1 -1.5)\lvec(1 -3)
\move(0 -2)\lvec(1 -2)
\move(0 -3)\lvec(1 -3)
\move(0 -2.5)\lvec(1 -2.5)
\htext(0.5 -2.75){$n$}
\htext(0.5 -2.25){$n$}
\htext(0.5 -1.75){$n\!\!-\!\!1$}
\lpatt(0.05 0.15)
\move(0 -1.5)\lvec(1 -1.5)
\esegment
\move(-2.4 5)
\bsegment
\move(0 -3)\lvec(0 -1.5)
\move(1.4 -1.5)\lvec(1.4 -3)
\move(0 -2)\lvec(1.4 -2)
\move(0 -3)\lvec(1.4 -3)
\move(0 -2.5)\lvec(1.4 -2.5)
\htext(0.7 -1.75){$\cdots$}
\htext(0.7 -2.25){$\cdots$}
\htext(0.7 -2.75){$\cdots$}
\lpatt(0.05 0.15)
\move(0 -1.5)\lvec(1.4 -1.5)
\esegment
\move(-1 5)
\bsegment
\move(0 -3)\lvec(0 -1.5)
\move(1 -1.5)\lvec(1 -3.75)
\move(0 -2)\lvec(1 -2)
\move(0 -3)\lvec(1 -3)
\move(0 -2.5)\lvec(1 -2.5)
\htext(0.5 -2.75){$n$}
\htext(0.5 -2.25){$n$}
\htext(0.5 -1.75){$n\!\!-\!\!1$}
\lpatt(0.05 0.15)
\move(0 -1.5)\lvec(1 -1.5)
\esegment%%%%%%%%%%%
\move(0 0)
\bsegment
\move(0 -0.75)\lvec(0 3)\lvec(1 3)\lvec(1 2)
\move(0 2.5)\lvec(1 2.5)
\move(0 2)\lvec(1 2)
\htext(0.5 2.25){$n$}
\htext(0.5 2.75){$n$}
\esegment
\move(1 0)
\bsegment
\move(0 2)\lvec(0 3)\lvec(1.4 3)\lvec(1.4 2)
\move(0 2.5)\lvec(1.4 2.5)
\move(0 2)\lvec(1.4 2)
\htext(0.7 2.25){$\cdots$}
\htext(0.7 2.75){$\cdots$}
\esegment
\move(2.4 0)
\bsegment
\move(0 2)\lvec(0 3)\lvec(1 3)\lvec(1 -0.5)
\move(0 2.5)\lvec(1 2.5)
\move(0 2)\lvec(1 2)
\htext(0.5 2.25){$n$}
\htext(0.5 2.75){$n$}
\esegment%%%%%%%%%%%
\move(3.4 0)
\bsegment
\move(0 -0.75)\lvec(0 2.5)\lvec(1 2.5)\lvec(1 1)
\move(0 1)\lvec(1 1)
\move(0 2)\lvec(1 2)
\htext(0.5 2.25){$n$}
\htext(0.5 1.5){$n\!\!-\!\!1$}
\esegment
\move(4.4 0)
\bsegment
\move(0 1)\lvec(0 2.5)\lvec(1.4 2.5)\lvec(1.4 1)
\move(0 1)\lvec(1.4 1)
\move(0 2)\lvec(1.4 2)
\htext(0.7 2.25){$\cdots$}
\htext(0.7 1.5){$\cdots$}
\esegment
\move(5.8 0)
\bsegment
\move(0 1)\lvec(0 2.5)\lvec(1 2.5)\lvec(1 -0.75)
\move(0 1)\lvec(1 1)
\move(0 2)\lvec(1 2)
\htext(0.5 2.25){$n$}
\htext(0.5 1.5){$n\!\!-\!\!1$}
\esegment%%%%%%%%%%%
\move(6.8 0)
\bsegment
\move(0 -0.75)\lvec(0 1.5)
\move(1 1.5)\lvec(1 0)
\move(0 1)\lvec(1 1)
\move(0 0)\lvec(1 0)
\htext(0.5 0.5){$n\!\!-\!\!2$}
\htext(0.5 1.25){$n\!\!-\!\!1$}
\lpatt(0.05 0.15)
\move(0 1.5)\lvec(1 1.5)
\esegment
\move(7.8 0)
\bsegment
\move(0 0)\lvec(0 1.5)
\move(1.4 1.5)\lvec(1.4 1)
\move(0 1)\lvec(1.4 1)
\move(0 0)\lvec(1.4 0)
\htext(0.7 0.5){$\cdots$}
\htext(0.7 1.25){$\cdots$}
\lpatt(0.05 0.15)
\move(0 1.5)\lvec(1.4 1.5)
\esegment
\move(9.2 0)
\bsegment
\move(0 0)\lvec(0 1.5)
\move(1 1.5)\lvec(1 -0.75)
\move(0 1)\lvec(1 1)
\move(0 0)\lvec(1 0)
\htext(0.5 0.5){$n\!\!-\!\!2$}
\htext(0.5 1.25){$n\!\!-\!\!1$}
\lpatt(0.05 0.15)
\move(0 1.5)\lvec(1 1.5)
\esegment%%%%%%%%%%%%%%%%%
\move(10.2 0)
\bsegment
\move(0 -0.75)\lvec(0 1)\lvec(3.4 1)\lvec(3.4 -0.75)
\esegment%%%%%%%%%%%%%%%%%%%%
\htext(1.7 -1.25){$\underbrace{\rule{5.7em}{0em}}$}
\htext(5.1 -1.25){$\underbrace{\rule{5.7em}{0em}}$}
\htext(8.5 -1.25){$\underbrace{\rule{5.7em}{0em}}$}
\htext(11.9 -1.25){$\underbrace{\rule{5.7em}{0em}}$}
\htext(-1.7 -1.25){$\underbrace{\rule{5.7em}{0em}}$}
\htext(-1.7 -2){$a_3$}
\htext(1.7 -2){$a_1-a_3$}
\htext(5.1 -2){$a_2$}
\htext(8.5 -2){$a_3$}
\htext(11.9 -2){$a_4$}
\move(0 -2.3)
%\drawbb
\end{texdraw}
}
\savebox{\tmpfigb}{
\begin{texdraw}
\fontsize{6}{6}\selectfont
\textref h:C v:C
\drawdim em
\setunitscale 1.7
\move(-3.4 5)
\bsegment
\move(0 -5.75)\lvec(0 -1.5)
\move(1 -1.5)\lvec(1 -3)
\move(0 -2)\lvec(1 -2)
\move(0 -3)\lvec(1 -3)
\move(0 -2.5)\lvec(1 -2.5)
\htext(0.5 -2.75){$n$}
\htext(0.5 -2.25){$n$}
\htext(0.5 -1.75){$n\!\!-\!\!1$}
\lpatt(0.05 0.15)
\move(0 -1.5)\lvec(1 -1.5)
\esegment
\move(-2.4 5)
\bsegment
\move(0 -3)\lvec(0 -1.5)
\move(1.4 -1.5)\lvec(1.4 -3)
\move(0 -2)\lvec(1.4 -2)
\move(0 -3)\lvec(1.4 -3)
\move(0 -2.5)\lvec(1.4 -2.5)
\htext(0.7 -1.75){$\cdots$}
\htext(0.7 -2.25){$\cdots$}
\htext(0.7 -2.75){$\cdots$}
\lpatt(0.05 0.15)
\move(0 -1.5)\lvec(1.4 -1.5)
\esegment
\move(-1 5)
\bsegment
\move(0 -3)\lvec(0 -1.5)
\move(1 -1.5)\lvec(1 -3.75)
\move(0 -2)\lvec(1 -2)
\move(0 -3)\lvec(1 -3)
\move(0 -2.5)\lvec(1 -2.5)
\htext(0.5 -2.75){$n$}
\htext(0.5 -2.25){$n$}
\htext(0.5 -1.75){$n\!\!-\!\!1$}
\lpatt(0.05 0.15)
\move(0 -1.5)\lvec(1 -1.5)
\esegment%%%%%%%%%%%
\move(0 0)
\bsegment
\move(0 -0.75)\lvec(0 2.5)\lvec(1 2.5)\lvec(1 1)
\move(0 1)\lvec(1 1)
\move(0 2)\lvec(1 2)
\htext(0.5 2.25){$n$}
\htext(0.5 1.5){$n\!\!-\!\!1$}
\esegment
\move(1 0)
\bsegment
\move(0 1)\lvec(0 2.5)\lvec(1.4 2.5)\lvec(1.4 1)
\move(0 1)\lvec(1.4 1)
\move(0 2)\lvec(1.4 2)
\htext(0.7 2.25){$\cdots$}
\htext(0.7 1.5){$\cdots$}
\esegment
\move(2.4 0)
\bsegment
\move(0 1)\lvec(0 2.5)\lvec(1 2.5)\lvec(1 -0.75)
\move(0 1)\lvec(1 1)
\move(0 2)\lvec(1 2)
\htext(0.5 2.25){$n$}
\htext(0.5 1.5){$n\!\!-\!\!1$}
\esegment%%%%%%%%%%%%%%%%
\move(3.4 0)
\bsegment
\move(0 -0.75)\lvec(0 2)\lvec(1 2)\lvec(1 0)
\move(0 1)\lvec(1 1)
\move(0 0)\lvec(1 0)
\htext(0.5 0.5){$n\!\!-\!\!2$}
\htext(0.5 1.5){$n\!\!-\!\!1$}
\esegment
\move(4.4 0)
\bsegment
\move(0 0)\lvec(0 2)\lvec(1.4 2)\lvec(1.4 0)
\move(0 1)\lvec(1.4 1)
\move(0 0)\lvec(1.4 0)
\htext(0.7 0.5){$\cdots$}
\htext(0.7 1.5){$\cdots$}
\esegment
\move(5.8 0)
\bsegment
\move(0 0)\lvec(0 2)\lvec(1 2)\lvec(1 -0.75)
\move(0 1)\lvec(1 1) \move(0 0)\lvec(1 0)
\htext(0.5 0.5){$n\!\!-\!\!2$}
\htext(0.5 1.5){$n\!\!-\!\!1$}
\esegment%%%%%%%%%%%%%
\move(6.8 0)
\bsegment
\move(0 -0.75)\lvec(0 1.5)
\move(1 1.5)\lvec(1 0)
\move(0 1)\lvec(1 1)
\move(0 0)\lvec(1 0)
\htext(0.5 0.5){$n\!\!-\!\!2$}
\htext(0.5 1.25){$n\!\!-\!\!1$}
\lpatt(0.05 0.15)
\move(0 1.5)\lvec(1 1.5)
\esegment
\move(7.8 0)
\bsegment
\move(0 0)\lvec(0 1.5)
\move(1.4 1.5)\lvec(1.4 1)
\move(0 1)\lvec(1.4 1)
\move(0 0)\lvec(1.4 0)
\htext(0.7 0.5){$\cdots$}
\htext(0.7 1.25){$\cdots$}
\lpatt(0.05 0.15)
\move(0 1.5)\lvec(1.4 1.5)
\esegment
\move(9.2 0)
\bsegment
\move(0 0)\lvec(0 1.5)
\move(1 1.5)\lvec(1 -0.75)
\move(0 1)\lvec(1 1)
\move(0 0)\lvec(1 0)
\htext(0.5 0.5){$n\!\!-\!\!2$}
\htext(0.5 1.25){$n\!\!-\!\!1$}
\lpatt(0.05 0.15)
\move(0 1.5)\lvec(1 1.5)
\esegment%%%%%%%%%%%%%%%%%
\move(10.2 0)
\bsegment
\move(0 -0.75)\lvec(0 1)\lvec(3.4 1)\lvec(3.4 -0.75)
\esegment%%%%%%%%%%%%%
\htext(1.7 -1.25){$\underbrace{\rule{5.7em}{0em}}$}
\htext(5.1 -1.25){$\underbrace{\rule{5.7em}{0em}}$}
\htext(8.5 -1.25){$\underbrace{\rule{5.7em}{0em}}$}
\htext(11.9 -1.25){$\underbrace{\rule{5.7em}{0em}}$}
\htext(-1.7 -1.25){$\underbrace{\rule{5.7em}{0em}}$}
\htext(-1.7 -2){$a_1$}
\htext(1.7 -2){$a_2$}
\htext(5.1 -2){$a_3-a_1$}
\htext(8.3 -2){$a_1$}
\htext(11.9 -2){$a_4$}
\move(0 -2.3)
%\drawbb
\end{texdraw}}
\noindent case \Na\ :
\begin{center}
%\begin{minipage}{0.68\textwidth}
\usebox{\tmpfiga}
%\end{minipage}
%\begin{minipage}{0.28\textwidth}
%\setlength{\tabcolsep}{0.8mm}
%\begin{tabular}{rcl}
%$\veps$ &:& $b_2+2(b_1-b_3)\ge 1$\\
%$\vphi$ &:& $a_2\ge 1$\\
%notation &:& L1
%\end{tabular}
%\end{minipage}\\[5mm]
\end{center}
\noindent case \Nb\ :
\begin{center}
\usebox{\tmpfigb}
\end{center}

Now, we will place two slices
%Young wall columns
side by side and also consider the corresponding pair of perfect
crystal elements. When we use the above notations \Na\ and \Nb\
for the \emph{right} of the two slices, we will take the number of
layers to be given by $b_i$ instead of $a_i$. For the left slices,
we will use $a_i$ as given in the figure.

Given any two slices, we may either add or remove finitely many
$\delta$'s to or from either of the two slices so that the two may
be considered as a part of a reduced proper Young wall. We shall
take the following convention. First fix the right slice and next
remove finitely many $\delta$'s to the left slice so that the two
may be considered as a part of a reduced proper Young wall. To
express exactly how many times this layer rotation(or
$\delta$-removal) has taken place, we first need to designate a
starting point for the two slices.

To make things easier later on, we choose the following starting
shapes and relative heights for the two columns. The shape of
right slice should be so that all layers with $n$-blocks at the
top are placed at the rear, when every block possible in it has
been split. Also the starting shape of the left slice should be so
that all layers with $n$-slots at the top are placed at the rear,
when every block possible in it has been split. Below is an
example showing right slice \Na and left slice \Nb. We have given
the figures with every possible block split,
for the readers' convenience.\\[2mm]
\savebox{\tmpfiga}{
\begin{texdraw}
\fontsize{6}{6}\selectfont
\textref h:C v:C
\drawdim em
\setunitscale 1.7
\move(0 0)
\bsegment
\move(0 -0.5)\lvec(0 3)\lvec(1 3)\lvec(1 2)
\move(0 2.5)\lvec(1 2.5)
\move(0 2)\lvec(1 2)
\htext(0.5 2.25){$n$}
\htext(0.5 2.75){$n$}
\move(0 -1)\lvec(0 -3.75)
\esegment
\move(1 0)
\bsegment
\move(0 2)\lvec(0 3)\lvec(1.4 3)\lvec(1.4 2)
\move(0 2.5)\lvec(1.4 2.5)
\move(0 2)\lvec(1.4 2)
\htext(0.7 2.25){$\cdots$}
\htext(0.7 2.75){$\cdots$}
\esegment
\move(2.4 0)
\bsegment
\move(0 2)\lvec(0 3)\lvec(1 3)\lvec(1 -0.5)
\move(0 2.5)\lvec(1 2.5)
\move(0 2)\lvec(1 2)
\htext(0.5 2.25){$n$}
\htext(0.5 2.75){$n$}
\esegment%%%%%%%%%%%
\move(3.4 0) \bsegment \move(0 -0.5)\lvec(0 2.5)\lvec(1
2.5)\lvec(1 1) \move(0 1)\lvec(1 1) \move(0 2)\lvec(1 2)
\htext(0.5 2.25){$n$} \htext(0.5 1.5){$n\!\!-\!\!1$} \move(0
-1)\lvec(0 -3.75) \esegment \move(4.4 0) \bsegment \move(0
1)\lvec(0 2.5)\lvec(1.4 2.5)\lvec(1.4 1) \move(0 1)\lvec(1.4 1)
\move(0 2)\lvec(1.4 2) \htext(0.7 2.25){$\cdots$} \htext(0.7
1.5){$\cdots$} \esegment \move(5.8 0) \bsegment \move(0 1)\lvec(0
2.5)\lvec(1 2.5)\lvec(1 -0.5) \move(0 1)\lvec(1 1) \move(0
2)\lvec(1 2) \htext(0.5 2.25){$n$} \htext(0.5 1.5){$n\!\!-\!\!1$}
\esegment%%%%%%%%%%%
\move(6.8 0)
\bsegment
\move(0 -0.5)\lvec(0 1.5)
\move(1 1.5)\lvec(1 0)
\move(0 1)\lvec(1 1)
\move(0 0)\lvec(1 0)
\htext(0.5 0.5){$n\!\!-\!\!2$}
\htext(0.5 1.25){$n\!\!-\!\!1$}
\move(0 -1)\lvec(0 -3.75)
\lpatt(0.05 0.15)
\move(0 1.5)\lvec(1 1.5)
\esegment
\move(7.8 0)
\bsegment
\move(0 0)\lvec(0 1.5)
\move(1.4 1.5)\lvec(1.4 1)
\move(0 1)\lvec(1.4 1)
\move(0 0)\lvec(1.4 0)
\htext(0.7 0.5){$\cdots$}
\htext(0.7 1.25){$\cdots$}
\lpatt(0.05 0.15)
\move(0 1.5)\lvec(1.4 1.5)
\esegment
\move(9.2 0)
\bsegment
\move(0 0)\lvec(0 1.5)
\move(1 1.5)\lvec(1 -0.5)
\move(0 1)\lvec(1 1)
\move(0 0)\lvec(1 0)
\htext(0.5 0.5){$n\!\!-\!\!2$}
\htext(0.5 1.25){$n\!\!-\!\!1$}
\lpatt(0.05 0.15)
\move(0 1.5)\lvec(1 1.5)
\esegment%%%%%%%%%%%%%%%%%
\move(10.2 0)
\bsegment
\move(0 -0.5)\lvec(0 1)\lvec(3.4
1)\lvec(3.4 -0.5) \move(0 -1)\lvec(0 -3.75) \move(3.4 -1)\lvec(3.4
-3.75)
\esegment%%%%%%%%%%%%%%%%%%%%
\move(13.6 0)
\bsegment
\move(0 -3.75)\lvec(0 -1.5)
\move(1 -1.5)\lvec(1 -3)
\move(0 -2)\lvec(1 -2)
\move(0 -3)\lvec(1 -3)
\move(0 -2.5)\lvec(1 -2.5)
\htext(0.5 -2.75){$n$}
\htext(0.5 -2.25){$n$}
\htext(0.5 -1.75){$n\!\!-\!\!1$}
\lpatt(0.05 0.15)
\move(0 -1.5)\lvec(1 -1.5)
\esegment
\move(14.6 0)
\bsegment
\move(0 -3)\lvec(0 -1.5)
\move(1.4 -1.5)\lvec(1.4 -3)
\move(0 -2)\lvec(1.4 -2)
\move(0 -3)\lvec(1.4 -3)
\move(0 -2.5)\lvec(1.4 -2.5)
\htext(0.7 -1.75){$\cdots$}
\htext(0.7 -2.25){$\cdots$}
\htext(0.7 -2.75){$\cdots$}
\lpatt(0.05 0.15)
\move(0 -1.5)\lvec(1.4 -1.5)
\esegment
\move(16 0)
\bsegment
\move(0 -3)\lvec(0 -1.5)
\move(1 -1.5)\lvec(1 -3.75)
\move(0 -2)\lvec(1 -2)
\move(0 -3)\lvec(1 -3)
\move(0 -2.5)\lvec(1 -2.5)
\htext(0.5 -2.75){$n$}
\htext(0.5 -2.25){$n$}
\htext(0.5 -1.75){$n\!\!-\!\!1$}
\lpatt(0.05 0.15)
\move(0 -1.5)\lvec(1 -1.5)
\esegment%%%%%%%%%%%
\htext(1.7 -4.25){$\underbrace{\rule{5.7em}{0em}}$}
\htext(5.1 -4.25){$\underbrace{\rule{5.7em}{0em}}$}
\htext(8.5 -4.25){$\underbrace{\rule{5.7em}{0em}}$}
\htext(11.9 -4.25){$\underbrace{\rule{5.7em}{0em}}$}
\htext(15.3 -4.25){$\underbrace{\rule{5.7em}{0em}}$}
\htext(15.3 -5){$b_3$}
\htext(1.7 -5){$b_1-b_3$}
\htext(5.1 -5){$b_2$}
\htext(8.5 -5){$b_3$}
\htext(11.9 -5){$b_4$}
\move(0 -5.3)
%\drawbb
\end{texdraw}}
\savebox{\tmpfigb}{\begin{texdraw} \fontsize{6}{6}\selectfont
\textref h:C v:C \drawdim em \setunitscale 1.7 \move(0 0)
\bsegment \move(0 -0.5)\lvec(0 2.5)\lvec(1 2.5)\lvec(1 1) \move(0
1)\lvec(1 1) \move(0 2)\lvec(1 2) \htext(0.5 2.25){$n$} \htext(0.5
1.5){$n\!\!-\!\!1$} \move(0 -1)\lvec(0 -3.75) \esegment \move(1 0)
\bsegment \move(0 1)\lvec(0 2.5)\lvec(1.4 2.5)\lvec(1.4 1) \move(0
1)\lvec(1.4 1) \move(0 2)\lvec(1.4 2) \htext(0.7 2.25){$\cdots$}
\htext(0.7 1.5){$\cdots$} \esegment \move(2.4 0) \bsegment \move(0
1)\lvec(0 2.5)\lvec(1 2.5)\lvec(1 -0.5) \move(0 1)\lvec(1 1)
\move(0 2)\lvec(1 2) \htext(0.5 2.25){$n$} \htext(0.5
1.5){$n\!\!-\!\!1$}
\esegment%%%%%%%%%%%%%%%%
\move(3.4 0) \bsegment \move(0 -0.5)\lvec(0 2)\lvec(1 2)\lvec(1 0)
\move(0 1)\lvec(1 1) \move(0 0)\lvec(1 0) \htext(0.5
0.5){$n\!\!-\!\!2$} \htext(0.5 1.5){$n\!\!-\!\!1$} \move(0
-1)\lvec(0 -3.75) \esegment \move(4.4 0) \bsegment \move(0
0)\lvec(0 2)\lvec(1.4 2)\lvec(1.4 0) \move(0 1)\lvec(1.4 1)
\move(0 0)\lvec(1.4 0) \htext(0.7 0.5){$\cdots$} \htext(0.7
1.5){$\cdots$} \esegment \move(5.8 0) \bsegment \move(0 0)\lvec(0
2)\lvec(1 2)\lvec(1 -0.5) \move(0 1)\lvec(1 1) \move(0 0)\lvec(1
0) \htext(0.5 0.5){$n\!\!-\!\!2$} \htext(0.5 1.5){$n\!\!-\!\!1$}
\esegment%%%%%%%%%%%%%
\move(6.8 0)
\bsegment
\move(0 -0.5)\lvec(0 1.5)
\move(1 1.5)\lvec(1 0)
\move(0 1)\lvec(1 1)
\move(0 0)\lvec(1 0)
\htext(0.5 0.5){$n\!\!-\!\!2$}
\htext(0.5 1.25){$n\!\!-\!\!1$}
\move(0 -1)\lvec(0 -3.75)
\lpatt(0.05 0.15)
\move(0 1.5)\lvec(1 1.5)
\esegment
\move(7.8 0)
\bsegment
\move(0 0)\lvec(0 1.5)
\move(1.4 1.5)\lvec(1.4 1)
\move(0 1)\lvec(1.4 1)
\move(0 0)\lvec(1.4 0)
\htext(0.7 0.5){$\cdots$}
\htext(0.7 1.25){$\cdots$}
\lpatt(0.05 0.15)
\move(0 1.5)\lvec(1.4 1.5)
\esegment
\move(9.2 0)
\bsegment
\move(0 0)\lvec(0 1.5)
\move(1 1.5)\lvec(1 -0.5)
\move(0 1)\lvec(1 1)
\move(0 0)\lvec(1 0)
\htext(0.5 0.5){$n\!\!-\!\!2$}
\htext(0.5 1.25){$n\!\!-\!\!1$}
\lpatt(0.05 0.15)
\move(0 1.5)\lvec(1 1.5)
\esegment%%%%%%%%%%%%%%%%%
\move(10.2 0)
\bsegment
\move(0 -0.5)\lvec(0 1)\lvec(3.4 1)\lvec(3.4 -0.5)
\move(0 -1)\lvec(0 -3.75)
\esegment%%%%%%%%%%%%%
\move(13.6 0)
\bsegment
\move(0 -3.75)\lvec(0 -1.5)
\move(1 -1.5)\lvec(1 -3)
\move(0 -2)\lvec(1 -2)
\move(0 -3)\lvec(1 -3)
\move(0 -2.5)\lvec(1 -2.5)
\htext(0.5 -2.75){$n$}
\htext(0.5 -2.25){$n$}
\htext(0.5 -1.75){$n\!\!-\!\!1$}
\move(0 -1)\lvec(0 -3.75)
\lpatt(0.05 0.15)
\move(0 -1.5)\lvec(1 -1.5)
\esegment
\move(14.6 0)
\bsegment
\move(0 -3)\lvec(0 -1.5)
\move(1.4 -1.5)\lvec(1.4 -3)
\move(0 -2)\lvec(1.4 -2)
\move(0 -3)\lvec(1.4 -3)
\move(0 -2.5)\lvec(1.4 -2.5)
\htext(0.7 -1.75){$\cdots$}
\htext(0.7 -2.25){$\cdots$}
\htext(0.7 -2.75){$\cdots$}
\lpatt(0.05 0.15)
\move(0 -1.5)\lvec(1.4 -1.5)
\esegment
\move(16 0)
\bsegment
\move(0 -3)\lvec(0 -1.5)
\move(1 -1.5)\lvec(1 -3.75)
\move(0 -2)\lvec(1 -2)
\move(0 -3)\lvec(1 -3)
\move(0 -2.5)\lvec(1 -2.5)
\htext(0.5 -2.75){$n$}
\htext(0.5 -2.25){$n$}
\htext(0.5 -1.75){$n\!\!-\!\!1$}
\lpatt(0.05 0.15)
\move(0 -1.5)\lvec(1 -1.5)
\esegment%%%%%%%%%%%
\htext(1.7 -4.25){$\underbrace{\rule{5.7em}{0em}}$}
\htext(5.1 -4.25){$\underbrace{\rule{5.7em}{0em}}$}
\htext(8.5 -4.25){$\underbrace{\rule{5.7em}{0em}}$}
\htext(11.9 -4.25){$\underbrace{\rule{5.7em}{0em}}$}
\htext(15.3 -4.25){$\underbrace{\rule{5.7em}{0em}}$}
\htext(15.3 -5){$a_1$}
\htext(1.7 -5){$a_2$}
\htext(5.1 -5){$a_3-a_1$}
\htext(8.3 -5){$a_1$}
\htext(11.9 -5){$a_4$}
\move(0 -5.3)
%\drawbb
\end{texdraw}}
\noindent right :
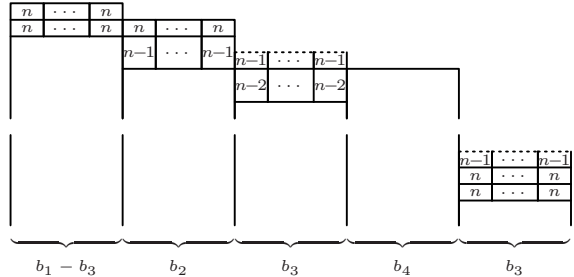
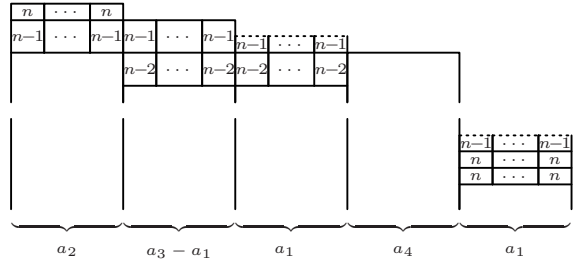
\begin{center}
\usebox{\tmpfiga}
\end{center}
\noindent left :
\begin{center}
\usebox{\tmpfigb}
\end{center}
Finally, join the two slices so that the highest layer of the
result forms a part of a level-$l$ reduced proper Young wall that has had
all its blocks that may be split, split.
This is taken to be the starting point.

Now, to bring this into a reduced proper form, we need to
\emph{remove} $\delta$'s from the \emph{left} slice. We denote
the number of $\delta$ removals needed by $k$.

Below, we list left-$\vphi$ and right-$\veps$ values for each
possible Young wall column pair. The line containing the bullet
lists the two column types in left-right order. Each case is again
separated into two cases according to the range $k$ falls into.
These numbers have been directly observed from (mental) drawings
of the two columns in reduced proper form.
\begin{itemize}
\item \Na\ \Na\\[2mm]
\begin{tabular}{rcl}
rotation & : & $0\le k\le a_2$\\
left-$\vphi$ & : & $k$\\
right-$\veps$ & : &$b_2+2(b_1-b_3)-(a_2-k)$\\[2mm]
rotation & : & $a_2 \leq k$\\
left-$\vphi$ & : & $a_2$\\
right-$\veps$ & : &$b_2+2(b_1-b_3)$\\
\end{tabular}\\[2mm]
\item \Na\ \Nb\\[2mm]
\begin{tabular}{rcl}
rotation & : & $0\le k\le a_2$\\
left-$\vphi$ & : & $k$\\
right-$\veps$ & : &$b_2-(a_2-k)$\\[2mm]
rotation & : & $a_2 \leq k$\\
left-$\vphi$ & : & $a_2$\\
right-$\veps$ & : &$b_2$\\
\end{tabular}\\[2mm]
\item \Nb\ \Na\\[2mm]
\begin{tabular}{rcl}
rotation & : & $0\le k\le a_2+(a_3-a_1)$\\
left-$\vphi$ & : & $k+(a_3-a_1)$\\
right-$\veps$ & : &$b_2+2(b_1-b_3)-(a_2+(a_3-a_1)-k)$\\[2mm]
rotation & : & $a_2+(a_3-a_1)\leq k$\\
left-$\vphi$ & : & $a_2+2(a_3-a_1)$\\
right-$\veps$ & : &$b_2+2(b_1-b_3)$\\
\end{tabular}\\[2mm]
\item \Nb\ \Nb\\[2mm]
\begin{tabular}{rcl}
rotation & : & $0\le k\le a_2+(a_3-a_1)$\\
left-$\vphi$ & : & $k+(a_3-a_1)$\\
right-$\veps$ & : &$b_2-(a_2+(a_3-a_1)-k)$\\[2mm]
rotation & : & $a_2+(a_3-a_1) \leq k$\\
left-$\vphi$ & : & $a_2+2(a_3-a_1)$\\
right-$\veps$ & : &$b_2$\\
\end{tabular}\\[2mm]
\end{itemize}
Similarly, the following gives the signatures of the path
description. The number in the list are the left-$\vphi$ and
right-$\veps$ values for the two corresponding crystal elements.
\begin{itemize}
\item \Na\ \Na\\[2mm]
  \begin{tabular}{rcl}
  left-$\vphi$ & : & $a_2$\\
  right-$\veps$ & : &$b_2+2(b_1-b_3)$\\
  \end{tabular}\\[2mm]
\item \Na\ \Nb\\[2mm]
  \begin{tabular}{rcl}
  left-$\vphi$ & : & $a_2$\\
  right-$\veps$ & : &$b_2$\\
  \end{tabular}\\[2mm]
\item \Nb\ \Na\\[2mm]
  \begin{tabular}{rcl}
  left-$\vphi$ & : & $a_2+2(a_3-a_1)$\\
  right-$\veps$ & : &$b_2+2(b_1-b_3)$\\
  \end{tabular}\\[2mm]
\item \Nb\ \Nb\\[2mm]
  \begin{tabular}{rcl}
  left-$\vphi$ & : & $a_2+2(a_3-a_1)$\\
  right-$\veps$ & : &$b_2$\\
  \end{tabular}\\[2mm]
\end{itemize}
We can easily see that the signatures agree with those of the
corresponding path description in all of the cases after
$(0,1)$-pair cancellations.

We have covered the remaining part of this proof for the $B_n^{(1)}$ case
in the Appendices.
\end{proof}

We close this paper with an example of Young wall realization of
irreducible highest weight crystals.

\begin{example} \hfill
\begin{enumerate}
\item The top part of the crystal graph
$\rpwspace(3\La_0)$ for $\uq(B_3^{(1)})$ is given below.
\vspace{2mm}
\newcommand{\filledbox}{%
\bsegment
\setsegscale 0.8
\move(0 0)\lvec(-1 0)\lvec(-1 0.133)\lvec(0 0.133)\lvec(0 0)
\ifill f:0.5
\esegment
}
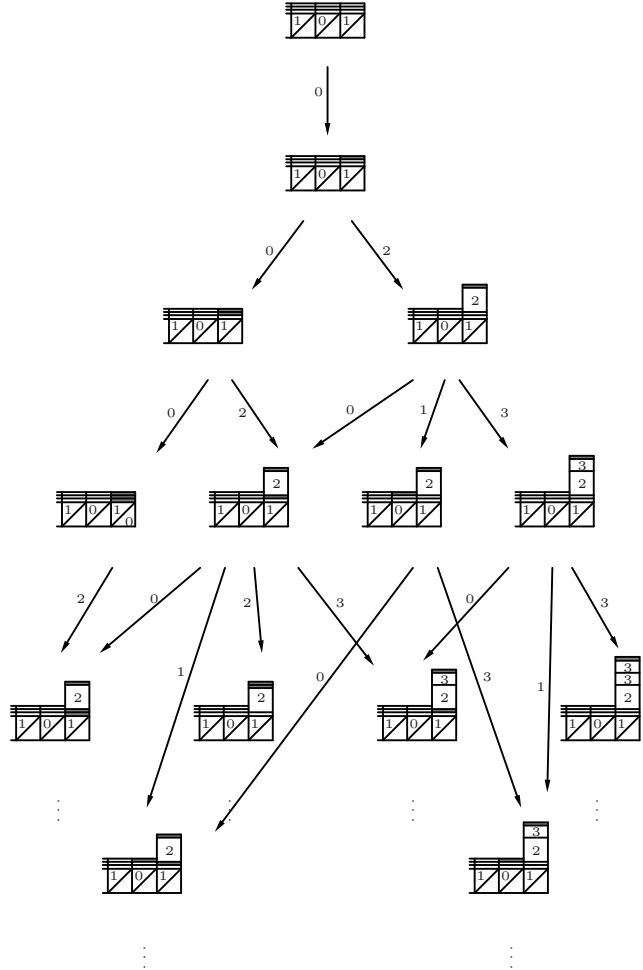
\begin{center}
\begin{texdraw}
\drawdim in
\arrowheadsize l:0.065 w:0.03
\arrowheadtype t:F
\fontsize{5}{5}\selectfont
\textref h:C v:C
\drawdim em
\setunitscale 1.7
\move(0 10)
\bsegment
\setsegscale 0.8
\move(0 0)\rlvec(-3.2 0)
\move(0 1)\rlvec(-3.2 0)
\move(0 1.133)\rlvec(-3.2 0)
\move(0 1.266)\rlvec(-3.2 0)
\move(0 1.399)\rlvec(-3.2 0)%
\move(0 0)\rlvec(0 1.399)
\move(-1 0)\rlvec(0 1.399)
\move(-2 0)\rlvec(0 1.399)
\move(-3 0)\rlvec(0 1.399)%
\move(-1 0)\rlvec(1 1)
\move(-2 0)\rlvec(1 1)
\move(-3 0)\rlvec(1 1)%
\htext(-0.75 0.72){$1$}\htext(-1.75 0.72){$0$}\htext(-2.75 0.72){$1$}
\esegment
\move(0 5)
\bsegment
\setsegscale 0.8
\move(0 1.266)\filledbox%
\move(0 0)\rlvec(-3.2 0)
\move(0 1)\rlvec(-3.2 0)
\move(0 1.133)\rlvec(-3.2 0)
\move(0 1.266)\rlvec(-3.2 0)
\move(0 1.399)\rlvec(-3.2 0)%
\move(0 0)\rlvec(0 1.399)
\move(-1 0)\rlvec(0 1.399)
\move(-2 0)\rlvec(0 1.399)
\move(-3 0)\rlvec(0 1.399)%
\move(-1 0)\rlvec(1 1)
\move(-2 0)\rlvec(1 1)
\move(-3 0)\rlvec(1 1)%
\htext(-0.75 0.72){$1$}\htext(-1.75 0.72){$0$}\htext(-2.75 0.72){$1$}
\esegment
\move(0 0)
\bsegment
\move(-4 0)
\bsegment
\setsegscale 0.8
\move(0 1.266)\filledbox
\move(0 1.133)\filledbox%
\move(0 0)\rlvec(-3.2 0)
\move(0 1)\rlvec(-3.2 0)
\move(0 1.133)\rlvec(-3.2 0)
\move(0 1.266)\rlvec(-3.2 0)
\move(0 1.399)\rlvec(-3.2 0)%
\move(0 0)\rlvec(0 1.399)
\move(-1 0)\rlvec(0 1.399)
\move(-2 0)\rlvec(0 1.399)
\move(-3 0)\rlvec(0 1.399)%
\move(-1 0)\rlvec(1 1)
\move(-2 0)\rlvec(1 1)
\move(-3 0)\rlvec(1 1)%
\htext(-0.75 0.72){$1$}\htext(-1.75 0.72){$0$}\htext(-2.75 0.72){$1$}
\esegment
\move(4 0)
\bsegment
\setsegscale 0.8
\move(0 2.266)\filledbox%
\move(0 0)\rlvec(-3.2 0)
\move(0 1)\rlvec(-3.2 0)
\move(0 1.133)\rlvec(-3.2 0)
\move(0 1.266)\rlvec(-3.2 0)
\move(-1 1.4)\rlvec(-2.2 0)
\move(0 2.266)\rlvec(-1 0)
\move(0 2.399)\rlvec(-1 0)
\move(0 0)\rlvec(0 2.399)
\move(-1 0)\rlvec(0 2.399)
\move(-2 0)\rlvec(0 1.399)
\move(-3 0)\rlvec(0 1.399)%
\move(-1 0)\rlvec(1 1)
\move(-2 0)\rlvec(1 1)
\move(-3 0)\rlvec(1 1)%
\htext(-0.75 0.72){$1$}\htext(-1.75 0.72){$0$}\htext(-2.75 0.72){$1$}
\htext(-0.5 1.766){$2$}
\esegment
\esegment
\move(0 -6)
\bsegment
\move(-7.5 0)
\bsegment
\setsegscale 0.8
\move(0 1.266)\filledbox
\move(0 1.133)\filledbox
\move(0 1)\filledbox%
\move(0 0)\rlvec(-3.2 0)
\move(0 1)\rlvec(-3.2 0)
\move(0 1.133)\rlvec(-3.2 0)
\move(0 1.266)\rlvec(-3.2 0)
\move(0 1.399)\rlvec(-3.2 0)%
\move(0 0)\rlvec(0 1.399)
\move(-1 0)\rlvec(0 1.399)
\move(-2 0)\rlvec(0 1.399)
\move(-3 0)\rlvec(0 1.399)%
\move(-1 0)\rlvec(1 1)
\move(-2 0)\rlvec(1 1)
\move(-3 0)\rlvec(1 1)%
\htext(-0.75 0.72){$1$}\htext(-1.75 0.72){$0$}\htext(-2.75 0.72){$1$}
\htext(-0.25 0.23){$0$}
\esegment
\move(-2.5 0)
\bsegment
\setsegscale 0.8
\move(0 2.266)\filledbox
\move(0 1.133)\filledbox%
\move(0 0)\rlvec(-3.2 0)
\move(0 1)\rlvec(-3.2 0)
\move(0 1.133)\rlvec(-3.2 0)
\move(0 1.266)\rlvec(-3.2 0)
\move(-1 1.399)\rlvec(-2.2 0)
\move(0 2.266)\rlvec(-1 0)
\move(0 2.399)\rlvec(-1 0)%
\move(0 0)\rlvec(0 2.399)
\move(-1 0)\rlvec(0 2.399)
\move(-2 0)\rlvec(0 1.399)
\move(-3 0)\rlvec(0 1.399)%
\move(-1 0)\rlvec(1 1)
\move(-2 0)\rlvec(1 1)
\move(-3 0)\rlvec(1 1)%
\htext(-0.75 0.72){$1$}\htext(-1.75 0.72){$0$}\htext(-2.75 0.72){$1$}
\htext(-0.5 1.766){$2$}
\esegment
\move(2.5 0)
\bsegment
\setsegscale 0.8
\move(0 2.266)\filledbox
\move(-1 1.266)\filledbox%
\move(0 0)\rlvec(-3.2 0)
\move(0 1)\rlvec(-3.2 0)
\move(0 1.133)\rlvec(-3.2 0)
\move(0 1.266)\rlvec(-3.2 0)
\move(-1 1.399)\rlvec(-2.2 0)
\move(0 2.266)\rlvec(-1 0)
\move(0 2.399)\rlvec(-1 0)%
\move(0 0)\rlvec(0 2.399)
\move(-1 0)\rlvec(0 2.399)
\move(-2 0)\rlvec(0 1.399)
\move(-3 0)\rlvec(0 1.399)%
\move(-1 0)\rlvec(1 1)
\move(-2 0)\rlvec(1 1)
\move(-3 0)\rlvec(1 1)%
\htext(-0.75 0.72){$1$}\htext(-1.75 0.72){$0$}\htext(-2.75 0.72){$1$}
\htext(-0.5 1.766){$2$}
\esegment
\move(7.5 0)
\bsegment
\setsegscale 0.8
\move(0 2.766)\filledbox%
\move(0 0)\rlvec(-3.2 0)
\move(0 1)\rlvec(-3.2 0)
\move(0 1.133)\rlvec(-3.2 0)
\move(0 1.266)\rlvec(-3.2 0)
\move(-1 1.399)\rlvec(-2.2 0)
\move(0 2.266)\rlvec(-1 0)
\move(0 2.766)\rlvec(-1 0)
\move(0 2.899)\rlvec(-1 0)%
\move(0 0)\rlvec(0 2.899)
\move(-1 0)\rlvec(0 2.899)
\move(-2 0)\rlvec(0 1.399)
\move(-3 0)\rlvec(0 1.399)%
\move(-1 0)\rlvec(1 1)
\move(-2 0)\rlvec(1 1)
\move(-3 0)\rlvec(1 1)%
\htext(-0.75 0.72){$1$}\htext(-1.75 0.72){$0$}\htext(-2.75 0.72){$1$}
\htext(-0.5 1.766){$2$}
\htext(-0.5 2.516){$3$}
\esegment
\esegment
\move(0 -13)
\bsegment
\move(-9 0)
\bsegment
\setsegscale 0.8
\move(0 2.266)\filledbox
\move(0 1.133)\filledbox
\move(0 1)\filledbox
\move(0 0)\rlvec(-3.2 0)
\move(0 1)\rlvec(-3.2 0)
\move(0 1.133)\rlvec(-3.2 0)
\move(0 1.266)\rlvec(-3.2 0)
\move(-1 1.4)\rlvec(-2.2 0)
\move(0 2.266)\rlvec(-1 0)
\move(0 2.4)\rlvec(-1 0)
\move(0 0)\rlvec(0 2.4)
\move(-1 0)\rlvec(0 2.4)
\move(-2 0)\rlvec(0 1.4)
\move(-3 0)\rlvec(0 1.4)
\move(-1 0)\rlvec(1 1)
\move(-2 0)\rlvec(1 1)
\move(-3 0)\rlvec(1 1)
\htext(-0.75 0.72){$1$}\htext(-1.75 0.72){$0$}\htext(-2.75 0.72){$1$}
\htext(-0.5 1.766){$2$}
\esegment
\move(-3 0)
\bsegment
\setsegscale 0.8
\move(0 2.133)\filledbox
\move(0 2.266)\filledbox
\move(0 0)\rlvec(-3.2 0)
\move(0 1)\rlvec(-3.2 0)
\move(0 1.133)\rlvec(-3.2 0)
\move(-1 1.266)\rlvec(-2.2 0)
\move(-1 1.399)\rlvec(-2.2 0)
\move(0 2.133)\rlvec(-1 0)
\move(0 2.266)\rlvec(-1 0)
\move(0 2.399)\rlvec(-1 0)
\move(0 0)\rlvec(0 2.399)
\move(-1 0)\rlvec(0 2.399)
\move(-2 0)\rlvec(0 1.399)
\move(-3 0)\rlvec(0 1.399)
\move(-1 0)\rlvec(1 1)
\move(-2 0)\rlvec(1 1)
\move(-3 0)\rlvec(1 1)
\htext(-0.75 0.72){$1$}\htext(-1.75 0.72){$0$}\htext(-2.75 0.72){$1$}
\htext(-0.5 1.766){$2$}
\esegment
\move(3 0)
\bsegment
\setsegscale 0.8
\move(0 2.766)\filledbox
\move(0 1.133)\filledbox
\move(0 0)\rlvec(-3.2 0)
\move(0 1)\rlvec(-3.2 0)
\move(0 1.133)\rlvec(-3.2 0)
\move(0 1.266)\rlvec(-3.2 0)
\move(-1 1.4)\rlvec(-2.2 0)
\move(0 2.266)\rlvec(-1 0)
\move(0 2.766)\rlvec(-1 0)
\move(0 2.9)\rlvec(-1 0)
\move(0 0)\rlvec(0 2.9)
\move(-1 0)\rlvec(0 2.9)
\move(-2 0)\rlvec(0 1.4)
\move(-3 0)\rlvec(0 1.4)
\move(-1 0)\rlvec(1 1)
\move(-2 0)\rlvec(1 1)
\move(-3 0)\rlvec(1 1)
\htext(-0.75 0.72){$1$}\htext(-1.75 0.72){$0$}\htext(-2.75 0.72){$1$}
\htext(-0.5 1.766){$2$}
\htext(-0.5 2.516){$3$}
\esegment
\move(9 0)
\bsegment
\setsegscale 0.8
\move(0 3.266)\filledbox
\move(0 0)\rlvec(-3.2 0)
\move(0 1)\rlvec(-3.2 0)
\move(0 1.133)\rlvec(-3.2 0)
\move(0 1.266)\rlvec(-3.2 0)
\move(-1 1.4)\rlvec(-2.2 0)
\move(0 2.266)\rlvec(-1 0)
\move(0 2.766)\rlvec(-1 0)
\move(0 3.266)\rlvec(-1 0)
\move(0 3.4)\rlvec(-1 0)
\move(0 0)\rlvec(0 3.4)
\move(-1 0)\rlvec(0 3.4)
\move(-2 0)\rlvec(0 1.4)
\move(-3 0)\rlvec(0 1.4)
\move(-1 0)\rlvec(1 1)
\move(-2 0)\rlvec(1 1)
\move(-3 0)\rlvec(1 1)
\htext(-0.75 0.72){$1$}\htext(-1.75 0.72){$0$}\htext(-2.75 0.72){$1$}
\htext(-0.5 1.766){$2$}
\htext(-0.5 2.516){$3$}
\htext(-0.5 3.016){$3$}
\esegment
\esegment
\move(0 -18)
\bsegment
\move(-6 0)
\bsegment
\setsegscale 0.8
\move(-1 1.266)\filledbox
\move(0 2.266)\filledbox
\move(0 1.133)\filledbox
\move(0 0)\rlvec(-3.2 0)
\move(0 1)\rlvec(-3.2 0)
\move(0 1.133)\rlvec(-3.2 0)
\move(0 1.266)\rlvec(-3.2 0)
\move(-1 1.4)\rlvec(-2.2 0)
\move(0 2.266)\rlvec(-1 0)
\move(0 2.4)\rlvec(-1 0)
\move(0 0)\rlvec(0 2.4)
\move(-1 0)\rlvec(0 2.4)
\move(-2 0)\rlvec(0 1.4)
\move(-3 0)\rlvec(0 1.4)
\move(-1 0)\rlvec(1 1)
\move(-2 0)\rlvec(1 1)
\move(-3 0)\rlvec(1 1)
\htext(-0.75 0.72){$1$}\htext(-1.75 0.72){$0$}\htext(-2.75 0.72){$1$}
\htext(-0.5 1.766){$2$}
\esegment
\move(6 0)
\bsegment
\setsegscale 0.8
\move(0 2.766)\filledbox
\move(-1 1.266)\filledbox
\move(0 0)\rlvec(-3.2 0)
\move(0 1)\rlvec(-3.2 0)
\move(0 1.133)\rlvec(-3.2 0)
\move(0 1.266)\rlvec(-3.2 0)
\move(-1 1.4)\rlvec(-2.2 0)
\move(0 2.266)\rlvec(-1 0)
\move(0 2.766)\rlvec(-1 0)
\move(0 2.9)\rlvec(-1 0)
\move(0 0)\rlvec(0 2.9)
\move(-1 0)\rlvec(0 2.9)
\move(-2 0)\rlvec(0 1.4)
\move(-3 0)\rlvec(0 1.4)
\move(-1 0)\rlvec(1 1)
\move(-2 0)\rlvec(1 1)
\move(-3 0)\rlvec(1 1)
\htext(-0.75 0.72){$1$}\htext(-1.75 0.72){$0$}\htext(-2.75 0.72){$1$}
\htext(-0.5 1.766){$2$}
\htext(-0.5 2.516){$3$}
\esegment
\esegment
%\drawbb
%
\move(-1.2 0)
\bsegment
\setsegscale 0.8
\move(0 11.3)\ravec(0 -2.8)
\htext(-0.4 10.3){$0$}
\move(-1 5)\ravec(-2.1 -2.8)
\move(1 5)\ravec(2.1 -2.8)
\htext(-2.4 3.8){$0$}
\htext(2.4 3.8){$2$}
\move(-4.9 -1.5)\ravec(-2.1 -3)
\move(-3.9 -1.5)\ravec(1.9 -2.8)
\move(3.5 -1.5)\ravec(-4.1 -2.8)
\move(4.8 -1.5)\ravec(-1 -2.8)
\move(5.4 -1.5)\ravec(2 -2.8)
\htext(-6.4 -2.8){$0$}
\htext(-3.5 -2.8){$2$}
\htext(0.9 -2.7){$0$}
\htext(3.9 -2.7){$1$}
\htext(7.2 -2.8){$3$}
\move(-8.8 -9.2)\ravec(-2.1 -3.5)
\move(-5.2 -9.2)\ravec(-4.1 -3.5)
\move(-4.2 -9.2)\ravec(-3.2 -9.8)
\move(-3 -9.2)\ravec(0.3 -3.5)
\move(-1.2 -9.2)\ravec(3.1 -4)
\move(3.5 -9.2)\ravec(-8.1 -10.8)
\move(4.5 -9.2)\ravec(3.4 -9.8)
\move(7.4 -9.2)\ravec(-3.5 -3.8)
\move(9.2 -9.2)\ravec(-0.2 -9.2)
\move(10 -9.2)\ravec(1.9 -3.3)
\htext(-10.1 -10.4){$2$}
\htext(-7.1 -10.4){$0$}
\htext(-6 -13.4){$1$}
\htext(-3.3 -10.6){$2$}
\htext(0.5 -10.6){$3$}
\htext(-0.3 -13.6){$0$}
\htext(6.5 -13.6){$3$}
\htext(5.8 -10.4){$0$}
\htext(8.7 -14){$1$}
\htext(11.3 -10.6){$3$}
\vtext(-11 -19){$\cdots$}
\vtext(-7.5 -25){$\cdots$}
\vtext(-4 -19){$\cdots$}
\vtext(3.5 -19){$\cdots$}
\vtext(7.5 -25){$\cdots$}
\vtext(11 -19){$\cdots$}
\esegment
\end{texdraw}
\end{center}

\item The next figure shows how the Kashiwara operators act on a
reduced proper Young wall in $\rpwspace(3\La_0)$.

\vskip 5mm
\begin{center}
\begin{texdraw}
\fontsize{5}{5}\selectfont
\drawdim in
\arrowheadsize l:0.065 w:0.03
\arrowheadtype t:F
\textref h:C v:C
\drawdim em \setunitscale 1.7
\nc{\dtri}{ \bsegment
%\move(-1 0)\lvec(0 1)\lvec(0 0)\lvec(-1 0)\ifill f:0.7
\esegment }
\move(-8 12)
\bsegment
%\move(-3 1.5) \rlvec(-1 0) \rlvec(0 0.25) \rlvec(1 0) \rlvec(0 -0.25) \ifill f:0.6
\move(-2 2.5) \rlvec(-1 0) \rlvec(0 0.25) \rlvec(1 0) \rlvec(0 -0.25) \ifill f:0.6
\move(-2 3.25) \rlvec(-1 0) \rlvec(0 0.25) \rlvec(1 0) \rlvec(0 -0.25) \ifill f:0.6
\move(-2 3.5) \rlvec(-1 0) \rlvec(0 0.25) \rlvec(1 0) \rlvec(0 -0.25) \ifill f:0.6
\move(-1 3) \rlvec(-1 0) \rlvec(0 0.25) \rlvec(1 0) \rlvec(0 -0.25) \ifill f:0.6
\move(-1 4.25) \rlvec(-1 0) \rlvec(0 0.25) \rlvec(1 0) \rlvec(0 -0.25) \ifill f:0.6
\move(-1 6.5) \rlvec(-1 0) \rlvec(0 0.25) \rlvec(1 0) \rlvec(0 -0.25) \ifill f:0.6
\move(0 4) \rlvec(-1 0) \rlvec(0 0.25) \rlvec(1 0) \rlvec(0 -0.25) \ifill f:0.6
\move(0 5.25) \rlvec(-1 0) \rlvec(0 0.25) \rlvec(1 0) \rlvec(0 -0.25) \ifill f:0.6
\move(0 7.5) \rlvec(-1 0) \rlvec(0 0.25) \rlvec(1 0) \rlvec(0 -0.25) \ifill f:0.6

\move(0 0)\dtri \move(-1 0)\dtri \move(-2 0)\dtri \move(-3 0)\dtri
\move(0 0)\rlvec(-4 0) \move(0 1)\rlvec(-4 0) \move(0 2)\rlvec(-3
0)

\move(0 1) \rlvec(-1 -1) \move (-1 1) \rlvec(-1 -1)
\move (-2 1) \rlvec(-1 -1) \move (-3 1) \rlvec(-1 -1)

\move(0 2.5)\rlvec(-3 0)  \move(0 3) \rlvec(-2 0)
\move(0 4)\rlvec(-1 0)

\move(-3 1.25) \rlvec(-1 0) \move(-3 1.5) \rlvec(-1 0)
\move (-3 1.75) \rlvec(-1 0) \move(-2 2.75) \rlvec(-1 0)
\move(-2 3.25) \rlvec(-1 0) \move(-2 3.5) \rlvec(-1 0)
\move(-2 3.75) \rlvec(-1 0) \move(-1 3.25) \rlvec(-1 0)
\move(-1 4.25) \rlvec(-1 0) \move(-1 4.5) \rlvec(-1 0)
\move(-1 5.5) \rlvec(-1 0) \move(-1 6.5) \rlvec(-1 0)
\move(-1 6.75) \rlvec(-1 0) \move(-1 5.5) \rlvec(-1 -1)

\move(0 0)\rlvec(0 7.75) \move(-1 0)\rlvec(0 7.75)
\move(-2 0)\rlvec(0 6.75) \move(-3 0)\rlvec(0 3.75)
\move(-4 0)\rlvec(0 1.75)

\move(0 4.25) \rlvec(-1 0) \move(0 5.25) \rlvec(-1 -1)
\move(0 5.25) \rlvec(-1 0) \move(0 5.5) \rlvec(-1 0)
\move(0 6.5) \rlvec(-1 0) \move(0 7) \rlvec(-1 0)
\move(0 7.5) \rlvec(-1 0) \move(0 7.75) \rlvec(-1 0)

\htext(-0.3 0.25){$0$} \htext(-0.75 0.75){$1$} \htext(-1.3 0.25){$1$} \htext(-1.75 0.75){$0$}
\htext(-2.3 0.25){$0$} \htext(-2.75 0.75){$1$}\htext(-3.75 0.75){$0$}
\htext(-0.5 1.5){$2$}\htext(-1.5 1.5){$2$} \htext(-2.5 1.5){$2$}
\htext(-0.5 2.25){$3$}\htext(-1.5 2.25){$3$} \htext(-2.5 2.25){$3$}
\htext(-0.5 2.75){$3$}\htext(-1.5 2.75){$3$} \htext(-2.5 3){$3$}
\htext(-0.5 3.5){$2$} \htext(-1.5 3.75){$2$}
\htext(-0.3 4.5){$0$} \htext(-0.75 5){$1$}
\htext(-1.3 4.75){$1$} \htext(-1.75 5.25){$0$}
\htext(-0.5 6){$2$} \htext(-1.5 6){$2$}
\htext(-0.5 6.75){$3$}
\htext(-0.5 7.25){$3$}

\move(-4 0)\rlvec(-0.2 0)
\move(-4 1)\rlvec(-0.2 0)
\move(-4 1.25)\rlvec(-0.2 0)
\move(-4 1.5)\rlvec(-0.2 0)
\move(-4 1.75)\rlvec(-0.2 0)
\esegment

\move(8 12)
\bsegment

\move(-3 1.5) \rlvec(-1 0) \rlvec(0 0.25) \rlvec(1 0) \rlvec(0 -0.25) \ifill f:0.6
\move(-2 2.5) \rlvec(-1 0) \rlvec(0 0.25) \rlvec(1 0) \rlvec(0 -0.25) \ifill f:0.6
\move(-2 3.25) \rlvec(-1 0) \rlvec(0 0.25) \rlvec(1 0) \rlvec(0 -0.25) \ifill f:0.6
\move(-2 3.5) \rlvec(-1 0) \rlvec(0 0.25) \rlvec(1 0) \rlvec(0 -0.25) \ifill f:0.6
\move(-1 3) \rlvec(-1 0) \rlvec(0 0.25) \rlvec(1 0) \rlvec(0 -0.25) \ifill f:0.6
\move(-1 4.25) \rlvec(-1 0) \rlvec(0 0.25) \rlvec(1 0) \rlvec(0 -0.25) \ifill f:0.6
\move(-1 6.5) \rlvec(-1 0) \rlvec(0 0.25) \rlvec(1 0) \rlvec(0 -0.25) \ifill f:0.6
\move(0 4) \rlvec(-1 0) \rlvec(0 0.25) \rlvec(1 0) \rlvec(0 -0.25) \ifill f:0.6
\move(0 5.25) \rlvec(-1 0) \rlvec(0 0.25) \rlvec(1 0) \rlvec(0 -0.25) \ifill f:0.6
\move(0 7) \rlvec(-1 0) \rlvec(0 0.25) \rlvec(1 0) \rlvec(0 -0.25) \ifill f:0.6

\move(0 0)\dtri \move(-1 0)\dtri \move(-2 0)\dtri \move(-3 0)\dtri
\move(0 0)\rlvec(-4 0) \move(0 1)\rlvec(-4 0) \move(0 2)\rlvec(-3
0)

\move(0 1) \rlvec(-1 -1) \move (-1 1) \rlvec(-1 -1)
\move (-2 1) \rlvec(-1 -1) \move (-3 1) \rlvec(-1 -1)

\move(0 2.5)\rlvec(-3 0)  \move(0 3) \rlvec(-2 0)
\move(0 4)\rlvec(-1 0)

\move(-3 1.25) \rlvec(-1 0) \move(-3 1.5) \rlvec(-1 0)
\move (-3 1.75) \rlvec(-1 0) \move(-2 2.75) \rlvec(-1 0)
\move(-2 3.25) \rlvec(-1 0) \move(-2 3.5) \rlvec(-1 0)
\move(-2 3.75) \rlvec(-1 0) \move(-1 3.25) \rlvec(-1 0)
\move(-1 4.25) \rlvec(-1 0) \move(-1 4.5) \rlvec(-1 0)
\move(-1 5.5) \rlvec(-1 0) \move(-1 6.5) \rlvec(-1 0)
\move(-1 6.75) \rlvec(-1 0) \move(-1 5.5) \rlvec(-1 -1)

\move(0 0)\rlvec(0 7.25) \move(-1 0)\rlvec(0 7.25)
\move(-2 0)\rlvec(0 6.75) \move(-3 0)\rlvec(0 3.75)
\move(-4 0)\rlvec(0 1.75)

\move(0 4.25) \rlvec(-1 0) \move(0 5.25) \rlvec(-1 -1)
\move(0 5.25) \rlvec(-1 0) \move(0 5.5) \rlvec(-1 0)
\move(0 6.5) \rlvec(-1 0) \move(0 7) \rlvec(-1 0)
\move(0 7.25) \rlvec(-1 0) %\move(0 7.75) \rlvec(-1 0)

\htext(-0.3 0.25){$0$} \htext(-0.75 0.75){$1$} \htext(-1.3 0.25){$1$} \htext(-1.75 0.75){$0$}
\htext(-2.3 0.25){$0$} \htext(-2.75 0.75){$1$}\htext(-3.75 0.75){$0$}
\htext(-0.5 1.5){$2$}\htext(-1.5 1.5){$2$} \htext(-2.5 1.5){$2$}
\htext(-0.5 2.25){$3$}\htext(-1.5 2.25){$3$} \htext(-2.5 2.25){$3$}
\htext(-0.5 2.75){$3$}\htext(-1.5 2.75){$3$} \htext(-2.5 3){$3$}
\htext(-0.5 3.5){$2$} \htext(-1.5 3.75){$2$}
\htext(-0.3 4.5){$0$} \htext(-0.75 5){$1$}
\htext(-1.3 4.75){$1$} \htext(-1.75 5.25){$0$}
\htext(-0.5 6){$2$} \htext(-1.5 6){$2$}
\htext(-0.5 6.75){$3$}
%\htext(-0.5 7.25){$3$}

\move(-4 0)\rlvec(-0.2 0)
\move(-4 1)\rlvec(-0.2 0)
\move(-4 1.25)\rlvec(-0.2 0)
\move(-4 1.5)\rlvec(-0.2 0)
\move(-4 1.75)\rlvec(-0.2 0)
\esegment

\move(0 0)

\bsegment

\move(-3 1.5) \rlvec(-1 0) \rlvec(0 0.25) \rlvec(1 0) \rlvec(0 -0.25) \ifill f:0.6
\move(-2 2.5) \rlvec(-1 0) \rlvec(0 0.25) \rlvec(1 0) \rlvec(0 -0.25) \ifill f:0.6
\move(-2 3.25) \rlvec(-1 0) \rlvec(0 0.25) \rlvec(1 0) \rlvec(0 -0.25) \ifill f:0.6
\move(-2 3.5) \rlvec(-1 0) \rlvec(0 0.25) \rlvec(1 0) \rlvec(0 -0.25) \ifill f:0.6
\move(-1 3) \rlvec(-1 0) \rlvec(0 0.25) \rlvec(1 0) \rlvec(0 -0.25) \ifill f:0.6
\move(-1 4.25) \rlvec(-1 0) \rlvec(0 0.25) \rlvec(1 0) \rlvec(0 -0.25) \ifill f:0.6
\move(-1 6.5) \rlvec(-1 0) \rlvec(0 0.25) \rlvec(1 0) \rlvec(0 -0.25) \ifill f:0.6
\move(0 4) \rlvec(-1 0) \rlvec(0 0.25) \rlvec(1 0) \rlvec(0 -0.25) \ifill f:0.6
\move(0 5.25) \rlvec(-1 0) \rlvec(0 0.25) \rlvec(1 0) \rlvec(0 -0.25) \ifill f:0.6
\move(0 7.5) \rlvec(-1 0) \rlvec(0 0.25) \rlvec(1 0) \rlvec(0 -0.25) \ifill f:0.6

\move(0 0)\dtri \move(-1 0)\dtri \move(-2 0)\dtri \move(-3 0)\dtri
\move(0 0)\rlvec(-4 0) \move(0 1)\rlvec(-4 0) \move(0 2)\rlvec(-3 0)

\move(0 1) \rlvec(-1 -1) \move (-1 1) \rlvec(-1 -1)
\move (-2 1) \rlvec(-1 -1) \move (-3 1) \rlvec(-1 -1)

\move(0 2.5)\rlvec(-3 0)  \move(0 3) \rlvec(-2 0)
\move(0 4)\rlvec(-1 0)

\move(-3 1.25) \rlvec(-1 0) \move(-3 1.5) \rlvec(-1 0)
\move (-3 1.75) \rlvec(-1 0) \move(-2 2.75) \rlvec(-1 0)
\move(-2 3.25) \rlvec(-1 0) \move(-2 3.5) \rlvec(-1 0)
\move(-2 3.75) \rlvec(-1 0) \move(-1 3.25) \rlvec(-1 0)
\move(-1 4.25) \rlvec(-1 0) \move(-1 4.5) \rlvec(-1 0)
\move(-1 5.5) \rlvec(-1 0) \move(-1 6.5) \rlvec(-1 0)
\move(-1 6.75) \rlvec(-1 0) \move(-1 5.5) \rlvec(-1 -1)

\move(0 0)\rlvec(0 7.75) \move(-1 0)\rlvec(0 7.75)
\move(-2 0)\rlvec(0 6.75) \move(-3 0)\rlvec(0 3.75)
\move(-4 0)\rlvec(0 1.75)

\move(0 4.25) \rlvec(-1 0) \move(0 5.25) \rlvec(-1 -1)
\move(0 5.25) \rlvec(-1 0) \move(0 5.5) \rlvec(-1 0)
\move(0 6.5) \rlvec(-1 0) \move(0 7) \rlvec(-1 0)
\move(0 7.5) \rlvec(-1 0) \move(0 7.75) \rlvec(-1 0)

\htext(-0.3 0.25){$0$} \htext(-0.75 0.75){$1$} \htext(-1.3 0.25){$1$} \htext(-1.75 0.75){$0$}
\htext(-2.3 0.25){$0$} \htext(-2.75 0.75){$1$}\htext(-3.75 0.75){$0$}
\htext(-0.5 1.5){$2$}\htext(-1.5 1.5){$2$} \htext(-2.5 1.5){$2$}
\htext(-0.5 2.25){$3$}\htext(-1.5 2.25){$3$} \htext(-2.5 2.25){$3$}
\htext(-0.5 2.75){$3$}\htext(-1.5 2.75){$3$} \htext(-2.5 3){$3$}
\htext(-0.5 3.5){$2$} \htext(-1.5 3.75){$2$}
\htext(-0.3 4.5){$0$} \htext(-0.75 5){$1$}
\htext(-1.3 4.75){$1$} \htext(-1.75 5.25){$0$}
\htext(-0.5 6){$2$} \htext(-1.5 6){$2$}
\htext(-0.5 6.75){$3$}
\htext(-0.5 7.25){$3$}

\move(-4 0)\rlvec(-0.2 0)
\move(-4 1)\rlvec(-0.2 0)
\move(-4 1.25)\rlvec(-0.2 0)
\move(-4 1.5)\rlvec(-0.2 0)
\move(-4 1.75)\rlvec(-0.2 0)
\esegment

\move(-10 -12)

\bsegment

\move(-3 1.5) \rlvec(-1 0) \rlvec(0 0.25) \rlvec(1 0) \rlvec(0 -0.25) \ifill f:0.6
\move(-2 2.5) \rlvec(-1 0) \rlvec(0 0.25) \rlvec(1 0) \rlvec(0 -0.25) \ifill f:0.6
\move(-2 3.25) \rlvec(-1 0) \rlvec(0 0.25) \rlvec(1 0) \rlvec(0 -0.25) \ifill f:0.6
\move(-2 3.5) \rlvec(-1 0) \rlvec(0 0.25) \rlvec(1 0) \rlvec(0 -0.25) \ifill f:0.6
\move(-1 3) \rlvec(-1 0) \rlvec(0 0.25) \rlvec(1 0) \rlvec(0 -0.25) \ifill f:0.6
\move(-1 4.25) \rlvec(-1 0) \rlvec(0 0.25) \rlvec(1 0) \rlvec(0 -0.25) \ifill f:0.6
\move(-1 6.5) \rlvec(-1 0) \rlvec(0 0.25) \rlvec(1 0) \rlvec(0 -0.25) \ifill f:0.6
%\move(0 4) \rlvec(-1 0) \rlvec(0 0.25) \rlvec(1 0) \rlvec(0 -0.25) \ifill f:0.6
\move(0 5.25) \rlvec(-1 0) \rlvec(0 0.25) \rlvec(1 0) \rlvec(0 -0.25) \ifill f:0.6
\move(0 7.5) \rlvec(-1 0) \rlvec(0 0.25) \rlvec(1 0) \rlvec(0 -0.25) \ifill f:0.6

\move(0 0)\dtri \move(-1 0)\dtri \move(-2 0)\dtri \move(-3 0)\dtri
\move(0 0)\rlvec(-4 0) \move(0 1)\rlvec(-4 0) \move(0 2)\rlvec(-3
0)

\move(0 1) \rlvec(-1 -1) \move (-1 1) \rlvec(-1 -1)
\move (-2 1) \rlvec(-1 -1) \move (-3 1) \rlvec(-1 -1)

\move(0 2.5)\rlvec(-3 0)  \move(0 3) \rlvec(-2 0)
\move(0 4)\rlvec(-1 0)

\move(-3 1.25) \rlvec(-1 0) \move(-3 1.5) \rlvec(-1 0)
\move (-3 1.75) \rlvec(-1 0) \move(-2 2.75) \rlvec(-1 0)
\move(-2 3.25) \rlvec(-1 0) \move(-2 3.5) \rlvec(-1 0)
\move(-2 3.75) \rlvec(-1 0) \move(-1 3.25) \rlvec(-1 0)
\move(-1 4.25) \rlvec(-1 0) \move(-1 4.5) \rlvec(-1 0)
\move(-1 5.5) \rlvec(-1 0) \move(-1 6.5) \rlvec(-1 0)
\move(-1 6.75) \rlvec(-1 0) \move(-1 5.5) \rlvec(-1 -1)

\move(0 0)\rlvec(0 7.75) \move(-1 0)\rlvec(0 7.75)
\move(-2 0)\rlvec(0 6.75) \move(-3 0)\rlvec(0 3.75)
\move(-4 0)\rlvec(0 1.75)

\move(0 5)\rlvec(-1 -1) \move(0 5)\rlvec(-1 0)
\move(0 5.25) \rlvec(-1 0) \move(0 5.5) \rlvec(-1 0)
\move(0 6.5) \rlvec(-1 0) \move(0 7) \rlvec(-1 0)
\move(0 7.5) \rlvec(-1 0) \move(0 7.75) \rlvec(-1 0)

\htext(-0.3 0.25){$0$} \htext(-0.75 0.75){$1$} \htext(-1.3 0.25){$1$} \htext(-1.75 0.75){$0$}
\htext(-2.3 0.25){$0$} \htext(-2.75 0.75){$1$}\htext(-3.75 0.75){$0$}
\htext(-0.5 1.5){$2$}\htext(-1.5 1.5){$2$} \htext(-2.5 1.5){$2$}
\htext(-0.5 2.25){$3$}\htext(-1.5 2.25){$3$} \htext(-2.5 2.25){$3$}
\htext(-0.5 2.75){$3$}\htext(-1.5 2.75){$3$} \htext(-2.5 3){$3$}
\htext(-0.5 3.5){$2$} \htext(-1.5 3.75){$2$}
\htext(-0.3 4.25){$0$} %\htext(-0.75 5){$1$}
\htext(-1.3 4.75){$1$} \htext(-1.75 5.25){$0$}
\htext(-0.5 6){$2$} \htext(-1.5 6){$2$}
\htext(-0.5 6.75){$3$}
\htext(-0.5 7.25){$3$}

\move(-4 0)\rlvec(-0.2 0)
\move(-4 1)\rlvec(-0.2 0)
\move(-4 1.25)\rlvec(-0.2 0)
\move(-4 1.5)\rlvec(-0.2 0)
\move(-4 1.75)\rlvec(-0.2 0)
\esegment

\move(0 -12)

\bsegment

\move(-3 1.25) \rlvec(-1 0) \rlvec(0 0.25) \rlvec(1 0) \rlvec(0 -0.25) \ifill f:0.6
\move(-3 1.5) \rlvec(-1 0) \rlvec(0 0.25) \rlvec(1 0) \rlvec(0 -0.25) \ifill f:0.6
\move(-2 2.5) \rlvec(-1 0) \rlvec(0 0.25) \rlvec(1 0) \rlvec(0 -0.25) \ifill f:0.6
\move(-2 3.25) \rlvec(-1 0) \rlvec(0 0.25) \rlvec(1 0) \rlvec(0 -0.25) \ifill f:0.6
\move(-2 3.5) \rlvec(-1 0) \rlvec(0 0.25) \rlvec(1 0) \rlvec(0 -0.25) \ifill f:0.6
\move(-1 3) \rlvec(-1 0) \rlvec(0 0.25) \rlvec(1 0) \rlvec(0 -0.25) \ifill f:0.6
\move(-1 4.25) \rlvec(-1 0) \rlvec(0 0.25) \rlvec(1 0) \rlvec(0 -0.25) \ifill f:0.6
\move(-1 6.5) \rlvec(-1 0) \rlvec(0 0.25) \rlvec(1 0) \rlvec(0 -0.25) \ifill f:0.6
\move(0 4) \rlvec(-1 0) \rlvec(0 0.25) \rlvec(1 0) \rlvec(0 -0.25) \ifill f:0.6
\move(0 5.25) \rlvec(-1 0) \rlvec(0 0.25) \rlvec(1 0) \rlvec(0 -0.25) \ifill f:0.6
\move(0 7.5) \rlvec(-1 0) \rlvec(0 0.25) \rlvec(1 0) \rlvec(0 -0.25) \ifill f:0.6

\move(0 0)\dtri \move(-1 0)\dtri \move(-2 0)\dtri \move(-3 0)\dtri
\move(0 0)\rlvec(-4 0) \move(0 1)\rlvec(-4 0) \move(0 2)\rlvec(-3
0)

\move(0 1) \rlvec(-1 -1) \move (-1 1) \rlvec(-1 -1)
\move (-2 1) \rlvec(-1 -1) \move (-3 1) \rlvec(-1 -1)

\move(0 2.5)\rlvec(-3 0)  \move(0 3) \rlvec(-2 0)
\move(0 4)\rlvec(-1 0)

\move(-3 1.25) \rlvec(-1 0) \move(-3 1.5) \rlvec(-1 0)
\move (-3 1.75) \rlvec(-1 0) \move(-2 2.75) \rlvec(-1 0)
\move(-2 3.25) \rlvec(-1 0) \move(-2 3.5) \rlvec(-1 0)
\move(-2 3.75) \rlvec(-1 0) \move(-1 3.25) \rlvec(-1 0)
\move(-1 4.25) \rlvec(-1 0) \move(-1 4.5) \rlvec(-1 0)
\move(-1 5.5) \rlvec(-1 0) \move(-1 6.5) \rlvec(-1 0)
\move(-1 6.75) \rlvec(-1 0) \move(-1 5.5) \rlvec(-1 -1)

\move(0 0)\rlvec(0 7.75) \move(-1 0)\rlvec(0 7.75)
\move(-2 0)\rlvec(0 6.75) \move(-3 0)\rlvec(0 3.75)
\move(-4 0)\rlvec(0 1.75)

\move(0 4.25)\rlvec(-1 0)
\move(0 5.25)\rlvec(-1 -1)\move(0 5.25) \rlvec(-1 0)
\move(0 5.5) \rlvec(-1 0)
\move(0 6.5) \rlvec(-1 0) \move(0 7) \rlvec(-1 0)
\move(0 7.5) \rlvec(-1 0) \move(0 7.75) \rlvec(-1 0)

\htext(-0.3 0.25){$0$} \htext(-0.75 0.75){$1$} \htext(-1.3 0.25){$1$} \htext(-1.75 0.75){$0$}
\htext(-2.3 0.25){$0$} \htext(-2.75 0.75){$1$}\htext(-3.75 0.75){$0$}
\htext(-0.5 1.5){$2$}\htext(-1.5 1.5){$2$} \htext(-2.5 1.5){$2$}
\htext(-0.5 2.25){$3$}\htext(-1.5 2.25){$3$} \htext(-2.5 2.25){$3$}
\htext(-0.5 2.75){$3$}\htext(-1.5 2.75){$3$} \htext(-2.5 3){$3$}
\htext(-0.5 3.5){$2$} \htext(-1.5 3.75){$2$}
\htext(-0.3 4.5){$0$} \htext(-0.75 5){$1$}
\htext(-1.3 4.75){$1$} \htext(-1.75 5.25){$0$}
\htext(-0.5 6){$2$} \htext(-1.5 6){$2$}
\htext(-0.5 6.75){$3$}
\htext(-0.5 7.25){$3$}

\move(-4 0)\rlvec(-0.2 0)
\move(-4 1)\rlvec(-0.2 0)
\move(-4 1.25)\rlvec(-0.2 0)
\move(-4 1.5)\rlvec(-0.2 0)
\move(-4 1.75)\rlvec(-0.2 0)
\esegment

\move(10 -12)

\bsegment

\move(-3 2.5) \rlvec(-1 0) \rlvec(0 0.25) \rlvec(1 0) \rlvec(0 -0.25) \ifill f:0.6
\move(-2 2.5) \rlvec(-1 0) \rlvec(0 0.25) \rlvec(1 0) \rlvec(0 -0.25) \ifill f:0.6
\move(-2 3.25) \rlvec(-1 0) \rlvec(0 0.25) \rlvec(1 0) \rlvec(0 -0.25) \ifill f:0.6
\move(-2 3.5) \rlvec(-1 0) \rlvec(0 0.25) \rlvec(1 0) \rlvec(0 -0.25) \ifill f:0.6
\move(-1 3) \rlvec(-1 0) \rlvec(0 0.25) \rlvec(1 0) \rlvec(0 -0.25) \ifill f:0.6
\move(-1 4.25) \rlvec(-1 0) \rlvec(0 0.25) \rlvec(1 0) \rlvec(0 -0.25) \ifill f:0.6
\move(-1 6.5) \rlvec(-1 0) \rlvec(0 0.25) \rlvec(1 0) \rlvec(0 -0.25) \ifill f:0.6
\move(0 4) \rlvec(-1 0) \rlvec(0 0.25) \rlvec(1 0) \rlvec(0 -0.25) \ifill f:0.6
\move(0 5.25) \rlvec(-1 0) \rlvec(0 0.25) \rlvec(1 0) \rlvec(0 -0.25) \ifill f:0.6
\move(0 7.5) \rlvec(-1 0) \rlvec(0 0.25) \rlvec(1 0) \rlvec(0 -0.25) \ifill f:0.6

\move(0 0)\dtri \move(-1 0)\dtri \move(-2 0)\dtri \move(-3 0)\dtri
\move(0 0)\rlvec(-4 0) \move(0 1)\rlvec(-4 0) \move(0 2)\rlvec(-3
0)

\move(0 1) \rlvec(-1 -1) \move (-1 1) \rlvec(-1 -1)
\move (-2 1) \rlvec(-1 -1) \move (-3 1) \rlvec(-1 -1)

\move(0 2.5)\rlvec(-3 0)  \move(0 3) \rlvec(-2 0)
\move(0 4)\rlvec(-1 0)

\move(-3 1.5) \rlvec(-1 0) \move(-4 1.5) \rlvec(0 1)
\move(-4 2.5) \rlvec(1 0)
\move(-3 2.5) \rlvec(-1 0) \rlvec(0 0.25) \rlvec(1 0)
\move(-3 1.25) \rlvec(-1 0) \move(-3 1.5) \rlvec(-1 0)
%\move (-3 1.75) \rlvec(-1 0)
\move(-2 2.75) \rlvec(-1 0)
\move(-2 3.25) \rlvec(-1 0) \move(-2 3.5) \rlvec(-1 0)
\move(-2 3.75) \rlvec(-1 0)
\move(-1 3.25) \rlvec(-1 0)
\move(-1 4.25) \rlvec(-1 0) \move(-1 4.5) \rlvec(-1 0)
\move(-1 5.5) \rlvec(-1 0) \move(-1 6.5) \rlvec(-1 0)
\move(-1 6.75) \rlvec(-1 0) \move(-1 5.5) \rlvec(-1 -1)

\move(0 0)\rlvec(0 7.75) \move(-1 0)\rlvec(0 7.75)
\move(-2 0)\rlvec(0 6.75) \move(-3 0)\rlvec(0 3.75)
\move(-4 0)\rlvec(0 1.75)

\move(0 4.25) \rlvec(-1 0) \move(0 5.25) \rlvec(-1 -1)
\move(0 5.25) \rlvec(-1 0) \move(0 5.5) \rlvec(-1 0)
\move(0 6.5) \rlvec(-1 0) \move(0 7) \rlvec(-1 0)
\move(0 7.5) \rlvec(-1 0) \move(0 7.75) \rlvec(-1 0)

\htext(-0.3 0.25){$0$} \htext(-0.75 0.75){$1$} \htext(-1.3 0.25){$1$}
\htext(-1.75 0.75){$0$}
\htext(-2.3 0.25){$0$} \htext(-2.75 0.75){$1$}\htext(-3.75 0.75){$0$}
\htext(-0.5 1.5){$2$}\htext(-1.5 1.5){$2$} \htext(-2.5 1.5){$2$} \htext(-3.5 2){$2$}
\htext(-0.5 2.25){$3$}\htext(-1.5 2.25){$3$} \htext(-2.5 2.25){$3$}
\htext(-0.5 2.75){$3$}\htext(-1.5 2.75){$3$} \htext(-2.5 3){$3$}
\htext(-0.5 3.5){$2$} \htext(-1.5 3.75){$2$}
\htext(-0.3 4.5){$0$} \htext(-0.75 5){$1$}
\htext(-1.3 4.75){$1$} \htext(-1.75 5.25){$0$}
\htext(-0.5 6){$2$} \htext(-1.5 6){$2$}
\htext(-0.5 6.75){$3$}
\htext(-0.5 7.25){$3$}

\move(-4 0)\rlvec(-0.2 0)
\move(-4 1)\rlvec(-0.2 0)
\move(-4 1.25)\rlvec(-0.2 0)
\move(-4 1.5)\rlvec(-0.2 0)
\move(-4 1.75)\rlvec(-0.2 0)
\esegment

\move(-9.5 21) \htext{\vdots}
%\move(-1.5 21) \htext{\vdots}
\move(6.5 21) \htext{\vdots}
\move(-11.5 -13) \htext{\vdots}
\move(-1.5 -13) \htext{\vdots}
\move(8.5 -13) \htext{\vdots}

%\move(-1.5 11.5) \ravec(0 -3) \move(-1.2 10) \htext{$1$}
\move(-1.5 -0.5) \ravec(0 -3) \move(-1.2 -2) \htext{$1$}
\move(-8 11.5) \ravec(5.5 -3) \move(-4.3 9.8) \htext{$1$}
\move(4 11.5) \ravec(-4.5 -3) \move(0.9 9.8) \htext{$3$}
\move(-4 -0.5) \ravec(-6.5 -3) \move(-7.2 -1.6) \htext{$0$}
\move(0 -0.5) \ravec(7.5 -3) \move(4 -1.6) \htext{$2$}

\end{texdraw}%
\end{center}

\end{enumerate}
\end{example}

%%%%%%%%%%%%%%%%%%%%%%%%%%%%%%%%%%%%%%%%%%%%%%%%%%%%%%%%%%%%%
\bibliographystyle{amsplain}
%\bibliography{cnone}

\def\cprime{$'$}
\providecommand{\bysame}{\leavevmode\hbox to3em{\hrulefill}\thinspace}
\providecommand{\MR}{\relax\ifhmode\unskip\space\fi MR }
\providecommand{\MRhref}[2]{%
  \href{http://www.ams.org/mathscinet-getitem?mr=#1}{#2}
}
\providecommand{\href}[2]{#2}

%%%%%%%%%%%%%%%%%%%%%%%%%%%%%%%%%%%%%%%%%%%%%%%%%

\appendix
\section{case $i=0,1$}
We are dealing with the $B_n^{(1)}$ case with $i=0 \;\text{or}\, 1$.
We will write the general level-$l$ Young wall column(or, slice)
in the following forms and fix the notations.\\
\savebox{\tmpfiga}{\begin{texdraw}
\fontsize{6}{6}\selectfont
\textref h:C v:C
\drawdim em
\setunitscale 1.7
\move(0 -0.75)\lvec(0 3)\lvec(3.4 3)\lvec(3.4 -0.75)
\move(3.9 3)\lvec(4.4 3)
\move(4.4 3)\lvec(5.8 3)
\move(6.3 3)\lvec(6.8 3)\lvec(6.8 -0.75)
\move(0 2)\lvec(10.2 2)\lvec(10.2 -0.75)
\move(3.4 1)\lvec(13.6 1)\lvec(13.6 -0.75)
\move(6.8 0)\lvec(10.2 0)
\move(1 3)\lvec(1 2)
\move(2.4 3)\lvec(2.4 2)
\move(4.4 3)\lvec(4.4 1)
\move(5.8 3)\lvec(5.8 1)
\move(7.8 2)\lvec(7.8 0)
\move(9.2 2)\lvec(9.2 0)
\move(0.5 3)\lvec(0.5 2)
\move(2.9 3)\lvec(2.9 2)
\move(3.9 3)\lvec(3.9 2)
\move(6.3 3)\lvec(6.3 2)
\htext(0.25 2.5){$0$}\htext(0.75 2.5){$1$}
\htext(1.7 2.5){$\cdots$}
\htext(2.65 2.5){$0$}\htext(3.15 2.5){$1$}
\htext(4.15 2.5){$1$}
\htext(5.1 2.5){$\cdots$}
\htext(6.55 2.5){$1$}
\htext(5.1 1.5){$\cdots$}
\htext(3.9 1.5){$2$}\htext(6.3 1.5){$2$}
\htext(7.3 1.5){$2$}\htext(9.7 1.5){$2$}
\htext(8.5 1.5){$\cdots$}
\htext(7.3 0.5){$3$}\htext(9.7 0.5){$3$}
\htext(8.5 0.5){$\cdots$}
\htext(1.7 -1.25){$\underbrace{\rule{5.6em}{0em}}$}
\htext(5.1 -1.25){$\underbrace{\rule{5.6em}{0em}}$}
\htext(8.5 -1.25){$\underbrace{\rule{5.6em}{0em}}$}
\htext(11.9 -1.25){$\underbrace{\rule{5.6em}{0em}}$}
\htext(1.7 -2){$a_1$}\htext(5.1 -2){$a_2\!-\!a_1$}
\htext(8.5 -2){$a_3$}\htext(11.9 -2){$a_4$}
\move(0 -2.3)
%\drawbb
\end{texdraw}}
\savebox{\tmpfigb}{\begin{texdraw}
\fontsize{6}{6}\selectfont
\textref h:C v:C
\drawdim em
\setunitscale 1.7
\move(0 -0.75)\lvec(0 3)\lvec(3.4 3)\lvec(3.4 -0.75)
\move(3.4 3)\lvec(3.9 3)
\move(4.4 3)\lvec(6.3 3)
\move(6.8 2)\lvec(6.8 -0.75)
\move(0 2)\lvec(10.2 2)\lvec(10.2 -0.75)
\move(3.4 1)\lvec(13.6 1)\lvec(13.6 -0.75)
\move(6.8 0)\lvec(10.2 0)
\move(1 3)\lvec(1 2)
\move(2.4 3)\lvec(2.4 2)
\move(4.4 3)\lvec(4.4 1)
\move(5.8 3)\lvec(5.8 1)
\move(7.8 2)\lvec(7.8 0)
\move(9.2 2)\lvec(9.2 0)
\move(0.5 3)\lvec(0.5 2)
\move(2.9 3)\lvec(2.9 2)
\move(3.9 3)\lvec(3.9 2)
\move(6.3 3)\lvec(6.3 2)
\htext(0.25 2.5){$0$}\htext(0.75 2.5){$1$}
\htext(1.7 2.5){$\cdots$}
\htext(2.65 2.5){$0$}\htext(3.15 2.5){$1$}
\htext(3.65 2.5){$0$}
\htext(5.1 2.5){$\cdots$}
\htext(6.05 2.5){$0$}
\htext(5.1 1.5){$\cdots$}
\htext(3.9 1.5){$2$}\htext(6.3 1.5){$2$}
\htext(7.3 1.5){$2$}\htext(9.7 1.5){$2$}
\htext(8.5 1.5){$\cdots$}
\htext(7.3 0.5){$3$}\htext(9.7 0.5){$3$}
\htext(8.5 0.5){$\cdots$}
\htext(1.7 -1.25){$\underbrace{\rule{5.6em}{0em}}$}
\htext(5.1 -1.25){$\underbrace{\rule{5.6em}{0em}}$}
\htext(8.5 -1.25){$\underbrace{\rule{5.6em}{0em}}$}
\htext(11.9 -1.25){$\underbrace{\rule{5.6em}{0em}}$}
\htext(1.7 -2){$a_2$}\htext(5.1 -2){$a_1\!-\!a_2$}
\htext(8.5 -2){$a_3$}\htext(11.9 -2){$a_4$}
\move(0 -2.3)
%\drawbb
\end{texdraw}}
\savebox{\tmpfigc}{\begin{texdraw}
\fontsize{6}{6}\selectfont
\textref h:C v:C
\drawdim em
\setunitscale 1.7
\move(0 -0.75)\lvec(0 3)\lvec(3.4 3)\lvec(3.4 -0.75)
\move(3.9 3)\lvec(4.4 3)
\move(4.4 3)\lvec(5.8 3)
\move(6.3 3)\lvec(6.8 3)\lvec(6.8 -0.75)
\move(0 2)\lvec(10.2 2)\lvec(10.2 -0.75)
\move(3.4 1)\lvec(13.6 1)\lvec(13.6 -0.75)
\move(6.8 0)\lvec(10.2 0)
\move(1 3)\lvec(1 2)
\move(2.4 3)\lvec(2.4 2)
\move(4.4 3)\lvec(4.4 1)
\move(5.8 3)\lvec(5.8 1)
\move(7.8 2)\lvec(7.8 0)
\move(9.2 2)\lvec(9.2 0)
\move(0.5 3)\lvec(0.5 2)
\move(2.9 3)\lvec(2.9 2)
\move(3.9 3)\lvec(3.9 2)
\move(6.3 3)\lvec(6.3 2)
\htext(0.25 2.5){$1$}\htext(0.75 2.5){$0$}
\htext(1.7 2.5){$\cdots$}
\htext(2.65 2.5){$1$}\htext(3.15 2.5){$0$}
\htext(4.15 2.5){$0$}
\htext(5.1 2.5){$\cdots$}
\htext(6.55 2.5){$0$}
\htext(5.1 1.5){$\cdots$}
\htext(3.9 1.5){$2$}\htext(6.3 1.5){$2$}
\htext(7.3 1.5){$2$}\htext(9.7 1.5){$2$}
\htext(8.5 1.5){$\cdots$}
\htext(7.3 0.5){$3$}\htext(9.7 0.5){$3$}
\htext(8.5 0.5){$\cdots$}
\htext(1.7 -1.25){$\underbrace{\rule{5.6em}{0em}}$}
\htext(5.1 -1.25){$\underbrace{\rule{5.6em}{0em}}$}
\htext(8.5 -1.25){$\underbrace{\rule{5.6em}{0em}}$}
\htext(11.9 -1.25){$\underbrace{\rule{5.6em}{0em}}$}
\htext(1.7 -2){$b_1$}\htext(5.1 -2){$b_2\!-\!b_1$}
\htext(8.5 -2){$b_3$}\htext(11.9 -2){$b_4$}
\move(0 -2.3)
%\drawbb
\end{texdraw}}
\savebox{\tmpfigd}{\begin{texdraw}
\fontsize{6}{6}\selectfont
\textref h:C v:C
\drawdim em
\setunitscale 1.7
\move(0 -0.75)\lvec(0 3)\lvec(3.4 3)\lvec(3.4 -0.75)
\move(3.4 3)\lvec(3.9 3)
\move(4.4 3)\lvec(6.3 3)
\move(6.8 2)\lvec(6.8 -0.75)
\move(0 2)\lvec(10.2 2)\lvec(10.2 -0.75)
\move(3.4 1)\lvec(13.6 1)\lvec(13.6 -0.75)
\move(6.8 0)\lvec(10.2 0)
\move(1 3)\lvec(1 2)
\move(2.4 3)\lvec(2.4 2)
\move(4.4 3)\lvec(4.4 1)
\move(5.8 3)\lvec(5.8 1)
\move(7.8 2)\lvec(7.8 0)
\move(9.2 2)\lvec(9.2 0)
\move(0.5 3)\lvec(0.5 2)
\move(2.9 3)\lvec(2.9 2)
\move(3.9 3)\lvec(3.9 2)
\move(6.3 3)\lvec(6.3 2)
\htext(0.25 2.5){$1$}\htext(0.75 2.5){$0$}
\htext(1.7 2.5){$\cdots$}
\htext(2.65 2.5){$1$}\htext(3.15 2.5){$0$}
\htext(3.65 2.5){$1$}
\htext(5.1 2.5){$\cdots$}
\htext(6.05 2.5){$1$}
\htext(5.1 1.5){$\cdots$}
\htext(3.9 1.5){$2$}\htext(6.3 1.5){$2$}
\htext(7.3 1.5){$2$}\htext(9.7 1.5){$2$}
\htext(8.5 1.5){$\cdots$}
\htext(7.3 0.5){$3$}\htext(9.7 0.5){$3$}
\htext(8.5 0.5){$\cdots$}
\htext(1.7 -1.25){$\underbrace{\rule{5.6em}{0em}}$}
\htext(5.1 -1.25){$\underbrace{\rule{5.6em}{0em}}$}
\htext(8.5 -1.25){$\underbrace{\rule{5.6em}{0em}}$}
\htext(11.9 -1.25){$\underbrace{\rule{5.6em}{0em}}$}
\htext(1.7 -2){$b_2$}\htext(5.1 -2){$b_1\!-\!b_2$}
\htext(8.5 -2){$b_3$}\htext(11.9 -2){$b_4$}
\move(0 -2.3)
%\drawbb
\end{texdraw}}
\begin{center}
\begin{minipage}{0.6\textwidth}
\usebox{\tmpfigb}
\end{minipage}
\quad
\begin{minipage}{0.2\textwidth}
\begin{tabular}{l}
$a_2\le a_1$\\[0.8mm]
$l=a_1+\sum_{i=3}^{4} a_i$\\[0.8mm]
notation : L$(0)_1$
\end{tabular}
\end{minipage}
\end{center}
\begin{center}
\begin{minipage}{0.6\textwidth}
\usebox{\tmpfiga}
\end{minipage}
\quad
\begin{minipage}{0.2\textwidth}
\begin{tabular}{l}
$a_1\le a_2$\\[0.8mm]
$l=\sum_{i=2}^{4} a_i$\\[0.8mm]
notation : L$(1)_0$
\end{tabular}
\end{minipage}
\end{center}
\begin{center}
\begin{minipage}{0.6\textwidth}
\usebox{\tmpfigd}
\end{minipage}
\quad
\begin{minipage}{0.2\textwidth}
\begin{tabular}{l}
$b_2\le b_1$\\[0.8mm]
$l=b_1+\sum_{i=3}^{4} b_i$\\[0.8mm]
notation : R$(1)_0$
\end{tabular}
\end{minipage}
\end{center}
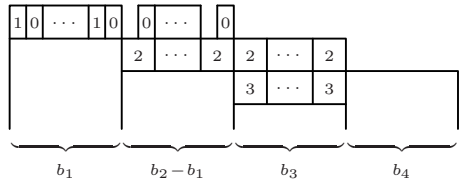
\begin{center}
\begin{minipage}{0.6\textwidth}
\usebox{\tmpfigc}
\end{minipage}
\quad
\begin{minipage}{0.2\textwidth}
\setlength{\tabcolsep}{0.8mm}
\begin{tabular}{l}
$b_1\le b_2$\\[0.8mm]
$l=\sum_{i=2}^{4} b_i$\\[0.8mm]
notation : R$(0)_1$
\end{tabular}
\end{minipage}
\end{center}
Here, $a_1,b_2$ denote the number of layers having $0$-blocks at
the top, $a_2,b_1$ denote the number of layers having $1$-blocks
at the top, and $a_3,b_3$ denote the number of layers having
covering $2$-blocks at the top, as given in the figure. Also
$a_4,b_4$ denote the number of all other layers, that is, any
layer having a top block that comes between the covering $3$-block
and the supporting $2$-block (inclusive) is counted in $a_4,b_4$.

Now, we break each of these into two cases and fix the notations
for each of the cases.
\begin{itemize}
\item $a_2\le a_1$ and $a_3\le a_2$ : L$(0)_1$$(01)_{\overline{2}}$
\item $a_2\le a_1$ and $a_2\le a_3$ : L$(0)_1$$(\overline{2})_{01}$\\
\item $a_1\le a_2$ and $a_3\le a_1$ : L$(1)_0$$(01)_{\overline{2}}$
\item $a_1\le a_2$ and $a_1\le a_3$ : L$(1)_0$$(\overline{2})_{01}$\\
\item $b_2\le b_1$ and $b_3\le b_2$ : R$(1)_0$$(10)_{\overline{2}}$
\item $b_2\le b_1$ and $b_2\le b_3$ : R$(1)_0$$(\overline{2})_{10}$\\
\item $b_1\le b_2$ and $b_3\le b_1$ : R$(0)_1$$(10)_{\overline{2}}$
\item $b_1\le b_2$ and $b_1\le b_3$ : R$(0)_1$$(\overline{2})_{10}$\\
\end{itemize}
The same notation will be used to denote any other
layer-rotation of these slices
and the perfect crystal elements corresponding to these
slices.

If we split every block possible from these slices, the
resulting slices will be of the following eight shapes.\\[2mm]
\savebox{\tmpfiga}{\begin{texdraw}
\fontsize{6}{6}\selectfont
\textref h:C v:C
\drawdim em
\setunitscale 1.7
\move(0 0)
\bsegment
\move(0 -0.75)\lvec(0 2)\lvec(3.4 2)\lvec(3.4 -0.75)
\move(0 1)\lvec(3.4 1)
\move(1 2)\lvec(1 1)
\move(2.4 2)\lvec(2.4 1)
\move(0.5 2)\lvec(0.5 1)
\move(2.9 2)\lvec(2.9 1)
\htext(0.25 1.5){$0$}\htext(2.65 1.5){$0$}
\htext(0.75 1.5){$1$}\htext(3.15 1.5){$1$}
\htext(1.7 1.5){$\cdots$}
\esegment
\move(3.4 0)
\bsegment
%\move(0 -0.75)\lvec(0 1.5)\lvec(3.4 1.5)\lvec(3.4 -0.75)
\move(0 1)\lvec(3.4 1)
\move(0.5 1.5)\lvec(0.5 1)
\move(1 1.5)\lvec(1 1)
\move(2.4 1.5)\lvec(2.4 1)
\move(2.9 1.5)\lvec(2.9 1)
\htext(0.25 1.25){$0$}\htext(2.65 1.25){$0$}
\htext(0.75 1.25){$1$}\htext(3.15 1.25){$1$}
\htext(1.7 1.25){$\cdots$}
\lpatt(0.05 0.15)
\move(0 1.5)\lvec(3.4 1.5)
\esegment
\move(6.8 0)
\bsegment
\move(0 -0.75)\lvec(0 2)\lvec(0.5 2)
\move(1 2)\lvec(2.9 2)
\move(0 1)\lvec(3.4 1)
\move(0 0)\lvec(3.4 0)
\move(1 2)\lvec(1 0)
\move(2.4 2)\lvec(2.4 0)
\move(0.5 2)\lvec(0.5 1)
\move(2.9 2)\lvec(2.9 1)
\htext(0.25 1.5){$0$}\htext(2.65 1.5){$0$}
\htext(0.5 0.5){$2$}\htext(2.9 0.5){$2$}
\htext(1.7 1.5){$\cdots$}
\htext(1.7 0.5){$\cdots$}
\esegment
\move(10.2 0)
\bsegment
\move(0 -0.75)\lvec(0 1.5)
\move(3.4 1.5)\lvec(3.4 -0.75)
\move(0 1)\lvec(3.4 1)
\move(0 0)\lvec(3.4 0)
\move(0.5 1.5)\lvec(0.5 1)
\move(1 1.5)\lvec(1 0)
\move(2.4 1.5)\lvec(2.4 0)
\move(2.9 1.5)\lvec(2.9 1)
\htext(0.25 1.25){$0$}\htext(2.65 1.25){$0$}
\htext(0.75 1.25){$1$}\htext(3.15 1.25){$1$}
\htext(0.5 0.5){$2$}\htext(2.9 0.5){$2$}
\htext(1.7 0.5){$\cdots$}
\htext(1.7 1.25){$\cdots$}
\lpatt(0.05 0.15)
\move(0 1.5)\lvec(3.4 1.5)
\esegment
\move(13.6 0)
\bsegment
\move(0 -0.75)\lvec(0 0)\lvec(3.4 0)\lvec(3.4 -0.75)
\esegment
\htext(1.7 -1.25){$\underbrace{\rule{5.6em}{0em}}$}
\htext(5.1 -1.25){$\underbrace{\rule{5.6em}{0em}}$}
\htext(8.5 -1.25){$\underbrace{\rule{5.6em}{0em}}$}
\htext(11.9 -1.25){$\underbrace{\rule{5.6em}{0em}}$}
\htext(15.3 -1.25){$\underbrace{\rule{5.6em}{0em}}$}
\htext(1.7 -2){$a_2-a_3$}\htext(5.1 -2){$a_3$}
\htext(8.5 -2){$a_1\!-\!a_2$}\htext(11.9 -2){$a_3$}
\htext(15.3 -2){$a_4$}
\move(0 -2.3)
%\drawbb
\end{texdraw}}
\savebox{\tmpfigb}{\begin{texdraw}
\fontsize{6}{6}\selectfont
\textref h:C v:C
\drawdim em
\setunitscale 1.7
\move(0 0)
\bsegment
\move(0 -0.75)\lvec(0 2.5)
\move(3.4 2.5)\lvec(3.4 -0.75)
\move(0 2)\lvec(3.4 2)
\move(0.5 2.5)\lvec(0.5 2)
\move(1 2.5)\lvec(1 2)
\move(2.4 2.5)\lvec(2.4 2)
\move(2.9 2.5)\lvec(2.9 2)
\htext(0.25 2.25){$0$}\htext(2.65 2.25){$0$}
\htext(0.75 2.25){$1$}\htext(3.15 2.25){$1$}
\htext(1.7 2.25){$\cdots$}
\lpatt(0.05 0.15)
\move(0 2.5)\lvec(3.4 2.5)
\esegment
\move(3.4 0)
\bsegment
\move(0 -0.75)\lvec(0 3)\lvec(0.5 3)
\move(1 3)\lvec(2.9 3)
%\move(0 -0.75)\lvec(0 3)\lvec(3.4 3)\lvec(3.4 -0.75)
\move(0 2)\lvec(3.4 2)
\move(0 1)\lvec(3.4 1)
\move(1 3)\lvec(1 1)
\move(2.4 3)\lvec(2.4 1)
\move(0.5 3)\lvec(0.5 2)
\move(2.9 3)\lvec(2.9 2)
\htext(0.25 2.5){$0$}\htext(2.65 2.5){$0$}
\htext(0.5 1.5){$2$}\htext(2.9 1.5){$2$}
\htext(1.7 2.5){$\cdots$}
\htext(1.7 1.5){$\cdots$}
\esegment
\move(6.8 0)
\bsegment
\move(0 -0.75)\lvec(0 2.5)
\move(3.4 2.5)\lvec(3.4 -0.75)
\move(0 2)\lvec(3.4 2)
\move(0 1)\lvec(3.4 1)
\move(0.5 2.5)\lvec(0.5 2)
\move(1 2.5)\lvec(1 1)
\move(2.4 2.5)\lvec(2.4 1)
\move(2.9 2.5)\lvec(2.9 2)
\htext(0.25 2.25){$0$}\htext(2.65 2.25){$0$}
\htext(0.75 2.25){$1$}\htext(3.15 2.25){$1$}
\htext(0.5 1.5){$2$}\htext(2.9 1.5){$2$}
\htext(1.7 1.5){$\cdots$}
\htext(1.7 2.25){$\cdots$}
\lpatt(0.05 0.15)
\move(0 2.5)\lvec(3.4 2.5)
\esegment
\move(10.2 0)
\bsegment
\move(0 -0.75)\lvec(0 2)\lvec(3.4 2)\lvec(3.4 -0.75)
\move(0 0)\lvec(3.4 0)
\move(0 1)\lvec(3.4 1)
\move(1 2)\lvec(1 0)
\move(2.4 2)\lvec(2.4 0)
\htext(1.7 1.5){$\cdots$}
\htext(0.5 1.5){$2$}\htext(2.9 1.5){$2$}
\htext(1.7 0.5){$\cdots$}
\htext(0.5 0.5){$3$}\htext(2.9 0.5){$3$}
\esegment
\move(13.6 0)
\bsegment
\move(0 -0.75)\lvec(0 1)\lvec(3.4 1)\lvec(3.4 -0.75)
\esegment
\htext(1.7 -1.25){$\underbrace{\rule{5.6em}{0em}}$}
\htext(5.1 -1.25){$\underbrace{\rule{5.6em}{0em}}$}
\htext(8.5 -1.25){$\underbrace{\rule{5.6em}{0em}}$}
\htext(11.9 -1.25){$\underbrace{\rule{5.6em}{0em}}$}
\htext(15.3 -1.25){$\underbrace{\rule{5.6em}{0em}}$}
\htext(1.7 -2){$a_2$}\htext(5.1 -2){$a_1\!-\!a_2$}
\htext(8.5 -2){$a_2$}\htext(11.9 -2){$a_3-a_2$}
\htext(15.3 -2){$a_4$}
\move(0 -2.3)
%\drawbb
\end{texdraw}}
\savebox{\tmpfigc}{\begin{texdraw}
\fontsize{6}{6}\selectfont
\textref h:C v:C
\drawdim em
\setunitscale 1.7
\move(0 0)
\bsegment
\move(0 -0.75)\lvec(0 2)\lvec(3.4 2)\lvec(3.4 -0.75)
\move(0 1)\lvec(3.4 1)
\move(1 2)\lvec(1 1)
\move(2.4 2)\lvec(2.4 1)
\move(0.5 2)\lvec(0.5 1)
\move(2.9 2)\lvec(2.9 1)
\htext(0.25 1.5){$0$}\htext(2.65 1.5){$0$}
\htext(0.75 1.5){$1$}\htext(3.15 1.5){$1$}
\htext(1.7 1.5){$\cdots$}
\esegment
\move(3.4 0)
\bsegment
\move(0 -0.75)\lvec(0 1.5)
\move(3.4 1.5)\lvec(3.4 -0.75)
\move(0 1)\lvec(3.4 1)
\move(0.5 1.5)\lvec(0.5 1)
\move(1 1.5)\lvec(1 1)
\move(2.4 1.5)\lvec(2.4 1)
\move(2.9 1.5)\lvec(2.9 1)
\htext(0.25 1.25){$0$}\htext(2.65 1.25){$0$}
\htext(0.75 1.25){$1$}\htext(3.15 1.25){$1$}
\htext(1.7 1.25){$\cdots$}
\lpatt(0.05 0.15)
\move(0 1.5)\lvec(3.4 1.5)
\esegment
\move(6.8 0)
\bsegment
\move(0.5 2)\lvec(1 2)
\move(1 2)\lvec(2.4 2)
\move(2.9 2)\lvec(3.4 2)\lvec(3.4 1.5)
\move(0 1)\lvec(3.4 1)
\move(0 0)\lvec(3.4 0)
\move(1 2)\lvec(1 0)
\move(2.4 2)\lvec(2.4 0)
\move(0.5 2)\lvec(0.5 1)
\move(2.9 2)\lvec(2.9 1)
\htext(0.75 1.5){$1$}\htext(3.15 1.5){$1$}
\htext(0.5 0.5){$2$}\htext(2.9 0.5){$2$}
\htext(1.7 1.5){$\cdots$}
\htext(1.7 0.5){$\cdots$}
\esegment
\move(10.2 0)
\bsegment
\move(0 -0.75)\lvec(0 1.5)
\move(3.4 1.5)\lvec(3.4 -0.75)
\move(0 1)\lvec(3.4 1)
\move(0 0)\lvec(3.4 0)
\move(1 1.5)\lvec(1 0)
\move(0.5 1.5)\lvec(0.5 1)
\move(2.4 1.5)\lvec(2.4 0)
\move(2.9 1.5)\lvec(2.9 1)
\htext(0.25 1.25){$0$}\htext(2.65 1.25){$0$}
\htext(0.75 1.25){$1$}\htext(3.15 1.25){$1$}
\htext(0.5 0.5){$2$}\htext(2.9 0.5){$2$}
\htext(1.7 0.5){$\cdots$}
\htext(1.7 1.25){$\cdots$}
\lpatt(0.05 0.15)
\move(0 1.5)\lvec(3.4 1.5)
\esegment
\move(13.6 0)
\bsegment
\move(0 -0.75)\lvec(0 0)\lvec(3.4 0)\lvec(3.4 -0.75)
\esegment
\htext(1.7 -1.25){$\underbrace{\rule{5.6em}{0em}}$}
\htext(5.1 -1.25){$\underbrace{\rule{5.6em}{0em}}$}
\htext(8.5 -1.25){$\underbrace{\rule{5.6em}{0em}}$}
\htext(11.9 -1.25){$\underbrace{\rule{5.6em}{0em}}$}
\htext(15.3 -1.25){$\underbrace{\rule{5.6em}{0em}}$}
\htext(1.7 -2){$a_1-a_3$}\htext(5.1 -2){$a_3$}
\htext(8.5 -2){$a_2\!-\!a_1$}\htext(11.9 -2){$a_3$}
\htext(15.3 -2){$a_4$}
\move(0 -2.3)
%\drawbb
\end{texdraw}}
\savebox{\tmpfigd}{\begin{texdraw}
\fontsize{6}{6}\selectfont
\textref h:C v:C
\drawdim em
\setunitscale 1.7
\move(0 0)
\bsegment
\move(0 -0.75)\lvec(0 2.5)
\move(3.4 2.5)\lvec(3.4 -0.75)
\move(0 2)\lvec(3.4 2)
\move(1 2.5)\lvec(1 2)
\move(0.5 2.5)\lvec(0.5 2)
\move(2.4 2.5)\lvec(2.4 2)
\move(2.9 2.5)\lvec(2.9 2)
\htext(0.25 2.25){$0$}\htext(2.65 2.25){$0$}
\htext(0.75 2.25){$1$}\htext(3.15 2.25){$1$}
\htext(1.7 2.25){$\cdots$}
\lpatt(0.05 0.15)
\move(0 2.5)\lvec(3.4 2.5)
\esegment
\move(3.4 0)
\bsegment
\move(0.5 3)\lvec(1 3)
\move(1 3)\lvec(2.4 3)
\move(2.9 3)\lvec(3.4 3)\lvec(3.4 2.5)
%\move(0 -0.75)\lvec(0 3)\lvec(3.4 3)\lvec(3.4 -0.75)
\move(0 2)\lvec(3.4 2)
\move(0 1)\lvec(3.4 1)
\move(1 3)\lvec(1 1)
\move(2.4 3)\lvec(2.4 1)
\move(0.5 3)\lvec(0.5 2)
\move(2.9 3)\lvec(2.9 2)
\htext(0.75 2.5){$1$}\htext(3.15 2.5){$1$}
\htext(0.5 1.5){$2$}\htext(2.9 1.5){$2$}
\htext(1.7 2.5){$\cdots$}
\htext(1.7 1.5){$\cdots$}
\esegment
\move(6.8 0)
\bsegment
\move(0 -0.75)\lvec(0 2.5)
\move(3.4 2.5)\lvec(3.4 -0.75)
\move(0 2)\lvec(3.4 2)
\move(0 1)\lvec(3.4 1)
\move(1 2.5)\lvec(1 1)
\move(0.5 2.5)\lvec(0.5 2)
\move(2.4 2.5)\lvec(2.4 1)
\move(2.9 2.5)\lvec(2.9 2)
\htext(0.25 2.25){$0$}\htext(2.65 2.25){$0$}
\htext(0.75 2.25){$1$}\htext(3.15 2.25){$1$}
\htext(0.5 1.5){$2$}\htext(2.9 1.5){$2$}
\htext(1.7 1.5){$\cdots$}
\htext(1.7 2.25){$\cdots$}
\lpatt(0.05 0.15)
\move(0 2.5)\lvec(3.4 2.5)
\esegment
\move(10.2 0)
\bsegment
\move(0 -0.75)\lvec(0 2)\lvec(3.4 2)\lvec(3.4 -0.75)
\move(0 0)\lvec(3.4 0)
\move(0 1)\lvec(3.4 1)
\move(1 2)\lvec(1 0)
\move(2.4 2)\lvec(2.4 0)
\htext(1.7 1.5){$\cdots$}
\htext(0.5 1.5){$2$}\htext(2.9 1.5){$2$}
\htext(1.7 0.5){$\cdots$}
\htext(0.5 0.5){$3$}\htext(2.9 0.5){$3$}
\esegment
\move(13.6 0)
\bsegment
\move(0 -0.75)\lvec(0 1)\lvec(3.4 1)\lvec(3.4 -0.75)
\esegment
\htext(1.7 -1.25){$\underbrace{\rule{5.6em}{0em}}$}
\htext(5.1 -1.25){$\underbrace{\rule{5.6em}{0em}}$}
\htext(8.5 -1.25){$\underbrace{\rule{5.6em}{0em}}$}
\htext(11.9 -1.25){$\underbrace{\rule{5.6em}{0em}}$}
\htext(15.3 -1.25){$\underbrace{\rule{5.6em}{0em}}$}
\htext(1.7 -2){$a_1$}\htext(5.1 -2){$a_2\!-\!a_1$}
\htext(8.5 -2){$a_1$}\htext(11.9 -2){$a_3-a_1$}
\htext(15.3 -2){$a_4$}
\move(0 -2.3)
%\drawbb
\end{texdraw}}
\savebox{\tmpfige}{\begin{texdraw}
\fontsize{6}{6}\selectfont
\textref h:C v:C
\drawdim em
\setunitscale 1.7
\move(0 0)
\bsegment
\move(0 -0.75)\lvec(0 2)\lvec(3.4 2)\lvec(3.4 -0.75)
\move(0 1)\lvec(3.4 1)
\move(1 2)\lvec(1 1)
\move(2.4 2)\lvec(2.4 1)
\move(0.5 2)\lvec(0.5 1)
\move(2.9 2)\lvec(2.9 1)
\htext(0.25 1.5){$1$}\htext(2.65 1.5){$1$}
\htext(0.75 1.5){$0$}\htext(3.15 1.5){$0$}
\htext(1.7 1.5){$\cdots$}
\esegment
\move(3.4 0)
\bsegment
\move(0 -0.75)\lvec(0 1.5)
\move(3.4 1.5)\lvec(3.4 -0.75)
\move(0 1)\lvec(3.4 1)
\move(0.5 1.5)\lvec(0.5 1)
\move(1 1.5)\lvec(1 1)
\move(2.4 1.5)\lvec(2.4 1)
\move(2.9 1.5)\lvec(2.9 1)
\htext(0.25 1.25){$1$}\htext(2.65 1.25){$1$}
\htext(0.75 1.25){$0$}\htext(3.15 1.25){$0$}
\htext(1.7 1.25){$\cdots$}
\lpatt(0.05 0.15)
\move(0 1.5)\lvec(3.4 1.5)
\esegment
\move(6.8 0)
\bsegment
%\move(0 -0.75)\lvec(0 2)\lvec(3.4 2)\lvec(3.4 -0.75)
\move(0 -0.75)\lvec(0 2)\lvec(0.5 2)
\move(1 2)\lvec(2.9 2)
\move(0 1)\lvec(3.4 1)
\move(0 0)\lvec(3.4 0)
\move(1 2)\lvec(1 0)
\move(2.4 2)\lvec(2.4 0)
\move(0.5 2)\lvec(0.5 1)
\move(2.9 2)\lvec(2.9 1)
\htext(0.25 1.5){$1$}\htext(2.65 1.5){$1$}
\htext(0.5 0.5){$2$}\htext(2.9 0.5){$2$}
\htext(1.7 1.5){$\cdots$}
\htext(1.7 0.5){$\cdots$}
\esegment
\move(10.2 0)
\bsegment
\move(0 -0.75)\lvec(0 1.5)
\move(3.4 1.5)\lvec(3.4 -0.75)
\move(0 1)\lvec(3.4 1)
\move(0 0)\lvec(3.4 0)
\move(0.5 1.5)\lvec(0.5 1)
\move(1 1.5)\lvec(1 0)
\move(2.4 1.5)\lvec(2.4 0)
\move(2.9 1.5)\lvec(2.9 1)
\htext(0.25 1.25){$1$}\htext(2.65 1.25){$1$}
\htext(0.75 1.25){$0$}\htext(3.15 1.25){$0$}
\htext(0.5 0.5){$2$}\htext(2.9 0.5){$2$}
\htext(1.7 0.5){$\cdots$}
\htext(1.7 1.25){$\cdots$}
\lpatt(0.05 0.15)
\move(0 1.5)\lvec(3.4 1.5)
\esegment
\move(13.6 0)
\bsegment
\move(0 -0.75)\lvec(0 0)\lvec(3.4 0)\lvec(3.4 -0.75)
\esegment
\htext(1.7 -1.25){$\underbrace{\rule{5.6em}{0em}}$}
\htext(5.1 -1.25){$\underbrace{\rule{5.6em}{0em}}$}
\htext(8.5 -1.25){$\underbrace{\rule{5.6em}{0em}}$}
\htext(11.9 -1.25){$\underbrace{\rule{5.6em}{0em}}$}
\htext(15.3 -1.25){$\underbrace{\rule{5.6em}{0em}}$}
\htext(1.7 -2){$b_2-b_3$}\htext(5.1 -2){$b_3$}
\htext(8.5 -2){$b_1\!-\!b_2$}\htext(11.9 -2){$b_3$}
\htext(15.3 -2){$b_4$}
\move(0 -2.3)
%\drawbb
\end{texdraw}}
\savebox{\tmpfigf}{\begin{texdraw}
\fontsize{6}{6}\selectfont
\textref h:C v:C
\drawdim em
\setunitscale 1.7
\move(0 0)
\bsegment
\move(0 -0.75)\lvec(0 2.5)
\move(3.4 2.5)\lvec(3.4 -0.75)
\move(0 2)\lvec(3.4 2)
\move(0.5 2.5)\lvec(0.5 2)
\move(1 2.5)\lvec(1 2)
\move(2.4 2.5)\lvec(2.4 2)
\move(2.9 2.5)\lvec(2.9 2)
\htext(0.25 2.25){$1$}\htext(2.65 2.25){$1$}
\htext(0.75 2.25){$0$}\htext(3.15 2.25){$0$}
\htext(1.7 2.25){$\cdots$}
\lpatt(0.05 0.15)
\move(0 2.5)\lvec(3.4 2.5)
\esegment
\move(3.4 0)
\bsegment
\move(0 -0.75)\lvec(0 3)\lvec(0.5 3)
\move(1 3)\lvec(2.9 3)
%\move(0 -0.75)\lvec(0 3)\lvec(3.4 3)\lvec(3.4 -0.75)
\move(0 2)\lvec(3.4 2)
\move(0 1)\lvec(3.4 1)
\move(1 3)\lvec(1 1)
\move(2.4 3)\lvec(2.4 1)
\move(0.5 3)\lvec(0.5 2)
\move(2.9 3)\lvec(2.9 2)
\htext(0.25 2.5){$1$}\htext(2.65 2.5){$1$}
\htext(0.5 1.5){$2$}\htext(2.9 1.5){$2$}
\htext(1.7 2.5){$\cdots$}
\htext(1.7 1.5){$\cdots$}
\esegment
\move(6.8 0)
\bsegment
\move(0 -0.75)\lvec(0 2.5)
\move(3.4 2.5)\lvec(3.4 -0.75)
\move(0 2)\lvec(3.4 2)
\move(0 1)\lvec(3.4 1)
\move(0.5 2.5)\lvec(0.5 2)
\move(1 2.5)\lvec(1 1)
\move(2.4 2.5)\lvec(2.4 1)
\move(2.9 2.5)\lvec(2.9 2)
\htext(0.25 2.25){$1$}\htext(2.65 2.25){$1$}
\htext(0.75 2.25){$0$}\htext(3.15 2.25){$0$}
\htext(0.5 1.5){$2$}\htext(2.9 1.5){$2$}
\htext(1.7 1.5){$\cdots$}
\htext(1.7 2.25){$\cdots$}
\lpatt(0.05 0.15)
\move(0 2.5)\lvec(3.4 2.5)
\esegment
\move(10.2 0)
\bsegment
\move(0 -0.75)\lvec(0 2)\lvec(3.4 2)\lvec(3.4 -0.75)
\move(0 0)\lvec(3.4 0)
\move(0 1)\lvec(3.4 1)
\move(1 2)\lvec(1 0)
\move(2.4 2)\lvec(2.4 0)
\htext(1.7 1.5){$\cdots$}
\htext(0.5 1.5){$2$}\htext(2.9 1.5){$2$}
\htext(1.7 0.5){$\cdots$}
\htext(0.5 0.5){$3$}\htext(2.9 0.5){$3$}
\esegment
\move(13.6 0)
\bsegment
\move(0 -0.75)\lvec(0 1)\lvec(3.4 1)\lvec(3.4 -0.75)
\esegment
\htext(1.7 -1.25){$\underbrace{\rule{5.6em}{0em}}$}
\htext(5.1 -1.25){$\underbrace{\rule{5.6em}{0em}}$}
\htext(8.5 -1.25){$\underbrace{\rule{5.6em}{0em}}$}
\htext(11.9 -1.25){$\underbrace{\rule{5.6em}{0em}}$}
\htext(15.3 -1.25){$\underbrace{\rule{5.6em}{0em}}$}
\htext(1.7 -2){$b_2$}\htext(5.1 -2){$b_1\!-\!b_2$}
\htext(8.5 -2){$b_2$}\htext(11.9 -2){$b_3-b_2$}
\htext(15.3 -2){$b_4$}
\move(0 -2.3)
%\drawbb
\end{texdraw}}
\savebox{\tmpfigg}{\begin{texdraw}
\fontsize{6}{6}\selectfont
\textref h:C v:C
\drawdim em
\setunitscale 1.7
\move(0 0)
\bsegment
\move(0 -0.75)\lvec(0 2)\lvec(3.4 2)\lvec(3.4 -0.75)
\move(0 1)\lvec(3.4 1)
\move(1 2)\lvec(1 1)
\move(2.4 2)\lvec(2.4 1)
\move(0.5 2)\lvec(0.5 1)
\move(2.9 2)\lvec(2.9 1)
\htext(0.25 1.5){$1$}\htext(2.65 1.5){$1$}
\htext(0.75 1.5){$0$}\htext(3.15 1.5){$0$}
\htext(1.7 1.5){$\cdots$}
\esegment
\move(3.4 0)
\bsegment
\move(0 -0.75)\lvec(0 1.5)
\move(3.4 1.5)\lvec(3.4 -0.75)
\move(0 1)\lvec(3.4 1)
\move(0.5 1.5)\lvec(0.5 1)
\move(1 1.5)\lvec(1 1)
\move(2.4 1.5)\lvec(2.4 1)
\move(2.9 1.5)\lvec(2.9 1)
\htext(0.25 1.25){$1$}\htext(2.65 1.25){$1$}
\htext(0.75 1.25){$0$}\htext(3.15 1.25){$0$}
\htext(1.7 1.25){$\cdots$}
\lpatt(0.05 0.15)
\move(0 1.5)\lvec(3.4 1.5)
\esegment
\move(6.8 0)
\bsegment
\move(0.5 2)\lvec(1 2)
\move(1 2)\lvec(2.4 2)
\move(2.9 2)\lvec(3.4 2)\lvec(3.4 1.5)
%\move(0 -0.75)\lvec(0 2)\lvec(3.4 2)\lvec(3.4 -0.75)
\move(0 1)\lvec(3.4 1)
\move(0 0)\lvec(3.4 0)
\move(1 2)\lvec(1 0)
\move(2.4 2)\lvec(2.4 0)
\move(0.5 2)\lvec(0.5 1)
\move(2.9 2)\lvec(2.9 1)
\htext(0.75 1.5){$0$}\htext(3.15 1.5){$0$}
\htext(0.5 0.5){$2$}\htext(2.9 0.5){$2$}
\htext(1.7 1.5){$\cdots$}
\htext(1.7 0.5){$\cdots$}
\esegment
\move(10.2 0)
\bsegment
\move(0 -0.75)\lvec(0 1.5)
\move(3.4 1.5)\lvec(3.4 -0.75)
\move(0 1)\lvec(3.4 1)
\move(0 0)\lvec(3.4 0)
\move(1 1.5)\lvec(1 0)
\move(0.5 1.5)\lvec(0.5 1)
\move(2.4 1.5)\lvec(2.4 0)
\move(2.9 1.5)\lvec(2.9 1)
\htext(0.25 1.25){$1$}\htext(2.65 1.25){$1$}
\htext(0.75 1.25){$0$}\htext(3.15 1.25){$0$}
\htext(0.5 0.5){$2$}\htext(2.9 0.5){$2$}
\htext(1.7 0.5){$\cdots$}
\htext(1.7 1.25){$\cdots$}
\lpatt(0.05 0.15)
\move(0 1.5)\lvec(3.4 1.5)
\esegment
\move(13.6 0)
\bsegment
\move(0 -0.75)\lvec(0 0)\lvec(3.4 0)\lvec(3.4 -0.75)
\esegment
\htext(1.7 -1.25){$\underbrace{\rule{5.6em}{0em}}$}
\htext(5.1 -1.25){$\underbrace{\rule{5.6em}{0em}}$}
\htext(8.5 -1.25){$\underbrace{\rule{5.6em}{0em}}$}
\htext(11.9 -1.25){$\underbrace{\rule{5.6em}{0em}}$}
\htext(15.3 -1.25){$\underbrace{\rule{5.6em}{0em}}$}
\htext(1.7 -2){$b_1-b_3$}\htext(5.1 -2){$b_3$}
\htext(8.5 -2){$b_2\!-\!b_1$}\htext(11.9 -2){$b_3$}
\htext(15.3 -2){$b_4$}
\move(0 -2.3)
%\drawbb
\end{texdraw}}
\savebox{\tmpfigh}{\begin{texdraw}
\fontsize{6}{6}\selectfont
\textref h:C v:C
\drawdim em
\setunitscale 1.7
\move(0 0)
\bsegment
\move(0 -0.75)\lvec(0 2.5)
\move(3.4 2.5)\lvec(3.4 -0.75)
\move(0 2)\lvec(3.4 2)
\move(1 2.5)\lvec(1 2)
\move(0.5 2.5)\lvec(0.5 2)
\move(2.4 2.5)\lvec(2.4 2)
\move(2.9 2.5)\lvec(2.9 2)
\htext(0.25 2.25){$1$}\htext(2.65 2.25){$1$}
\htext(0.75 2.25){$0$}\htext(3.15 2.25){$0$}
\htext(1.7 2.25){$\cdots$}
\lpatt(0.05 0.15)
\move(0 2.5)\lvec(3.4 2.5)
\esegment
\move(3.4 0)
\bsegment
%\move(0 -0.75)\lvec(0 3)\lvec(3.4 3)\lvec(3.4 -0.75)
\move(0.5 3)\lvec(1 3)
\move(1 3)\lvec(2.4 3)
\move(2.9 3)\lvec(3.4 3)\lvec(3.4 2.5)
\move(0 2)\lvec(3.4 2)
\move(0 1)\lvec(3.4 1)
\move(1 3)\lvec(1 1)
\move(2.4 3)\lvec(2.4 1)
\move(0.5 3)\lvec(0.5 2)
\move(2.9 3)\lvec(2.9 2)
\htext(0.75 2.5){$0$}\htext(3.15 2.5){$0$}
\htext(0.5 1.5){$2$}\htext(2.9 1.5){$2$}
\htext(1.7 2.5){$\cdots$}
\htext(1.7 1.5){$\cdots$}
\esegment
\move(6.8 0)
\bsegment
\move(0 -0.75)\lvec(0 2.5)
\move(3.4 2.5)\lvec(3.4 -0.75)
\move(0 2)\lvec(3.4 2)
\move(0 1)\lvec(3.4 1)
\move(1 2.5)\lvec(1 1)
\move(0.5 2.5)\lvec(0.5 2)
\move(2.4 2.5)\lvec(2.4 1)
\move(2.9 2.5)\lvec(2.9 2)
\htext(0.25 2.25){$1$}\htext(2.65 2.25){$1$}
\htext(0.75 2.25){$0$}\htext(3.15 2.25){$0$}
\htext(0.5 1.5){$2$}\htext(2.9 1.5){$2$}
\htext(1.7 1.5){$\cdots$}
\htext(1.7 2.25){$\cdots$}
\lpatt(0.05 0.15)
\move(0 2.5)\lvec(3.4 2.5)
\esegment
\move(10.2 0)
\bsegment
\move(0 -0.75)\lvec(0 2)\lvec(3.4 2)\lvec(3.4 -0.75)
\move(0 0)\lvec(3.4 0)
\move(0 1)\lvec(3.4 1)
\move(1 2)\lvec(1 0)
\move(2.4 2)\lvec(2.4 0)
\htext(1.7 1.5){$\cdots$}
\htext(0.5 1.5){$2$}\htext(2.9 1.5){$2$}
\htext(1.7 0.5){$\cdots$}
\htext(0.5 0.5){$3$}\htext(2.9 0.5){$3$}
\esegment
\move(13.6 0)
\bsegment
\move(0 -0.75)\lvec(0 1)\lvec(3.4 1)\lvec(3.4 -0.75)
\esegment
\htext(1.7 -1.25){$\underbrace{\rule{5.6em}{0em}}$}
\htext(5.1 -1.25){$\underbrace{\rule{5.6em}{0em}}$}
\htext(8.5 -1.25){$\underbrace{\rule{5.6em}{0em}}$}
\htext(11.9 -1.25){$\underbrace{\rule{5.6em}{0em}}$}
\htext(15.3 -1.25){$\underbrace{\rule{5.6em}{0em}}$}
\htext(1.7 -2){$b_1$}\htext(5.1 -2){$b_2\!-\!b_1$}
\htext(8.5 -2){$b_1$}\htext(11.9 -2){$b_3-b_1$}
\htext(15.3 -2){$b_4$}
\move(0 -2.3)
%\drawbb
\end{texdraw}}
\noindent
case L$(0)_1$$(01)_{\overline{2}}$ :
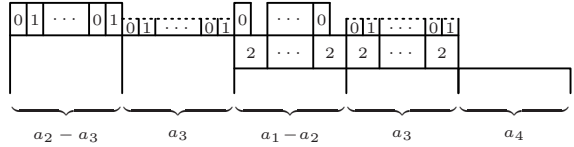
\begin{center}
\usebox{\tmpfiga}
\end{center}
\noindent
case L$(0)_1$$(\overline{2})_{01}$ :
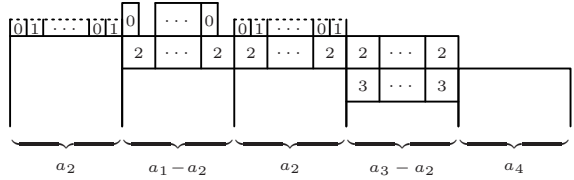
\begin{center}
\usebox{\tmpfigb}
\end{center}
\noindent
case L$(1)_0$$(01)_{\overline{2}}$ :
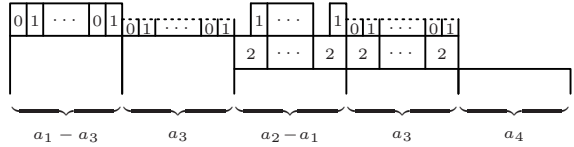
\begin{center}
\usebox{\tmpfigc}
\end{center}
\noindent
case L$(1)_0$$(\overline{2})_{01}$ :
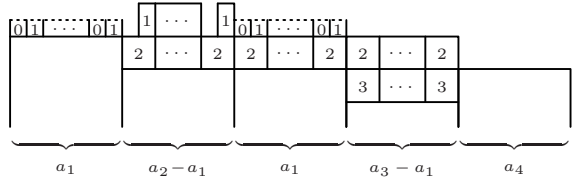
\begin{center}
\usebox{\tmpfigd}
\end{center}
%%%%%%%%%%%
case R$(1)_0$$(10)_{\overline{2}}$ :
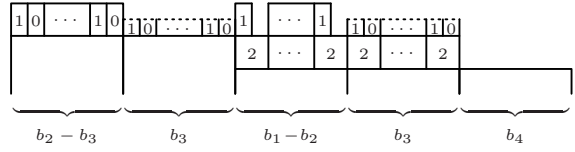
\begin{center}
\usebox{\tmpfige}
\end{center}
\noindent
case R$(1)_0$$(\overline{2})_{10}$ :
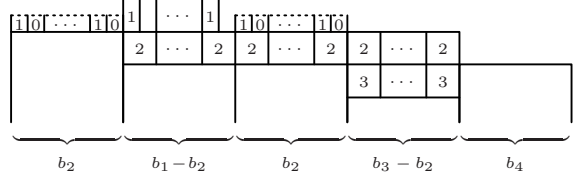
\begin{center}
\usebox{\tmpfigf}
\end{center}
\noindent
case R$(0)_1$$(10)_{\overline{2}}$ :
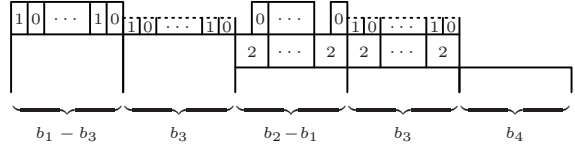
\begin{center}
\usebox{\tmpfigg}
\end{center}
\noindent
case R$(0)_1$$(\overline{2})_{10}$ :
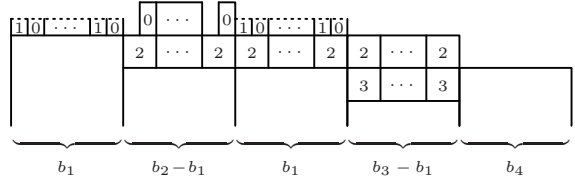
\begin{center}
\usebox{\tmpfigh}
\end{center}

Now, we will place two slices side by side and also consider the
corresponding pair of perfect crystal elements. For the left
slice, we will use L$(0)_1$$(01)_{\overline{2}}$,
L$(0)_1$$(\overline{2})_{01}$, L$(1)_0$$(01)_{\overline{2}}$, or
L$(1)_0$$(\overline{2})_{01}$ and for the right column, use
R$(1)_0$$(10)_{\overline{2}}$, R$(1)_0$$(\overline{2})_{10}$,
R$(0)_1$$(10)_{\overline{2}}$, or R$(0)_1$$(\overline{2})_{10}$,
as given in the above figures.

First, we fix the right slice and next remove finitely many
$\delta$'s from the left slice. We give the following starting
shapes and relative heights. Shape of the right slice is to be
given in the form R$(0)_1$ or R$(1)_0$. And the starting shape of
the left slice should be so that when every possible block is
split, all bottom halves of $01$-blocks appear at the rear layers.
The following is an example showing right slice
R$(1)_0$$(10)_{\overline{2}}$
and left slice L$(0)_1$$(\overline{2})_{01}$.\\[2mm]
\savebox{\tmpfige}{\begin{texdraw}
\fontsize{6}{6}\selectfont
\textref h:C v:C
\drawdim em
\setunitscale 1.7
\move(0 0)
\bsegment
\move(0 -0.75)\lvec(0 2)\lvec(3.4 2)\lvec(3.4 -0.75)
\move(0 1)\lvec(3.4 1)
\move(1 2)\lvec(1 1)
\move(2.4 2)\lvec(2.4 1)
\move(0.5 2)\lvec(0.5 1)
\move(2.9 2)\lvec(2.9 1)
\htext(0.25 1.5){$1$}\htext(2.65 1.5){$1$}
\htext(0.75 1.5){$0$}\htext(3.15 1.5){$0$}
\htext(1.7 1.5){$\cdots$}
\esegment
\move(3.4 0)
\bsegment
\move(0 -0.75)\lvec(0 1.5)
\move(3.4 1.5)\lvec(3.4 -0.75)
\move(0 1)\lvec(3.4 1)
\move(0.5 1.5)\lvec(0.5 1)
\move(1 1.5)\lvec(1 1)
\move(2.4 1.5)\lvec(2.4 1)
\move(2.9 1.5)\lvec(2.9 1)
\htext(0.25 1.25){$1$}\htext(2.65 1.25){$1$}
\htext(0.75 1.25){$0$}\htext(3.15 1.25){$0$}
\htext(1.7 1.25){$\cdots$}
\lpatt(0.05 0.15)
\move(0 1.5)\lvec(3.4 1.5)
\esegment
\move(6.8 0)
\bsegment
\move(0 -0.75)\lvec(0 2)\lvec(0.5 2)
\move(1 2)\lvec(2.9 2)
%\move(0 -0.75)\lvec(0 2)\lvec(3.4 2)\lvec(3.4 -0.75)
\move(0 1)\lvec(3.4 1)
\move(0 0)\lvec(3.4 0)
\move(1 2)\lvec(1 0)
\move(2.4 2)\lvec(2.4 0)
\move(0.5 2)\lvec(0.5 1)
\move(2.9 2)\lvec(2.9 1)
\htext(0.25 1.5){$1$}\htext(2.65 1.5){$1$}
\htext(0.5 0.5){$2$}\htext(2.9 0.5){$2$}
\htext(1.7 1.5){$\cdots$}
\htext(1.7 0.5){$\cdots$}
\esegment
\move(10.2 0)
\bsegment
\move(0 -0.75)\lvec(0 1.5)
\move(3.4 1.5)\lvec(3.4 -0.75)
\move(0 1)\lvec(3.4 1)
\move(0 0)\lvec(3.4 0)
\move(0.5 1.5)\lvec(0.5 1)
\move(1 1.5)\lvec(1 0)
\move(2.4 1.5)\lvec(2.4 0)
\move(2.9 1.5)\lvec(2.9 1)
\htext(0.25 1.25){$1$}\htext(2.65 1.25){$1$}
\htext(0.75 1.25){$0$}\htext(3.15 1.25){$0$}
\htext(0.5 0.5){$2$}\htext(2.9 0.5){$2$}
\htext(1.7 0.5){$\cdots$}
\htext(1.7 1.25){$\cdots$}
\lpatt(0.05 0.15)
\move(0 1.5)\lvec(3.4 1.5)
\esegment
\move(13.6 0)
\bsegment
\move(0 -0.75)\lvec(0 0)\lvec(3.4 0)\lvec(3.4 -0.75)
\esegment
\htext(1.7 -1.25){$\underbrace{\rule{5.6em}{0em}}$}
\htext(5.1 -1.25){$\underbrace{\rule{5.6em}{0em}}$}
\htext(8.5 -1.25){$\underbrace{\rule{5.6em}{0em}}$}
\htext(11.9 -1.25){$\underbrace{\rule{5.6em}{0em}}$}
\htext(15.3 -1.25){$\underbrace{\rule{5.6em}{0em}}$}
\htext(1.7 -2){$b_2-b_3$}\htext(5.1 -2){$b_3$}
\htext(8.5 -2){$b_1\!-\!b_2$}\htext(11.9 -2){$b_3$}
\htext(15.3 -2){$b_4$}
\move(0 -2.3)
%\drawbb
\end{texdraw}}
\savebox{\tmpfigb}{\begin{texdraw}
\fontsize{6}{6}\selectfont
\textref h:C v:C
\drawdim em
\setunitscale 1.7
\move(0 0)
\bsegment
\move(0 -0.75)\lvec(0 2.5)
\move(3.4 2.5)\lvec(3.4 -0.75)
\move(0 2)\lvec(3.4 2)
\move(0.5 2.5)\lvec(0.5 2)
\move(1 2.5)\lvec(1 2)
\move(2.4 2.5)\lvec(2.4 2)
\move(2.9 2.5)\lvec(2.9 2)
\htext(0.25 2.25){$0$}\htext(2.65 2.25){$0$}
\htext(0.75 2.25){$1$}\htext(3.15 2.25){$1$}
\htext(1.7 2.25){$\cdots$}
\lpatt(0.05 0.15)
\move(0 2.5)\lvec(3.4 2.5)
\esegment
\move(3.4 0)
\bsegment
\move(0 -0.75)\lvec(0 3)\lvec(0.5 3)
\move(1 3)\lvec(2.9 3)
%\move(0 -0.75)\lvec(0 3)\lvec(3.4 3)\lvec(3.4 -0.75)
\move(0 2)\lvec(3.4 2)
\move(0 1)\lvec(3.4 1)
\move(1 3)\lvec(1 1)
\move(2.4 3)\lvec(2.4 1)
\move(0.5 3)\lvec(0.5 2)
\move(2.9 3)\lvec(2.9 2)
\htext(0.25 2.5){$0$}\htext(2.65 2.5){$0$}
\htext(0.5 1.5){$2$}\htext(2.9 1.5){$2$}
\htext(1.7 2.5){$\cdots$}
\htext(1.7 1.5){$\cdots$}
\esegment
\move(6.8 0)
\bsegment
\move(0 -0.75)\lvec(0 2.5)
\move(3.4 2.5)\lvec(3.4 -0.75)
\move(0 2)\lvec(3.4 2)
\move(0 1)\lvec(3.4 1)
\move(0.5 2.5)\lvec(0.5 2)
\move(1 2.5)\lvec(1 1)
\move(2.4 2.5)\lvec(2.4 1)
\move(2.9 2.5)\lvec(2.9 2)
\htext(0.25 2.25){$0$}\htext(2.65 2.25){$0$}
\htext(0.75 2.25){$1$}\htext(3.15 2.25){$1$}
\htext(0.5 1.5){$2$}\htext(2.9 1.5){$2$}
\htext(1.7 1.5){$\cdots$}
\htext(1.7 2.25){$\cdots$}
\lpatt(0.05 0.15)
\move(0 2.5)\lvec(3.4 2.5)
\esegment
\move(10.2 0)
\bsegment
\move(0 -0.75)\lvec(0 2)\lvec(3.4 2)\lvec(3.4 -0.75)
\move(0 0)\lvec(3.4 0)
\move(0 1)\lvec(3.4 1)
\move(1 2)\lvec(1 0)
\move(2.4 2)\lvec(2.4 0)
\htext(1.7 1.5){$\cdots$}
\htext(0.5 1.5){$2$}\htext(2.9 1.5){$2$}
\htext(1.7 0.5){$\cdots$}
\htext(0.5 0.5){$3$}\htext(2.9 0.5){$3$}
\esegment
\move(13.6 0)
\bsegment
\move(0 -0.75)\lvec(0 1)\lvec(3.4 1)\lvec(3.4 -0.75)
\esegment
\htext(1.7 -1.25){$\underbrace{\rule{5.6em}{0em}}$}
\htext(5.1 -1.25){$\underbrace{\rule{5.6em}{0em}}$}
\htext(8.5 -1.25){$\underbrace{\rule{5.6em}{0em}}$}
\htext(11.9 -1.25){$\underbrace{\rule{5.6em}{0em}}$}
\htext(15.3 -1.25){$\underbrace{\rule{5.6em}{0em}}$}
\htext(1.7 -2){$a_2$}\htext(5.1 -2){$a_1\!-\!a_2$}
\htext(8.5 -2){$a_2$}\htext(11.9 -2){$a_3-a_2$}
\htext(15.3 -2){$a_4$}
\move(0 -2.3)
%\drawbb
\end{texdraw}}
\noindent
right :
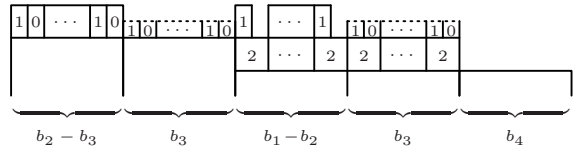
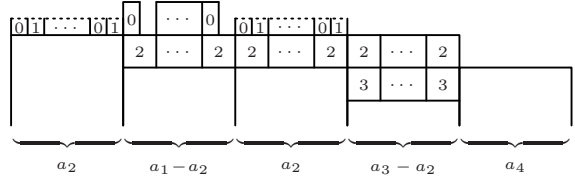
\begin{center}
\usebox{\tmpfige}
\end{center}
\noindent
left :
\begin{center}
\usebox{\tmpfigb}
\end{center}
Finally, join the two slices in such a way that the highest layer
of the result forms a part of a level-$l$ reduced proper Young wall that has
had all its blocks that may be split, split.

Now, to bring this into a reduced proper form, we need to
\emph{remove} $\delta$'s from the \emph{left} slice.
We denote the number of $\delta$ removals needed by $k$.

Below, we list left-$\vphi$ and right-$\veps$ values
for only four of the sixteen possible cases of Young wall column pairs.
The remaining cases are less complicated.
The line containing the bullet lists the two column
types in left-right order.
\begin{itemize}
\item L$(0)_1$$(\overline{2})_{01}$ R$(1)_0$$(10)_{\overline{2}}$ or
      L$(0)_1$$({\overline{2}})_{01}$ R$(0)_1$$(10)_{\overline{2}}$\\[2mm]
  \begin{enumerate}
  \item $a_3-a_2\le b_1-a_1$\\[2mm]
  \begin{tabular}{rcl}
  rotation & : & $0\le k\le a_2$\\
  left-$\vphi$ & : & $k+a_3-a_2$\\
  right-$\veps$ & : &$k+b_2-a_2$\\[2mm]
  rotation & : & $a_2\leq k$\\
  left-$\vphi$ & : & $a_3$\\
  right-$\veps$ & : &$b_2$\\
  \end{tabular}\\[2mm]
  \item $b_1-a_1\le a_3-a_2$\\[2mm]
  \begin{tabular}{rcl}
  rotation & : & $0\le k\le a_1+a_3-b_1$\\
  left-$\vphi$ & : & $k+b_1-a_1$\\
  right-$\veps$ & : &$k+b_2+b_1-a_1-a_3$\\[2mm]
  rotation & : & $a_1+a_3-b_1\leq k$\\
  left-$\vphi$ & : & $a_3$\\
  right-$\veps$ & : &$b_2$\\
  \end{tabular}\\[2mm]
  \end{enumerate}
%\item L$(0)_1$$(2)_{01}$ R$(0)_1$$(10)_2$\\[2mm]
%  \begin{enumerate}
%  \item $a_3-a_2\le b_1-a_1$\\[2mm]
%  \begin{tabular}{rcl}
%  rotation & : & $0\le k\le a_2$\\
%  left-$\vphi$ & : & $a_3-a_2+k$\\
%  right-$\veps$ & : &$b_2-a_2+k$\\[2mm]
%  rotation & : & $a_2\leq k$\\
%  left-$\vphi$ & : & $a_3$\\
%  right-$\veps$ & : &$b_2$\\
%  \end{tabular}\\[2mm]
%  \item $b_1-a_1\le a_3-a_2$\\[2mm]
%  \begin{tabular}{rcl}
%  rotation & : & $0\le k\le a_1+a_3-b_1$\\
%  left-$\vphi$ & : & $b_1-a_1+k$\\
%  right-$\veps$ & : &$b_2-a_1+b_1-a_3+k$\\[2mm]
%  rotation & : & $a_1+a_3-b_1\leq k$\\
%  left-$\vphi$ & : & $a_3$\\
%  right-$\veps$ & : &$b_2$\\
%  \end{tabular}\\[2mm]
%  \end{enumerate}
\item L$(1)_0$$({\overline{2}})_{01}$ R$(1)_0$$(10)_{\overline{2}}$ or
      L$(1)_0$$({\overline{2}})_{01}$ R$(0)_1$$(10)_{\overline{2}}$\\[2mm]
  \begin{enumerate}
  \item $b_1\le a_3$\\[2mm]
  \begin{tabular}{rcl}
  rotation & : & $0\le k\le a_2+a_3-b_1$\\
  left-$\vphi$ & : & $k+b_1-a_1$\\
  right-$\veps$ & : &$k+b_2+b_1-a_2-a_3$\\[2mm]
  rotation & : & $a_2+a_3-b_1\leq k$\\
  left-$\vphi$ & : & $a_2+a_3-a_1$\\
  right-$\veps$ & : &$b_2$\\
  \end{tabular}\\[2mm]
  \item $a_3\le b_1$\\[2mm]
  \begin{tabular}{rcl}
  rotation & : & $0\le k\le a_2$\\
  left-$\vphi$ & : & $k+a_3-a_1$\\
  right-$\veps$ & : &$k+b_2-a_2$\\[2mm]
  rotation & : & $a_2\leq k$\\
  left-$\vphi$ & : & $a_2+a_3-a_1$\\
  right-$\veps$ & : &$b_2$\\
  \end{tabular}\\[2mm]
  \end{enumerate}
%\item L$(1)_0$$(2)_{01}$ R$(0)_1$$(10)_2$\\[2mm]
%  \begin{enumerate}
%  \item $b_1\le a_3$\\[2mm]
%  \begin{tabular}{rcl}
%  rotation & : & $0\le k\le a_2+a_3-b_1$\\
%  left-$\vphi$ & : & $b_1-a_1+k$\\
%  right-$\veps$ & : &$b_2-a_2+b_1-a_3+k$\\[2mm]
%  rotation & : & $a_2+a_3-b_1\leq k$\\
%  left-$\vphi$ & : & $a_2+a_3-a_1$\\
%  right-$\veps$ & : &$b_2$\\
%  \end{tabular}\\[2mm]
%  \item $a_3\le b_1$\\[2mm]
% \begin{tabular}{rcl}
%  rotation & : & $0\le k\le a_2$\\
%  left-$\vphi$ & : & $a_3-a_1+k$\\
%  right-$\veps$ & : &$b_2-a_2+k$\\[2mm]
%  rotation & : & $a_2\leq k$\\
%  left-$\vphi$ & : & $a_2+a_3-a_1$\\
%  right-$\veps$ & : &$b_2$\\
%  \end{tabular}\\[2mm]
%  \end{enumerate}
\end{itemize}
Similarly, the following gives the signatures of the
path description.
The number in the list are the
left-$\vphi$ and right-$\veps$ values for the two
corresponding crystal elements.
\begin{itemize}
\item Any pair with L$(0)_1$ in the left column\\[2mm]
\begin{tabular}{rcl}
left-$\vphi$ & : & $a_3$\\
right-$\veps$ & : & $b_2$\\
\end{tabular}\\[2mm]
\item Any pair with L$(1)_0$ in the left column\\[2mm]
\begin{tabular}{rcl}
left-$\vphi$ & : & $a_3+(a_2-a_1)$\\
right-$\veps$ & : & $b_2$\\
\end{tabular}\\[2mm]
\end{itemize}
We can easily see that the signatures agree with those of the
corresponding path description in all of the cases after
$(0,1)$-pair cancellations.

Results for $i=1$ can be obtained
if we apply the following substitutions
to the parts giving $\vphi$ and $\veps$ values,
appearing just above.
\begin{itemize}
\item L$(0)_1$$(01)_{{\overline{2}}}$
      $\leftrightarrow$ L$(1)_0$$(01)_{{\overline{2}}}$
\item L$(0)_1$$({\overline{2}})_{01}$
      $\leftrightarrow$ L$(1)_0$$({\overline{2}})_{01}$
\item R$(1)_0$$(10)_{{\overline{2}}}$
      $\leftrightarrow$ R$(0)_1$$(10)_{{\overline{2}}}$
\item R$(1)_0$$({\overline{2}})_{10}$
      $\leftrightarrow$ R$(0)_1$$({\overline{2}})_{10}$
\item $a_1$ $\leftrightarrow$ $a_2$
\item $b_1$ $\leftrightarrow$ $b_2$
\end{itemize}

%%%%%%%%%%%%%%%%%%%%%%%%%%%%%%%%%%%%%%%%%%%%%%%%%%%
\section{case $2 \le i \le n-1$}

Finally, we are dealing with the $2 \le i \le n-1$ case.
Let us write the general level-$l$ Young wall column in the
following form.\\
\savebox{\tmpfiga}{\begin{texdraw}
\fontsize{6}{6}\selectfont
\textref h:C v:C
\drawdim em
\setunitscale 1.7
\move(0 0)
\bsegment
\move(0 -0.5)\lvec(0 3)\lvec(3.4 3)\lvec(3.4 -0.5)
\move(0 2)\lvec(3.4 2)
\move(0 1)\lvec(3.4 1)
\move(1 3)\lvec(1 1)
\move(2.4 3)\lvec(2.4 1)
\htext(0.5 2.5){$i$}\htext(2.9 2.5){$i$}
\htext(1.7 2.5){$\cdots$}
\htext(0.5 1.5){$i\!\!-\!\!1$}\htext(2.9 1.5){$i\!\!-\!\!1$}
\htext(1.7 1.5){$\cdots$}
\esegment
\move(3.4 -1)
\bsegment
\move(0 0.5)\lvec(0 3)\lvec(3.4 3)\lvec(3.4 0.5)
\move(0 2)\lvec(3.4 2)
\move(0 1)\lvec(3.4 1)
\move(1 3)\lvec(1 1)
\move(2.4 3)\lvec(2.4 1)
\htext(0.5 2.5){$i\!\!-\!\!1$}\htext(2.9 2.5){$i\!\!-\!\!1$}
\htext(1.7 2.5){$\cdots$}
\htext(0.5 1.5){$i\!\!-\!\!2$}\htext(2.9 1.5){$i\!\!-\!\!2$}
\htext(1.7 1.5){$\cdots$}
\esegment
\move(6.8 0)
\bsegment
\move(0 -0.5)\lvec(0 1)\lvec(3.4 1)\lvec(3.4 -0.5)
\esegment
%%%%%
\move(0 -1)\lvec(0 -5.25)
\move(3.4 -1)\lvec(3.4 -5.25)
\move(6.8 -1)\lvec(6.8 -5.25)
\move(10.2 -1)\lvec(10.2 -1.5)
%%%%%
\move(10.2 -4.5)
\bsegment
\move(0 0)
\bsegment
\move(0 -0.75)\lvec(0 3)\lvec(3.4 3)\lvec(3.4 -0.75)
\move(0 2)\lvec(3.4 2)
\move(0 1)\lvec(3.4 1)
\move(1 3)\lvec(1 1)
\move(2.4 3)\lvec(2.4 1)
\htext(0.5 2.5){$i$}\htext(2.9 2.5){$i$}
\htext(1.7 2.5){$\cdots$}
\htext(0.5 1.5){$i\!\!+\!\!1$}\htext(2.9 1.5){$i\!\!+\!\!1$}
\htext(1.7 1.5){$\cdots$}
\esegment
\move(3.4 -1)
\bsegment
\move(0 0.75)\lvec(0 3)\lvec(3.4 3)\lvec(3.4 0.75)
\move(0 2)\lvec(3.4 2)
\move(0 1)\lvec(3.4 1)
\move(1 3)\lvec(1 1)
\move(2.4 3)\lvec(2.4 1)
\htext(0.5 2.5){$i\!\!+\!\!1$}\htext(2.9 2.5){$i\!\!+\!\!1$}
\htext(1.7 2.5){$\cdots$}
\htext(0.5 1.5){$i\!\!+\!\!2$}\htext(2.9 1.5){$i\!\!+\!\!2$}
\htext(1.7 1.5){$\cdots$}
\esegment
\move(6.8 0)
\bsegment
\move(0 -0.75)\lvec(0 1)\lvec(3.4 1)\lvec(3.4 -0.75)
\esegment
\esegment
\htext(1.7 -5.75){$\underbrace{\rule{5.6em}{0em}}$}
\htext(5.1 -5.75){$\underbrace{\rule{5.6em}{0em}}$}
\htext(8.5 -5.75){$\underbrace{\rule{5.6em}{0em}}$}
\htext(11.9 -5.75){$\underbrace{\rule{5.6em}{0em}}$}
\htext(15.3 -5.75){$\underbrace{\rule{5.6em}{0em}}$}
\htext(18.7 -5.75){$\underbrace{\rule{5.6em}{0em}}$}
\htext(1.7 -6.5){$a_1$}\htext(5.1 -6.5){$a_2$}
\htext(8.5 -6.5){$a_3$}\htext(11.9 -6.5){$a_4$}
\htext(15.3 -6.5){$a_5$}\htext(18.7 -6.5){$a_6$}
\move(0 -6.8)
%\drawbb
\end{texdraw}}
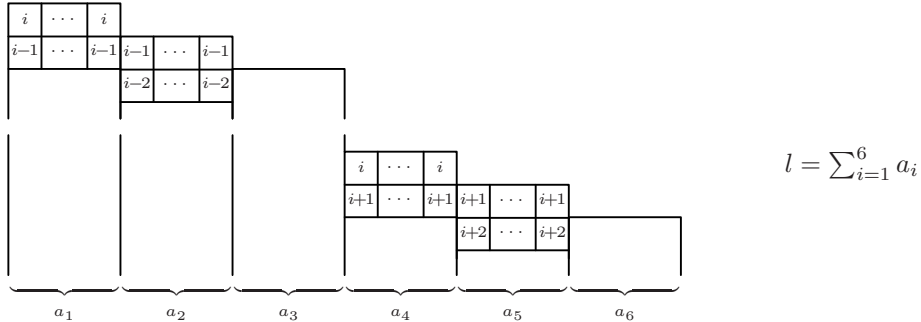
\begin{center}
\begin{minipage}{0.8\textwidth}
\usebox{\tmpfiga}
\end{minipage}
\begin{minipage}{0.15\textwidth}
\setlength{\tabcolsep}{0.8mm}
\begin{tabular}{rcl}
$l=\sum_{i=1}^6 a_i$
\end{tabular}
\end{minipage}
\end{center}
Here, $a_1, a_2, a_4, a_5$ denote the number of layers having top
blocks of the form given in the figure. Also $a_3, a_6$ denote the
number of all other layers, that is, any layer having a top block
that comes between the supporting $(i-2)$-block and the covering
$(i-1)$-block (inclusive) is counted in $a_3$ and any layer having
a top block that comes between the covering $(i+2)$-block and the
supporting $(i+1)$-block (inclusive) is counted in $a_6$.

We break this into four cases and fix notations for the cases.
\begin{itemize}
\item $a_1\geq a_5$ and $a_2\geq a_4$
      : $(\underline{i})_{\overline{i+1}}$$(\underline{i-1})_{\overline{i}}$
\item $a_1\geq a_5$ and $a_2\leq a_4$
      : $(\underline{i})_{\overline{i+1}}$$(\overline{i})_{\underline{i-1}}$
\item $a_1\leq a_5$ and $a_2\geq a_4$
      : $(\overline{i+1})_{\underline{i}}$$(\underline{i-1})_{\overline{i}}$
\item $a_1\leq a_5$ and $a_2\leq a_4$
      : $(\overline{i+1})_{\underline{i}}$$(\overline{i})_{\underline{i-1}}$
\end{itemize}
%The same notation will be used to denote any other
%layer-rotation of these columns.
%Also the perfect crystal elements corresponding to these
%columns will be denoted with the same notations.

If we split every block possible from these slices, the
result will be of the following four shapes.\\[2mm]
\savebox{\tmpfiga}{\begin{texdraw}
\fontsize{6}{6}\selectfont
\textref h:C v:C
\drawdim em
\setunitscale 1.7
\move(0 3)
\bsegment
\move(0 -0.5)\lvec(0 3)\lvec(3.4 3)\lvec(3.4 -0.5)
\move(0 2)\lvec(3.4 2)
\move(0 1)\lvec(3.4 1)
\move(1 3)\lvec(1 1)
\move(2.4 3)\lvec(2.4 1)
\htext(1.7 2.5){$\cdots$}
\htext(0.5 2.5){$i$}\htext(2.9 2.5){$i$}
\htext(1.7 1.5){$\cdots$}
\htext(0.5 1.5){$i\!\!-\!\!1$}\htext(2.9 1.5){$i\!\!-\!\!1$}
\move(0 -1)\rlvec(0 -4.75)
\esegment
\move(3.4 2.5)
\bsegment
\move(0 0)\lvec(0 3)
\move(3.4 3)\lvec(3.4 0)
\move(0 2.5)\lvec(3.4 2.5)
\move(0 1.5)\lvec(3.4 1.5)
\move(1 3)\lvec(1 1.5)
\move(2.4 3)\lvec(2.4 1.5)
\htext(1.7 2.75){$\cdots$}
\htext(0.5 2.75){$i$}\htext(2.9 2.75){$i$}
\htext(1.7 2){$\cdots$}
\htext(0.5 2){$i\!\!-\!\!1$}\htext(2.9 2){$i\!\!-\!\!1$}
\move(0 -0.5)\rlvec(0 -4.75)
\lpatt(0.05 0.15)
\move(0 3)\lvec(3.4 3)
\esegment
\move(6.8 2)
\bsegment
\move(0 0.5)\lvec(0 3)\lvec(3.4 3)\lvec(3.4 0.5)
\move(0 2)\lvec(3.4 2)
\move(0 1)\lvec(3.4 1)
\move(1 3)\lvec(1 1)
\move(2.4 3)\lvec(2.4 1)
\htext(1.7 2.5){$\cdots$}
\htext(0.5 2.5){$i\!\!-\!\!1$}\htext(2.9 2.5){$i\!\!-\!\!1$}
\htext(1.7 1.5){$\cdots$}
\htext(0.5 1.5){$i\!\!-\!\!2$}\htext(2.9 1.5){$i\!\!-\!\!2$}
\move(0 0)\rlvec(0 -4.75)
\esegment
\move(10.2 1.5)
\bsegment
\move(0 1)\lvec(0 3)
\move(3.4 3)\lvec(3.4 1)
\move(0 2.5)\lvec(3.4 2.5)
\move(0 1.5)\lvec(3.4 1.5)
\move(1 3)\lvec(1 1.5)
\move(2.4 3)\lvec(2.4 1.5)
\htext(1.7 2.75){$\cdots$}
\htext(0.5 2.75){$i\!\!-\!\!1$}\htext(2.9 2.75){$i\!\!-\!\!1$}
\htext(1.7 2){$\cdots$}
\htext(0.5 2){$i\!\!-\!\!2$}\htext(2.9 2){$i\!\!-\!\!2$}
\move(0 0.5)\rlvec(0 -4.75)
\lpatt(0.05 0.15)
\move(0 3)\lvec(3.4 3)
\esegment
\move(13.6 1)
\bsegment
\move(0 1.5)\lvec(0 3)\lvec(3.4 3)\lvec(3.4 1.5)
\move(0 1)\rlvec(0 -4.75)
\esegment
\move(17 -1.5)
\bsegment
\move(0 -1.25)\lvec(0 3)
\move(3.4 3)\lvec(3.4 -1.25)
\move(0 2.5)\lvec(3.4 2.5)
\move(0 1.5)\lvec(3.4 1.5)
\move(1 3)\lvec(1 1.5)
\move(2.4 3)\lvec(2.4 1.5)
\htext(1.7 2.75){$\cdots$}
\htext(0.5 2.75){$i\!\!-\!\!1$}\htext(2.9 2.75){$i\!\!-\!\!1$}
\htext(1.7 2){$\cdots$}
\htext(0.5 2){$i$}\htext(2.9 2){$i$}
\move(0 3)\rlvec(0 0.5)
\lpatt(0.05 0.15)
\move(0 3)\lvec(3.4 3)
\esegment
\move(20.4 -2.5)
\bsegment
\move(0 -0.25)\lvec(0 3)
\move(3.4 3)\lvec(3.4 -0.25)
\move(0 2.5)\lvec(3.4 2.5)
\move(0 1.5)\lvec(3.4 1.5)
\move(1 3)\lvec(1 1.5)
\move(2.4 3)\lvec(2.4 1.5)
\htext(1.7 2.75){$\cdots$}
\htext(0.5 2.75){$i$}\htext(2.9 2.75){$i$}
\htext(1.7 2){$\cdots$}
\htext(0.5 2){$i\!\!+\!\!1$}\htext(2.9 2){$i\!\!+\!\!1$}
\lpatt(0.05 0.15)
\move(0 3)\lvec(3.4 3)
\esegment
\move(23.8 -4)
\bsegment
\move(0 1.25)\lvec(0 3)\lvec(3.4 3)\lvec(3.4 1.25)
\esegment
\htext(1.7 -3.25){$\underbrace{\rule{5.6em}{0em}}$}
\htext(5.1 -3.25){$\underbrace{\rule{5.6em}{0em}}$}
\htext(8.5 -3.25){$\underbrace{\rule{5.6em}{0em}}$}
\htext(11.9 -3.25){$\underbrace{\rule{5.6em}{0em}}$}
\htext(15.3 -3.25){$\underbrace{\rule{5.6em}{0em}}$}
\htext(18.7 -3.25){$\underbrace{\rule{5.6em}{0em}}$}
\htext(22.1 -3.25){$\underbrace{\rule{5.6em}{0em}}$}
\htext(25.5 -3.25){$\underbrace{\rule{5.6em}{0em}}$}
\htext(1.7 -4){$a_1-a_5$}\htext(5.1 -4){$a_5$}
\htext(8.5 -4){$a_2-a_4$}\htext(11.9 -4){$a_4$}
\htext(15.3 -4){$a_3$}\htext(18.7 -4){$a_4$}
\htext(22.1 -4){$a_5$}\htext(25.5 -4){$a_6$}
\move(0 -4.3)
%\drawbb
\end{texdraw}}
\savebox{\tmpfigb}{\begin{texdraw}
\fontsize{6}{6}\selectfont
\textref h:C v:C
\drawdim em
\setunitscale 1.7
\move(0 3)
\bsegment
\move(0 -0.5)\lvec(0 3)\lvec(3.4 3)\lvec(3.4 -0.5)
\move(0 2)\lvec(3.4 2)
\move(0 1)\lvec(3.4 1)
\move(1 3)\lvec(1 1)
\move(2.4 3)\lvec(2.4 1)
\htext(1.7 2.5){$\cdots$}
\htext(0.5 2.5){$i$}\htext(2.9 2.5){$i$}
\htext(1.7 1.5){$\cdots$}
\htext(0.5 1.5){$i\!\!-\!\!1$}\htext(2.9 1.5){$i\!\!-\!\!1$}
\move(0 -1)\rlvec(0 -4.75)
\esegment
\move(3.4 2.5)
\bsegment
\move(0 0)\lvec(0 3)
\move(3.4 3)\lvec(3.4 0)
\move(0 2.5)\lvec(3.4 2.5)
\move(0 1.5)\lvec(3.4 1.5)
\move(1 3)\lvec(1 1.5)
\move(2.4 3)\lvec(2.4 1.5)
\htext(1.7 2.75){$\cdots$}
\htext(0.5 2.75){$i$}\htext(2.9 2.75){$i$}
\htext(1.7 2){$\cdots$}
\htext(0.5 2){$i\!\!-\!\!1$}\htext(2.9 2){$i\!\!-\!\!1$}
\move(0 -0.5)\rlvec(0 -4.75)
\lpatt(0.05 0.15)
\move(0 3)\lvec(3.4 3)
\esegment
\move(6.8 1.5)
\bsegment
\move(0 1)\lvec(0 3)
\move(3.4 3)\lvec(3.4 1)
\move(0 2.5)\lvec(3.4 2.5)
\move(0 1.5)\lvec(3.4 1.5)
\move(1 3)\lvec(1 1.5)
\move(2.4 3)\lvec(2.4 1.5)
\htext(1.7 2.75){$\cdots$}
\htext(0.5 2.75){$i\!\!-\!\!1$}\htext(2.9 2.75){$i\!\!-\!\!1$}
\htext(1.7 2){$\cdots$}
\htext(0.5 2){$i\!\!-\!\!2$}\htext(2.9 2){$i\!\!-\!\!2$}
\move(0 0.5)\rlvec(0 -4.75)
\lpatt(0.05 0.15)
\move(0 3)\lvec(3.4 3)
\esegment
\move(10.2 1)
\bsegment
\move(0 1.5)\lvec(0 3)\lvec(3.4 3)\lvec(3.4 1.5)
\move(0 1)\rlvec(0 -4.75)
\esegment
\move(13.6 -1.5)
\bsegment
\move(0 -1.25)\lvec(0 3)
\move(3.4 3)\lvec(3.4 -1.25)
\move(0 2.5)\lvec(3.4 2.5)
\move(0 1.5)\lvec(3.4 1.5)
\move(1 3)\lvec(1 1.5)
\move(2.4 3)\lvec(2.4 1.5)
\htext(1.7 2.75){$\cdots$}
\htext(0.5 2.75){$i\!\!-\!\!1$}\htext(2.9 2.75){$i\!\!-\!\!1$}
\htext(1.7 2){$\cdots$}
\htext(0.5 2){$i$}\htext(2.9 2){$i$}
\move(0 3)\rlvec(0 0.5)
\lpatt(0.05 0.15)
\move(0 3)\lvec(3.4 3)
\esegment
\move(17 -2)
\bsegment
\move(0 -0.5)\lvec(0 3)\lvec(3.4 3)\lvec(3.4 -0.5)
\move(0 2)\lvec(3.4 2)
\move(0 1)\lvec(3.4 1)
\move(1 3)\lvec(1 1)
\move(2.4 3)\lvec(2.4 1)
\htext(1.7 2.5){$\cdots$}
\htext(0.5 2.5){$i$}\htext(2.9 2.5){$i$}
\htext(1.7 1.5){$\cdots$}
\htext(0.5 1.5){$i\!\!+\!\!1$}\htext(2.9 1.5){$i\!\!+\!\!1$}
\esegment
\move(20.4 -2.5)
\bsegment
\move(0 -0.25)\lvec(0 3)
\move(3.4 3)\lvec(3.4 -0.25)
\move(0 2.5)\lvec(3.4 2.5)
\move(0 1.5)\lvec(3.4 1.5)
\move(1 3)\lvec(1 1.5)
\move(2.4 3)\lvec(2.4 1.5)
\htext(1.7 2.75){$\cdots$}
\htext(0.5 2.75){$i$}\htext(2.9 2.75){$i$}
\htext(1.7 2){$\cdots$}
\htext(0.5 2){$i\!\!+\!\!1$}\htext(2.9 2){$i\!\!+\!\!1$}
\lpatt(0.05 0.15)
\move(0 3)\lvec(3.4 3)
\esegment
\move(23.8 -4)
\bsegment
\move(0 1.25)\lvec(0 3)\lvec(3.4 3)\lvec(3.4 1.25)
\esegment
\htext(1.7 -3.25){$\underbrace{\rule{5.6em}{0em}}$}
\htext(5.1 -3.25){$\underbrace{\rule{5.6em}{0em}}$}
\htext(8.5 -3.25){$\underbrace{\rule{5.6em}{0em}}$}
\htext(11.9 -3.25){$\underbrace{\rule{5.6em}{0em}}$}
\htext(15.3 -3.25){$\underbrace{\rule{5.6em}{0em}}$}
\htext(18.7 -3.25){$\underbrace{\rule{5.6em}{0em}}$}
\htext(22.1 -3.25){$\underbrace{\rule{5.6em}{0em}}$}
\htext(25.5 -3.25){$\underbrace{\rule{5.6em}{0em}}$}
\htext(1.7 -4){$a_1-a_5$}\htext(5.1 -4){$a_5$}
\htext(8.5 -4){$a_2$}\htext(11.9 -4){$a_3$}
\htext(15.3 -4){$a_2$}\htext(18.7 -4){$a_4-a_2$}
\htext(22.1 -4){$a_5$}\htext(25.5 -4){$a_6$}
\move(0 -4.3)
%\drawbb
\end{texdraw}}
\savebox{\tmpfigc}{\begin{texdraw}
\fontsize{6}{6}\selectfont
\textref h:C v:C
\drawdim em
\setunitscale 1.7
\move(0 2.5)
\bsegment
\move(0 0)\lvec(0 3)
\move(3.4 3)\lvec(3.4 0)
\move(0 2.5)\lvec(3.4 2.5)
\move(0 1.5)\lvec(3.4 1.5)
\move(1 3)\lvec(1 1.5)
\move(2.4 3)\lvec(2.4 1.5)
\htext(1.7 2.75){$\cdots$}
\htext(0.5 2.75){$i$}\htext(2.9 2.75){$i$}
\htext(1.7 2){$\cdots$}
\htext(0.5 2){$i\!\!-\!\!1$}\htext(2.9 2){$i\!\!-\!\!1$}
\move(0 -0.5)\rlvec(0 -4.75)
\lpatt(0.05 0.15)
\move(0 3)\lvec(3.4 3)
\esegment
\move(3.4 2)
\bsegment
\move(0 0.5)\lvec(0 3)\lvec(3.4 3)\lvec(3.4 0.5)
\move(0 2)\lvec(3.4 2)
\move(0 1)\lvec(3.4 1)
\move(1 3)\lvec(1 1)
\move(2.4 3)\lvec(2.4 1)
\htext(1.7 2.5){$\cdots$}
\htext(0.5 2.5){$i\!\!-\!\!1$}\htext(2.9 2.5){$i\!\!-\!\!1$}
\htext(1.7 1.5){$\cdots$}
\htext(0.5 1.5){$i\!\!-\!\!2$}\htext(2.9 1.5){$i\!\!-\!\!2$}
\move(0 0)\rlvec(0 -4.75)
\esegment
\move(6.8 1.5)
\bsegment
\move(0 1)\lvec(0 3)
\move(3.4 3)\lvec(3.4 1)
\move(0 2.5)\lvec(3.4 2.5)
\move(0 1.5)\lvec(3.4 1.5)
\move(1 3)\lvec(1 1.5)
\move(2.4 3)\lvec(2.4 1.5)
\htext(1.7 2.75){$\cdots$}
\htext(0.5 2.75){$i\!\!-\!\!1$}\htext(2.9 2.75){$i\!\!-\!\!1$}
\htext(1.7 2){$\cdots$}
\htext(0.5 2){$i\!\!-\!\!2$}\htext(2.9 2){$i\!\!-\!\!2$}
\move(0 0.5)\rlvec(0 -4.75)
\lpatt(0.05 0.15)
\move(0 3)\lvec(3.4 3)
\esegment
\move(10.2 1)
\bsegment
\move(0 1.5)\lvec(0 3)\lvec(3.4 3)\lvec(3.4 1.5)
\move(0 1)\rlvec(0 -4.75)
\esegment
\move(13.6 -1.5)
\bsegment
\move(0 -1.25)\lvec(0 3)
\move(3.4 3)\lvec(3.4 -1.25)
\move(0 2.5)\lvec(3.4 2.5)
\move(0 1.5)\lvec(3.4 1.5)
\move(1 3)\lvec(1 1.5)
\move(2.4 3)\lvec(2.4 1.5)
\htext(1.7 2.75){$\cdots$}
\htext(0.5 2.75){$i\!\!-\!\!1$}\htext(2.9 2.75){$i\!\!-\!\!1$}
\htext(1.7 2){$\cdots$}
\htext(0.5 2){$i$}\htext(2.9 2){$i$}
\move(0 3)\rlvec(0 0.5)
\lpatt(0.05 0.15)
\move(0 3)\lvec(3.4 3)
\esegment
\move(17 -2.5)
\bsegment
\move(0 -0.25)\lvec(0 3)
\move(3.4 3)\lvec(3.4 -0.25)
\move(0 2.5)\lvec(3.4 2.5)
\move(0 1.5)\lvec(3.4 1.5)
\move(1 3)\lvec(1 1.5)
\move(2.4 3)\lvec(2.4 1.5)
\htext(1.7 2.75){$\cdots$}
\htext(0.5 2.75){$i$}\htext(2.9 2.75){$i$}
\htext(1.7 2){$\cdots$}
\htext(0.5 2){$i\!\!+\!\!1$}\htext(2.9 2){$i\!\!+\!\!1$}
\lpatt(0.05 0.15)
\move(0 3)\lvec(3.4 3)
\esegment
\move(20.4 -3)
\bsegment
\move(0 0.25)\lvec(0 3)\lvec(3.4 3)\lvec(3.4 0.25)
\move(0 2)\lvec(3.4 2)
\move(0 1)\lvec(3.4 1)
\move(1 3)\lvec(1 1)
\move(2.4 3)\lvec(2.4 1)
\htext(1.7 2.5){$\cdots$}
\htext(0.5 2.5){$i\!\!+\!\!1$}\htext(2.9 2.5){$i\!\!+\!\!1$}
\htext(1.7 1.5){$\cdots$}
\htext(0.5 1.5){$i\!\!+\!\!2$}\htext(2.9 1.5){$i\!\!+\!\!2$}
\esegment
\move(23.8 -4)
\bsegment
\move(0 1.25)\lvec(0 3)\lvec(3.4 3)\lvec(3.4 1.25)
\esegment
\htext(1.7 -3.25){$\underbrace{\rule{5.6em}{0em}}$}
\htext(5.1 -3.25){$\underbrace{\rule{5.6em}{0em}}$}
\htext(8.5 -3.25){$\underbrace{\rule{5.6em}{0em}}$}
\htext(11.9 -3.25){$\underbrace{\rule{5.6em}{0em}}$}
\htext(15.3 -3.25){$\underbrace{\rule{5.6em}{0em}}$}
\htext(18.7 -3.25){$\underbrace{\rule{5.6em}{0em}}$}
\htext(22.1 -3.25){$\underbrace{\rule{5.6em}{0em}}$}
\htext(25.5 -3.25){$\underbrace{\rule{5.6em}{0em}}$}
\htext(1.7 -4){$a_1$}\htext(5.1 -4){$a_2-a_4$}
\htext(8.5 -4){$a_4$}\htext(11.9 -4){$a_3$}
\htext(15.3 -4){$a_4$}\htext(18.7 -4){$a_1$}
\htext(22.1 -4){$a_5-a_1$}\htext(25.5 -4){$a_6$}
\move(0 -4.3)
%\drawbb
\end{texdraw}}
\savebox{\tmpfigd}{\begin{texdraw}
\fontsize{6}{6}\selectfont
\textref h:C v:C
\drawdim em
\setunitscale 1.7
\move(0 2.5)
\bsegment
\move(0 0)\lvec(0 3)
\move(3.4 3)\lvec(3.4 0)
\move(0 2.5)\lvec(3.4 2.5)
\move(0 1.5)\lvec(3.4 1.5)
\move(1 3)\lvec(1 1.5)
\move(2.4 3)\lvec(2.4 1.5)
\htext(1.7 2.75){$\cdots$}
\htext(0.5 2.75){$i$}\htext(2.9 2.75){$i$}
\htext(1.7 2){$\cdots$}
\htext(0.5 2){$i\!\!-\!\!1$}\htext(2.9 2){$i\!\!-\!\!1$}
\move(0 -0.5)\rlvec(0 -4.75)
\lpatt(0.05 0.15)
\move(0 3)\lvec(3.4 3)
\esegment
\move(3.4 1.5)
\bsegment
\move(0 1)\lvec(0 3)
\move(3.4 3)\lvec(3.4 1)
\move(0 2.5)\lvec(3.4 2.5)
\move(0 1.5)\lvec(3.4 1.5)
\move(1 3)\lvec(1 1.5)
\move(2.4 3)\lvec(2.4 1.5)
\htext(1.7 2.75){$\cdots$}
\htext(0.5 2.75){$i\!\!-\!\!1$}\htext(2.9 2.75){$i\!\!-\!\!1$}
\htext(1.7 2){$\cdots$}
\htext(0.5 2){$i\!\!-\!\!2$}\htext(2.9 2){$i\!\!-\!\!2$}
\move(0 0.5)\rlvec(0 -4.75)
\lpatt(0.05 0.15)
\move(0 3)\lvec(3.4 3)
\esegment
\move(6.8 1)
\bsegment
\move(0 1.5)\lvec(0 3)\lvec(3.4 3)\lvec(3.4 1.5)
\move(0 1)\rlvec(0 -4.75)
\esegment
\move(10.2 -1.5)
\bsegment
\move(0 -1.25)\lvec(0 3)
\move(3.4 3)\lvec(3.4 -1.25)
\move(0 2.5)\lvec(3.4 2.5)
\move(0 1.5)\lvec(3.4 1.5)
\move(1 3)\lvec(1 1.5)
\move(2.4 3)\lvec(2.4 1.5)
\htext(1.7 2.75){$\cdots$}
\htext(0.5 2.75){$i\!\!-\!\!1$}\htext(2.9 2.75){$i\!\!-\!\!1$}
\htext(1.7 2){$\cdots$}
\htext(0.5 2){$i$}\htext(2.9 2){$i$}
\move(0 3)\rlvec(0 0.5)
\lpatt(0.05 0.15)
\move(0 3)\lvec(3.4 3)
\esegment
\move(13.6 -2)
\bsegment
\move(0 0.5)\lvec(0 3)\lvec(3.4 3)\lvec(3.4 0.5)
\move(0 2)\lvec(3.4 2)
\move(0 1)\lvec(3.4 1)
\move(1 3)\lvec(1 1)
\move(2.4 3)\lvec(2.4 1)
\htext(1.7 2.5){$\cdots$}
\htext(0.5 2.5){$i$}\htext(2.9 2.5){$i$}
\htext(1.7 1.5){$\cdots$}
\htext(0.5 1.5){$i\!\!+\!\!1$}\htext(2.9 1.5){$i\!\!+\!\!1$}
\esegment
\move(17 -2.5)
\bsegment
\move(0 -0.25)\lvec(0 3)
\move(3.4 3)\lvec(3.4 -0.25)
\move(0 2.5)\lvec(3.4 2.5)
\move(0 1.5)\lvec(3.4 1.5)
\move(1 3)\lvec(1 1.5)
\move(2.4 3)\lvec(2.4 1.5)
\htext(1.7 2.75){$\cdots$}
\htext(0.5 2.75){$i$}\htext(2.9 2.75){$i$}
\htext(1.7 2){$\cdots$}
\htext(0.5 2){$i\!\!+\!\!1$}\htext(2.9 2){$i\!\!+\!\!1$}
\lpatt(0.05 0.15)
\move(0 3)\lvec(3.4 3)
\esegment
\move(20.4 -3)
\bsegment
\move(0 0.25)\lvec(0 3)\lvec(3.4 3)\lvec(3.4 0.25)
\move(0 2)\lvec(3.4 2)
\move(0 1)\lvec(3.4 1)
\move(1 3)\lvec(1 1)
\move(2.4 3)\lvec(2.4 1)
\htext(1.7 2.5){$\cdots$}
\htext(0.5 2.5){$i\!\!+\!\!1$}\htext(2.9 2.5){$i\!\!+\!\!1$}
\htext(1.7 1.5){$\cdots$}
\htext(0.5 1.5){$i\!\!+\!\!2$}\htext(2.9 1.5){$i\!\!+\!\!2$}
\esegment
\move(23.8 -4)
\bsegment
\move(0 1.25)\lvec(0 3)\lvec(3.4 3)\lvec(3.4 1.25)
\esegment
\htext(1.7 -3.25){$\underbrace{\rule{5.6em}{0em}}$}
\htext(5.1 -3.25){$\underbrace{\rule{5.6em}{0em}}$}
\htext(8.5 -3.25){$\underbrace{\rule{5.6em}{0em}}$}
\htext(11.9 -3.25){$\underbrace{\rule{5.6em}{0em}}$}
\htext(15.3 -3.25){$\underbrace{\rule{5.6em}{0em}}$}
\htext(18.7 -3.25){$\underbrace{\rule{5.6em}{0em}}$}
\htext(22.1 -3.25){$\underbrace{\rule{5.6em}{0em}}$}
\htext(25.5 -3.25){$\underbrace{\rule{5.6em}{0em}}$}
\htext(1.7 -4){$a_1$}\htext(5.1 -4){$a_2$}
\htext(8.5 -4){$a_3$}\htext(11.9 -4){$a_2$}
\htext(15.3 -4){$a_4-a_2$}\htext(18.7 -4){$a_1$}
\htext(22.1 -4){$a_5-a_1$}\htext(25.5 -4){$a_6$}
\move(0 -4.3)
%\drawbb
\end{texdraw}}
\noindent
case $(\underline{i})_{\overline{i+1}}$$(\underline{i-1})_{\overline{i}}$ :\\
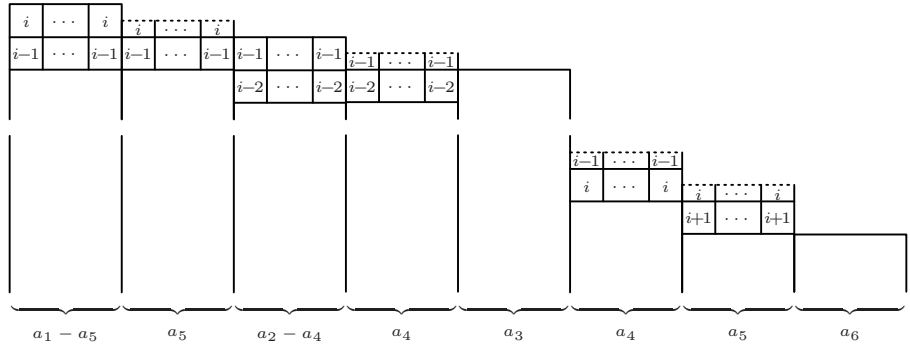
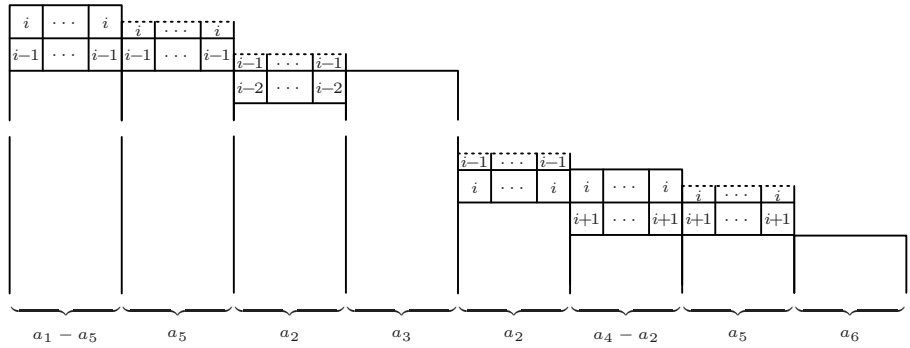
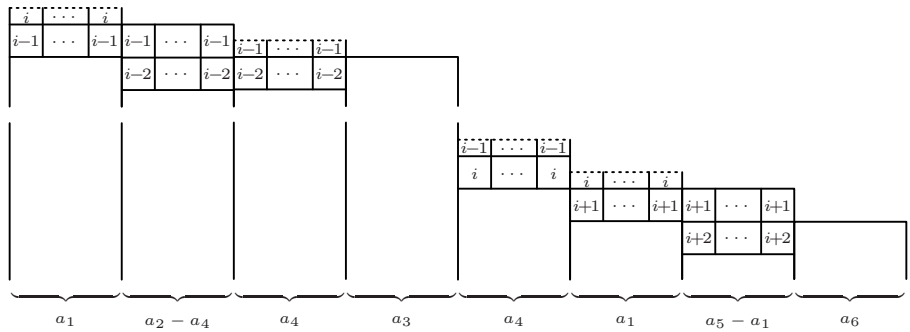
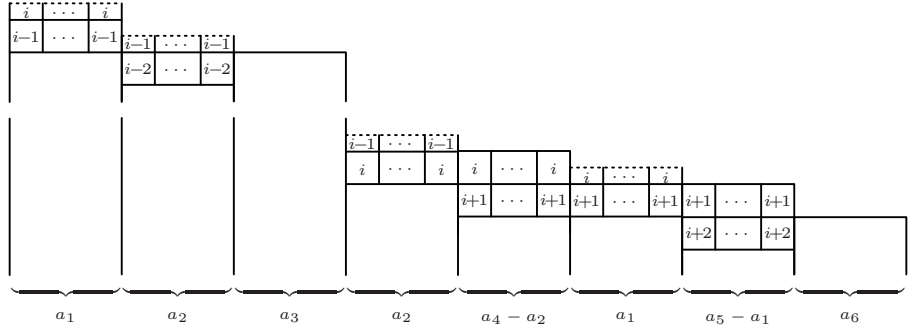
\begin{center}
\usebox{\tmpfiga}
\end{center}
case $(\underline{i})_{\overline{i+1}}$$(\overline{i})_{\underline{i-1}}$ :\\
\begin{center}
\usebox{\tmpfigb}
\end{center}
case $(\overline{i+1})_{\underline{i}}$$(\underline{i-1})_{\overline{i}}$ :\\
\begin{center}
\usebox{\tmpfigc}
\end{center}
case $(\overline{i+1})_{\underline{i}}$$(\overline{i})_{\underline{i-1}}$ :\\
\begin{center}
\usebox{\tmpfigd}
\end{center}

Now, we will place two slices side by side and also
consider the corresponding pair of perfect crystal elements.
When we use the above notations for the \emph{right}
of the two slices,
we will take the number of layers to be given by
$b_i$ instead of $a_i$.
For the left slice, we will use $a_i$ as given in the figure.

%First we fix the right column
%and next remove finitely many $\delta$'s to the left column
%so that the two may be considered as a part of a reduced
%proper Young wall.
We give the following starting shapes and relative height. We fix
the shape of the right slice to be one of the forms
$(\underline{i})_{\overline{i+1}}$$(\underline{i-1})_{\overline{i}}$,
$(\underline{i})_{\overline{i+1}}$$(\overline{i})_{\underline{i-1}}$,
$(\overline{i+1})_{\underline{i}}$$(\underline{i-1})_{\overline{i}}$,
$(\overline{i+1})_{\underline{i}}$$(\overline{i})_{\underline{i-1}}$.
The starting shape of the left slice should be so that when every
possible block is split, all bottom halves of $i$-blocks(stacked
in supporting place) appear at the rear layers.
%Note that there exist group of layers related to
%covering $i$-blocks and slots at the top
%when split every block possible in them.
The following is an example showing
right slice $(\underline{i})_{\overline{i+1}}$$(\overline{i})_{\underline{i-1}}$
and left slice
$(\overline{i+1})_{\underline{i}}$$(\underline{i-1})_{\overline{i}}$. \\[2mm]
\savebox{\tmpfigc}{\begin{texdraw}
\fontsize{6}{6}\selectfont
\textref h:C v:C
\drawdim em
\setunitscale 1.7
\move(0 2.5)
\bsegment
\move(0 0)\lvec(0 3)
\move(3.4 3)\lvec(3.4 0)
\move(0 2.5)\lvec(3.4 2.5)
\move(0 1.5)\lvec(3.4 1.5)
\move(1 3)\lvec(1 1.5)
\move(2.4 3)\lvec(2.4 1.5)
\htext(1.7 2.75){$\cdots$}
\htext(0.5 2.75){$i$}\htext(2.9 2.75){$i$}
\htext(1.7 2){$\cdots$}
\htext(0.5 2){$i\!\!-\!\!1$}\htext(2.9 2){$i\!\!-\!\!1$}
\move(0 -0.5)\rlvec(0 -4.75)
\lpatt(0.05 0.15)
\move(0 3)\lvec(3.4 3)
\esegment
\move(3.4 2)
\bsegment
\move(0 0.5)\lvec(0 3)\lvec(3.4 3)\lvec(3.4 0.5)
\move(0 2)\lvec(3.4 2)
\move(0 1)\lvec(3.4 1)
\move(1 3)\lvec(1 1)
\move(2.4 3)\lvec(2.4 1)
\htext(1.7 2.5){$\cdots$}
\htext(0.5 2.5){$i\!\!-\!\!1$}\htext(2.9 2.5){$i\!\!-\!\!1$}
\htext(1.7 1.5){$\cdots$}
\htext(0.5 1.5){$i\!\!-\!\!2$}\htext(2.9 1.5){$i\!\!-\!\!2$}
\move(0 0)\rlvec(0 -4.75)
\esegment
\move(6.8 1.5)
\bsegment
\move(0 1)\lvec(0 3)
\move(3.4 3)\lvec(3.4 1)
\move(0 2.5)\lvec(3.4 2.5)
\move(0 1.5)\lvec(3.4 1.5)
\move(1 3)\lvec(1 1.5)
\move(2.4 3)\lvec(2.4 1.5)
\htext(1.7 2.75){$\cdots$}
\htext(0.5 2.75){$i\!\!-\!\!1$}\htext(2.9 2.75){$i\!\!-\!\!1$}
\htext(1.7 2){$\cdots$}
\htext(0.5 2){$i\!\!-\!\!2$}\htext(2.9 2){$i\!\!-\!\!2$}
\move(0 0.5)\rlvec(0 -4.75)
\lpatt(0.05 0.15)
\move(0 3)\lvec(3.4 3)
\esegment
\move(10.2 1)
\bsegment
\move(0 1.5)\lvec(0 3)\lvec(3.4 3)\lvec(3.4 1.5)
\move(0 1)\rlvec(0 -4.75)
\esegment
\move(13.6 -1.5)
\bsegment
\move(0 -1.25)\lvec(0 3)
\move(3.4 3)\lvec(3.4 -1.25)
\move(0 2.5)\lvec(3.4 2.5)
\move(0 1.5)\lvec(3.4 1.5)
\move(1 3)\lvec(1 1.5)
\move(2.4 3)\lvec(2.4 1.5)
\htext(1.7 2.75){$\cdots$}
\htext(0.5 2.75){$i\!\!-\!\!1$}\htext(2.9 2.75){$i\!\!-\!\!1$}
\htext(1.7 2){$\cdots$}
\htext(0.5 2){$i$}\htext(2.9 2){$i$}
\move(0 3)\rlvec(0 0.5)
\lpatt(0.05 0.15)
\move(0 3)\lvec(3.4 3)
\esegment
\move(17 -2.5)
\bsegment
\move(0 -0.25)\lvec(0 3)
\move(3.4 3)\lvec(3.4 -0.25)
\move(0 2.5)\lvec(3.4 2.5)
\move(0 1.5)\lvec(3.4 1.5)
\move(1 3)\lvec(1 1.5)
\move(2.4 3)\lvec(2.4 1.5)
\htext(1.7 2.75){$\cdots$}
\htext(0.5 2.75){$i$}\htext(2.9 2.75){$i$}
\htext(1.7 2){$\cdots$}
\htext(0.5 2){$i\!\!+\!\!1$}\htext(2.9 2){$i\!\!+\!\!1$}
\lpatt(0.05 0.15)
\move(0 3)\lvec(3.4 3)
\esegment
\move(20.4 -3)
\bsegment
\move(0 0.25)\lvec(0 3)\lvec(3.4 3)\lvec(3.4 0.25)
\move(0 2)\lvec(3.4 2)
\move(0 1)\lvec(3.4 1)
\move(1 3)\lvec(1 1)
\move(2.4 3)\lvec(2.4 1)
\htext(1.7 2.5){$\cdots$}
\htext(0.5 2.5){$i\!\!+\!\!1$}\htext(2.9 2.5){$i\!\!+\!\!1$}
\htext(1.7 1.5){$\cdots$}
\htext(0.5 1.5){$i\!\!+\!\!2$}\htext(2.9 1.5){$i\!\!+\!\!2$}
\esegment
\move(23.8 -4)
\bsegment
\move(0 1.25)\lvec(0 3)\lvec(3.4 3)\lvec(3.4 1.25)
\esegment
\htext(1.7 -3.25){$\underbrace{\rule{5.6em}{0em}}$}
\htext(5.1 -3.25){$\underbrace{\rule{5.6em}{0em}}$}
\htext(8.5 -3.25){$\underbrace{\rule{5.6em}{0em}}$}
\htext(11.9 -3.25){$\underbrace{\rule{5.6em}{0em}}$}
\htext(15.3 -3.25){$\underbrace{\rule{5.6em}{0em}}$}
\htext(18.7 -3.25){$\underbrace{\rule{5.6em}{0em}}$}
\htext(22.1 -3.25){$\underbrace{\rule{5.6em}{0em}}$}
\htext(25.5 -3.25){$\underbrace{\rule{5.6em}{0em}}$}
\htext(1.7 -4){$a_1$}\htext(5.1 -4){$a_2-a_4$}
\htext(8.5 -4){$a_4$}\htext(11.9 -4){$a_3$}
\htext(15.3 -4){$a_4$}\htext(18.7 -4){$a_1$}
\htext(22.1 -4){$a_5-a_1$}\htext(25.5 -4){$a_6$}
\move(0 -4.3)
%\drawbb
\end{texdraw}}
\savebox{\tmpfigb}{\begin{texdraw}
\fontsize{6}{6}\selectfont
\textref h:C v:C
\drawdim em
\setunitscale 1.7
\move(0 3)
\bsegment
\move(0 -0.5)\lvec(0 3)\lvec(3.4 3)\lvec(3.4 -0.5)
\move(0 2)\lvec(3.4 2)
\move(0 1)\lvec(3.4 1)
\move(1 3)\lvec(1 1)
\move(2.4 3)\lvec(2.4 1)
\htext(1.7 2.5){$\cdots$}
\htext(0.5 2.5){$i$}\htext(2.9 2.5){$i$}
\htext(1.7 1.5){$\cdots$}
\htext(0.5 1.5){$i\!\!-\!\!1$}\htext(2.9 1.5){$i\!\!-\!\!1$}
\move(0 -1)\rlvec(0 -4.75)
\esegment
\move(3.4 2.5)
\bsegment
\move(0 0)\lvec(0 3)
\move(3.4 3)\lvec(3.4 0)
\move(0 2.5)\lvec(3.4 2.5)
\move(0 1.5)\lvec(3.4 1.5)
\move(1 3)\lvec(1 1.5)
\move(2.4 3)\lvec(2.4 1.5)
\htext(1.7 2.75){$\cdots$}
\htext(0.5 2.75){$i$}\htext(2.9 2.75){$i$}
\htext(1.7 2){$\cdots$}
\htext(0.5 2){$i\!\!-\!\!1$}\htext(2.9 2){$i\!\!-\!\!1$}
\move(0 -0.5)\rlvec(0 -4.75)
\lpatt(0.05 0.15)
\move(0 3)\lvec(3.4 3)
\esegment
\move(6.8 1.5)
\bsegment
\move(0 1)\lvec(0 3)
\move(3.4 3)\lvec(3.4 1)
\move(0 2.5)\lvec(3.4 2.5)
\move(0 1.5)\lvec(3.4 1.5)
\move(1 3)\lvec(1 1.5)
\move(2.4 3)\lvec(2.4 1.5)
\htext(1.7 2.75){$\cdots$}
\htext(0.5 2.75){$i\!\!-\!\!1$}\htext(2.9 2.75){$i\!\!-\!\!1$}
\htext(1.7 2){$\cdots$}
\htext(0.5 2){$i\!\!-\!\!2$}\htext(2.9 2){$i\!\!-\!\!2$}
\move(0 0.5)\rlvec(0 -4.75)
\lpatt(0.05 0.15)
\move(0 3)\lvec(3.4 3)
\esegment
\move(10.2 1)
\bsegment
\move(0 1.5)\lvec(0 3)\lvec(3.4 3)\lvec(3.4 1.5)
\move(0 1)\rlvec(0 -4.75)
\esegment
\move(13.6 -1.5)
\bsegment
\move(0 -1.25)\lvec(0 3)\move(3.4 3)\lvec(3.4 -1.25)
\move(0 2.5)\lvec(3.4 2.5)
\move(0 1.5)\lvec(3.4 1.5)
\move(1 3)\lvec(1 1.5)
\move(2.4 3)\lvec(2.4 1.5)
\htext(1.7 2.75){$\cdots$}
\htext(0.5 2.75){$i\!\!-\!\!1$}\htext(2.9 2.75){$i\!\!-\!\!1$}
\htext(1.7 2){$\cdots$}
\htext(0.5 2){$i$}\htext(2.9 2){$i$}
\move(0 3)\rlvec(0 0.5)
\lpatt(0.05 0.15)
\move(0 3)\lvec(3.4 3)
\esegment
\move(17 -2)
\bsegment
\move(0 -0.5)\lvec(0 3)\lvec(3.4 3)\lvec(3.4 -0.5)
\move(0 2)\lvec(3.4 2)
\move(0 1)\lvec(3.4 1)
\move(1 3)\lvec(1 1)
\move(2.4 3)\lvec(2.4 1)
\htext(1.7 2.5){$\cdots$}
\htext(0.5 2.5){$i$}\htext(2.9 2.5){$i$}
\htext(1.7 1.5){$\cdots$}
\htext(0.5 1.5){$i\!\!+\!\!1$}\htext(2.9 1.5){$i\!\!+\!\!1$}
\esegment
\move(20.4 -2.5)
\bsegment
\move(0 -0.25)\lvec(0 3)\move(3.4 3)\lvec(3.4 -0.25)
\move(0 2.5)\lvec(3.4 2.5)
\move(0 1.5)\lvec(3.4 1.5)
\move(1 3)\lvec(1 1.5)
\move(2.4 3)\lvec(2.4 1.5)
\htext(1.7 2.75){$\cdots$}
\htext(0.5 2.75){$i$}\htext(2.9 2.75){$i$}
\htext(1.7 2){$\cdots$}
\htext(0.5 2){$i\!\!+\!\!1$}\htext(2.9 2){$i\!\!+\!\!1$}
\lpatt(0.05 0.15)
\move(0 3)\lvec(3.4 3)
\esegment
\move(23.8 -4)
\bsegment
\move(0 1.25)\lvec(0 3)\lvec(3.4 3)\lvec(3.4 1.25)
\esegment
\htext(1.7 -3.25){$\underbrace{\rule{5.6em}{0em}}$}
\htext(5.1 -3.25){$\underbrace{\rule{5.6em}{0em}}$}
\htext(8.5 -3.25){$\underbrace{\rule{5.6em}{0em}}$}
\htext(11.9 -3.25){$\underbrace{\rule{5.6em}{0em}}$}
\htext(15.3 -3.25){$\underbrace{\rule{5.6em}{0em}}$}
\htext(18.7 -3.25){$\underbrace{\rule{5.6em}{0em}}$}
\htext(22.1 -3.25){$\underbrace{\rule{5.6em}{0em}}$}
\htext(25.5 -3.25){$\underbrace{\rule{5.6em}{0em}}$}
\htext(1.7 -4){$b_1-b_5$}\htext(5.1 -4){$b_5$}
\htext(8.5 -4){$b_2$}\htext(11.9 -4){$b_3$}
\htext(15.3 -4){$b_2$}\htext(18.7 -4){$b_4-b_2$}
\htext(22.1 -4){$b_5$}\htext(25.5 -4){$b_6$}
\move(0 -4.3)
%\drawbb
\end{texdraw}}
right :\\
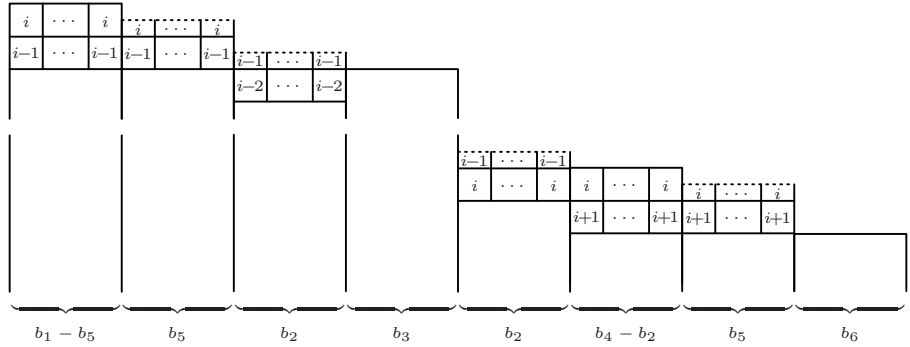
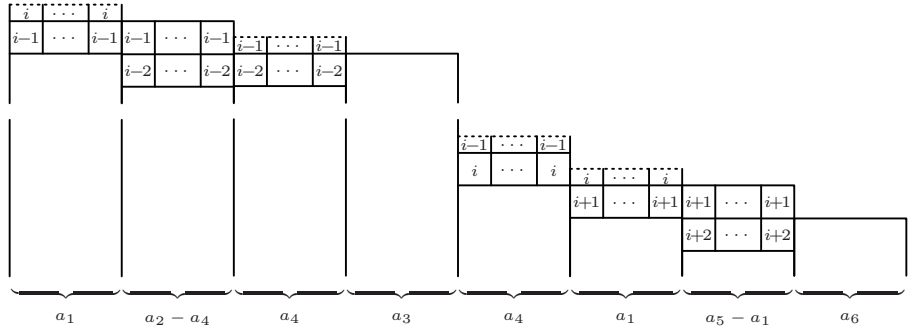
\begin{center}
\usebox{\tmpfigb}
\end{center}
left :\\
\begin{center}
\usebox{\tmpfigc}
\end{center}
Finally, join the two slices in such a way that the highest layer
of the result forms a part of a level-$l$ reduced proper Young wall that has
had all its blocks that may be split, split.

Now, to bring this into a reduced proper form, we need to
\emph{remove} $\delta$'s from the \emph{left} slice.
We denote the number of $\delta$ removals needed by $k$.

Below, we list left-$\vphi$ and right-$\veps$ values
for only the most complicated case among the sixteen possible cases of
Young wall column pairs.
The remaining cases are less complicated.
The line containing the bullet lists the two column
types in left-right order.
Unlike the $i=0,1$ or $i=n$ case, in this type,
we know that there exist two groups of layers with $i$-slots and blocks
at the top, that are apart from each other.
But we shall determine $k$ based only on the shape of one of the two groups
(the group of layers with the supporting part at the top).
So, for some $k$, left-$\vphi$ and right-$\veps$ values
will have  several possibilities.
\begin{itemize}
\item $(\overline{i+1})_{\underline{i}}$$(\underline{i-1})_{\overline{i}}$
      $(\underline{i})_{\overline{i+1}}$$(\overline{i})_{\underline{i-1}}$\\
  \begin{enumerate}
  \item $a_5-a_1\le a_2-a_4$, $b_4-b_2\le b_1-b_5$, $a_5\le (b_4-b_2)+b_5$
  \ or\quad
  $a_5-a_1\le a_2-a_4$, $b_4-b_2\ge b_1-b_5$, $a_5\le b_1$\\[2mm]
  \begin{tabular}{rcl}
  rotation & : & $0\le k\le (a_2-a_4)-(a_5-a_1)$\\
  left-$\vphi$ & : & $(a_5-a_1)+k$\\
  right-$\veps$ & : & $(b_4-b_2)+b_1-(a_1+(a_2-a_4)-k)$\\[2mm]
  rotation & : & $(a_2-a_4)-(a_5-a_1)\le k\le a_1+(a_2-a_4)$\\
  left-$\vphi$ & : & $(a_2-a_4)$\\
  right-$\veps$ & : & $(b_4-b_2)+b_1-a_5$\\
  left-$\vphi$ & : & $(a_2-a_4)+1$\\
  right-$\veps$ & : & $(b_4-b_2)+b_1-a_5+1$\\
  \quad\ $\vdots$ & & \quad\ $\vdots$\\
  left-$\vphi$ & : & $(a_5-a_1)+k$\\
  right-$\veps$ & : & $(b_4-b_2)+b_1-(a_1+(a_2-a_4)-k)$\\[2mm]
  rotation & : & $a_1+(a_2-a_4)\le k$\\
  left-$\vphi$ & : & $(a_2-a_4)$\\
  right-$\veps$ & : & $(b_4-b_2)+b_1-a_5$\\
  left-$\vphi$ & : & $(a_2-a_4)+1$\\
  right-$\veps$ & : & $(b_4-b_2)+b_1-a_5+1$\\
  \quad\ $\vdots$ & & \quad\ $\vdots$\\
  left-$\vphi$ & : & $a_5+(a_2-a_4)$\\
  right-$\veps$ & : & $(b_4-b_2)+b_1$\\
  \end{tabular}\\[2mm]
  \item $a_5-a_1\le a_2-a_4$, $b_4-b_2\le b_1-b_5$, $a_5\ge (b_4-b_2)+b_5$\\[2mm]
  \begin{tabular}{rcl}
  rotation & : & $0\le k\le (a_2-a_4)+a_1-(b_4-b_2)-b_5$\\
  left-$\vphi$ & : & $(a_5-a_1)+k$\\
  right-$\veps$ & : & $(b_4-b_2)+b_1-(a_1+(a_2-a_4)-k)$\\[2mm]
  rotation & : & $(a_2-a_4)+a_1-(b_4-b_2)-b_5\le k\le a_1+(a_2-a_4)$\\
  left-$\vphi$ & : & $(a_2-a_4)+a_5-(b_4-b_2)-b_5$\\
  right-$\veps$ & : & $b_1-b_5$\\
  left-$\vphi$ & : & $(a_2-a_4)+a_5-(b_4-b_2)-b_5+1$\\
  right-$\veps$ & : & $b_1-b_5+1$\\
  \quad\ $\vdots$ & & \quad\ $\vdots$\\
  left-$\vphi$ & : & $(a_5-a_1)+k$\\
  right-$\veps$ & : & $(b_4-b_2)+b_1-(a_1+(a_2-a_4)-k)$\\[2mm]
  rotation & : & $a_1+(a_2-a_4)\le k$\\
  left-$\vphi$ & : & $(a_2-a_4)+a_5-(b_4-b_2)-b_5$\\
  right-$\veps$ & : & $b_1-b_5$\\
  left-$\vphi$ & : & $(a_2-a_4)+a_5-(b_4-b_2)-b_5+1$\\
  right-$\veps$ & : & $b_1-b_5+1$\\
  \quad\ $\vdots$ & & \quad\ $\vdots$\\
  left-$\vphi$ & : & $a_5+(a_2-a_4)$\\
  right-$\veps$ & : & $(b_4-b_2)+b_1$\\
  \end{tabular}\\[2mm]
  \item $a_5-a_1\le a_2-a_4$, $b_4-b_2\ge b_1-b_5$, $a_5\ge b_1$\\[2mm]
  \begin{tabular}{rcl}
  rotation & : & $0\le k\le (a_2-a_4)+a_1-b_1$\\
  left-$\vphi$ & : & $(a_5-a_1)+k$\\
  right-$\veps$ & : & $(b_4-b_2)+b_1-(a_1+(a_2-a_4)-k)$\\[2mm]
  rotation & : & $(a_2-a_4)+a_1-b_1\le k\le a_1+(a_2-a_4)$\\
  left-$\vphi$ & : & $(a_2-a_4)+a_5-b_1$\\
  right-$\veps$ & : & $b_4-b_2$\\
  left-$\vphi$ & : & $(a_2-a_4)+a_5-b_1+1$\\
  right-$\veps$ & : & $b_4-a_2+1$\\
  \quad\ $\vdots$ & & \quad\ $\vdots$\\
  left-$\vphi$ & : & $(a_5-a_1)+k$\\
  right-$\veps$ & : & $(b_4-b_2)+b_1-(a_1+(a_2-a_4)-k)$\\[2mm]
  rotation & : & $a_1+(a_2-a_4)\le k$\\
  left-$\vphi$ & : & $(a_2-a_4)+a_5-b_1$\\
  right-$\veps$ & : & $b_4-b_2$\\
  left-$\vphi$ & : & $(a_2-a_4)+a_5-b_1+1$\\
  right-$\veps$ & : & $b_4-b_2+1$\\
  \quad\ $\vdots$ & & \quad\ $\vdots$\\
  left-$\vphi$ & : & $a_5+(a_2-a_4)$\\
  right-$\veps$ & : & $(b_4-b_2)+b_1$\\
  \end{tabular}\\[2mm]
  \item $a_5-a_1\ge a_2-a_4$, $b_4-b_2\le b_1-b_5$,
        $a_1+(a_2-a_4)\le (b_4-b_2)+b_5$
        \ or \quad
        $a_5-a_1\ge a_2-a_4$, $b_4-b_2\ge b_1-b_5$,
        $a_1+(a_2-a_4)\le b_1$\\[2mm]
  \begin{tabular}{rcl}
  rotation & : & $0\le k\le (a_5-a_1)-(a_2-a_4)$\\
  left-$\vphi$ & : & $(a_2-a_4)+k$\\
  right-$\veps$ & : & $(b_4-b_2)+b_1-(a_5-k)$\\[2mm]
  rotation & : & $(a_5-a_1)-(a_2-a_4)\le k\le a_5$\\
  left-$\vphi$ & : & $(a_5-a_1)$\\
  right-$\veps$ & : & $(b_4-b_2)+b_1-(a_1+(a_2-a_4))$\\
  left-$\vphi$ & : & $(a_5-a_1)+1$\\
  right-$\veps$ & : & $(b_4-b_2)+b_1-(a_1+(a_2-a_4))+1$\\
  \quad\ $\vdots$ & & \quad\ $\vdots$\\
  left-$\vphi$ & : & $(a_2-a_4)+k$\\
  right-$\veps$ & : & $(b_4-b_2)+b_1-a_5+k$\\[2mm]
  rotation & : & $a_5\le k$\\
  left-$\vphi$ & : & $(a_5-a_1)$\\
  right-$\veps$ & : & $(b_4-b_2)+b_1-(a_1+(a_2-a_4))$\\
  left-$\vphi$ & : & $(a_5-a_1)+1$\\
  right-$\veps$ & : & $(b_4-b_2)+b_1-(a_1+(a_2-a_4))+1$\\
  \quad\ $\vdots$ & & \quad\ $\vdots$\\
  left-$\vphi$ & : & $(a_2-a_4)+a_5$\\
  right-$\veps$ & : & $(b_4-b_2)+b_1$\\
  \end{tabular}\\[2mm]
  \item $a_5-a_1\ge a_2-a_4$, $b_4-b_2\le b_1-b_5$,
        $a_1+(a_2-a_4)\ge (b_4-b_2)+b_5$\\[2mm]
  \begin{tabular}{rcl}
  rotation & : & $0\le k\le a_5-b_5-(b_4-b_2)$\\
  left-$\vphi$ & : & $(a_2-a_4)+k$\\
  right-$\veps$ & : & $(b_4-b_2)+b_1-(a_5-k)$\\[2mm]
  rotation & : & $a_5-b_5-(b_4-b_2)\le k\le a_5$\\
  left-$\vphi$ & : & $(a_2-a_4)+a_5-b_5-(b_4-b_2)$\\
  right-$\veps$ & : & $b_1-b_5$\\
  left-$\vphi$ & : & $(a_2-a_4)+a_5-b_5-(b_4-b_2)+1$\\
  right-$\veps$ & : & $b_1-b_5+1$\\
  \quad\ $\vdots$ & & \quad\ $\vdots$\\
  left-$\vphi$ & : & $(a_2-a_4)+k$\\
  right-$\veps$ & : & $(b_4-b_2)+b_1-a_5+k$\\[2mm]
  rotation & : & $a_5\le k$\\
  left-$\vphi$ & : & $(a_2-a_4)+a_5-b_5-(b_4-b_2)$\\
  right-$\veps$ & : & $b_1-b_5$\\
  left-$\vphi$ & : & $(a_2-a_4)+a_5-b_5-(b_4-b_2)+1$\\
  right-$\veps$ & : & $b_1-b_5+1$\\
  \quad\ $\vdots$ & & \quad\ $\vdots$\\
  left-$\vphi$ & : & $(a_2-a_4)+a_5$\\
  right-$\veps$ & : & $(b_4-b_2)+b_1$\\
  \end{tabular}\\[2mm]
  \item $a_5-a_1\ge a_2-a_4$, $b_4-b_2\ge b_1-b_5$,
        $a_1+(a_2-a_4)\ge b_1$\\[2mm]
  \begin{tabular}{rcl}
  rotation & : & $0\le k\le a_5-b_1$\\
  left-$\vphi$ & : & $(a_2-a_4)+k$\\
  right-$\veps$ & : & $(b_4-b_2)+b_1-(a_5-k)$\\[2mm]
  rotation & : & $a_5-b_1\le k\le a_5$\\
  left-$\vphi$ & : & $(a_2-a_4)+a_5-b_1$\\
  right-$\veps$ & : & $(b_4-b_2)$\\
  left-$\vphi$ & : & $(a_2-a_4)+a_5-b_1+1$\\
  right-$\veps$ & : & $(b_4-b_2)+1$\\
  \quad\ $\vdots$ & & \quad\ $\vdots$\\
  left-$\vphi$ & : & $(a_2-a_4)+k$\\
  right-$\veps$ & : & $(b_4-b_2)+b_1-a_5+k$\\[2mm]
  rotation & : & $a_5\le k$\\
  left-$\vphi$ & : & $(a_2-a_4)+a_5-b_1$\\
  right-$\veps$ & : & $(b_4-b_2)$\\
  left-$\vphi$ & : & $(a_2-a_4)+a_5-b_1+1$\\
  right-$\veps$ & : & $(b_4-b_2)+1$\\
  \quad\ $\vdots$ & & \quad\ $\vdots$\\
  left-$\vphi$ & : & $(a_2-a_4)+a_5$\\
  right-$\veps$ & : & $(b_4-b_2)+b_1$\\
  \end{tabular}\\
  \end{enumerate}
\end{itemize}
Similarly, the following gives the signatures of the
path description.
The numbers in the list are the
left-$\vphi$ and right-$\veps$ values for the two
corresponding crystal elements.
\begin{itemize}
\item Any pair with
      $(\underline{i})_{\overline{i+1}}$$(\underline{i-1})_{\overline{i}}$ or
      $(\overline{i+1})_{\underline{i}}$$(\underline{i-1})_{\overline{i}}$
      in both the left and right columns\\[2mm]
 \begin{tabular}{rcl}
 left-$\vphi$ & : & $a_5+(a_2-a_4)$\\
 right-$\veps$ & : & $b_1$\\
 \end{tabular}\\[2mm]
\item Any pair with
      $(\underline{i})_{\overline{i+1}}$$(\overline{i})_{\underline{i-1}}$ or
      $(\overline{i+1})_{\underline{i}}$$(\overline{i})_{\underline{i-1}}$
      in both the left and right columns\\[2mm]
 \begin{tabular}{rcl}
 left-$\vphi$ & : & $a_5$\\
 right-$\veps$ & : & $b_1+(b_4-b_2)$\\
 \end{tabular}\\[2mm]
\item Any pair with
      $(\underline{i})_{\overline{i+1}}$$(\underline{i-1})_{\overline{i}}$ or
      $(\overline{i+1})_{\underline{i}}$$(\underline{i-1})_{\overline{i}}$
      in the left column,
      $(\underline{i})_{\overline{i+1}}$$(\overline{i})_{\underline{i-1}}$
      or
      $(\overline{i+1})_{\underline{i}}$$(\overline{i})_{\underline{i-1}}$
      in the right column\\[2mm]
 \begin{tabular}{rcl}
 left-$\vphi$ & : & $a_5+(a_2-a_4)$\\
 right-$\veps$ & : & $b_1+(b_4-b_2)$\\
 \end{tabular}\\[2mm]
\item Any pair with
      $(\underline{i})_{\overline{i+1}}$$(\overline{i})_{\underline{i-1}}$
      or
      $(\overline{i+1})_{\underline{i}}$$(\overline{i})_{\underline{i-1}}$
      in the left column,
      $(\underline{i})_{\overline{i+1}}$$(\underline{i-1})_{\overline{i}}$
      or
      $(\overline{i+1})_{\underline{i}}$$(\underline{i-1})_{\overline{i}}$
      in the right column\\[2mm]
 \begin{tabular}{rcl}
 left-$\vphi$ & : & $a_5$\\
 right-$\veps$ & : & $b_1$\\
 \end{tabular}\\[2mm]
\end{itemize}
We can easily see that the signatures agree with those of the
corresponding path description in all of the cases after
$(0,1)$-pair cancellations.

This completes the proof of Lemma~5.5.

\end{document}